\documentclass[reqno, 8pt]{amsart}

\usepackage{amsmath,amssymb,amsthm,amsfonts,mathrsfs}
\usepackage{fullpage}
\usepackage{xcolor}
\usepackage{graphicx}
\usepackage{makecell}
\usepackage{listings}
\usepackage{bbm}

\definecolor{codegreen}{rgb}{0,0.6,0}
\definecolor{codegray}{rgb}{0.5,0.5,0.5}
\definecolor{codepurple}{rgb}{0.58,0,0.82}
\definecolor{backcolour}{rgb}{0.95,0.95,0.92}

\lstdefinestyle{mystyle}{
  backgroundcolor=\color{backcolour},   commentstyle=\color{codegreen},
  keywordstyle=\color{magenta},
  numberstyle=\tiny\color{codegray},
  stringstyle=\color{codepurple},
  basicstyle=\ttfamily\footnotesize,
  breakatwhitespace=false,         
  breaklines=true,                 
  captionpos=b,                    
  keepspaces=true,                 
  numbers=left,                    
  numbersep=5pt,                  
  showspaces=false,                
  showstringspaces=false,
  showtabs=false,                  
  tabsize=2
}

\usepackage{listings}
\usepackage[colorlinks=true,
linkcolor=blue,
anchorcolor=blue,
citecolor=violet
]{hyperref}
\allowdisplaybreaks 

\newtheorem{theorem}{Theorem}[section]
\newtheorem{lemma}[theorem]{Lemma}

\newtheorem*{conjecture*}{Conjecture}
\newtheorem*{theorem*}{Theorem}
\theoremstyle{definition}

\theoremstyle{remark}

\newtheorem*{remark*}{Remark}

\newcommand{\on}{\operatorname}

\renewcommand{\mod}[1]{\on{mod}#1}

\usepackage{colonequals}

\author{\textbf{Runbo Li}}
\address{International Curriculum Center, The High School Affiliated to Renmin University of China, Beijing, China}
\email{runbo.li.carey@gmail.com}

\makeatletter
\@namedef{subjclassname@2020}{\textup{}2020 Mathematics Subject Classification}
\makeatother

\title[]{Primes in arithmetic progressions to large moduli and refinements of Harman's sieve}

\begin{document}

\begin{abstract}
We study the average distribution of primes of size $x$ in arithmetic progressions to moduli larger than $x^{\frac{1}{2}}$. Using arithmetic information from the works of many authors together with different variants of the original Harman's sieve, we construct suitable majorants and minorants for the prime indicator function $\mathbbm{1}_{p}(n)$ that satisfy Bombieri--Vinogradov type mean value theorems with different types of moduli. Specifically, we obtain some mean value theorems for primes with bilinear forms of moduli up to $x^{\frac{9}{17}}$ or with trilinear forms of moduli up to $x^{\frac{17}{32}}$. As a by-product, we obtain new upper and lower bounds for $\pi(x; q, a)$ that hold for almost all moduli $q$.
\end{abstract}

\maketitle

\tableofcontents

\section{Introduction}
One of the famous topics in prime number theory is the distribution of primes in arithmetic progressions. The Siegel--Walfisz Theorem states that
\begin{equation}
\left| \pi(x; q, a) - \frac{\pi(x)}{\varphi(q)} \right| \ll \frac{x}{(\log x)^A}
\end{equation}
uniformly for $q \leqslant (\log x)^A$ and $(a, q) = 1$. Under the Generalized Riemann Hypothesis, the range of $q$ such that (1) holds can be extended to $q \leqslant x^{\frac{1}{2}} (\log x)^{-B}$, where $B = B(A)$. In \cite{MontgomeryEH}, Montgomery conjectured that (1) should hold for all $q \leqslant x^{1 - \varepsilon}$.

In many applications, however, mathematicians only need an average distribution result of primes in arithmetic progressions rather than an individual estimate (1). In this case, a substitute of the Generalized Riemann Hypothesis is the Bombieri--Vinogradov Theorem \cite{BombieriBV} \cite{VinogradovBV}. This theorem states that
\begin{equation}
\sum_{q \leqslant x^{\frac{1}{2}} (\log x)^{-B} } \max_{(a, q) = 1} \left| \pi(x; q, a) - \frac{\pi(x)}{\varphi(q)} \right| \ll \frac{x}{(\log x)^A}.
\end{equation}
By (2), we also know that (1) holds for almost all $q \leqslant x^{\frac{1}{2}} (\log x)^{-B}$ with $(a, q) = 1$. As well as (1), Elliott and Halberstam \cite{ElliottHalberstam} conjectured that the range of $q$ in (2) can be extended to $q \leqslant x^{1 - \varepsilon}$. However, one still cannot prove (2) even with $q \leqslant x^{\frac{1}{2} + \varepsilon}$ now.

From here, we suppose that $a \in \mathbb{Z} \backslash \{0\}$ is fixed. In 2013, Zhang \cite{ZhangYitang} proved that (2) is valid for $q \leqslant x^{\frac{293}{584} - \varepsilon}$ when the moduli $q$ is square-free and only has small prime factors:
\begin{equation}
\sum_{\substack{q \leqslant x^{\frac{293}{584} - \varepsilon} \\ q \mid P(x^{\delta}) \\ (q, a) = 1}} \left| \pi(x; q, a) - \frac{\pi(x)}{\varphi(q)} \right| \ll \frac{x}{(\log x)^A},
\end{equation}
where $\delta = \delta(\varepsilon)$ is a small positive number. Polymath8a \cite{Polymath8a} and Stadlmann \cite{Stadlmann525} further extended the range of $q$ in (3) to $q \leqslant x^{\frac{157}{300} - \varepsilon}$ and $q \leqslant x^{\frac{21}{40} - \varepsilon}$.

In 2025, Maynard \cite{MaynardLargeModuliI} proved that
\begin{equation}
\sum_{\substack{q_1 \leqslant Q_1 \\ q_2 \leqslant Q_2 \\ (q_1 q_2, a) = 1}} \left| \pi(x; q_1 q_2, a) - \frac{\pi(x)}{\varphi(q_1 q_2)} \right| \ll \frac{x}{(\log x)^A}
\end{equation}
if
$$
Q_1^2 Q_2 < x^{1 - \varepsilon}, \quad Q_1^7 Q_2^{12} < x^{4 - \varepsilon}, \quad Q_1^{19} Q_2^{20} < x^{10 - \varepsilon}.
$$
In another paper, he \cite{MaynardLargeModuliIII} also proved that
\begin{equation}
\sum_{\substack{q_1 \leqslant Q_1 \\ q_2 \leqslant Q_2 \\ q_3 \leqslant Q_3 \\ (q_1 q_2 q_3, a) = 1}} \left| \pi(x; q_1 q_2 q_3, a) - \frac{\pi(x)}{\varphi(q_1 q_2 q_3)} \right| \ll \frac{x}{(\log x)^A}
\end{equation}
if $\delta \in (0,\ 0.001)$ and
$$
Q_1 Q_2 Q_3 = x^{\frac{1}{2} + \delta}, \quad x^{40 \delta} < Q_2 < x^{\frac{1}{20} - 7 \delta}, \quad x^{\frac{1}{10} + 12 \delta} Q_2^{-1} < Q_3 < x^{\frac{1}{10} - 4 \delta} Q_2^{-\frac{3}{5}}.
$$

The problem of bounding sums of the form
$$
\sum_{q \leqslant Q } \left| \pi(x; q, a) - \frac{\pi(x)}{\varphi(q)} \right|
$$
can be seen as a problem equivalent to bounding 
$$
\sum_{q \leqslant Q } \lambda_{q} \left( \pi(x; q, a) - \frac{\pi(x)}{\varphi(q)} \right),
$$
for an arbitrary divisor-bounded weight $\lambda_{q}$. Naturally, we can try to make the weight $\lambda_{q}$ more ``flexible'' to extend the range of $Q$ beyond $x^{\frac{1}{2}}$. One example is to introduce the ``bilinear weights'' and prove results of the following form:
\begin{equation}
\sum_{\substack{q_1 \leqslant Q_1 \\ q_2 \leqslant Q_2 \\ (q_1 q_2, a) = 1}} \lambda_{1, q_1} \lambda_{2, q_2} \left( \pi(x; q_1 q_2, a) - \frac{\pi(x)}{\varphi(q_1 q_2)} \right) \ll \frac{x}{(\log x)^A}.
\end{equation}
In 1986, Bombieri, Friedlander and Iwaniec \cite{BFI} first showed that (6) holds if
$$
Q_1 < x^{\frac{1}{3}}, \quad Q_2 < x^{\frac{1}{5}}, \quad Q_1^{5} Q_2^{2} < x^{2}, \quad Q_1 Q_2 < x^{\frac{29}{56}}.
$$
In 1987, Fouvry \cite{FouvryA2} proved that (6) holds if either
$$
Q_1 Q_2^{3} < x, \quad Q_1 Q_2 < x^{\frac{29}{56}}, \quad Q_1 < \max\left(x^{\frac{1}{2}} Q_2^{-1}, x^{\frac{2}{5}} Q_2^{-\frac{2}{5}} \right)
$$
or
$$
Q_1 Q_2^{3} < x, \quad Q_1 Q_2 < x^{\frac{29}{56}}, \quad Q_1^{4} Q_2 < x^{\frac{403}{266}}, \quad Q_1^{\frac{7}{4}} Q_2 < x^{\frac{403}{532}}.
$$
In 1998, Baker and Harman \cite{677} improved the result of Bombieri, Friedlander and Iwaniec \cite{BFI} by removing the condition $Q_1^{5} Q_2^{2} < x^{2}$ in their result above. They showed that (6) holds if
$$
Q_1 < x^{\frac{1}{3}}, \quad Q_2 < x^{\frac{1}{5}}, \quad Q_1 Q_2 < x^{\frac{29}{56}}.
$$
Those results extended the range of $Q = Q_1 Q_2$ up to $x^{\frac{29}{56}}$ in some special cases. The main result of Maynard \cite{MaynardLargeModuliI} can also be seen as a stronger ``bilinear'' result, with $Q = Q_1 Q_2$ up to $x^{\frac{11}{21}}$. In 2022, Lichtman \cite{Lichtman2} considered a more ``flexible'' case with quadrilinear weights and extended the moduli up to $x^{\frac{17}{32}}$. He showed that
\begin{equation}
\sum_{\substack{q_1 \leqslant Q_1 \\ q_2 \leqslant Q_2 \\ q_3 \leqslant Q_3 \\ q_4 \leqslant Q_3 \\ (q_1 q_2 q_3 q_4, a) = 1}} \lambda_{1, q_1} \lambda_{2, q_2} \lambda_{3, q_3} \lambda_{4, q_4} \left( \pi(x; q_1 q_2 q_3 q_4, a) - \frac{\pi(x)}{\varphi(q_1 q_2 q_3 q_4)} \right) \ll \frac{x}{(\log x)^A}
\end{equation}
holds if
$$
Q_1 Q_2 < x^{\frac{1}{2} + \varepsilon}, \quad Q_1 Q_3^{2} < x^{\frac{1}{2} - 2 \varepsilon}, \quad Q_3^{2} < Q_2 < x^{\frac{1}{32} - \varepsilon}.
$$

Weights that are more ``flexible'' than above weights are the well-factorable weights, which means that for any $Q_1, Q_2$ such that $Q_1 Q_2 = Q$, we can ``split'' a well-factorable function $\lambda_{q}$ to $\lambda_{q} = \lambda_{1, q_1} * \lambda_{2, q_2}$ supported on $[1, Q_1]$ and $[1, Q_2]$ respectively. There are also lots of works on this topic, and we refer the readers to \cite{FouvryIwaniec}, \cite{FouvryA1}, \cite{BFI}, \cite{MaynardLargeModuliII}, \cite{Lichtman}, \cite{Lichtman3}, \cite{Pascadi2} and \cite{ZhengZongkun}.

In 1996, Baker and Harman \cite{676} considered a different variant of (2): They used Harman's sieve \cite{HarmanBOOK} to construct majorants and minorants for the prime indicator function $\mathbbm{1}_{p}(n)$ and studied their behaviors in fixed residue classes with moduli $q \geqslant x^{\frac{1}{2}}$. Write $q \asymp x^{\theta}$. Baker and Harman \cite{676} constructed majorants $\rho_1(n) \geqslant \mathbbm{1}_{p}(n)$ for $0.5 \leqslant \theta \leqslant 0.56$ and minorants $\rho_0(n) \leqslant \mathbbm{1}_{p}(n)$ for $0.5 \leqslant \theta \leqslant 0.52$. They used their majorants and minorants to study the average Brun--Titchmarsh theorem and the shifted prime problem, and similar results before them are obtained by Motohashi \cite{Motohashi1970}, Hooley \cite{HooleyBT1} \cite{Hooleypa}, Iwaniec \cite{IwaniecBT}, Deshouillers and Iwaniec \cite{DI1984}, Fouvry \cite{Fouvry1984} \cite{FouvryFermat} and Rousselet \cite{Rousselet} respectively. In 2001, Mikawa \cite{Mikawa} further constructed minorants $\rho_0(n) \leqslant \mathbbm{1}_{p}(n)$ for $0.5 \leqslant \theta < \frac{17}{32}$ using a different sieve method.

In this paper, we refine the methods developed by Baker and Harman \cite{676} and Mikawa \cite{Mikawa} to construct majorants and minorants for the prime indicator function $\mathbbm{1}_{p}(n)$, and we study the average distributions of them in fixed residue classes with different types of moduli. One result we obtain in this paper is the following theorem.
\begin{theorem}\label{t11}
Let $Q_1$ and $Q_2$ satisfy
$$
Q_1^2 Q_2 < x^{1 - \varepsilon} \qquad \text{and} \qquad Q_1^7 Q_2^{12} < x^{4 - \varepsilon}. 
$$
Let $\lambda_{1, q_1}$ and $\lambda_{2, q_2}$ be divisor-bounded complex sequences. Then, for any fixed $a \in \mathbb{Z} \backslash \{0\}$ and any $A>0$, we have
$$
\sum_{\substack{q_1 \leqslant Q_1 \\ q_2 \leqslant Q_2 \\ (q_1 q_2, a) = 1}} \lambda_{1, q_1} \lambda_{2, q_2} \left( \pi(x; q_1 q_2, a) - \frac{\pi(x)}{\varphi(q_1 q_2)} \right) \ll \frac{x}{(\log x)^A}.
$$
\end{theorem}

Throughout this paper, we always suppose that $x$ is sufficiently large, $\varepsilon > 0$ is sufficiently small, $A, B > 0$ (may depend on other variables), and $Q, Q_i, M_i \geqslant 1$. We shall ignore the presence of $\varepsilon$ in many places below (mainly in the definitions of regions $\boldsymbol{A}$, $\boldsymbol{B}$, $\boldsymbol{C}$, $\boldsymbol{E}$, $\boldsymbol{F}$ and their subregions) for clarity. Let $\theta, \theta_i \in (0,\ 1)$ and $\delta = 10^{-100}$. The letters $p$ and $\beta$, with or without subscript, are reserved for primes and almost-primes respectively. We put $p_i \asymp x^{\alpha_i}$ and write $\boldsymbol{\alpha}_{n}$ to denote $(\alpha_1, \ldots, \alpha_n)$. We use $P^{+}(n)$ and $P^{-}(n)$ to denote the largest and smallest prime factor of $n$. We shall use the terms \textit{partition} and \textit{exactly partition} many times in the rest of our paper, and one can see [\cite{HarmanBOOK}, Page 162] for a definition. Put
$$
P(z)=\prod_{p<z} p, \quad \psi(n, z) = \mathbbm{1}_{(n, P(z))=1}, \quad \Psi(n, z) = \mathbbm{1}_{P^{+}(n) < z}.
$$
Let $\mathcal{C}$ denote a finite set of positive integers and put
$$
\mathcal{C}_d=\left\{n: n d \in \mathcal{C}\right\}, \quad S(\mathcal{C}, z) = \sum_{n \in \mathcal{C}} \psi(n, z) = \sum_{\substack{n \in \mathcal{C} \\ (n, P(z))=1}} 1.
$$

\textit{Buchstab's identity} is the equation
$$
S\left(\mathcal{C}, z\right) = S\left(\mathcal{C}, w\right) - \sum_{w \leqslant p < z} S\left(\mathcal{C}_{p}, p\right),
$$
where $2 \leqslant w < z$.

Let $\omega(u)$ denote the Buchstab function determined by the following differential-difference equation
\begin{align*}
\begin{cases}
\omega(u)=\frac{1}{u}, & \quad 1 \leqslant u \leqslant 2, \\
(u \omega(u))^{\prime}= \omega(u-1), & \quad u \geqslant 2 .
\end{cases}
\end{align*}
Moreover, we have the upper and lower bounds for $\omega(u)$:
\begin{align*}
\omega(u) \geqslant \omega_{0}(u) =
\begin{cases}
\frac{1}{u}, & \quad 1 \leqslant u < 2, \\
\frac{1+\log(u-1)}{u}, & \quad 2 \leqslant u < 3, \\
\frac{1+\log(u-1)}{u} + \frac{1}{u} \int_{2}^{u-1}\frac{\log(t-1)}{t} d t \geqslant 0.5607, & \quad 3 \leqslant u < 4, \\
0.5612, & \quad u \geqslant 4, \\
\end{cases}
\end{align*}
\begin{align*}
\omega(u) \leqslant \omega_{1}(u) =
\begin{cases}
\frac{1}{u}, & \quad 1 \leqslant u < 2, \\
\frac{1+\log(u-1)}{u}, & \quad 2 \leqslant u < 3, \\
\frac{1+\log(u-1)}{u} + \frac{1}{u} \int_{2}^{u-1}\frac{\log(t-1)}{t} d t \leqslant 0.5644, & \quad 3 \leqslant u < 4, \\
0.5617, & \quad u \geqslant 4. \\
\end{cases}
\end{align*}
Sometimes we also need this simple equality in our numerical calculations:
$$
\int_{1 - \theta}^{\frac{1}{2}} \frac{\omega\left(\frac{1-t}{t} \right)}{t^2} d t = \int_{1 - \theta}^{\frac{1}{2}} \frac{1}{t (1-t)} d t = \log\left(\frac{\theta}{1 - \theta} \right) \quad \text{for} \quad \frac{1}{2} < \theta < \frac{2}{3}.
$$

In this paper, a ``Type-I$_{j}$'' sum refers to a sum of type
$$
\sum_{m_0, m_1, \ldots, m_j} a_{0, m_0},
$$
and a ``Type-II$_{j}$'' sum refers to a sum of type
$$
\sum_{m_1, \ldots, m_j} a_{1, m_1} \cdots a_{j, m_j},
$$
where $a_{i, m_i}$ ($0 \leqslant i \leqslant j$) are divisor-bounded complex sequences. For the sake of simplicity, we often write ``Type-I$_{1}$'' as ``Type-I'' and ``Type-II$_{2}$'' as ``Type-II''.

\section{General Moduli}
In this section we focus on the general case, where the moduli $q \sim Q = x^{\theta}$. We put
$$
\mathcal{A}^q =\left\{n : n \sim x,\ n \equiv a (\bmod q) \right\} \quad \text{and} \quad \mathcal{B}^q =\left\{n : n \sim x,\ (n, q) = 1 \right\}.
$$
By the definitions of the sieved sets $\mathcal{A}^{q}$, $\mathcal{B}^{q}$ and the sieve function $S\left(\mathcal{C}, z\right)$, using Prime Number Theorem, we have
\begin{equation}
\pi(2 x; q, a) - \pi(x; q, a) = \sum_{p \in \mathcal{A}^{q}} 1 = S\left(\mathcal{A}^{q}, (2x)^{\frac{1}{2}} \right) \quad \text{and} \quad S\left(\mathcal{B}^{q}, (2x)^{\frac{1}{2}} \right) = (1+o(1)) \frac{x}{\log x}.
\end{equation}
Our aim is to show that the sparser set $\mathcal{A}^{q}$ contains the expected proportion of primes compared to the larger set $\mathcal{B}^{q}$, which requires us to decompose $S\left(\mathcal{A}^{q}, (2x)^{\frac{1}{2}} \right)$ and prove ``asymptotic formulas'' of the form
\begin{equation}
\sum_{\substack{q \sim Q \\ (q, a) = 1}} \left| S\left(\mathcal{A}^{q}, z \right) - \frac{1}{\varphi(q)} S\left(\mathcal{B}^{q}, z \right) \right| \ll \frac{x}{(\log x)^A}.
\end{equation}
for some parts of it, and drop the remaining parts to construct a suitable majorant or minorant. For the majorant case we can only drop negative parts, while for the minorant case we can only drop positive parts. After the final decompositions, we can get the following result with some $0 < C_0(\theta) \leqslant 1$ and $C_1(\theta) \geqslant 1$:
\begin{theorem}\label{t21}
There exist functions $\rho_0$ and $\rho_1$ which satisfies the following properties:

(Majorant / Minorant). $\rho_0(n)$ is a minorant for the prime indicator function $\mathbbm{1}_{p}(n)$, and $\rho_1(n)$ is a majorant for the prime indicator function $\mathbbm{1}_{p}(n)$. That is, we have
$$
\rho_0(n) \leqslant \mathbbm{1}_{p}(n) \leqslant \rho_1(n).
$$

(Upper and Lower bounds). We have
$$
\sum_{n \leqslant x} \rho_0(n) \geqslant (1+o(1))\frac{C_0(\theta) x}{\log x} \quad \text{and} \quad \sum_{n \leqslant x} \rho_1(n) \leqslant (1+o(1))\frac{C_1(\theta) x}{\log x}
$$
for two functions $C_0(\theta)$ and $C_1(\theta)$ satisfy $0 < C_0(\theta) \leqslant 1$ and $C_1(\theta) \geqslant 1$.

(Distributions in Arithmetic Progressions). For any $a \in \mathbb{Z} \backslash \{0\}$ and any $A>0$, we have
$$
\sum_{\substack{q \sim Q \\ (q, a) = 1}} \left| \sum_{\substack{n \leqslant x \\ n \equiv a (\bmod q)}} \rho_j(n) - \frac{1}{\varphi(q)} \sum_{\substack{n \leqslant x \\ (n, q) = 1}} \rho_j(n) \right| \ll \frac{x}{(\log x)^A}
$$
for $j = 0, 1$.
\end{theorem}

The ``asymptotic formula'' (9) yields the following asymptotic formula for almost all $q \sim Q$:
\begin{equation}
S\left(\mathcal{A}^{q}, z \right) = (1+o(1)) \frac{1}{\varphi(q)} S\left(\mathcal{B}^{q}, z \right).
\end{equation}
By a similar decomposing process, one can deduce the following result using asymptotic formulas of the form (10):
\begin{theorem}\label{t22}
For almost all $q \leqslant Q$, we have
$$
\frac{C_0^{*}(\theta) x}{\varphi(q) \log x} \leqslant \pi(x; q, a) \leqslant \frac{C_1^{*}(\theta) x}{\varphi(q) \log x}
$$
for two functions $C_0^{*}(\theta)$ and $C_1^{*}(\theta)$ satisfy $0 < C_0(\theta) \leqslant C_0^{*}(\theta) \leqslant 1$ and $1 \leqslant C_1^{*}(\theta) \leqslant C_1(\theta)$.
\end{theorem}

In order to give asymptotic formulas (9) and (10) for sieve functions $S\left(\mathcal{A}^{q}, z \right)$, we need results of the form
\begin{equation}
\sum_{\substack{q \sim Q \\ (q, a) = 1}} \left| \sum_{\substack{n \sim x \\ n \equiv a (\bmod q)}} f(n) - \frac{1}{\varphi(q)} \sum_{\substack{n \sim x \\ (n, q) = 1}} f(n) \right| \ll \frac{x}{(\log x)^A}.
\end{equation}

There are two conditions that we may want the coefficients to satisfy. We shall use a divisor-bounded coefficient sequence $\lambda_l$ as an example. The first one is the Siegel--Walfisz condition, which demonstrate that at least one of the coefficient sequences is well-distributed in arithmetic progressions having small moduli. This condition is necessary in the dispersion estimates.

(\textbf{Condition A}: Siegel--Walfisz condition) For any $f \geqslant 1$, $k \geqslant 1$, $b \neq 0$ and $(k, b) = 1$, we have
$$
\sum_{\substack{l \sim L \\ l \equiv b (\mod k) \\ (l, f) = 1 }} \lambda_l = \frac{1}{\varphi(k)} \sum_{\substack{l \sim L \\ (l, f k) = 1 }} \lambda_l + O \left(\frac{L (d(f))^B }{(\log L)^{A}} \right).
$$
We note that $\lambda_l$ certainly satisfies the Siegel--Walfisz condition if $\lambda_l = 1$, if $\lambda_l = \mu(n)$, or if
$$
\lambda_l = \sum_{\substack{p_1 \cdots p_j = l \\ p_j \sim P_j}} 1
$$
by the Siegel--Walfisz theorem.

The next condition ensures that $\lambda_l$ is supported on ``almost-primes'': integers with all prime factors larger than $\exp\left(\log x (\log \log x)^{-3} \right)$.

(\textbf{Condition B}: No small prime factors) We have $\lambda_l = 0$ whenever $l$ has a prime factor smaller than $\exp\left(\log x (\log \log x)^{-3} \right)$.

\subsection{Preliminary Lemmas}
Before constructing the majorant and minorant, we need estimate results of the form (11). Note that many of them are still useful in the later sections.

\subsubsection{Type-II estimate}
The first lemma comes from \cite{FouvryA1}, and it served as one of the most important Type-II information inputs in previous works \cite{676}, \cite{677}, \cite{LRB679} and \cite{Mikawa}.
\begin{lemma}\label{l23} ([\cite{FouvryA1}, Théorème 1]).
Let $M_1 M_2 \asymp x$ and $M_2 \geqslant x^{\varepsilon}$. Let $a_{1, m_1}$ and $a_{2, m_2}$ be divisor-bounded complex sequences. Suppose that $a_{2, m_2}$ satisfies \textbf{Condition A}. If we have
$$
Q^2 x^{- 1 + \varepsilon} \leqslant M_2 \leqslant Q^{-\frac{4}{3}} x^{\frac{5}{6} - \varepsilon},
$$
then
$$
\sum_{\substack{q \sim Q \\ (q, a) = 1}} \left| \sum_{\substack{m_1 \sim M_1 \\ m_2 \sim M_2 \\ m_1 m_2 \equiv a (\bmod q)}} a_{1, m_1} a_{2, m_2} - \frac{1}{\varphi(q)} \sum_{\substack{m_1 \sim M_1 \\ m_2 \sim M_2 \\ (m_1 m_2, q) = 1}} a_{1, m_1} a_{2, m_2} \right| \ll \frac{x}{(\log x)^A}.
$$
Note that this lemma is nontrivial when $Q < x^{\frac{11}{20}}$.
\end{lemma}

The second lemma comes from [\cite{Wright26}, Corollary 2.2], and it will play a vital role in our final decomposition.
\begin{lemma}\label{l24} ([\cite{Wright26}, Corollary 2.2(i)]).
Let $M_1 M_2 \asymp x$. Let $a_{1, m_1}$ and $a_{2, m_2}$ be divisor-bounded complex sequences. Suppose that $a_{2, m_2}$ satisfies \textbf{Condition A}. If we have
$$
\exp\left((\log x)^{\varepsilon}\right) \leqslant M_2 \leqslant Q^{-\frac{33}{28}} x^{\frac{17}{28} - \varepsilon},
$$
then
$$
\sum_{\substack{q \sim Q \\ (q, a) = 1}} \left| \sum_{\substack{m_1 \sim M_1 \\ m_2 \sim M_2 \\ m_1 m_2 \equiv a (\bmod q)}} a_{1, m_1} a_{2, m_2} - \frac{1}{\varphi(q)} \sum_{\substack{m_1 \sim M_1 \\ m_2 \sim M_2 \\ (m_1 m_2, q) = 1}} a_{1, m_1} a_{2, m_2} \right| \ll \frac{x}{(\log x)^A}.
$$
\end{lemma}

Many other estimate results, such as [\cite{FouvryA2}, Corollaire 1], [\cite{FouvryRadziwill2}, Theorem 1.1], [\cite{Fouvry1984}, Lemme 3] and [\cite{FouvryRadziwill1}, Corollary 1.1], are also applicable in this problem. However, all of them can be deduced by Lemma~\ref{l23} or Lemma~\ref{l24} when $a$ is a fixed nonzero integer. Combining Lemma~\ref{l23} and Lemma~\ref{l24}, we can deduce [\cite{Wright26}, Corollary 2.2(ii)(iii)] for a fixed nonzero integer $a$.

\subsubsection{Type-II$_3$ estimate}
Most of the next 5 lemmas were used in previous works \cite{676}, \cite{677} and \cite{LRB679}, and we still need them in this section and later sections. Note that Lemma~\ref{l28} gives the main Type-I information in this and later (except for the last two) sections.
\begin{lemma}\label{l25} ([\cite{BFI2}, Theorem 3]).
Let $M_1 M_2 M_3 \asymp x$ and $\min(M_1, M_2, M_3) > x^{\varepsilon}$. Let $a_{1, m_1}$, $a_{2, m_2}$ and $a_{3, m_3}$ be divisor-bounded complex sequences. Suppose that $a_{1, m_1}$, $a_{2, m_2}$ and $a_{3, m_3}$ satisfy \textbf{Condition B}, and $a_{2, m_2}$ also satisfies \textbf{Condition A}. If we have
$$
Q x^{\varepsilon} < M_1 M_2, \quad M_1^2 M_2^3 < Q x^{1 - \varepsilon}, \quad M_1^5 M_2^2 < x^{2 - \varepsilon}, \quad M_1^4 M_2^3 < x^{2 - \varepsilon},
$$
then
$$
\sum_{\substack{q \sim Q \\ (q, a) = 1}} \left| \sum_{\substack{m_1 \sim M_1 \\ m_2 \sim M_2 \\ m_3 \sim M_3 \\ m_1 m_2 m_3 \equiv a (\bmod q)}} a_{1, m_1} a_{2, m_2} a_{3, m_3} - \frac{1}{\varphi(q)} \sum_{\substack{m_1 \sim M_1 \\ m_2 \sim M_2 \\ m_3 \sim M_3 \\ (m_1 m_2 m_3, q) = 1}} a_{1, m_1} a_{2, m_2} a_{3, m_3} \right| \ll \frac{x}{(\log x)^A}.
$$
\end{lemma}
\begin{lemma}\label{l26} ([\cite{676}, Lemma 5]).
Let $M_1 M_2 M_3 \asymp x$ and $\min(M_1, M_2, M_3) > x^{\varepsilon}$. Let $a_{1, m_1}$, $a_{2, m_2}$ and $a_{3, m_3}$ be divisor-bounded complex sequences. Suppose that $a_{1, m_1}$, $a_{2, m_2}$ and $a_{3, m_3}$ satisfy \textbf{Condition B}, and $a_{2, m_2}$ also satisfies \textbf{Condition A}. If we have
$$
Q x^{\varepsilon} < M_1 M_2, \quad M_1 M_2^2 Q^2 < x^{2 - 3 \varepsilon}, \quad M_1^5 M_2^2 < x^{2 - 3 \varepsilon},
$$
then
$$
\sum_{\substack{q \sim Q \\ (q, a) = 1}} \left| \sum_{\substack{m_1 \sim M_1 \\ m_2 \sim M_2 \\ m_3 \sim M_3 \\ m_1 m_2 m_3 \equiv a (\bmod q)}} a_{1, m_1} a_{2, m_2} a_{3, m_3} - \frac{1}{\varphi(q)} \sum_{\substack{m_1 \sim M_1 \\ m_2 \sim M_2 \\ m_3 \sim M_3 \\ (m_1 m_2 m_3, q) = 1}} a_{1, m_1} a_{2, m_2} a_{3, m_3} \right| \ll \frac{x}{(\log x)^A}.
$$
\end{lemma}
\begin{lemma}\label{l27} ([\cite{MaynardLargeModuliI}, Proposition 8.3]).
Let $M_1 M_2 M_3 \asymp x$ and $\min(M_1, M_2, M_3) > x^{\varepsilon}$. Let $a_{1, m_1}$, $a_{2, m_2}$ and $a_{3, m_3}$ be divisor-bounded complex sequences. Suppose that $a_{1, m_1}$, $a_{2, m_2}$ and $a_{3, m_3}$ satisfy \textbf{Condition B}, and $a_{2, m_2}$ also satisfies \textbf{Condition A}. If we have
$$
Q < x^{0.7 - \varepsilon}, \quad Q x^{\varepsilon} < M_1 M_2, \quad M_2 < Q^{-1} x^{1 - 2 \varepsilon}, \quad M_1 M_2 < Q^{-\frac{1}{7}} x^{\frac{153}{224} - 10 \varepsilon}, \quad M_1^4 M_2 < Q^{-1} x^{\frac{57}{32} - 10 \varepsilon},
$$
then
$$
\sum_{\substack{q \sim Q \\ (q, a) = 1}} \left| \sum_{\substack{m_1 \sim M_1 \\ m_2 \sim M_2 \\ m_3 \sim M_3 \\ m_1 m_2 m_3 \equiv a (\bmod q)}} a_{1, m_1} a_{2, m_2} a_{3, m_3} - \frac{1}{\varphi(q)} \sum_{\substack{m_1 \sim M_1 \\ m_2 \sim M_2 \\ m_3 \sim M_3 \\ (m_1 m_2 m_3, q) = 1}} a_{1, m_1} a_{2, m_2} a_{3, m_3} \right| \ll \frac{x}{(\log x)^A}.
$$
\end{lemma}
\begin{lemma}\label{l28} ([\cite{MaynardLargeModuliI}, Proposition 8.6]).
Let $M_1 M_2 M_3 \asymp x$ and $\min(M_1, M_2, M_3) > x^{\varepsilon}$. Let $a_{1, m_1}$, $a_{2, m_2}$ and $a_{3, m_3}$ be divisor-bounded complex sequences. Suppose that $a_{1, m_1}$, $a_{2, m_2}$ and $a_{3, m_3}$ satisfy \textbf{Condition B}, and $a_{2, m_2}$ also satisfies \textbf{Condition A}. If we have
$$
Q < x^{\frac{127}{224} - \varepsilon}, \quad Q x^{\varepsilon} < M_1 M_2 \leqslant x^{\frac{4}{7} - \varepsilon}, \quad M_1 \leqslant M_2 \leqslant (M_1 M_2)^{\frac{3}{4}},
$$
then
$$
\sum_{\substack{q \sim Q \\ (q, a) = 1}} \left| \sum_{\substack{m_1 \sim M_1 \\ m_2 \sim M_2 \\ m_3 \sim M_3 \\ m_1 m_2 m_3 \equiv a (\bmod q)}} a_{1, m_1} a_{2, m_2} a_{3, m_3} - \frac{1}{\varphi(q)} \sum_{\substack{m_1 \sim M_1 \\ m_2 \sim M_2 \\ m_3 \sim M_3 \\ (m_1 m_2 m_3, q) = 1}} a_{1, m_1} a_{2, m_2} a_{3, m_3} \right| \ll \frac{x}{(\log x)^A}.
$$
Note that this lemma can be deduced from Lemma~\ref{l25} and Lemma~\ref{l27}.
\end{lemma}

\subsubsection{Type-I/II estimate}
\begin{lemma}\label{l29} ([\cite{BFI2}, Theorems 5 and 5$^{*}$]).
Let $M_1 M_2 M_3 \asymp x$ and $z \ll \exp\left(\log x (\log \log x)^{-1} \right)$. Let $a_{1, m_1}$, $a_{2, m_2}$ and $a_{3, m_3}$ be divisor-bounded complex sequences. Suppose that
$$
a_{2, m_2} = \mathbbm{1}_{m_2 \in \mathbf{M}} \qquad \text{or} \qquad a_{2, m_2} = \mathbbm{1}_{\substack{m_2 \in \mathbf{M} \\ \left(m_2, P(z)\right) = 1}}
$$
for some interval $\mathbf{M} \subseteq [M_2, 2 M_2]$. If we have
$$
M_3 Q < x^{1 - \varepsilon}, \quad M_1^4 M_3 Q < x^{2 - \varepsilon}, \quad M_1^2 M_3 Q^2 < x^{2 - \varepsilon},
$$
then
$$
\sum_{\substack{q \sim Q \\ (q, a) = 1}} \left| \sum_{\substack{m_1 \sim M_1 \\ m_2 \in \mathbf{M} \\ m_3 \sim M_3 \\ m_1 m_2 m_3 \equiv a (\bmod q)}} a_{1, m_1} a_{2, m_2} a_{3, m_3} - \frac{1}{\varphi(q)} \sum_{\substack{m_1 \sim M_1 \\ m_2 \in \mathbf{M} \\ m_3 \sim M_3 \\ (m_1 m_2 m_3, q) = 1}} a_{1, m_1} a_{2, m_2} a_{3, m_3} \right| \ll \frac{x}{(\log x)^A}.
$$
\end{lemma}

\subsubsection{Another Type-II estimate}
The next estimate is a new type of Type-II estimate for convolutions, and it will be useful in proving a better asymptotic formula of the form (9) than the corresponding result used in \cite{676}, \cite{677} and \cite{LRB679}. The proof of this lemma requires Lemma~\ref{l28} and Lemma~\ref{l29}.
\begin{lemma}\label{l210} ([\cite{MaynardLargeModuliI}, Lemma 8.12]).
Let $Q \leqslant x^{\frac{127}{224} - \varepsilon}$ and $M_1, M_2, \ldots, M_r \geqslant 1$ be such that $\prod_{1 \leqslant i \leqslant r} M_i \asymp x$. Let $a_{j, m_j}$ ($1 \leqslant j \leqslant r$) be divisor-bounded complex sequences. Suppose that $a_{j, m_j}$ ($1 \leqslant j \leqslant r$) satisfy \textbf{Conditions A and B} and
$$
a_{j, m_j} = \mathbbm{1}_{\left(m_j, P(z)\right) = 1}, \qquad z = \exp\left(\log x (\log \log x)^{-3} \right)
$$
for all $m_j > x^{\frac{1}{15}}$. Let $\mathbf{M}_i$ ($1 \leqslant i \leqslant r$) be intervals such that $\mathbf{M}_i \subseteq [M_i, 2 M_i]$. If we have
$$
x^{\frac{3}{7} + \varepsilon} < \prod_{j \in \mathcal{J}} M_j < Q^{-1} x^{1 - \varepsilon}
$$
for some set $\mathcal{J} \subseteq \{1, \ldots, r\}$, then
$$
\sum_{\substack{q \sim Q \\ (q, a) = 1}} \left| \sum_{\substack{m_i \in \mathbf{M}_i \\ 1 \leqslant i \leqslant r}} \left( \prod_{1 \leqslant j \leqslant r} a_{j, m_j} \right) \left( \mathbbm{1}_{m_1 \cdots m_r \equiv a (\bmod q)} - \frac{\mathbbm{1}_{(m_1 \cdots m_r, q) = 1}}{\varphi(q)} \right) \right| \ll \frac{x}{(\log x)^A}.
$$
\end{lemma}

\subsection{Sieve Asymptotic Formulas}
In this subsection we give asymptotic formulas (9) for sums of sieve functions $S\left(\mathcal{A}^{q}_{p_1 \cdots p_n}, p_n\right)$ and $S\left(\mathcal{A}^{q}_{p_1 \cdots p_n}, x^{\kappa_0}\right)$ with $\kappa_0 = \kappa$ or $\kappa^{\prime}$ or other values, where
\begin{align*}
\kappa = \kappa(\theta) =
\begin{cases}
\frac{5-8\theta}{6} - \varepsilon, & \theta \leqslant \frac{17}{32} - \varepsilon, \\
\frac{5-8\theta}{12} - 3 \varepsilon, & \frac{17}{32} - \varepsilon < \theta \leqslant \frac{7}{13} - \varepsilon, \\
\frac{3-5\theta}{7} - 2 \varepsilon, & \frac{7}{13} - \varepsilon < \theta \leqslant \frac{4}{7} - \varepsilon, \\
\end{cases}
\end{align*}
and
\begin{align*}
\kappa^{\prime} = \kappa^{\prime}(\theta) =
\begin{cases}
\frac{11-20\theta}{6} - 2 \varepsilon, & \frac{7}{13} - \varepsilon < \theta \leqslant \frac{11}{20} - \varepsilon, \\
\kappa, & \text{otherwise}. \\
\end{cases}
\end{align*}
We also write
\begin{align*}
\tau = \tau(\theta) =
\begin{cases}
\frac{3(1-\theta)}{5} - \varepsilon, & \theta \leqslant \frac{11}{21}, \\
\frac{2}{7} - \varepsilon, & \frac{11}{21} < \theta \leqslant \frac{6}{11} - \varepsilon, \\
\frac{5-6\theta}{7} - \varepsilon, & \frac{6}{11} - \varepsilon < \theta, \\
\end{cases}
\end{align*}
and
\begin{align*}
\tau^{\prime} = \tau^{\prime}(\theta) =
\begin{cases}
\frac{5-6\theta}{7}, & \frac{7}{13} - \varepsilon < \theta \leqslant \frac{11}{20} - \varepsilon, \\
\tau, & \text{otherwise}. \\
\end{cases}
\end{align*}

\begin{lemma}\label{l211}
Let $\frac{1}{2} \leqslant \theta \leqslant \frac{4}{7} - \varepsilon$. Define
\begin{align}
\nonumber \boldsymbol{g}_{1} =&\ \boldsymbol{g}_{1}(\theta) = \left\{(s, t): 2 \theta - 1 + \varepsilon \leqslant s \leqslant \frac{5-8\theta}{6} - \varepsilon \right\}, \\
\nonumber \boldsymbol{g}_{2} =&\ \boldsymbol{g}_{2}(\theta) = \left\{(s, t): s \leqslant \frac{17-33\theta}{28} - \varepsilon \right\}, \\
\nonumber \boldsymbol{g}_{3} =&\ \boldsymbol{g}_{3}(\theta) = \left\{(s, t): s + t \geqslant \theta + \varepsilon ,\ 2 s + 3 t \leqslant 1 + \theta - \varepsilon,\ 5 s + 2 t \leqslant 2 - \varepsilon,\ 4 s + 3 t \leqslant 2 - \varepsilon \right\}, \\
\nonumber \boldsymbol{g}_{4} =&\ \boldsymbol{g}_{4}(\theta) = \left\{(s, t): s + t \geqslant \theta + \varepsilon,\ s + 2 t \leqslant 2 - 2 \theta - \varepsilon,\ 5 s + 2 t \leqslant 2 - \varepsilon \right\}, \\
\nonumber \boldsymbol{g}_{5} =&\ \boldsymbol{g}_{5}(\theta) = \left\{(s, t): s + t \geqslant \theta + \varepsilon,\ t \leqslant 1 - \theta - \varepsilon,\ s + t \leqslant \frac{153}{224} - \frac{1}{7} \theta - \varepsilon,\ 4 s + t \leqslant \frac{57}{32} - \theta - \varepsilon \right\}, \\
\nonumber \boldsymbol{\mathcal{G}}_{j} =&\ \boldsymbol{\mathcal{G}}_{j}(\theta) = \left\{\boldsymbol{\alpha}_{j}: \boldsymbol{\alpha}_{j} \text{ partitions into } \boldsymbol{g}_{1} \cup \boldsymbol{g}_{2} \cup \boldsymbol{g}_{3} \cup \boldsymbol{g}_{4} \cup \boldsymbol{g}_{5} \right\}.
\end{align}
Suppose that $\min \boldsymbol{\alpha}_{j} \geqslant (\log \log x)^{-3}$. Then
$$
\sum_{\boldsymbol{\alpha}_{j} \in \boldsymbol{\mathcal{G}}_{j}} S\left(\mathcal{A}^{q}_{p_1 \cdots p_j}, p_j\right)
$$
has an asymptotic formula of the form (9).
\end{lemma}
\begin{proof}
This lemma follows easily from Lemmas~\ref{l23}--\ref{l27}. Note that this lemma also holds for more general Type-II sums.
\end{proof}

\begin{lemma}\label{l212}
Let $\frac{1}{2} \leqslant \theta \leqslant \frac{17}{32} - \varepsilon$ and $M < x^{1 - \varepsilon} Q^{-1}$. Then, for a divisor-bounded sequence $a_m$,
$$
\sum_{m \sim M} a_m S\left(\mathcal{A}^{q}_{m}, x^{\kappa}\right)
$$
has an asymptotic formula of the form (9).

Let $\frac{1}{2} \leqslant \theta < \frac{11}{20} - \varepsilon$ and $M < x^{1 - \varepsilon} Q^{-1}$. Then, for a divisor-bounded sequence $a_m$,
$$
\sum_{m \sim M} a_m S\left(\mathcal{A}^{q}_{m}, x^{\frac{11 - 20 \theta}{6} - \varepsilon}\right)
$$
has an asymptotic formula of the form (9).
\end{lemma}
\begin{proof}
Using Buchstab's identity many times (or a Möbius inversion), we can replace the sum
$$
\sum_{m \sim M} a_m S\left(\mathcal{A}^{q}_{m}, x^{\frac{11 - 20 \theta}{6} - \varepsilon}\right)
$$
with $\ll \log x (\log \log x)^{-1}$ sums of the form
$$
\sum_{m \sim M} a_m \sum_{\substack{p_1 p_2 \cdots p_k m n \in \mathcal{A}^{q} \\ p_k < \ldots < p_1 < x^{\frac{11 - 20 \theta}{6} - \varepsilon} }} 1
$$
together with an outer summation $\sum_{k \geqslant 0}(-1)^k$. We write $p_1 p_2 \cdots p_k = d \asymp D$. When $D > x^{2 \theta - 1}$, we know that there must be a product of some prime variables that lies in $\left[2 \theta - 1 + \varepsilon, \frac{5 - 8 \theta}{6} - \varepsilon \right]$ since $p_k < \ldots < p_1 < x^{\frac{11 - 20 \theta}{6} - \varepsilon}$. We split the prime variables, starting from the largest $p_1$, and group the remaining ones together to form a Type-II sum after removing cross conditions.

When $D \leqslant x^{2 \theta - 1}$, we can use Lemma~\ref{l29} with $M_1 = D$ and $M_3 = M$. The conditions $D \leqslant x^{2 \theta - 1}$ and $M < x^{1 - \theta - \varepsilon}$ ensure the required conditions in Lemma~\ref{l29}: we have $M_3 < x^{1 - \theta - \varepsilon}$, $M_1^4 M_3 < x^{7 \theta - 3 - \varepsilon} < x^{2 - \theta - \varepsilon}$ (since $\theta < \frac{5}{8}$) and $M_1^2 M_3 < x^{3 \theta - 1 - \varepsilon} < x^{2 - 2 \theta - \varepsilon}$ (since $\theta < \frac{3}{5}$).

Now we can give an asymptotic formula of the form (9) for the whole sum
$$
\sum_{m \sim M} a_m S\left(\mathcal{A}^{q}_{m}, x^{\frac{11 - 20 \theta}{6} - \varepsilon}\right).
$$
The proof of an asymptotic formula of the form (9) for the sum
$$
\sum_{m \sim M} a_m S\left(\mathcal{A}^{q}_{m}, x^{\kappa}\right)
$$
can thus be done by applying Buchstab's identity once as in [\cite{HarmanBOOK}, Lemma 8.15]. Note that the two ranges $\left[\frac{1}{2},\ \frac{11}{20} - \varepsilon \right)$ and $\left[\frac{1}{2},\ \frac{17}{32} - \varepsilon \right]$ cannot be enlarged because of the restriction on the width of the Type-II range $\left[2 \theta - 1 + \varepsilon, \frac{5 - 8 \theta}{6} - \varepsilon \right]$ here.
\end{proof}
\begin{remark*}
When $\theta < \frac{4}{7} - 3 \varepsilon$, we have $x^{\frac{3}{7} + \varepsilon} < x^{1 - \varepsilon} Q^{-1}$. In [\cite{HarmanBOOK}, Lemmas 8.10 and 8.15], this lemma is stated with a condition $M < x^{2 - \varepsilon} Q^{-3}$ instead of $M < x^{1 - \varepsilon} Q^{-1}$.
\end{remark*}

\begin{lemma}\label{l213}
Let $\frac{1}{2} \leqslant \theta \leqslant \frac{4}{7} - \varepsilon$. Define
\begin{align}
\nonumber \boldsymbol{S} =&\ \boldsymbol{S}(\theta) = \left\{(s, t): s \leqslant 1 - \theta - \varepsilon,\ s + 2t \leqslant 2 - 2 \theta - \varepsilon,\ s + 4t \leqslant 2 - \theta - \varepsilon \right\}, \\
\nonumber \boldsymbol{S}_{j} =&\ \boldsymbol{S}_{j}(\theta) = \left\{\boldsymbol{\alpha}_{j}: \boldsymbol{\alpha}_{j} \text{ partitions exactly into } \boldsymbol{S} \right\}, \\
\nonumber \boldsymbol{A}_{j} =&\ \boldsymbol{A}_{j}(\theta) = \left\{\boldsymbol{\alpha}_{j}: (\log \log x)^{-3} \leqslant \alpha_j < \cdots < \alpha_1 < \tau,\ \alpha_1 + \cdots + \alpha_j \leqslant 1 \right\}, \\
\nonumber \boldsymbol{T}^{*} =&\ \boldsymbol{T}^{*}(\theta) = \left\{(s, t): 0 \leqslant s \leqslant \frac{8\theta-2}{7},\ 0 \leqslant t \leqslant \frac{5-6\theta}{7} \right\}, \\
\nonumber \boldsymbol{T}_{j}^{*} =&\ \boldsymbol{T}_{j}^{*}(\theta) = \left\{\boldsymbol{\alpha}_{j}: \boldsymbol{\alpha}_{j} \text{ partitions exactly into } \boldsymbol{T}^{*} \right\}, \\
\nonumber \boldsymbol{U}_{j}^{\prime} =&\ \boldsymbol{U}_{j}^{\prime}(\theta) = \left\{\boldsymbol{\alpha}_{j}: \boldsymbol{\alpha}_{j} \in \boldsymbol{A}_{j},\ (\alpha_1, \ldots, \alpha_j, 2\theta-1+\varepsilon) \in \boldsymbol{S}_{j+1} \right\}, \\
\nonumber \boldsymbol{U}_{j} =&\ \boldsymbol{U}_{j}(\theta) =
\begin{cases}
\boldsymbol{U}_{j}^{\prime}(\theta), & \theta < \frac{7}{13}, \\
\left\{\boldsymbol{\alpha}_{j}: \boldsymbol{\alpha}_{j} \in \boldsymbol{A}_{j},\ \boldsymbol{\alpha}_{j} \in \boldsymbol{T}_{j}^{*} \right\}, & \theta \geqslant \frac{7}{13}.
\end{cases}
\end{align}

Let $\frac{1}{2} \leqslant \theta \leqslant \frac{4}{7} - \varepsilon$. Then
$$
\sum_{\boldsymbol{\alpha}_{j} \in \boldsymbol{U}_{j}} S\left(\mathcal{A}^{q}_{p_1 \cdots p_j}, x^{\kappa}\right)
$$
has an asymptotic formula of the form (9).

Let $\frac{1}{2} \leqslant \theta < \frac{11}{20} - \varepsilon$. Then
$$
\sum_{\boldsymbol{\alpha}_{j} \in \boldsymbol{U}_{j}^{\prime}} S\left(\mathcal{A}^{q}_{p_1 \cdots p_j}, x^{\kappa^{\prime}}\right)
$$
has an asymptotic formula of the form (9).
\end{lemma}
\begin{proof}
The proof follows by similar arguments as in [\cite{676}, Lemmas 15 and 16].

We first prove the second case of Lemma~\ref{l213}, which also implies the first case when $\theta \leqslant \frac{7}{13} - \varepsilon$ by the definitions of $\nonumber \boldsymbol{U}_{j}$ and $\nonumber \boldsymbol{U}_{j}^{\prime}$. Let $\frac{1}{2} \leqslant \theta < \frac{11}{20} - \varepsilon$. Using Buchstab's identity repeatedly, we have
\begin{align}
\nonumber \sum_{\boldsymbol{\alpha}_{j} \in \boldsymbol{U}_{j}^{\prime}} S\left(\mathcal{A}^{q}_{p_1 \cdots p_j}, x^{\kappa^{\prime}}\right) =&\ \sum_{\boldsymbol{\alpha}_{j} \in \boldsymbol{U}_{j}^{\prime}} S\left(\mathcal{A}^{q}_{p_1 \cdots p_j}, \exp\left(\log x (\log \log x)^{-3} \right)\right) \\
\nonumber & - \sum_{\substack{\boldsymbol{\alpha}_{j} \in \boldsymbol{U}_{j}^{\prime} \\ \boldsymbol{\alpha}_{j+1} \in \boldsymbol{A}_{j+1} \\ 2 \theta - 1 + \varepsilon \leqslant \alpha_{j+1} < \kappa^{\prime} }} S\left(\mathcal{A}^{q}_{p_1 \cdots p_{j+1}}, p_{j+1}\right) \\
\nonumber & - \sum_{\substack{\boldsymbol{\alpha}_{j} \in \boldsymbol{U}_{j}^{\prime} \\ \boldsymbol{\alpha}_{j+1} \in \boldsymbol{A}_{j+1} \\ \alpha_{j+1} < 2 \theta - 1 + \varepsilon }} S\left(\mathcal{A}^{q}_{p_1 \cdots p_{j+1}}, p_{j+1}\right) \\
\nonumber =&\ \sum_{\boldsymbol{\alpha}_{j} \in \boldsymbol{U}_{j}^{\prime}} S\left(\mathcal{A}^{q}_{p_1 \cdots p_j}, \exp\left(\log x (\log \log x)^{-3} \right)\right) \\
\nonumber & - \sum_{\substack{\boldsymbol{\alpha}_{j} \in \boldsymbol{U}_{j}^{\prime} \\ \boldsymbol{\alpha}_{j+1} \in \boldsymbol{A}_{j+1} \\ 2 \theta - 1 + \varepsilon \leqslant \alpha_{j+1} < \kappa^{\prime} }} S\left(\mathcal{A}^{q}_{p_1 \cdots p_{j+1}}, p_{j+1}\right) \\
\nonumber & - \sum_{\substack{\boldsymbol{\alpha}_{j} \in \boldsymbol{U}_{j}^{\prime} \\ \boldsymbol{\alpha}_{j+1} \in \boldsymbol{A}_{j+1} \\ \alpha_{j+1} < 2 \theta - 1 + \varepsilon }} S\left(\mathcal{A}^{q}_{p_1 \cdots p_{j+1}}, \exp\left(\log x (\log \log x)^{-3} \right)\right) \\
\nonumber & + \sum_{\substack{\boldsymbol{\alpha}_{j} \in \boldsymbol{U}_{j}^{\prime} \\ \boldsymbol{\alpha}_{j+2} \in \boldsymbol{A}_{j+2} \\ \alpha_{j+1} < \kappa^{\prime} \\ \alpha_{j+1} < 2 \theta - 1 + \varepsilon \leqslant \alpha_{j+1} + \alpha_{j+2} }} S\left(\mathcal{A}^{q}_{p_1 \cdots p_{j+2}}, p_{j+2}\right) \\
\nonumber & + \sum_{\substack{\boldsymbol{\alpha}_{j} \in \boldsymbol{U}_{j}^{\prime} \\ \boldsymbol{\alpha}_{j+2} \in \boldsymbol{A}_{j+2} \\ \alpha_{j+1} < \kappa^{\prime} \\ \alpha_{j+1} + \alpha_{j+2} < 2 \theta - 1 + \varepsilon }} S\left(\mathcal{A}^{q}_{p_1 \cdots p_{j+2}}, p_{j+2}\right) \\
\nonumber =&\ \sum_{\boldsymbol{\alpha}_{j} \in \boldsymbol{U}_{j}^{\prime}} S\left(\mathcal{A}^{q}_{p_1 \cdots p_j}, \exp\left(\log x (\log \log x)^{-3} \right)\right) \\
\nonumber & - \sum_{\substack{\boldsymbol{\alpha}_{j} \in \boldsymbol{U}_{j}^{\prime} \\ \boldsymbol{\alpha}_{j+1} \in \boldsymbol{A}_{j+1} \\ 2 \theta - 1 + \varepsilon \leqslant \alpha_{j+1} < \kappa^{\prime} }} S\left(\mathcal{A}^{q}_{p_1 \cdots p_{j+1}}, p_{j+1}\right) \\
\nonumber & - \sum_{\substack{\boldsymbol{\alpha}_{j} \in \boldsymbol{U}_{j}^{\prime} \\ \boldsymbol{\alpha}_{j+1} \in \boldsymbol{A}_{j+1} \\ \alpha_{j+1} < 2 \theta - 1 + \varepsilon }} S\left(\mathcal{A}^{q}_{p_1 \cdots p_{j+1}}, \exp\left(\log x (\log \log x)^{-3} \right)\right) \\
\nonumber & + \sum_{\substack{\boldsymbol{\alpha}_{j} \in \boldsymbol{U}_{j}^{\prime} \\ \boldsymbol{\alpha}_{j+2} \in \boldsymbol{A}_{j+2} \\ \alpha_{j+1} < \kappa^{\prime} \\ \alpha_{j+1} < 2 \theta - 1 + \varepsilon \leqslant \alpha_{j+1} + \alpha_{j+2} }} S\left(\mathcal{A}^{q}_{p_1 \cdots p_{j+2}}, p_{j+2}\right) \\
\nonumber & + \sum_{\substack{\boldsymbol{\alpha}_{j} \in \boldsymbol{U}_{j}^{\prime} \\ \boldsymbol{\alpha}_{j+2} \in \boldsymbol{A}_{j+2} \\ \alpha_{j+1} < \kappa^{\prime} \\ \alpha_{j+1} + \alpha_{j+2} < 2 \theta - 1 + \varepsilon }} S\left(\mathcal{A}^{q}_{p_1 \cdots p_{j+2}}, \exp\left(\log x (\log \log x)^{-3} \right)\right) \\
\nonumber & - \sum_{\substack{\boldsymbol{\alpha}_{j} \in \boldsymbol{U}_{j}^{\prime} \\ \boldsymbol{\alpha}_{j+3} \in \boldsymbol{A}_{j+3} \\ \alpha_{j+1} < \kappa^{\prime} \\ \alpha_{j+1} + \alpha_{j+2} < 2 \theta - 1 + \varepsilon \leqslant \alpha_{j+1} + \alpha_{j+2} + \alpha_{j+3} }} S\left(\mathcal{A}^{q}_{p_1 \cdots p_{j+3}}, p_{j+3}\right) - \cdots \\
\nonumber =&\ \sum_{\boldsymbol{\alpha}_{j} \in \boldsymbol{U}_{j}^{\prime}} S\left(\mathcal{A}^{q}_{p_1 \cdots p_j}, \exp\left(\log x (\log \log x)^{-3} \right)\right) \\
\nonumber & - \sum_{\substack{\boldsymbol{\alpha}_{j} \in \boldsymbol{U}_{j}^{\prime} \\ \boldsymbol{\alpha}_{j+1} \in \boldsymbol{A}_{j+1} \\ 2 \theta - 1 + \varepsilon \leqslant \alpha_{j+1} < \kappa^{\prime} }} S\left(\mathcal{A}^{q}_{p_1 \cdots p_{j+1}}, p_{j+1}\right) \\
\nonumber & + \sum_{k \geqslant j+1} (-1)^{k-j} \sum_{\substack{\boldsymbol{\alpha}_{j} \in \boldsymbol{U}_{j}^{\prime} \\ \boldsymbol{\alpha}_{k} \in \boldsymbol{A}_{k} \\ \alpha_{j+1} < \kappa^{\prime} \\ \alpha_{j+1} + \cdots + \alpha_{k} < 2 \theta - 1 + \varepsilon }} S\left(\mathcal{A}^{q}_{p_1 \cdots p_{k}}, \exp\left(\log x (\log \log x)^{-3} \right)\right) \\
\nonumber & + \sum_{k \geqslant j+2} (-1)^{k-j} \sum_{\substack{\boldsymbol{\alpha}_{j} \in \boldsymbol{U}_{j}^{\prime} \\ \boldsymbol{\alpha}_{k} \in \boldsymbol{A}_{k} \\ \alpha_{j+1} < \kappa^{\prime} \\ \alpha_{j+1} + \cdots + \alpha_{k-1} < 2 \theta - 1 + \varepsilon \leqslant \alpha_{j+1} + \cdots + \alpha_{k} }} S\left(\mathcal{A}^{q}_{p_1 \cdots p_{k}}, p_{k}\right) \\
\nonumber =&\ \Sigma_{21321} - \Sigma_{21322} + \Sigma_{21323} + \Sigma_{21324}.
\end{align}
We can give an asymptotic formula of the form (9) for $\Sigma_{21322}$ by Lemma~\ref{l211}. We can give asymptotic formulas of the form (9) for $\Sigma_{21321}$ and $\Sigma_{21323}$ by the condition $(\alpha_1, \ldots, \alpha_j, 2\theta-1+\varepsilon) \in \boldsymbol{S}_{j+1}$ and an application of Lemma~\ref{l29}. We can give an asymptotic formula of the form (9) for $\Sigma_{21324}$ by [\cite{676}, Lemma 11] and Lemma~\ref{l211}. Now the proof of the second case is completed.

Next, we prove the first case of Lemma~\ref{l213} when $\frac{7}{13} - \varepsilon < \theta \leqslant \frac{4}{7} - \varepsilon$. Similarly, we have
\begin{align}
\nonumber \sum_{\boldsymbol{\alpha}_{j} \in \boldsymbol{U}_{j}} S\left(\mathcal{A}^{q}_{p_1 \cdots p_j}, x^{\kappa}\right) =&\ \sum_{\boldsymbol{\alpha}_{j} \in \boldsymbol{U}_{j}} S\left(\mathcal{A}^{q}_{p_1 \cdots p_j}, \exp\left(\log x (\log \log x)^{-3} \right)\right) \\
\nonumber & - \sum_{\substack{\boldsymbol{\alpha}_{j} \in \boldsymbol{U}_{j} \\ \boldsymbol{\alpha}_{j+1} \in \boldsymbol{A}_{j+1} \\ \alpha_{j+1} < \kappa \\ \theta + \varepsilon \leqslant \alpha_{1} + \cdots + \alpha_{j+1} }} S\left(\mathcal{A}^{q}_{p_1 \cdots p_{j+1}}, p_{j+1}\right) \\
\nonumber & - \sum_{\substack{\boldsymbol{\alpha}_{j} \in \boldsymbol{U}_{j} \\ \boldsymbol{\alpha}_{j+1} \in \boldsymbol{A}_{j+1} \\ \alpha_{j+1} < \kappa \\ \alpha_{1} + \cdots + \alpha_{j+1} < \theta + \varepsilon }} S\left(\mathcal{A}^{q}_{p_1 \cdots p_{j+1}}, p_{j+1}\right) \\
\nonumber =&\ \sum_{\boldsymbol{\alpha}_{j} \in \boldsymbol{U}_{j}} S\left(\mathcal{A}^{q}_{p_1 \cdots p_j}, \exp\left(\log x (\log \log x)^{-3} \right)\right) \\
\nonumber & - \sum_{\substack{\boldsymbol{\alpha}_{j} \in \boldsymbol{U}_{j} \\ \boldsymbol{\alpha}_{j+1} \in \boldsymbol{A}_{j+1} \\ \alpha_{j+1} < \kappa \\ \theta + \varepsilon \leqslant \alpha_{1} + \cdots + \alpha_{j+1} }} S\left(\mathcal{A}^{q}_{p_1 \cdots p_{j+1}}, p_{j+1}\right) \\
\nonumber & - \sum_{\substack{\boldsymbol{\alpha}_{j} \in \boldsymbol{U}_{j} \\ \boldsymbol{\alpha}_{j+1} \in \boldsymbol{A}_{j+1} \\ \alpha_{j+1} < \kappa \\ \alpha_{1} + \cdots + \alpha_{j+1} < \theta + \varepsilon }} S\left(\mathcal{A}^{q}_{p_1 \cdots p_{j+1}}, \exp\left(\log x (\log \log x)^{-3} \right)\right) \\
\nonumber & + \sum_{\substack{\boldsymbol{\alpha}_{j} \in \boldsymbol{U}_{j} \\ \boldsymbol{\alpha}_{j+2} \in \boldsymbol{A}_{j+2} \\ \alpha_{j+1} < \kappa \\ \alpha_{1} + \cdots + \alpha_{j+1} < \theta + \varepsilon \leqslant \alpha_{1} + \cdots + \alpha_{j+1} + \alpha_{j+2} }} S\left(\mathcal{A}^{q}_{p_1 \cdots p_{j+2}}, p_{j+2}\right) \\
\nonumber & + \sum_{\substack{\boldsymbol{\alpha}_{j} \in \boldsymbol{U}_{j} \\ \boldsymbol{\alpha}_{j+2} \in \boldsymbol{A}_{j+2} \\ \alpha_{j+1} < \kappa \\ \alpha_{1} + \cdots + \alpha_{j+1} + \alpha_{j+2} < \theta + \varepsilon }} S\left(\mathcal{A}^{q}_{p_1 \cdots p_{j+2}}, p_{j+2}\right) \\
\nonumber =&\ \sum_{\boldsymbol{\alpha}_{j} \in \boldsymbol{U}_{j}} S\left(\mathcal{A}^{q}_{p_1 \cdots p_j}, \exp\left(\log x (\log \log x)^{-3} \right)\right) \\
\nonumber & - \sum_{\substack{\boldsymbol{\alpha}_{j} \in \boldsymbol{U}_{j} \\ \boldsymbol{\alpha}_{j+1} \in \boldsymbol{A}_{j+1} \\ \alpha_{j+1} < \kappa \\ \theta + \varepsilon \leqslant \alpha_{1} + \cdots + \alpha_{j+1} }} S\left(\mathcal{A}^{q}_{p_1 \cdots p_{j+1}}, p_{j+1}\right) \\
\nonumber & - \sum_{\substack{\boldsymbol{\alpha}_{j} \in \boldsymbol{U}_{j} \\ \boldsymbol{\alpha}_{j+1} \in \boldsymbol{A}_{j+1} \\ \alpha_{j+1} < \kappa \\ \alpha_{1} + \cdots + \alpha_{j+1} < \theta + \varepsilon }} S\left(\mathcal{A}^{q}_{p_1 \cdots p_{j+1}}, \exp\left(\log x (\log \log x)^{-3} \right)\right) \\
\nonumber & + \sum_{\substack{\boldsymbol{\alpha}_{j} \in \boldsymbol{U}_{j} \\ \boldsymbol{\alpha}_{j+2} \in \boldsymbol{A}_{j+2} \\ \alpha_{j+1} < \kappa \\ \alpha_{1} + \cdots + \alpha_{j+1} < \theta + \varepsilon \leqslant \alpha_{1} + \cdots + \alpha_{j+1} + \alpha_{j+2} }} S\left(\mathcal{A}^{q}_{p_1 \cdots p_{j+2}}, p_{j+2}\right) \\
\nonumber & + \sum_{\substack{\boldsymbol{\alpha}_{j} \in \boldsymbol{U}_{j} \\ \boldsymbol{\alpha}_{j+2} \in \boldsymbol{A}_{j+2} \\ \alpha_{j+1} < \kappa \\ \alpha_{1} + \cdots + \alpha_{j+1} + \alpha_{j+2} < \theta + \varepsilon }} S\left(\mathcal{A}^{q}_{p_1 \cdots p_{j+2}}, \exp\left(\log x (\log \log x)^{-3} \right)\right) \\
\nonumber & - \sum_{\substack{\boldsymbol{\alpha}_{j} \in \boldsymbol{U}_{j} \\ \boldsymbol{\alpha}_{j+3} \in \boldsymbol{A}_{j+3} \\ \alpha_{j+1} < \kappa \\ \alpha_{1} + \cdots + \alpha_{j+1} + \alpha_{j+2} < \theta + \varepsilon \leqslant \alpha_{1} + \cdots + \alpha_{j+1} + \alpha_{j+2} + \alpha_{j+3} }} S\left(\mathcal{A}^{q}_{p_1 \cdots p_{j+3}}, p_{j+3}\right) - \cdots \\
\nonumber =&\ \sum_{\boldsymbol{\alpha}_{j} \in \boldsymbol{U}_{j}} S\left(\mathcal{A}^{q}_{p_1 \cdots p_j}, \exp\left(\log x (\log \log x)^{-3} \right)\right) \\
\nonumber & - \sum_{\substack{\boldsymbol{\alpha}_{j} \in \boldsymbol{U}_{j} \\ \boldsymbol{\alpha}_{j+1} \in \boldsymbol{A}_{j+1} \\ \alpha_{j+1} < \kappa \\ \theta + \varepsilon \leqslant \alpha_{1} + \cdots + \alpha_{j+1} }} S\left(\mathcal{A}^{q}_{p_1 \cdots p_{j+1}}, p_{j+1}\right) \\
\nonumber & + \sum_{k \geqslant j+1} (-1)^{k-j} \sum_{\substack{\boldsymbol{\alpha}_{j} \in \boldsymbol{U}_{j} \\ \boldsymbol{\alpha}_{k} \in \boldsymbol{A}_{k} \\ \alpha_{j+1} < \kappa \\ \alpha_{1} + \cdots + \alpha_{k} < \theta + \varepsilon }} S\left(\mathcal{A}^{q}_{p_1 \cdots p_{k}}, \exp\left(\log x (\log \log x)^{-3} \right)\right) \\
\nonumber & + \sum_{k \geqslant j+2} (-1)^{k-j} \sum_{\substack{\boldsymbol{\alpha}_{j} \in \boldsymbol{U}_{j} \\ \boldsymbol{\alpha}_{k} \in \boldsymbol{A}_{k} \\ \alpha_{j+1} < \kappa \\ \alpha_{1} + \cdots + \alpha_{k-1} < \theta + \varepsilon \leqslant \alpha_{1} + \cdots + \alpha_{k} }} S\left(\mathcal{A}^{q}_{p_1 \cdots p_{k}}, p_{k}\right) \\
\nonumber =&\ \Sigma_{21311} - \Sigma_{21312} + \Sigma_{21313} + \Sigma_{21314}.
\end{align}
We can give asymptotic formulas of the form (9) for $\Sigma_{21311}$ and $\Sigma_{21323}$ by [\cite{676}, Lemma 12] and an application of Lemma~\ref{l29}. We can give asymptotic formulas of the form (9) for $\Sigma_{21312}$ and $\Sigma_{21324}$ by [\cite{676}, Lemma 13] and Lemma~\ref{l211}. Now the proof of the first case when $\frac{7}{13} - \varepsilon < \theta \leqslant \frac{4}{7} - \varepsilon$ is completed.

Combining the above 2 cases, the proof of Lemma~\ref{l213} is completed.
\end{proof}
\begin{remark*}
When $\theta < \frac{17}{32} - \varepsilon$, the proof above can be simplified a lot using [\cite{MaynardLargeModuliI}, Lemma 10.2]: we take $z_1 = x^{\kappa}$, $z_2 = \exp\left(\log x (\log \log x)^{-3} \right)$ and $y = x^{2 \theta - 1 + \varepsilon}$. Let $\Sigma_{1}$ and $\Sigma_{2}$ denote the sums correspond to the first and the second terms on the right hand side of [\cite{MaynardLargeModuliI}, Lemma 10.2]. We can give asymptotic formulas of the form (9) for $\Sigma_{1}$ by an application of Lemma~\ref{l29}, and $\Sigma_{2}$ is a Type-II sum. Now we have $d \leqslant y < z_1$, $p \leqslant z_1$ and $d p > y$ in $\Sigma_{2}$. Suppose that $y^2 < z_1$. If $d, p \notin (y,\ z_1)$, then $y < d p \leqslant y^2 < z_1$. Since $2 (2 \theta - 1 + \varepsilon) < \kappa$ when $\theta < \frac{17}{32} - \varepsilon$, we can give asymptotic formulas of the form (9) for $\Sigma_{2}$ by Lemma~\ref{l211}. The proof is thus completed.

From the above proof, one can easily find that we can remove the condition $\alpha_1 < \tau$ in the definition of $\boldsymbol{U}_{j}$ when $\theta < \frac{17}{32} - \varepsilon$. This removal is very important in our final decompositions, especially in the lower bound case. Note that when $\theta \geqslant \frac{17}{32}$, we need to discard sums with $p_1 \in \left(\tau,\ \frac{3}{7} + \varepsilon \right)$ in the upper bound case, and we usually cannot get a nontrivial lower bound in the lower bound case. Hence $\theta < \frac{17}{32} - \varepsilon$ is enough in our applications.
\end{remark*}

For $\theta$ very close to $\frac{1}{2}$, we can obtain the following better result than Lemma~\ref{l213}. The next lemma is also useful when performing role-reversals in the final decompositions.
\begin{lemma}\label{l214}
Let $\frac{1}{2} \leqslant \theta \leqslant \frac{45}{89} - \varepsilon$, $m_1 \sim M_1$ and $m_2 \sim M_2$. Then
$$
\sum_{\left(\frac{\log M_1}{\log x}, \frac{\log M_2}{\log x} \right) \in \boldsymbol{S}} a_{1, m_1} a_{2, m_2} S\left(\mathcal{A}^{q}_{m_1 m_2}, x^{\kappa}\right)
$$
has an asymptotic formula of the form (9).
\end{lemma}
\begin{proof}
Using Buchstab's identity, we have
\begin{align}
\nonumber \sum_{\left(\frac{\log M_1}{\log x}, \frac{\log M_2}{\log x} \right) \in \boldsymbol{S}} a_{1, m_1} a_{2, m_2} S\left(\mathcal{A}^{q}_{m_1 m_2}, x^{\kappa}\right) =&\ \sum_{\left(\frac{\log M_1}{\log x}, \frac{\log M_2}{\log x} \right) \in \boldsymbol{S}} a_{1, m_1} a_{2, m_2} S\left(\mathcal{A}^{q}_{m_1 m_2}, \exp\left(\log x (\log \log x)^{-3} \right)\right) \\
\nonumber & - \sum_{\substack{\left(\frac{\log M_1}{\log x}, \frac{\log M_2}{\log x} \right) \in \boldsymbol{S} \\ (\log \log x)^{-3} \leqslant \alpha_{1} < \kappa }} a_{1, m_1} a_{2, m_2} S\left(\mathcal{A}^{q}_{m_1 m_2 p_1}, p_{1}\right) \\
\nonumber =&\ \Sigma_{2141} - \Sigma_{2142}.
\end{align}
We can give an asymptotic formula of the form (9) for $\Sigma_{2141}$ by an application of Lemma~\ref{l29}. Since $\exp\left(\log x (\log \log x)^{-3} \right) \gg \exp\left((\log x)^{\varepsilon}\right)$, we can give an asymptotic formula of the form (9) for $\Sigma_{2142}$ by Lemma~\ref{l211}. Now the proof of Lemma~\ref{l214} is completed.
\end{proof}

The next lemma is one of the most important asymptotic formulas used in \cite{676} and \cite{LRB679}, and we shall use this lemma only when $\theta$ is large.
\begin{lemma}\label{l215}
Let $\frac{1}{2} \leqslant \theta \leqslant \frac{4}{7} - \varepsilon$. Define
\begin{align}
\nonumber \boldsymbol{T}^{**} =&\ \boldsymbol{T}^{**}(\theta) = \left\{(s, t): \frac{3}{7} + \varepsilon \leqslant s \leqslant 1 - \theta - \varepsilon,\ 0 \leqslant t < \frac{1-s}{2} \right\}, \\
\nonumber \boldsymbol{U}_{j}^{*} =&\ \boldsymbol{U}_{j}^{*}(\theta) = \left\{\boldsymbol{\alpha}_{j}: \boldsymbol{\alpha}_{j} \text{ partitions exactly into } \boldsymbol{T}^{**} \right\}.
\end{align}
Let $\xi(\boldsymbol{\alpha_j})$ be a continuous function with $\varepsilon^2 \leqslant \xi(\boldsymbol{\alpha_j}) \leqslant \frac{1}{2}$. Then
$$
\sum_{\boldsymbol{\alpha}_{j} \in \boldsymbol{U}_{j}^{*}} S\left(\mathcal{A}^{q}_{p_1 \cdots p_j}, x^{\xi(\boldsymbol{\alpha_j})}\right)
$$
has an asymptotic formula of the form (9).
\end{lemma}
\begin{proof}
This lemma follows by similar arguments as in the proof of Lemma~\ref{l213} and the discussions in [\cite{HarmanBOOK}, Lemma 8.14].
\end{proof}
Lemma~\ref{l215} gives asymptotic formulas of the form (9) for sums of sieve functions $S\left(\mathcal{A}^{q}_{p_1 \cdots p_n}, p_n\right)$ when some of the variables can be grouped to lie in the interval $\left[\frac{3}{7} + \varepsilon,\ 1 - \theta - \varepsilon \right]$. However, we will not use this lemma in the final decompositions when $\theta$ is not too large. Instead, we will prove a stronger result that gives asymptotic formulas of the form (9) and does not require the condition $t < \frac{1}{2}(1 - s)$.
\begin{lemma}\label{l216}
Let $\frac{1}{2} \leqslant \theta \leqslant \frac{127}{224} - \varepsilon$. Define
\begin{align}
\nonumber \boldsymbol{T}^{***} =&\ \boldsymbol{T}^{***}(\theta) = \left\{(s, t): \frac{3}{7} + \varepsilon \leqslant s \leqslant 1 - \theta - \varepsilon \text{ or } \theta + \varepsilon \leqslant s \leqslant \frac{4}{7} - \varepsilon \right\}, \\
\nonumber \boldsymbol{V}_{j} =&\ \boldsymbol{V}_{j}(\theta) = \left\{\boldsymbol{\alpha}_{j}: \boldsymbol{\alpha}_{j} \text{ partitions into } \boldsymbol{T}^{***} \right\}.
\end{align}
Suppose that $\min \boldsymbol{\alpha}_{j} \geqslant \frac{1}{1000} + 10 \varepsilon$. Then
$$
\sum_{\boldsymbol{\alpha}_{j} \in \boldsymbol{V}_{j}} S\left(\mathcal{A}^{q}_{p_1 \cdots p_j}, p_j\right)
$$
has an asymptotic formula of the form (9).
\end{lemma}
\begin{proof}
This lemma can be proved by the same method used in the proof of [\cite{MaynardLargeModuliI}, Lemma 8.13]. The proof requires Heath-Brown identity, Lemma~\ref{l210}, a discussion on the contribution from higher prime-powers counted in the von Mangoldt function (see [\cite{MaynardLargeModuliI}, Lemma 8.10]), and a removal of the dependencies between variables (see [\cite{MaynardLargeModuliI}, Lemma 8.11]). Note that this lemma also holds for more general Type-II sums.
\end{proof}

Finally, we define the whole Type-II region $\boldsymbol{G}_{j}$ and the two-dimensional region $U_2$. We split $U_2 \backslash \boldsymbol{G}_2$ into three parts $A, B$ and $C$. Regions $A, B, C$ and $\boldsymbol{G}_{j}$ will also be used in the next sections.
\begin{lemma}\label{l217}
Define
\begin{align}
\nonumber \boldsymbol{G}_{j} =&\ \boldsymbol{G}_{j}(\theta) = \left\{\boldsymbol{\alpha}_{j}: \boldsymbol{\alpha}_{j} \in \boldsymbol{\mathcal{G}}_{j} \cup \boldsymbol{V}_{j} \right\}, \\
\nonumber U_2 =&\ U_2(\theta) = \left\{\boldsymbol{\alpha}_{2}: \kappa \leqslant \alpha_1 < \frac{3}{7} + \varepsilon,\ \kappa \leqslant \alpha_2 < \min\left(\alpha_1, \frac{1}{2}(1-\alpha_1) \right) \right\}, \\
\nonumber A =&\ A(\theta) = \left\{\boldsymbol{\alpha}_{2}: \boldsymbol{\alpha}_{2} \in U_2,\ \boldsymbol{\alpha}_{2} \notin \boldsymbol{G}_{2},\ \alpha_1 + \alpha_2 < \theta + \varepsilon \right\}, \\
\nonumber B =&\ B(\theta) = \left\{\boldsymbol{\alpha}_{2}: \boldsymbol{\alpha}_{2} \in U_2,\ \boldsymbol{\alpha}_{2} \notin \boldsymbol{G}_{2} \cup A,\ \alpha_1 + 4 \alpha_2 < 3 - 3 \theta - \varepsilon \right\}, \\
\nonumber C =&\ C(\theta) = \left\{\boldsymbol{\alpha}_{2}: \boldsymbol{\alpha}_{2} \in U_2,\ \boldsymbol{\alpha}_{2} \notin \boldsymbol{G}_{2} \cup A \cup B \right\}.
\end{align}
Let $\frac{1}{2} \leqslant \theta \leqslant \frac{11}{21} - \varepsilon$. Then
$$
\sum_{\boldsymbol{\alpha}_{2} \in A \cup B} S\left(\mathcal{A}^{q}_{p_1 p_2}, x^{\kappa}\right)
$$
has an asymptotic formula of the form (9).
\end{lemma}
\begin{proof}
This lemma follows by similar arguments as in the proof of Lemma~\ref{l213} and the discussions in [\cite{HarmanBOOK}, Lemma 8.16].
\end{proof}

\subsection{High-dimensional Sieves}
In this subsection, we mention several results regarding the upper and lower bounds for sums of sieve functions. Results in this subsection are only applicable for Theorem~\ref{t22} but not Theorem~\ref{t21}. The first two lemmas are proved using a two-dimensional Harman's sieve, and we shall use them in the final decomposition for both $C_0^{*}(\theta)$ and $C_1^{*}(\theta)$.

\begin{lemma}\label{l218} ([\cite{676}, Lemma 20]).
Let $\frac{1}{2} \leqslant \theta < \frac{17}{32}$. Then we have, for almost all $q \sim x^{\theta}$,
$$
- \sum_{1-\theta - \varepsilon < \alpha_1 < \frac{1}{2}} S\left(\mathcal{A}^{q}_{p_1}, p_1 \right) \leqslant (1+o(1)) \frac{1}{\varphi(q)} \left( - \sum_{1-\theta - \varepsilon < \alpha_1 < \frac{1}{2}} S\left(\mathcal{B}^{q}_{p_1}, p_1 \right) + \sum_{\substack{\kappa \leqslant \alpha_3 < \alpha_1 < \frac{\theta + \varepsilon}{2} \\ 1 - \theta - \varepsilon < \alpha_1 + v \leqslant \theta + \varepsilon \\ (\alpha_1, v) \notin \boldsymbol{G}_{2} \\ \alpha_3 < \min\left(v, \frac{1}{2}(1-\alpha_1-v) \right) \\ (\alpha_1, v, \alpha_3) \notin \boldsymbol{G}_{3} }} S\left(\mathcal{B}^{q}_{p_1 m p_3}, p_3 \right) \right),
$$
where $m = x^{v}$ and $\left(m, P(p_1)\right) = 1$.
\end{lemma}
\begin{lemma}\label{l219} ([\cite{676}, Lemma 21]).
Let $\frac{1}{2} \leqslant \theta < \frac{17}{32}$. Then we have, for almost all $q \sim x^{\theta}$,
$$
- \sum_{1-\theta - \varepsilon < \alpha_1 < \frac{1}{2}} S\left(\mathcal{A}^{q}_{p_1}, p_1 \right) \geqslant (1+o(1)) \frac{1}{\varphi(q)} \left( - \sum_{1-\theta - \varepsilon < \alpha_1 < \frac{1}{2}} S\left(\mathcal{B}^{q}_{p_1}, p_1 \right) - \sum_{\substack{\kappa \leqslant \alpha_1 < v \\ 1 - \theta - \varepsilon < \alpha_1 + v \leqslant \theta + \varepsilon \\ (\alpha_1, v) \notin \boldsymbol{G}_{2} }} S\left(\mathcal{B}^{q}_{p_1 m}, p_1 \right) \right),
$$
where $m = x^{v}$ and $\left(m, P(p_1)\right) = 1$.
\end{lemma}

The next two lemmas are proved using a three-dimensional Harman's sieve, and we shall use them in the final decomposition only for $C_1^{*}(\theta)$. Note that the second one is only applicable for prime moduli $q \sim Q$.
\begin{lemma}\label{l220} ([\cite{676}, Lemma 24]).
Define
\begin{align}
\nonumber \boldsymbol{R} =&\ \boldsymbol{R}(\theta) =  \left\{\boldsymbol{\alpha}_{2}: \alpha_2 \leqslant \alpha_1,\ \alpha_1 + 2 \alpha_2 \leqslant 1,\ \alpha_1 + 4 \alpha_2 \geqslant 3 - 3 \theta - \varepsilon,\ 3 \alpha_2 \geqslant 2 \alpha_1, \right.\\
\nonumber & \left. \qquad \qquad \quad \max\left(\tau,\ \frac{31 \theta - 15}{3} + 5 \varepsilon \right) \leqslant \alpha_1 \leqslant \min\left(\frac{3}{7} + \varepsilon,\ 4 - 7 \theta - 3 \varepsilon \right) \right\}, \\
\nonumber \boldsymbol{D}_1 =&\ \boldsymbol{D}_1(\theta) = \left\{\boldsymbol{\alpha}_{3}: \alpha_1 \geqslant \kappa,\ \alpha_2 \geqslant \kappa,\ \alpha_3 \geqslant \kappa,\ \boldsymbol{\alpha}_{3} \notin \boldsymbol{G}_3, \right.\\
\nonumber & \left. \qquad \qquad \quad (\alpha_1, \alpha_2 + \alpha_3) \in \boldsymbol{R} \text{ with } \alpha_2 \geqslant \alpha_3 \text{ or } (\alpha_1 + \alpha_2, \alpha_3) \in \boldsymbol{R} \text{ with } \alpha_1 \geqslant \alpha_2 \right\}, \\
\nonumber \boldsymbol{D}_2 =&\ \boldsymbol{D}_2(\theta) = \left\{\boldsymbol{\alpha}_{4}: \alpha_1 \geqslant \kappa,\ \alpha_2 \geqslant \kappa,\ \alpha_3 \geqslant \kappa,\ \alpha_4 \geqslant \kappa, \right.\\
\nonumber & \left. \qquad \qquad \quad \boldsymbol{\alpha}_{4} \notin \boldsymbol{G}_4,\ (\alpha_1 + \alpha_2, \alpha_3 + \alpha_4) \in \boldsymbol{R} \right\}.
\end{align}
Let $\frac{1}{2} \leqslant \theta < \frac{16}{31} - \varepsilon$. Then we have, for almost all $q \sim x^{\theta}$,
$$
\sum_{\boldsymbol{\alpha}_{2} \in \boldsymbol{R}} S\left(\mathcal{A}^{q}_{p_1 p_2}, x^{\kappa}\right) \leqslant (1+o(1)) \frac{1}{\varphi(q)} \left( \sum_{\boldsymbol{\alpha}_{2} \in \boldsymbol{R}} S\left(\mathcal{B}^{q}_{p_1 p_2}, x^{\kappa}\right) + \frac{x}{\log x} (I_1 + I_2) \right),
$$
where
$$
I_1 = \frac{1}{\kappa} \int_{(t_1, t_2, t_3) \in \boldsymbol{D}_{1}} \frac{\omega\left(\frac{1 - t_1 - t_2 - t_3}{\kappa}\right)}{t_1 t_2 t_3} d t_3 d t_2 d t_1,
$$
$$
I_2 = \frac{1}{\kappa} \int_{(t_1, t_2, t_3, t_4) \in \boldsymbol{D}_{2}} \frac{\omega\left(\frac{1 - t_1 - t_2 - t_3 - t_4}{\kappa}\right)}{t_1 t_2 t_3 t_4} d t_4 d t_3 d t_2 d t_1.
$$
\end{lemma}
\begin{lemma}\label{l221} ([\cite{677}, Lemma 24]).
Define
\begin{align}
\nonumber \boldsymbol{R}_0 =&\ \boldsymbol{R}_0(\theta) = \left\{\boldsymbol{\alpha}_{2}: \alpha_1 + 2 \alpha_2 \leqslant 1,\ \alpha_1 + 4 \alpha_2 \geqslant 3 - 3 \theta - \varepsilon,\ 3 \alpha_2 \geqslant 2 \alpha_1 - 3 \varepsilon, \right.\\
\nonumber & \left. \qquad \qquad \quad \max\left(\frac{19 \theta - 7}{7},\ \frac{50 \theta - 19}{17} \right) + 24 \varepsilon \leqslant \alpha_1 \leqslant \frac{3}{7} + \varepsilon \right\}, \\
\nonumber \boldsymbol{D}_3 =&\ \boldsymbol{D}_3(\theta) = \left\{\boldsymbol{\alpha}_{3}: \alpha_1 \geqslant \kappa,\ \alpha_2 \geqslant \kappa,\ \alpha_3 \geqslant \kappa,\ \boldsymbol{\alpha}_{3} \notin \boldsymbol{G}_3, \right.\\
\nonumber & \left. \qquad \qquad \quad (\alpha_1, \alpha_2 + \alpha_3) \in \boldsymbol{R}_0 \text{ with } \alpha_2 \geqslant \alpha_3 \text{ or } (\alpha_1 + \alpha_2, \alpha_3) \in \boldsymbol{R}_0 \text{ with } \alpha_1 \geqslant \alpha_2 \right\}, \\
\nonumber \boldsymbol{D}_4 =&\ \boldsymbol{D}_4(\theta) = \left\{\boldsymbol{\alpha}_{4}: \alpha_1 \geqslant \kappa,\ \alpha_2 \geqslant \kappa,\ \alpha_3 \geqslant \kappa,\ \alpha_4 \geqslant \kappa, \right.\\
\nonumber & \left. \qquad \qquad \quad \boldsymbol{\alpha}_{4} \notin \boldsymbol{G}_4,\ (\alpha_1 + \alpha_2, \alpha_3 + \alpha_4) \in \boldsymbol{R}_0 \right\}.
\end{align}
Let $\frac{25}{49} - \varepsilon \leqslant \theta \leqslant \frac{92}{175} - \varepsilon$. Then we have, for almost all prime $q \sim x^{\theta}$,
$$
\sum_{\boldsymbol{\alpha}_{2} \in \boldsymbol{R}_0} S\left(\mathcal{A}^{q}_{p_1 p_2}, x^{\kappa}\right) \leqslant (1+o(1)) \frac{1}{\varphi(q)} \left( \sum_{\boldsymbol{\alpha}_{2} \in \boldsymbol{R}_0} S\left(\mathcal{B}^{q}_{p_1 p_2}, x^{\kappa}\right) + \frac{x}{\log x} (I_3 + I_4) \right),
$$
where
$$
I_3 = \frac{1}{\kappa} \int_{(t_1, t_2, t_3) \in \boldsymbol{D}_{3}} \frac{\omega\left(\frac{1 - t_1 - t_2 - t_3}{\kappa}\right)}{t_1 t_2 t_3} d t_3 d t_2 d t_1,
$$
$$
I_4 = \frac{1}{\kappa} \int_{(t_1, t_2, t_3, t_4) \in \boldsymbol{D}_{4}} \frac{\omega\left(\frac{1 - t_1 - t_2 - t_3 - t_4}{\kappa}\right)}{t_1 t_2 t_3 t_4} d t_4 d t_3 d t_2 d t_1.
$$
\end{lemma}

\subsection{Upper Bounds}
We shall construct the majorant $\rho_1(n)$ and prove upper bounds for $C_1(\theta)$ and $C_1^{*}(\theta)$ in this subsection. Before our final decompositions, we first mention some existing results of $C_1(\theta)$ and $C_1^{*}(\theta)$.
\begin{theorem}\label{t222}
The functions $C_1(\theta)$ and $C_1^{*}(\theta)$ satisfy the following conditions:

(1). $C_1(\theta) = C_1^{*}(\theta) = 1$ for all $\theta \leqslant 0.5 - \varepsilon$;

(2). $C_1^{*}(\theta) \leqslant 1 + \varepsilon$ for $\theta = 0.5$;

(3). $C_1^{*}(\theta)$ is monotonic increasing for $0.5 \leqslant \theta \leqslant 0.6$.
\end{theorem}
\begin{proof}
This theorem follows easily from the Bombieri--Vinogradov Theorem and [\cite{676}, Theorem 1(i)(ii)].
\end{proof}

Now we split the range $\theta \in \left[\frac{1}{2},\ 1 \right)$ to several subranges and use different methods to treat them and obtain good bounds for $C_1(\theta)$ and $C_1^{*}(\theta)$. Note that when $\theta \in \left[\frac{1}{2},\ \frac{127}{224} - \varepsilon \right]$, we can remove the condition $t < \frac{1}{2}(1 - s)$ when applying the Type-II range $\left[\frac{3}{7} + \varepsilon,\ 1 - \theta - \varepsilon \right]$. We shall not repeatedly state this range in the decompositions for the sake of simplicity.

\subsubsection{Case 1. $\frac{1}{2} \leqslant \theta \leqslant \frac{45}{89} - \varepsilon$}
In this case, our Type-II range becomes
\begin{equation}
\left((\log x)^{\varepsilon - 1},\ \frac{1}{6}(5 - 8 \theta) - \varepsilon \right].
\end{equation}
Using Buchstab's identity, we have
\begin{align}
\nonumber S\left(\mathcal{A}^{q}, (2x)^{\frac{1}{2}} \right) =&\ S\left(\mathcal{A}^{q}, x^{\kappa} \right) - \sum_{\kappa \leqslant \alpha_1 < \frac{1}{2}} S\left(\mathcal{A}^{q}_{p_1}, p_1 \right) \\
\nonumber =&\ S\left(\mathcal{A}^{q}, x^{\kappa} \right) - \sum_{\kappa \leqslant \alpha_1 < \frac{3}{7} + \varepsilon} S\left(\mathcal{A}^{q}_{p_1}, p_1 \right) - \sum_{\frac{3}{7} + \varepsilon \leqslant \alpha_1 \leqslant 1 - \theta - \varepsilon} S\left(\mathcal{A}^{q}_{p_1}, p_1 \right) - \sum_{1 - \theta - \varepsilon < \alpha_1 < \frac{1}{2}} S\left(\mathcal{A}^{q}_{p_1}, p_1 \right) \\
=&\ S^{q}_{2411} - S^{q}_{2412} - S^{q}_{2413} - S^{q}_{2414}.
\end{align}
By Lemma~\ref{l212} and Lemma~\ref{l216}, we can give asymptotic formulas of the form (9) for $S^{q}_{2411}$ and $S^{q}_{2413}$. We need to discard the whole of $S^{q}_{2414}$ in the case $C_1(\theta)$. In the case $C_1^{*}(\theta)$, however, we can use Lemma~\ref{l218} to give an upper bound for $S^{q}_{2414}$ with a loss of
\begin{equation}
\sum_{\substack{\kappa \leqslant \alpha_3 < \alpha_1 < \frac{\theta + \varepsilon}{2} \\ 1 - \theta - \varepsilon \leqslant \alpha_1 + v \leqslant \theta + \varepsilon \\ (\alpha_1, v) \notin \boldsymbol{G}_{2} \\ \alpha_3 < \min\left(v, \frac{1}{2}(1-\alpha_1-v) \right) \\ (\alpha_1, v, \alpha_3) \notin \boldsymbol{G}_{3} }} S\left(\mathcal{A}^{q}_{p_1 m p_3}, p_3 \right).
\end{equation}

For $S^{q}_{2412}$ in the case $C_1(\theta)$, we first split it into two parts:
\begin{align}
\nonumber S^{q}_{2412} =&\ \sum_{\kappa \leqslant \alpha_1 < \frac{3}{7} + \varepsilon} S\left(\mathcal{A}^{q}_{p_1}, p_1 \right) \\
\nonumber =&\ \sum_{\kappa \leqslant \alpha_1 < \tau} S\left(\mathcal{A}^{q}_{p_1}, p_1 \right) + \sum_{\tau \leqslant \alpha_1 < \frac{3}{7} + \varepsilon} S\left(\mathcal{A}^{q}_{p_1}, p_1 \right) \\
=&\ S^{q}_{24121} + S^{q}_{24122}.
\end{align}
Since $\theta \leqslant \frac{11}{21} - \varepsilon$, we have $\boldsymbol{\alpha}_2 \in \boldsymbol{G}_{2} \cup A \cup B$ for any $\kappa \leqslant \alpha_2 < \alpha_1 < \tau$. Hence we can apply Buchstab's identity twice again on $S^{q}_{24121}$ to get
\begin{align}
\nonumber S^{q}_{24121} =&\ \sum_{\kappa \leqslant \alpha_1 < \tau} S\left(\mathcal{A}^{q}_{p_1}, p_1 \right) \\
\nonumber =&\ \sum_{\kappa \leqslant \alpha_1 < \tau} S\left(\mathcal{A}^{q}_{p_1}, x^{\kappa} \right) - \sum_{\substack{\kappa \leqslant \alpha_1 < \tau \\ \kappa \leqslant \alpha_2 < \min\left(\alpha_1, \frac{1}{2}(1-\alpha_1)\right) \\ \boldsymbol{\alpha}_2 \in \boldsymbol{G}_{2} }} S\left(\mathcal{A}^{q}_{p_1 p_2}, p_2 \right) \\
\nonumber &- \sum_{\substack{\kappa \leqslant \alpha_1 < \tau \\ \kappa \leqslant \alpha_2 < \min\left(\alpha_1, \frac{1}{2}(1-\alpha_1)\right) \\ \boldsymbol{\alpha}_2 \in A \cup B }} S\left(\mathcal{A}^{q}_{p_1 p_2}, p_2 \right) \\
\nonumber =&\ \sum_{\kappa \leqslant \alpha_1 < \tau} S\left(\mathcal{A}^{q}_{p_1}, x^{\kappa} \right) - \sum_{\substack{\kappa \leqslant \alpha_1 < \tau \\ \kappa \leqslant \alpha_2 < \min\left(\alpha_1, \frac{1}{2}(1-\alpha_1)\right) \\ \boldsymbol{\alpha}_2 \in \boldsymbol{G}_{2} }} S\left(\mathcal{A}^{q}_{p_1 p_2}, p_2 \right) \\
\nonumber &- \sum_{\substack{\kappa \leqslant \alpha_1 < \tau \\ \kappa \leqslant \alpha_2 < \min\left(\alpha_1, \frac{1}{2}(1-\alpha_1)\right) \\ \boldsymbol{\alpha}_2 \in A \cup B }} S\left(\mathcal{A}^{q}_{p_1 p_2}, x^{\kappa} \right) + \sum_{\substack{\kappa \leqslant \alpha_1 < \tau \\ \kappa \leqslant \alpha_2 < \min\left(\alpha_1, \frac{1}{2}(1-\alpha_1)\right) \\ \boldsymbol{\alpha}_2 \in A \cup B \\ \kappa \leqslant \alpha_3 < \min\left(\alpha_2, \frac{1}{2}(1-\alpha_1 -\alpha_2)\right) \\ \boldsymbol{\alpha}_{3} \in \boldsymbol{G}_{3} }} S\left(\mathcal{A}^{q}_{p_1 p_2 p_3}, p_3 \right) \\
\nonumber &+ \sum_{\substack{\kappa \leqslant \alpha_1 < \tau \\ \kappa \leqslant \alpha_2 < \min\left(\alpha_1, \frac{1}{2}(1-\alpha_1)\right) \\ \boldsymbol{\alpha}_2 \in A \cup B \\ \kappa \leqslant \alpha_3 < \min\left(\alpha_2, \frac{1}{2}(1-\alpha_1 -\alpha_2)\right) \\ \boldsymbol{\alpha}_{3} \notin \boldsymbol{G}_{3} }} S\left(\mathcal{A}^{q}_{p_1 p_2 p_3}, p_3 \right) \\
=&\ S^{q}_{241211} - S^{q}_{241212} - S^{q}_{241213} + S^{q}_{241214} + S^{q}_{241215}.
\end{align}
We can give an asymptotic formula of the form (9) for $S^{q}_{241211}$ by Lemma~\ref{l212}. We can give an asymptotic formula of the form (9) for $S^{q}_{241213}$ by Lemma~\ref{l217}. We can give asymptotic formulas of the form (9) for $S^{q}_{241212}$ and $S^{q}_{241214}$ by Lemma~\ref{l211} or Lemma~\ref{l216}. For $S^{q}_{241215}$ we can either discard it with a three-dimensional loss, or perform a further decomposition using Buchstab's identity and reversed Buchstab's identity if we can. A difference between here and the decomposing process in the corresponding parts of [\cite{LRB679}, Section 6.1] is that we can replace the required condition $(\alpha_1, \alpha_2, \alpha_3, \alpha_3) \in \boldsymbol{U}_{4}$ with $(\alpha_1, \alpha_2, \alpha_3, \alpha_3) \in \boldsymbol{S}_{4}$. This replacement can also be used in similar conditions involving $\boldsymbol{U}_{j}$. Again, the condition $t < \frac{1}{2}(1 - s)$ in the Type-II range $\left[\frac{3}{7} + \varepsilon,\ 1 - \theta - \varepsilon \right]$ can be removed. For the remaining $S^{q}_{24122}$, we need to discard the whole of it since there are some $\boldsymbol{\alpha}_2 \in C$ that we cannot give an asymptotic formula of the form (9) under the conditions $\kappa \leqslant \alpha_2 < \alpha_1$ and $\tau \leqslant \alpha_1 < \frac{3}{7} + \varepsilon$.

For $S^{q}_{2412}$ in the case $C_1^{*}(\theta)$, we can use Lemma~\ref{l220} to deal with those parts such that $\boldsymbol{\alpha}_2 \in C$. Here, we use Buchstab's identity again directly on $S^{q}_{2412}$ and split the resulting sum into three subsums corresponding to $\boldsymbol{\alpha}_2 \in \boldsymbol{G}_{2}$, $\boldsymbol{\alpha}_2 \in A \cup B$ and $\boldsymbol{\alpha}_2 \in C$. The next steps are almost the same as the decomposing process in [\cite{LRB679}, Section 6.1]: we give asymptotic formulas of the form (10) when $\boldsymbol{\alpha}_2 \in \boldsymbol{G}_{2}$, perform a straightforward decomposition when $\boldsymbol{\alpha}_2 \in A \cup B$, and use Lemma~\ref{l220} when $\boldsymbol{\alpha}_2 \in C$.

Another important device that can be used here is the \textit{role-reversal}, which can be seen as a trick to ``use Type-I information more effective'' by changing the roles of a larger explicit variable and a smaller implicit variable in the sum. A more detailed discussion of a role-reversal can be found in [\cite{HarmanBOOK}, Chapter 5]. In either case $C_1(\theta)$ or $C_1^{*}(\theta)$, we can perform a role-reversal on the three-dimensional sum
$$
\sum_{\substack{\boldsymbol{\alpha}_{2} \in A \cup B \\ \kappa \leqslant \alpha_3 < \min\left(\alpha_2, \frac{1}{2}(1-\alpha_1 -\alpha_2)\right) \\ \boldsymbol{\alpha}_{3} \notin \boldsymbol{G}_{3} \\ (\alpha_1, \alpha_2, \alpha_3, \alpha_3) \notin \boldsymbol{S}_4 }} S\left(\mathcal{A}^{q}_{p_1 p_2 p_3}, p_3 \right)
$$
if $\boldsymbol{\alpha}_3 \in \boldsymbol{S}_3$ and $(1 - \alpha_1 - \alpha_2 - \alpha_3, \alpha_2, \alpha_3) \in \boldsymbol{S}_3$. That is, we have
\begin{align}
\nonumber & \sum_{\substack{\boldsymbol{\alpha}_{2} \in A \cup B \\ \kappa \leqslant \alpha_3 < \min\left(\alpha_2, \frac{1}{2}(1-\alpha_1 -\alpha_2)\right) \\ \boldsymbol{\alpha}_{3} \notin \boldsymbol{G}_{3} \\ (\alpha_1, \alpha_2, \alpha_3, \alpha_3) \notin \boldsymbol{S}_4 \\ \boldsymbol{\alpha}_3 \in \boldsymbol{S}_3 \\ (1 - \alpha_1 - \alpha_2 - \alpha_3, \alpha_2, \alpha_3) \in \boldsymbol{S}_3 }} S\left(\mathcal{A}^{q}_{p_1 p_2 p_3}, p_3 \right) \\
\nonumber =&\ \sum_{\substack{\boldsymbol{\alpha}_{2} \in A \cup B \\ \kappa \leqslant \alpha_3 < \min\left(\alpha_2, \frac{1}{2}(1-\alpha_1 -\alpha_2)\right) \\ \boldsymbol{\alpha}_{3} \notin \boldsymbol{G}_{3} \\ (\alpha_1, \alpha_2, \alpha_3, \alpha_3) \notin \boldsymbol{S}_4 \\ \boldsymbol{\alpha}_3 \in \boldsymbol{S}_3 \\ (1 - \alpha_1 - \alpha_2 - \alpha_3, \alpha_2, \alpha_3) \in \boldsymbol{S}_3 }} S\left(\mathcal{A}^{q}_{p_1 p_2 p_3}, x^{\kappa} \right) - \sum_{\substack{\boldsymbol{\alpha}_{2} \in A \cup B \\ \kappa \leqslant \alpha_3 < \min\left(\alpha_2, \frac{1}{2}(1-\alpha_1 -\alpha_2)\right) \\ \boldsymbol{\alpha}_{3} \notin \boldsymbol{G}_{3} \\ (\alpha_1, \alpha_2, \alpha_3, \alpha_3) \notin \boldsymbol{S}_4 \\ \boldsymbol{\alpha}_3 \in \boldsymbol{S}_3 \\ (1 - \alpha_1 - \alpha_2 - \alpha_3, \alpha_2, \alpha_3) \in \boldsymbol{S}_3 \\ \kappa \leqslant \alpha_4 < \min\left(\alpha_3, \frac{1}{2}(1-\alpha_1 -\alpha_2 -\alpha_3)\right) \\ \boldsymbol{\alpha}_{4} \in \boldsymbol{G}_{4} }} S\left(\mathcal{A}^{q}_{p_1 p_2 p_3 p_4}, p_4 \right) \\
\nonumber & - \sum_{\substack{\boldsymbol{\alpha}_{2} \in A \cup B \\ \kappa \leqslant \alpha_3 < \min\left(\alpha_2, \frac{1}{2}(1-\alpha_1 -\alpha_2)\right) \\ \boldsymbol{\alpha}_{3} \notin \boldsymbol{G}_{3} \\ (\alpha_1, \alpha_2, \alpha_3, \alpha_3) \notin \boldsymbol{S}_4 \\ \boldsymbol{\alpha}_3 \in \boldsymbol{S}_3 \\ (1 - \alpha_1 - \alpha_2 - \alpha_3, \alpha_2, \alpha_3) \in \boldsymbol{S}_3 \\ \kappa \leqslant \alpha_4 < \min\left(\alpha_3, \frac{1}{2}(1-\alpha_1 -\alpha_2 -\alpha_3)\right) \\ \boldsymbol{\alpha}_{4} \notin \boldsymbol{G}_{4} }} S\left(\mathcal{A}^{q}_{p_1 p_2 p_3 p_4}, p_4 \right) \\
\nonumber =&\ \sum_{\substack{\boldsymbol{\alpha}_{2} \in A \cup B \\ \kappa \leqslant \alpha_3 < \min\left(\alpha_2, \frac{1}{2}(1-\alpha_1 -\alpha_2)\right) \\ \boldsymbol{\alpha}_{3} \notin \boldsymbol{G}_{3} \\ (\alpha_1, \alpha_2, \alpha_3, \alpha_3) \notin \boldsymbol{S}_4 \\ \boldsymbol{\alpha}_3 \in \boldsymbol{S}_3 \\ (1 - \alpha_1 - \alpha_2 - \alpha_3, \alpha_2, \alpha_3) \in \boldsymbol{S}_3 }} S\left(\mathcal{A}^{q}_{p_1 p_2 p_3}, x^{\kappa} \right) - \sum_{\substack{\boldsymbol{\alpha}_{2} \in A \cup B \\ \kappa \leqslant \alpha_3 < \min\left(\alpha_2, \frac{1}{2}(1-\alpha_1 -\alpha_2)\right) \\ \boldsymbol{\alpha}_{3} \notin \boldsymbol{G}_{3} \\ (\alpha_1, \alpha_2, \alpha_3, \alpha_3) \notin \boldsymbol{S}_4 \\ \boldsymbol{\alpha}_3 \in \boldsymbol{S}_3 \\ (1 - \alpha_1 - \alpha_2 - \alpha_3, \alpha_2, \alpha_3) \in \boldsymbol{S}_3 \\ \kappa \leqslant \alpha_4 < \min\left(\alpha_3, \frac{1}{2}(1-\alpha_1 -\alpha_2 -\alpha_3)\right) \\ \boldsymbol{\alpha}_{4} \in \boldsymbol{G}_{4} }} S\left(\mathcal{A}^{q}_{p_1 p_2 p_3 p_4}, p_4 \right) \\
\nonumber & - \sum_{\substack{\boldsymbol{\alpha}_{2} \in A \cup B \\ \kappa \leqslant \alpha_3 < \min\left(\alpha_2, \frac{1}{2}(1-\alpha_1 -\alpha_2)\right) \\ \boldsymbol{\alpha}_{3} \notin \boldsymbol{G}_{3} \\ (\alpha_1, \alpha_2, \alpha_3, \alpha_3) \notin \boldsymbol{S}_4 \\ \boldsymbol{\alpha}_3 \in \boldsymbol{S}_3 \\ (1 - \alpha_1 - \alpha_2 - \alpha_3, \alpha_2, \alpha_3) \in \boldsymbol{S}_3 \\ \kappa \leqslant \alpha_4 < \min\left(\alpha_3, \frac{1}{2}(1-\alpha_1 -\alpha_2 -\alpha_3)\right) \\ \boldsymbol{\alpha}_{4} \notin \boldsymbol{G}_{4} }} S\left(\mathcal{A}^{q}_{\beta_1 p_2 p_3 p_4}, \left(\frac{2x}{\beta_1 p_2 p_3 p_4}\right)^{\frac{1}{2}} \right) \\
\nonumber =&\ \sum_{\substack{\boldsymbol{\alpha}_{2} \in A \cup B \\ \kappa \leqslant \alpha_3 < \min\left(\alpha_2, \frac{1}{2}(1-\alpha_1 -\alpha_2)\right) \\ \boldsymbol{\alpha}_{3} \notin \boldsymbol{G}_{3} \\ (\alpha_1, \alpha_2, \alpha_3, \alpha_3) \notin \boldsymbol{S}_4 \\ \boldsymbol{\alpha}_3 \in \boldsymbol{S}_3 \\ (1 - \alpha_1 - \alpha_2 - \alpha_3, \alpha_2, \alpha_3) \in \boldsymbol{S}_3 }} S\left(\mathcal{A}^{q}_{p_1 p_2 p_3}, x^{\kappa} \right) - \sum_{\substack{\boldsymbol{\alpha}_{2} \in A \cup B \\ \kappa \leqslant \alpha_3 < \min\left(\alpha_2, \frac{1}{2}(1-\alpha_1 -\alpha_2)\right) \\ \boldsymbol{\alpha}_{3} \notin \boldsymbol{G}_{3} \\ (\alpha_1, \alpha_2, \alpha_3, \alpha_3) \notin \boldsymbol{S}_4 \\ \boldsymbol{\alpha}_3 \in \boldsymbol{S}_3 \\ (1 - \alpha_1 - \alpha_2 - \alpha_3, \alpha_2, \alpha_3) \in \boldsymbol{S}_3 \\ \kappa \leqslant \alpha_4 < \min\left(\alpha_3, \frac{1}{2}(1-\alpha_1 -\alpha_2 -\alpha_3)\right) \\ \boldsymbol{\alpha}_{4} \in \boldsymbol{G}_{4} }} S\left(\mathcal{A}^{q}_{p_1 p_2 p_3 p_4}, p_4 \right) \\
\nonumber & - \sum_{\substack{\boldsymbol{\alpha}_{2} \in A \cup B \\ \kappa \leqslant \alpha_3 < \min\left(\alpha_2, \frac{1}{2}(1-\alpha_1 -\alpha_2)\right) \\ \boldsymbol{\alpha}_{3} \notin \boldsymbol{G}_{3} \\ (\alpha_1, \alpha_2, \alpha_3, \alpha_3) \notin \boldsymbol{S}_4 \\ \boldsymbol{\alpha}_3 \in \boldsymbol{S}_3 \\ (1 - \alpha_1 - \alpha_2 - \alpha_3, \alpha_2, \alpha_3) \in \boldsymbol{S}_3 \\ \kappa \leqslant \alpha_4 < \min\left(\alpha_3, \frac{1}{2}(1-\alpha_1 -\alpha_2 -\alpha_3)\right) \\ \boldsymbol{\alpha}_{4} \notin \boldsymbol{G}_{4} }} S\left(\mathcal{A}^{q}_{\beta_1 p_2 p_3 p_4}, x^{\kappa} \right) + \sum_{\substack{\boldsymbol{\alpha}_{2} \in A \cup B \\ \kappa \leqslant \alpha_3 < \min\left(\alpha_2, \frac{1}{2}(1-\alpha_1 -\alpha_2)\right) \\ \boldsymbol{\alpha}_{3} \notin \boldsymbol{G}_{3} \\ (\alpha_1, \alpha_2, \alpha_3, \alpha_3) \notin \boldsymbol{S}_4 \\ \boldsymbol{\alpha}_3 \in \boldsymbol{S}_3 \\ (1 - \alpha_1 - \alpha_2 - \alpha_3, \alpha_2, \alpha_3) \in \boldsymbol{S}_3 \\ \kappa \leqslant \alpha_4 < \min\left(\alpha_3, \frac{1}{2}(1-\alpha_1 -\alpha_2 -\alpha_3)\right) \\ \boldsymbol{\alpha}_{4} \notin \boldsymbol{G}_{4} \\ \kappa \leqslant \alpha_5 < \frac{1}{2} \alpha_1 \\ (1 -\alpha_1 -\alpha_2 -\alpha_3 -\alpha_4, \alpha_2, \alpha_3, \alpha_4, \alpha_5) \in \boldsymbol{G}_{5} }} S\left(\mathcal{A}^{q}_{\beta_1 p_2 p_3 p_4 p_5}, p_5 \right) \\
\nonumber & + \sum_{\substack{\boldsymbol{\alpha}_{2} \in A \cup B \\ \kappa \leqslant \alpha_3 < \min\left(\alpha_2, \frac{1}{2}(1-\alpha_1 -\alpha_2)\right) \\ \boldsymbol{\alpha}_{3} \notin \boldsymbol{G}_{3} \\ (\alpha_1, \alpha_2, \alpha_3, \alpha_3) \notin \boldsymbol{S}_4 \\ \boldsymbol{\alpha}_3 \in \boldsymbol{S}_3 \\ (1 - \alpha_1 - \alpha_2 - \alpha_3, \alpha_2, \alpha_3) \in \boldsymbol{S}_3 \\ \kappa \leqslant \alpha_4 < \min\left(\alpha_3, \frac{1}{2}(1-\alpha_1 -\alpha_2 -\alpha_3)\right) \\ \boldsymbol{\alpha}_{4} \notin \boldsymbol{G}_{4} \\ \kappa \leqslant \alpha_5 < \frac{1}{2} \alpha_1 \\ (1 -\alpha_1 -\alpha_2 -\alpha_3 -\alpha_4, \alpha_2, \alpha_3, \alpha_4, \alpha_5) \notin \boldsymbol{G}_{5} }} S\left(\mathcal{A}^{q}_{\beta_1 p_2 p_3 p_4 p_5}, p_5 \right) \\
=&\ S^{q}_{24101} - S^{q}_{24102} - S^{q}_{24103} + S^{q}_{24104} + S^{q}_{24105}.
\end{align}
We can give asymptotic formulas of the form (9) for $S^{q}_{24101}$ and $S^{q}_{24103}$ by Lemma~\ref{l214}. We can give asymptotic formulas of the form (9) for $S^{q}_{24102}$ and $S^{q}_{24104}$ by Lemma~\ref{l211} or Lemma~\ref{l216}. We can simply discard $S^{q}_{24105}$ with a small five-dimensional loss. This process replaces a three-dimensional loss integral by a five-dimensional loss integral. However, if any $\boldsymbol{S}_3$ in the condition is replaced by $\boldsymbol{U}_3$, we cannot use the above role-reversal technique since we need $\boldsymbol{\alpha}_{j} \in \boldsymbol{A}_{j}$ in order to use the combinatorial lemmas such as [\cite{676}, Lemma 11] (see the proof of Lemma~\ref{l213}, for example). Thus, we shall only use role-reversals when $\frac{1}{2} \leqslant \theta \leqslant \frac{45}{89} - \varepsilon$.

\subsubsection{Case 2. $\frac{45}{89} - \varepsilon < \theta \leqslant \frac{17}{33} - \varepsilon$}
In this case, our Type-II range becomes
\begin{equation}
\left((\log x)^{\varepsilon - 1},\ \frac{17 - 33\theta}{28} - \varepsilon \right] \cup \left[2 \theta - 1 + \varepsilon,\ \frac{1}{6}(5 - 8 \theta) - \varepsilon \right].
\end{equation}
By a similar method as in the proof of Lemma~\ref{l214}, we can give an asymptotic formula of the form (9) for
$$
\sum_{\boldsymbol{\alpha}_{j} \in \boldsymbol{S}_{j}} S\left(\mathcal{A}^{q}_{p_1 \cdots p_j}, x^{\frac{17 - 33\theta}{28} - \varepsilon}\right).
$$
However, we must handle sums with some prime variables lie in $\left(x^{\frac{17 - 33\theta}{28} - \varepsilon}, x^{2 \theta - 1 + \varepsilon}\right)$ if we perform a straightforward decomposition to get the above sum. When $(\alpha_1, \alpha_2, \alpha_3, \alpha_3) \in \boldsymbol{U}_{4}$, performing a straightforward decomposition is good. But when $(\alpha_1, \alpha_2, \alpha_3, \alpha_3) \in \boldsymbol{S}_{4}$ but $(\alpha_1, \alpha_2, \alpha_3, \alpha_3) \notin \boldsymbol{U}_{4}$, doing a straightforward decomposition is not that good anymore. When $\theta = \frac{45}{89}$, we have $\frac{17 - 33\theta}{28} = \frac{1}{89}$, which means that the corresponding sums may count numbers with about $90$ prime factors. After numerical calculations we found that the sizes of resulting sums exceed the sizes of the original ones, hence we decide not to use this asymptotic formula in this case.

In the case $C_1(\theta)$, the decompositions are similar to the decompositions in \textit{Case 1}, with the condition $(\alpha_1, \alpha_2, \alpha_3, \alpha_3) \in \boldsymbol{S}_{4}$ and similar conditions involving $\boldsymbol{S}_{j}$ are replaced with $(\alpha_1, \alpha_2, \alpha_3, \alpha_3) \in \boldsymbol{U}_{4}$ and corresponding ones with $\boldsymbol{U}_{j}$. Unlike \textit{Case 1}, we do not perform any role-reversals when $\theta > \frac{45}{89} - \varepsilon$. We still need to discard two one-dimensional sums that summing over $p_1 \in \left[\tau,\ \frac{3}{7} + \varepsilon \right)$ (see $S^{q}_{24122}$) and $p_1 \in \left(1 - \theta - \varepsilon,\ \frac{1}{2} \right)$ (see $S^{q}_{2414}$).

In the case $C_1^{*}(\theta)$, the decomposition details are very similar to those in \cite{LRB679}, with several modifications. We can remove the condition $t < \frac{1}{2}(1 - s)$ when applying the Type-II range $\left[\frac{3}{7} + \varepsilon,\ 1 - \theta - \varepsilon \right]$. Since $\theta < \frac{17}{32} - \varepsilon$, the condition $\alpha_1 < \tau$ in the definition of $\boldsymbol{U}_{j}$ can be removed by the Remark of Lemma~\ref{l213}. Since $\frac{17}{33} < \frac{16}{31} - \varepsilon$, Lemma~\ref{l220} is still applicable. We cannot use Lemma~\ref{l221} except for prime $q$, and we need to discard a one-dimensional sum that summing over
$$
p_1 \in \left[\tau,\ \max\left(\tau,\ \frac{31 \theta - 15}{3} + 5 \varepsilon \right) \right) \cup \left(\min\left(\frac{3}{7} + \varepsilon,\ 4 - 7 \theta - 3 \varepsilon \right),\ \frac{3}{7} + \varepsilon \right).
$$
For prime moduli $q \sim Q$, we can use Lemma~\ref{l221} and the details are the same as those in \cite{LRB679}.

\subsubsection{Case 3. $\frac{17}{33} - \varepsilon < \theta \leqslant \frac{17}{32} - \varepsilon$}
From here, our Type-II range becomes
\begin{equation}
\left[2 \theta - 1 + \varepsilon,\ \frac{1}{6}(5 - 8 \theta) - \varepsilon \right].
\end{equation}

In the case $C_1(\theta)$, the decompositions are similar to the decompositions in \textit{Case 2}. We discard two one-dimensional sums that summing over $p_1 \in \left[\tau,\ \frac{3}{7} + \varepsilon \right)$ and $p_1 \in \left(1 - \theta - \varepsilon, \frac{1}{2} \right)$. More details can be found in \cite{LRB679}.

In the case $C_1^{*}(\theta)$, the decompositions are similar to the decompositions in \textit{Case 2} when $\theta < \frac{16}{31} - \varepsilon$. When $\frac{16}{31} - \varepsilon \leqslant \theta \leqslant \frac{11}{21} - \varepsilon$, we cannot use Lemma~\ref{l220}, and we need to discard a one-dimensional sum that summing over $p_1 \in \left[\tau,\ \frac{3}{7} + \varepsilon \right)$ unless $q$ is prime. For prime moduli $q \sim Q$, the decompositions are almost the same as those in \cite{LRB679} when $\frac{16}{31} - \varepsilon \leqslant \theta \leqslant \frac{11}{21} - \varepsilon$. When $\frac{11}{21} - \varepsilon < \theta \leqslant \frac{17}{32} - \varepsilon$, we use the simple bound
\begin{equation}
C_1^{*}(\theta) \leqslant 1.725 + 13.125(\theta - 0.523)
\end{equation}
from \cite{676}. Note that this bound does not require Lemma~\ref{l221}. In fact, the only difference between upper bounds for all moduli and for prime moduli occurs when $\theta \in \left[\frac{25}{49} - \varepsilon,\ \frac{92}{175} - \varepsilon \right]$, the range that we can apply Lemma~\ref{l221} to cover parts of $C$ that are not covered by Lemma~\ref{l220}.

\subsubsection{Case 4. $\frac{17}{32} - \varepsilon < \theta \leqslant \frac{4}{7} - \varepsilon$}
In this case, our Type-II range remains the same as in \textit{Case 3} when $\frac{17}{32} - \varepsilon < \theta < \frac{11}{20} - \varepsilon$ and becomes trivial when $\theta \geqslant \frac{11}{20}$. Note that all high-dimensional sieve results (see Lemmas~\ref{l218}--\ref{l221}) are unavailable from here, and we can only use those available results in Subsection 2.2 that gives asymptotic formulas of the form (9). Thus, the decompositions are the same in the cases $C_1(\theta)$ and $C_1^{*}(\theta)$. We only need to remove the restriction $t < \frac{1}{2}(1 - s)$ in Lemma~\ref{l215} and do the same decomposing process as in \cite{LRB679} when $\frac{17}{32} - \varepsilon < \theta \leqslant \frac{127}{224} - \varepsilon$. For $\frac{127}{224} - \varepsilon < \theta \leqslant \frac{4}{7} - \varepsilon$ we cannot do that removal, and the decompositions are exactly the same as in \cite{LRB679}. Note that we have better upper bounds for $C_1^{*}(\theta)$ when $\theta > 0.565$, and we shall discuss them in the next case.

\subsubsection{Case 5. $\frac{4}{7} - \varepsilon < \theta < 1$}
From here, almost all of the arithmetic information in this section becomes unavailable, and we can only use the standard upper bound sieve to construct the majorant $\rho_{1}(n)$ in the case $C_1(\theta)$. Following Fouvry--Radziwiłł \cite{FouvryRadziwill1}, we pick the sieve weights $\lambda_{d,z}$ as in Selberg's upper bound sieve of dimension 1:
$$
\sum_{\substack{d \mid n \\ d \leqslant z}} \lambda_{d,z} = \left(\sum_{\substack{d \mid n \\ d \leqslant \sqrt{z}}} \mu(d) \left(1 - \frac{\log d}{\log \sqrt{z}} \right) \right)^2
$$
and $\lambda_{d,z} = 0$ for $d > z$. It is clear that we have
$$
\mathbbm{1}_{p}(n) \leqslant \sum_{d \mid n} \lambda_{d,z}
$$
for $n > \sqrt{z}$ and
$$
\sum_{n \sim x} \sum_{d \mid n} \lambda_{d,z} = 2 (1+o(1)) \frac{x}{\log z}.
$$
Let $\rho_{1}(n) = \sum_{d \mid n} \lambda_{d,z}$. Assume first that $\theta > \frac{529}{630} + \varepsilon$. By our choice of $\rho_{1}(n)$ and [\cite{Wright26}, Corollary 2.4], we have the following bounds for $C_1(\theta)$:
\begin{align}
C_1(\theta) \leqslant
\begin{cases}
\frac{178}{45} + \varepsilon, & \frac{529}{630} + \varepsilon < \theta \leqslant \frac{177}{178}, \\
\frac{66}{28 \theta - 11} + \varepsilon, & \frac{177}{178} < \theta < 1.
\end{cases}
\end{align}
For $\frac{4}{7} - \varepsilon < \theta \leqslant \frac{529}{630} + \varepsilon$, we can use [\cite{BFI}, Theorem 0(b)] to prove an analog of [\cite{FouvryRadziwill1}, Corollary 1.5] that holds for $z \leqslant x^{\frac{1}{2} - 2 \varepsilon}$ and $x^{2 \varepsilon} < Q < x^{1 - 2 \varepsilon}$. This implies
\begin{equation}
C_1(\theta) \leqslant 4 + \varepsilon
\end{equation}
for $\theta$ in this interval.

In the case $C_1^{*}(\theta)$, we consider the range $0.565 < \theta < 1$. For $\frac{88}{89} < \theta < 1$ we can use the above upper bounds for $C_1(\theta)$ to bound $C_1^{*}(\theta)$, while for $0.565 < \theta \leqslant \frac{88}{89}$ we can instead use the existing better results of Baker and Harman \cite{676} and Fouvry \cite{FouvryFermat} based on the Rosser--Iwaniec sieve with a bilinear error term:
\begin{align}
C_1^{*}(\theta) \leqslant
\begin{cases}
\frac{14}{12 - 13\theta} - \log\left(\frac{4(1-\theta)}{3\theta} \right) + \varepsilon, & 0.565 < \theta < \frac{4}{7}, \\
\frac{14}{12 - 13\theta} + \varepsilon, & \frac{4}{7} \leqslant \theta < 0.6, \\
\frac{8}{3 - \theta} + \varepsilon, & 0.6 \leqslant \theta < \frac{5}{7}, \\
\frac{6}{1 + \theta} + \varepsilon, & \frac{5}{7} \leqslant \theta < \frac{3}{4}, \\
\frac{12}{5 - 2\theta} + \varepsilon, & \frac{3}{4} \leqslant \theta < \frac{5}{6}, \\
\frac{48}{15 - 2\theta} + \varepsilon, & \frac{5}{6} \leqslant \theta < \frac{9}{10}, \\
\frac{4}{2 - \theta} + \varepsilon, & \frac{9}{10} \leqslant \theta \leqslant \frac{88}{89}.
\end{cases}
\end{align}
Note that the bounds in the first two cases come from \cite{676}, and the remaining bounds come from \cite{FouvryFermat}.

The following tables give the values of the upper bounds for $C_1(\theta)$ and $C_1^{*}(\theta)$. Note that for $C_1^{*}(\theta)$, the only difference between upper bound values for all moduli and for prime moduli occurs when $\theta \in \left[\frac{25}{49} - \varepsilon,\ \frac{11}{21} - \varepsilon \right]$.
\begin{center}
\begin{tabular}{|c|c|c|c|c|c|}
\hline \boldmath{$\theta$} & \boldmath{$C_1(\theta)$} & \boldmath{$\theta$} & \boldmath{$C_1(\theta)$} & \boldmath{$\theta$} & \boldmath{$C_1(\theta)$} \\
\hline $0.500$ & $1.5823$ & $0.524$ & $1.7805$ & $0.548$ & $2.3725$ \\
\hline $0.501$ & $1.5900$ & $0.525$ & $1.7894$ & $0.549$ & $2.3738$ \\
\hline $0.502$ & $1.5977$ & $0.526$ & $1.7966$ & $0.550$ & $2.3996$ \\
\hline $0.503$ & $1.6055$ & $0.527$ & $1.8041$ & $0.551$ & $2.4318$ \\
\hline $0.504$ & $1.6132$ & $0.528$ & $1.8133$ & $0.552$ & $2.4393$ \\
\hline $0.505$ & $1.6214$ & $0.529$ & $1.8239$ & $0.553$ & $2.4404$ \\
\hline $0.506$ & $1.6292$ & $0.530$ & $1.8338$ & $0.554$ & $2.4426$ \\
\hline $0.507$ & $1.6368$ & $0.531$ & $1.8451$ & $0.555$ & $2.4540$ \\
\hline $0.508$ & $1.6446$ & $0.532$ & $1.8619$ & $0.556$ & $2.4776$ \\
\hline $0.509$ & $1.6532$ & $0.533$ & $1.8658$ & $0.557$ & $2.4987$ \\
\hline $0.510$ & $1.6614$ & $0.534$ & $1.8780$ & $0.558$ & $2.5120$ \\
\hline $0.511$ & $1.6692$ & $0.535$ & $1.9028$ & $0.559$ & $2.5327$ \\
\hline $0.512$ & $1.6776$ & $0.536$ & $1.9200$ & $0.560$ & $2.5345$ \\
\hline $0.513$ & $1.6862$ & $0.537$ & $1.9277$ & $0.561$ & $2.6895$ \\
\hline $0.514$ & $1.6944$ & $0.538$ & $1.9466$ & $0.562$ & $2.7144$ \\
\hline $0.515$ & $1.7031$ & $0.539$ & $1.9568$ & $0.563$ & $2.8054$ \\
\hline $0.516$ & $1.7104$ & $0.540$ & $1.9993$ & $0.564$ & $2.8788$ \\
\hline $0.517$ & $1.7193$ & $0.541$ & $2.0582$ & $0.565$ & $2.9548$ \\
\hline $0.518$ & $1.7283$ & $0.542$ & $2.0793$ & $0.566$ & $3.0791$ \\
\hline $0.519$ & $1.7342$ & $0.543$ & $2.1223$ & $0.567$ & $3.3065$ \\
\hline $0.520$ & $1.7394$ & $0.544$ & $2.1604$ & $0.568$ & $3.4768$ \\
\hline $0.521$ & $1.7493$ & $0.545$ & $2.2179$ & $0.569$ & $3.6696$ \\
\hline $0.522$ & $1.7611$ & $0.546$ & $2.3122$ & $0.570$ & $3.7733$ \\
\hline $0.523$ & $1.7710$ & $0.547$ & $2.3310$ & $0.571$ & $3.8666$ \\
\hline
\end{tabular} \\
\textbf{Table 2.1: Upper Bounds for }\boldmath{$C_1(\theta)$} \textbf{(}\boldmath{$0.5 \leqslant \theta \leqslant \frac{4}{7} - \varepsilon$}\textbf{)}
\end{center}
\begin{center}
\begin{tabular}{|c|c|c|c|c|c|}
\hline \boldmath{$\theta$} & \boldmath{$C_1^{*}(\theta)$} & \boldmath{$\theta$} & \boldmath{$C_1^{*}(\theta)$} & \boldmath{$\theta$} & \boldmath{$C_1^{*}(\theta)$} \\
\hline $0.500$ & $1 + \varepsilon$ & $0.522$ & $1.7320$ & $0.544$ & $2.1604$ \\
\hline $0.501$ & $1.0001$ & $0.523$ & $1.7382$ & $0.545$ & $2.2179$ \\
\hline $0.502$ & $1.0003$ & $0.524$ & $1.7382$ & $0.546$ & $2.3122$ \\
\hline $0.503$ & $1.0013$ & $0.525$ & $1.7500$ & $0.547$ & $2.3310$ \\
\hline $0.504$ & $1.0019$ & $0.526$ & $1.7644$ & $0.548$ & $2.3725$ \\
\hline $0.505$ & $1.0030$ & $0.527$ & $1.7775$ & $0.549$ & $2.3738$ \\
\hline $0.506$ & $1.0043$ & $0.528$ & $1.7907$ & $0.550$ & $2.3996$ \\
\hline $0.507$ & $1.0071$ & $0.529$ & $1.8038$ & $0.551$ & $2.4318$ \\
\hline $0.508$ & $1.0091$ & $0.530$ & $1.8169$ & $0.552$ & $2.4393$ \\
\hline $0.509$ & $1.0111$ & $0.531$ & $1.8300$ & $0.553$ & $2.4404$ \\
\hline $0.510$ & $1.0150$ & $0.532$ & $1.8619$ & $0.554$ & $2.4426$ \\
\hline $0.511$ & $1.0466$ & $0.533$ & $1.8658$ & $0.555$ & $2.4540$ \\
\hline $0.512$ & $1.0830$ & $0.534$ & $1.8780$ & $0.556$ & $2.4776$ \\
\hline $0.513$ & $1.1593$ & $0.535$ & $1.9028$ & $0.557$ & $2.4987$ \\
\hline $0.514$ & $1.2525$ & $0.536$ & $1.9200$ & $0.558$ & $2.5120$ \\
\hline $0.515$ & $1.3466$ & $0.537$ & $1.9277$ & $0.559$ & $2.5327$ \\
\hline $0.516$ & $1.4353$ & $0.538$ & $1.9466$ & $0.560$ & $2.5345$ \\
\hline $0.517$ & $1.6781$ & $0.539$ & $1.9568$ & $0.561$ & $2.6895$ \\
\hline $0.518$ & $1.6876$ & $0.540$ & $1.9993$ & $0.562$ & $2.7144$ \\
\hline $0.519$ & $1.6947$ & $0.541$ & $2.0582$ & $0.563$ & $2.8054$ \\
\hline $0.520$ & $1.7028$ & $0.542$ & $2.0793$ & $0.564$ & $2.8788$ \\
\hline $0.521$ & $1.7160$ & $0.543$ & $2.1223$ & $0.565$ & $2.9548$ \\
\hline
\end{tabular} \\
\textbf{Table 2.2: Upper Bounds for }\boldmath{$C_1^{*}(\theta)$} \textbf{(}\boldmath{$0.5 \leqslant \theta \leqslant 0.565$}\textbf{)}
\end{center}
\begin{center}
\begin{tabular}{|c|c|c|c|c|c|}
\hline \boldmath{$\theta$} & \boldmath{$C_1^{*}(\theta)$} & \boldmath{$\theta$} & \boldmath{$C_1^{*}(\theta)$} & \boldmath{$\theta$} & \boldmath{$C_1^{*}(\theta)$} \\
\hline $0.500$ & $1 + \varepsilon$ & $0.522$ & $1.6873$ & $0.544$ & $2.1604$ \\
\hline $0.501$ & $1.0001$ & $0.523$ & $1.7138$ & $0.545$ & $2.2179$ \\
\hline $0.502$ & $1.0003$ & $0.524$ & $1.7382$ & $0.546$ & $2.3122$ \\
\hline $0.503$ & $1.0013$ & $0.525$ & $1.7500$ & $0.547$ & $2.3310$ \\
\hline $0.504$ & $1.0019$ & $0.526$ & $1.7644$ & $0.548$ & $2.3725$ \\
\hline $0.505$ & $1.0030$ & $0.527$ & $1.7775$ & $0.549$ & $2.3738$ \\
\hline $0.506$ & $1.0043$ & $0.528$ & $1.7907$ & $0.550$ & $2.3996$ \\
\hline $0.507$ & $1.0071$ & $0.529$ & $1.8038$ & $0.551$ & $2.4318$ \\
\hline $0.508$ & $1.0091$ & $0.530$ & $1.8169$ & $0.552$ & $2.4393$ \\
\hline $0.509$ & $1.0111$ & $0.531$ & $1.8300$ & $0.553$ & $2.4404$ \\
\hline $0.510$ & $1.0150$ & $0.532$ & $1.8619$ & $0.554$ & $2.4426$ \\
\hline $0.511$ & $1.0237$ & $0.533$ & $1.8658$ & $0.555$ & $2.4540$ \\
\hline $0.512$ & $1.0315$ & $0.534$ & $1.8780$ & $0.556$ & $2.4776$ \\
\hline $0.513$ & $1.0789$ & $0.535$ & $1.9028$ & $0.557$ & $2.4987$ \\
\hline $0.514$ & $1.1431$ & $0.536$ & $1.9200$ & $0.558$ & $2.5120$ \\
\hline $0.515$ & $1.2199$ & $0.537$ & $1.9277$ & $0.559$ & $2.5327$ \\
\hline $0.516$ & $1.3199$ & $0.538$ & $1.9466$ & $0.560$ & $2.5345$ \\
\hline $0.517$ & $1.5740$ & $0.539$ & $1.9568$ & $0.561$ & $2.6895$ \\
\hline $0.518$ & $1.5948$ & $0.540$ & $1.9993$ & $0.562$ & $2.7144$ \\
\hline $0.519$ & $1.6136$ & $0.541$ & $2.0582$ & $0.563$ & $2.8054$ \\
\hline $0.520$ & $1.6320$ & $0.542$ & $2.0793$ & $0.564$ & $2.8788$ \\
\hline $0.521$ & $1.6591$ & $0.543$ & $2.1223$ & $0.565$ & $2.9548$ \\
\hline
\end{tabular} \\
\textbf{Table 2.3: Upper Bounds for }\boldmath{$C_1^{*}(\theta)$} \textbf{(}\boldmath{$0.5 \leqslant \theta \leqslant 0.565$}\textbf{), Prime Moduli}
\end{center}

Combining various bounds in this subsection, we can recover the estimate
\begin{equation}
\int_{0.5}^{0.679} C_1^{*}(\theta) d \theta < 0.5
\end{equation}
proved in \cite{LRB679} under the assumption that $q$ is prime. However, we are still unable to show that
\begin{equation}
\int_{0.5}^{0.68} C_1^{*}(\theta) d \theta < 0.5.
\end{equation}

\subsection{Lower Bounds}
We shall construct the minorant $\rho_0(n)$ and prove lower bounds for $C_0(\theta)$ and $C_0^{*}(\theta)$ in this subsection. Before our final decompositions, we first mention some existing results of $C_0(\theta)$ and $C_0^{*}(\theta)$.
\begin{theorem}\label{t223}
The function $C_0(\theta)$ satisfies the following conditions:

(1). $C_0(\theta) = C_0^{*}(\theta) = 1$ for all $\theta \leqslant 0.5 - \varepsilon$;

(2). $C_0^{*}(\theta) \geqslant 1 - \varepsilon$ for $\theta = 0.5$;

(3). $C_0^{*}(\theta)$ is monotonic decreasing for $0.5 \leqslant \theta < \frac{17}{32}$.
\end{theorem}
\begin{proof}
This theorem follows easily from the Bombieri--Vinogradov Theorem, [\cite{676}, Theorem 1(i)(iii)] and the main results proved in \cite{Mikawa}.
\end{proof}

We first consider the case $C_0(\theta)$. Let $\frac{1}{2} \leqslant \theta \leqslant \frac{17}{32} - \varepsilon$. Using Buchstab's identity, we have
\begin{align}
\nonumber S\left(\mathcal{A}^{q}, (2x)^{\frac{1}{2}} \right) =&\ S\left(\mathcal{A}^{q}, x^{\kappa} \right) - \sum_{\kappa \leqslant \alpha_1 < \frac{1}{2}} S\left(\mathcal{A}^{q}_{p_1}, p_1 \right) \\
\nonumber =&\ S\left(\mathcal{A}^{q}, x^{\kappa} \right) - \sum_{\kappa \leqslant \alpha_1 \leqslant 1 - \theta - \varepsilon} S\left(\mathcal{A}^{q}_{p_1}, p_1 \right) - \sum_{1 - \theta - \varepsilon < \alpha_1 < \frac{1}{2}} S\left(\mathcal{A}^{q}_{p_1}, p_1 \right).
\end{align}
However, we cannot give an asymptotic formula of the form (9) for the last sum above. Hence we cannot get any nontrivial lower bound for $C_0(\theta)$ for $\theta \geqslant \frac{1}{2}$.

Now we focus on the case $C_0^{*}(\theta)$. We shall use two different methods to give lower bounds for $C_0^{*}(\theta)$. The first method comes from Harman's sieve (see \cite{676}), while the second method comes from the idea of Mikawa \cite{Mikawa}.

\subsubsection{First Method}
The first method is to use Harman's sieve as in \cite{676}. Unlike the upper bound case, in this case we can only discard positive terms that do not have asymptotic formulas of the form (10). Now we focus on the range $\theta \in \left[\frac{1}{2},\ \frac{17}{32} - \varepsilon \right]$ and split it into two subranges.

\noindent 2.5.1.1. \textit{Case 1. $\frac{1}{2} \leqslant \theta \leqslant \frac{45}{89} - \varepsilon$.}
Just as in Subsubsection 2.4.1, our Type-II range in this case is still (12), and we can also replace $\boldsymbol{U}_{j}$ with $\boldsymbol{S}_{j}$ in various conditions. By Buchstab's identity, we have
\begin{align}
\nonumber S\left(\mathcal{A}^{q}, (2x)^{\frac{1}{2}} \right) =&\ S\left(\mathcal{A}^{q}, x^{\kappa} \right) - \sum_{\kappa \leqslant \alpha_1 < \frac{1}{2}} S\left(\mathcal{A}^{q}_{p_1}, p_1 \right) \\
\nonumber =&\ S\left(\mathcal{A}^{q}, x^{\kappa} \right) - \sum_{\kappa \leqslant \alpha_1 < \frac{3}{7} + \varepsilon} S\left(\mathcal{A}^{q}_{p_1}, p_1 \right) - \sum_{\frac{3}{7} + \varepsilon \leqslant \alpha_1 \leqslant 1 - \theta - \varepsilon} S\left(\mathcal{A}^{q}_{p_1}, p_1 \right) - \sum_{1 - \theta - \varepsilon < \alpha_1 < \frac{1}{2}} S\left(\mathcal{A}^{q}_{p_1}, p_1 \right) \\
\nonumber =&\ S\left(\mathcal{A}^{q}, x^{\kappa} \right) - \sum_{\kappa \leqslant \alpha_1 < \frac{3}{7} + \varepsilon} S\left(\mathcal{A}^{q}_{p_1}, x^{\kappa} \right) - \sum_{\frac{3}{7} + \varepsilon \leqslant \alpha_1 \leqslant 1 - \theta - \varepsilon} S\left(\mathcal{A}^{q}_{p_1}, p_1 \right) - \sum_{1 - \theta - \varepsilon < \alpha_1 < \frac{1}{2}} S\left(\mathcal{A}^{q}_{p_1}, p_1 \right) \\
\nonumber & + \sum_{\substack{\kappa \leqslant \alpha_1 < \frac{3}{7} + \varepsilon \\ \kappa \leqslant \alpha_2 < \min\left(\alpha_1, \frac{1}{2}(1-\alpha_1)\right)}} S\left(\mathcal{A}^{q}_{p_1 p_2}, p_2 \right) \\
=&\ S^{q}_{25111} - S^{q}_{25112} - S^{q}_{25113} - S^{q}_{25114} + S^{q}_{25115}.
\end{align}
We can give asymptotic formulas of the form (10) for $S^{q}_{25111}$--$S^{q}_{25113}$ by Lemma~\ref{l212} and Lemma~\ref{l216}. We use Lemma~\ref{l219} to give a lower bound for $- S^{q}_{25114}$ with a two-dimensional loss. For $S^{q}_{25115}$ we split it into three subsums:
\begin{align}
\nonumber S^{q}_{25115} =&\ \sum_{\substack{\kappa \leqslant \alpha_1 < \frac{3}{7} + \varepsilon \\ \kappa \leqslant \alpha_2 < \min\left(\alpha_1, \frac{1}{2}(1-\alpha_1)\right)}} S\left(\mathcal{A}^{q}_{p_1 p_2}, p_2 \right) \\
\nonumber =&\ \sum_{\substack{\kappa \leqslant \alpha_1 < \frac{3}{7} + \varepsilon \\ \kappa \leqslant \alpha_2 < \min\left(\alpha_1, \frac{1}{2}(1-\alpha_1)\right) \\ \boldsymbol{\alpha}_{2} \in \boldsymbol{G}_{2} }} S\left(\mathcal{A}^{q}_{p_1 p_2}, p_2 \right) + \sum_{\substack{\kappa \leqslant \alpha_1 < \frac{3}{7} + \varepsilon \\ \kappa \leqslant \alpha_2 < \min\left(\alpha_1, \frac{1}{2}(1-\alpha_1)\right) \\ \boldsymbol{\alpha}_{2} \in A \cup B }} S\left(\mathcal{A}^{q}_{p_1 p_2}, p_2 \right) + \sum_{\substack{\kappa \leqslant \alpha_1 < \frac{3}{7} + \varepsilon \\ \kappa \leqslant \alpha_2 < \min\left(\alpha_1, \frac{1}{2}(1-\alpha_1)\right) \\ \boldsymbol{\alpha}_{2} \in C }} S\left(\mathcal{A}^{q}_{p_1 p_2}, p_2 \right) \\
=&\ S^{q}_{251151} + S^{q}_{251152} + S^{q}_{251153}.
\end{align}
For $S^{q}_{251151}$ we can use Lemma~\ref{l211} or Lemma~\ref{l216} to give an asymptotic formula of the form (10). We discard the whole of $S^{q}_{251153}$, leading to a two-dimensional loss. For the remaining $S^{q}_{251152}$, we can perform further straightforward decompositions if $(\alpha_1, \alpha_2, \alpha_2) \in \boldsymbol{S}_{3}$. Note that the conditions and details of further decompositions are similar to the upper bound case. Also, role-reversals can be applied if $\boldsymbol{\alpha}_{2} \in \boldsymbol{S}_{2}$ and $(1 - \alpha_1 - \alpha_2, \alpha_2) \in \boldsymbol{S}_{2}$. In this way we obtain the following bounds of $C_0^{*}(\theta)$:
\begin{center}
\begin{tabular}{|c|c|}
\hline \boldmath{$\theta$} & \boldmath{$C_0^{*}(\theta)$} \\
\hline $0.501$ & $0.8636$ \\
\hline $0.502$ & $0.8455$ \\
\hline $0.503$ & $0.8298$ \\
\hline $0.504$ & $0.8135$ \\
\hline $0.505$ & $0.7972$ \\
\hline
\end{tabular} \\
\textbf{Table 2.4: Lower Bounds for }\boldmath{$C_0^{*}(\theta)$} \textbf{(First Method, }\boldmath{$\frac{1}{2} < \theta \leqslant \frac{45}{89} - \varepsilon$}\textbf{)}
\end{center}

\noindent 2.5.1.2. \textit{Case 2. $\frac{45}{89} - \varepsilon < \theta \leqslant \frac{17}{32} - \varepsilon$.}
The decompositions in this case are similar to the first case. We need to replace $\boldsymbol{S}_{j}$ with $\boldsymbol{U}_{j}$ in various sums, and we do not use the role-reversals (even it is possible to use in some regions). Again, we can remove the condition $\alpha_1 < \tau$ in the definition of $\boldsymbol{U}_{j}$ by the Remark of Lemma~\ref{l213}. Working like the above case we get
\begin{center}
\begin{tabular}{|c|c|c|c|}
\hline \boldmath{$\theta$} & \boldmath{$C_0^{*}(\theta)$} & \boldmath{$\theta$} & \boldmath{$C_0^{*}(\theta)$} \\
\hline $0.506$ & $0.7299$ & $0.518$ & $0.4086$ \\
\hline $0.507$ & $0.7082$ & $0.519$ & $0.3744$ \\
\hline $0.508$ & $0.6894$ & $0.520$ & $0.3363$ \\
\hline $0.509$ & $0.6606$ & $0.521$ & $0.2973$ \\
\hline $0.510$ & $0.6369$ & $0.522$ & $0.2546$ \\
\hline $0.511$ & $0.6128$ & $0.523$ & $0.2148$ \\
\hline $0.512$ & $0.5907$ & $0.524$ & $0.1706$ \\
\hline $0.513$ & $0.5567$ & $0.525$ & $0.1218$ \\
\hline $0.514$ & $0.5332$ & $0.526$ & $0.0778$ \\
\hline $0.515$ & $0.5040$ & $0.527$ & $0.0267$ \\
\hline $0.516$ & $0.4762$ & $0.528$ & $-0.03$ \\
\hline $0.517$ & $0.4455$ & & \\
\hline
\end{tabular} \\
\textbf{Table 2.5: Lower Bounds for }\boldmath{$C_0^{*}(\theta)$} \textbf{(First Method, }\boldmath{$\frac{45}{89} - \varepsilon < \theta \leqslant \frac{17}{32} - \varepsilon$}\textbf{)}
\end{center}
Note that the lower bound becomes trivial when $\theta \geqslant 0.528$.

\subsubsection{Second Method}
The second method is to use a modified Eratosthenes sieve (or a modified Harman's sieve) developed by Mikawa \cite{Mikawa}. The main idea of Mikawa \cite{Mikawa} is to introduce a generalized sieve function that involves a sum over Möbius function $\mu(d)$, ``split'' this sieve function in two different ways to get two sums
$$
\sum_{x - \frac{x}{(\log x)^4} < p \leqslant x} 1 \qquad \text{and} \qquad \sum_{\substack{x - \frac{x}{(\log x)^4} < p \leqslant x \\ p \equiv a (\bmod q) }} 1
$$
together with ``other sums'' that count almost primes, and use a combinatorial way to measure the contribution from ``other sums''. Before the decomposition, we need the following 3 combinatorial lemmas. The first 2 lemmas focus on the case $\Omega(n) \neq 5$:

\begin{lemma}\label{l224} ([\cite{Mikawa}, Lemma 4]).
For square-free $n$, we have
\begin{align*}
\sum_{\substack{d \mid n \\ d < \sqrt{n}}} \mu(d) =
\begin{cases}
0, & \mu(n) = 1, \\
1, & n = p, \\
0, & p \mid n \text{ with } \sqrt{n} < p < n, \\
-2, & n = p_1 p_2 p_3 \text{ with } p_3 < p_2 < p_1 < \sqrt{n}, \\
-20, & \Omega(n) = 7 \text{ with } P^{-}(n) > n^{\frac{1}{8}}.
\end{cases}
\end{align*}
\end{lemma}
\begin{proof}
Mikawa did not give a proof of this lemma in \cite{Mikawa}, so we provide a proof here for a reference. Let $n = p_1 p_2 \cdots p_k$ where $p_1 > p_2 > \ldots > p_k$, and write $\mathscr{K} = \{1, 2, \ldots, k\}$.

If $\mu(n) = 1$, then $k$ is even. For any $\mathscr{I} \subseteq \mathscr{K}$, we know that $\left|\mathscr{I}\right|$ and $\left|\mathscr{I}^{\complement} \right|$ must have the same parity. Thus, we have
$$
\mu\left(\prod_{i \in \mathscr{I}}p_i\right) = \mu\left(\prod_{i \in \mathscr{I}^{\complement}}p_i\right).
$$
Let $\mathscr{I}$ be a set that satisfies
$$
\prod_{i \in \mathscr{I}}p_i < \sqrt{n},
$$
then we must have
$$
\prod_{i \in \mathscr{I}^{\complement}}p_i > \sqrt{n}.
$$
Now, we have
$$
\sum_{\substack{d \mid n \\ d < \sqrt{n}}} \mu(d) = \sum_{\mathscr{I} \subseteq \mathscr{K}} \mu\left(\prod_{i \in \mathscr{I}}p_i\right) = \sum_{\mathscr{I} \subseteq \mathscr{K}} \mu\left(\prod_{i \in \mathscr{I}^{\complement}}p_i\right) = \sum_{\substack{d \mid n \\ d > \sqrt{n}}} \mu(d)
$$
and
$$
\sum_{\mathscr{I} \subseteq \mathscr{K}} \mu\left(\prod_{i \in \mathscr{I}}p_i\right) + \sum_{\mathscr{I} \subseteq \mathscr{K}} \mu\left(\prod_{i \in \mathscr{I}^{\complement}}p_i\right) = \sum_{\substack{d \mid n \\ d < \sqrt{n}}} \mu(d) + \sum_{\substack{d \mid n \\ d > \sqrt{n}}} \mu(d) = \sum_{d \mid n} \mu(d) = 0.
$$
Of course these imply 
$$
\sum_{\substack{d \mid n \\ d < \sqrt{n}}} \mu(d) = \sum_{\substack{d \mid n \\ d > \sqrt{n}}} \mu(d) = 0.
$$

If $n = p_1$, then the only $d$ that satisfies the condition $d < \sqrt{n} = \sqrt{p_1}$ is $d = 1$. Hence
$$
\sum_{\substack{d \mid n \\ d < \sqrt{n}}} \mu(d) = \mu(1) = 1.
$$

If we have $p_1 \mid n$ with $\sqrt{n} < p_1 < n$, we must need $p_1 \nmid d$ in the sum to ensure $d < \sqrt{n}$. Since $p_2 p_3 \cdots p_k = \frac{n}{p_1} < \sqrt{n}$, any possible choices of prime factors of $d$ among $p_2, p_3, \cdots, p_k$ are acceptable. Write $\mathscr{K}^{\prime} = \{2, \ldots, k\}$ and $\mathscr{I}^{\prime} \subseteq \mathscr{K}^{\prime}$, we know that 
$$
\sum_{\substack{d \mid n \\ d < \sqrt{n}}} \mu(d) = \sum_{\mathscr{I}^{\prime} \subseteq \mathscr{K}^{\prime}} \mu\left(\prod_{i \in \mathscr{I}^{\prime}} p_i\right) = \sum_{\substack{d \mid p_2 p_3 \cdots p_k}} \mu(d) = 0.
$$

If $n = p_1 p_2 p_3$ with $p_3 < p_2 < p_1 < \sqrt{n}$, we have $p_1 p_2 > p_1 p_3 > p_2 p_3 = \frac{n}{p_1} > \sqrt{n}$. Thus, the only possible choices of $d$ are $1$, $p_1$, $p_2$ and $p_3$. Now,
$$
\sum_{\substack{d \mid n \\ d < \sqrt{n}}} \mu(d) = \mu(1) + \mu(p_1) + \mu(p_2) + \mu(p_3) = 1 - 3 = -2.
$$

If $\Omega(n) = 7$ with $P^{-}(n) > n^{\frac{1}{8}}$, we know that $p_4 p_5 p_6 p_7 > \left(n^{\frac{1}{8}}\right)^4 = n^{\frac{1}{2}}$ and $p_1 p_2 p_3 = \frac{n}{p_4 p_5 p_6 p_7} < n^{\frac{1}{2}}$. Since we have $p_1 > p_2 > \ldots > p_k$, we know that $d < \sqrt{n}$ implies $\Omega(d) < 3$. Now,
$$
\sum_{\substack{d \mid n \\ d < \sqrt{n}}} \mu(d) = \sum_{0 \leqslant i \leqslant 3} (-1)^i \binom{7}{i} = -20.
$$

Combining the above 5 cases, the proof of Lemma~\ref{l224} is completed.
\end{proof}

\begin{lemma}\label{l225} ([\cite{Mikawa}, Lemma 6]).
For square-free $n$, we have
\begin{align*}
\sum_{\substack{d \mid n \\ \mu^2(d) = 1 \\ \Omega(d) = 3 \\ \sqrt{n} < d < \sqrt{n P^{-}(d)} }} 2 =
\begin{cases}
2, & n = p_1 p_2 p_3 p_4 (p_4 < p_3 < p_2 < p_1) \text{ with } p_1 < p_2 p_3 p_4 \text{ and } p_2 p_3 < p_1, \\
0, & \Omega(n) \leqslant 3, \\
0, & n = p_1 p_2 p_3 p_4 (p_4 < p_3 < p_2 < p_1) \text{ with } p_1 > p_2 p_3 p_4 \text{ or } p_2 p_3 > p_1, \\
\leqslant 20, & \Omega(n) = 6, \\
0, & \Omega(n) = 7 \text{ with } P^{-}(n) > n^{\frac{1}{8}}.
\end{cases}
\end{align*}
\end{lemma}
\begin{remark*}
The conditions $p_1 < p_2 p_3 p_4$ and $p_2 p_3 < p_1$ correspond to the condition $\sqrt{n} < d < \sqrt{n P^{-}(d)}$ with $d = p_2 p_3 p_4$, since
$$
\sqrt{p_1 p_2 p_3 p_4} < p_2 p_3 p_4 \Longrightarrow p_1 p_2 p_3 p_4 < p_2^2 p_3^2 p_4^2 \Longrightarrow p_1 < p_2 p_3 p_4
$$
and
$$
p_2 p_3 p_4 < \sqrt{p_1 p_2 p_3 p_4^2} \Longrightarrow p_2^2 p_3^2 < p_1 p_2 p_3 \Longrightarrow p_2 p_3 < p_1.
$$
For the case $\Omega(n) = 6$, we can improve the upper bound $20$ when $\boldsymbol{\alpha}_6$ lies in some special regions. We shall explain this improvement in the end of this subsubsection.
\end{remark*}

For the case $\Omega(n) = 5$, we use the following third lemma instead:
\begin{lemma}\label{l226} ([\cite{Mikawa}, Lemma 5]).
For $n = p_1 p_2 p_3 p_4 p_5$ with $p_5 < p_4 < p_3 < p_2 < p_1$, we have
\begin{align*}
0 \leqslant \sum_{\substack{d \mid n \\ \mu^2(d) = 1 \\ \Omega(d) = 3 \\ \sqrt{n} < d < \sqrt{n P^{-}(d)} }} 2 - \sum_{\substack{d \mid n \\ d < \sqrt{n}}} \mu(d) \leqslant
\begin{cases}
2, & p_2 p_3 < p_1 p_5, \\
0, & \text{otherwise}.
\end{cases}
\end{align*}
\end{lemma}

Now we start our final decomposition. Let $\frac{1}{2} \leqslant \theta \leqslant \frac{17}{32} - \varepsilon$. Write $y = \frac{x}{(\log x)^4}$ and put
$$
\mathcal{A} =\left\{n : x - y < n \leqslant x,\ n \equiv a (\bmod q) \right\} \quad \text{and} \quad \mathcal{B} =\left\{n : x - y < n \leqslant x \right\}.
$$
We begin with the generalized sieve functions
\begin{equation}
S_{252A} = \sum_{n \in \mathcal{A}} \left(\sum_{\substack{d \mid n \\ d < \sqrt{x} }} \mu(d) \right) \psi\left(n, x^{\kappa} \right) \quad \text{and} \quad S_{252B} = \sum_{n \in \mathcal{B}} \left(\sum_{\substack{d \mid n \\ d < \sqrt{x} }} \mu(d) \right) \psi\left(n, x^{\kappa} \right).
\end{equation}
For $S_{252A}$, we restrict $n$ to square-free integers and replace the condition $d < \sqrt{x}$ with $d < \sqrt{n}$, leading to a cost error term of
\begin{equation}
\sum_{n \in \mathcal{A}} \left(\sum_{\substack{p^2 \mid n \\ p \geqslant x^{\kappa} }} 1 + \sum_{\substack{d \mid n \\ d^2 \in \mathcal{B} }} 1 \right).
\end{equation}
The first term can be bounded by
\begin{equation}
\sum_{n \in \mathcal{A}} \sum_{\substack{p^2 \mid n \\ p \geqslant x^{\kappa} }} 1 \ll \sum_{n \in \mathcal{B}} \sum_{\substack{p^2 \mid n \\ p \geqslant x^{\kappa} }} d(n - a) \ll x^{\varepsilon} \sum_{x^{\kappa} \leqslant p \leqslant \sqrt{x}} \left(\frac{x}{p^2} + 1 \right) \ll x^{\frac{7}{8} + \varepsilon}
\end{equation}
since $\theta < \frac{17}{32}$. By the method of Hooley \cite{HooleyBOOK} and bounds for incomplete Kloosterman sums, we get
\begin{equation}
\sum_{n \in \mathcal{A}} \sum_{\substack{d \mid n \\ d^2 \in \mathcal{B} }} 1 \ll \sum_{d^2 \in \mathcal{B}} \frac{y}{q d} + q^{\frac{1}{2} + \varepsilon} + \frac{x^{1 - \varepsilon}}{q} \ll \frac{y}{q (\log x)^{4}},
\end{equation}
since $\theta < \frac{17}{32}$. Combining (28)--(31) and by Lemma~\ref{l224}, for almost all $q \sim Q$ we have
\begin{align}
\nonumber S_{252A} =&\ \sum_{n \in \mathcal{A}} \left(\sum_{\substack{d \mid n \\ d < \sqrt{x} }} \mu(d) \right) \psi\left(n, x^{\kappa} \right) \\
\nonumber =&\ \sum_{n \in \mathcal{A}} \mu^2(n) \left(\sum_{\substack{d \mid n \\ d < \sqrt{n} }} \mu(d) \right) \psi\left(n, x^{\kappa} \right) + O\left(\frac{y}{q (\log x)^{4}}\right) \\
\nonumber =&\ \sum_{p \in \mathcal{A}} 1 - \sum_{\substack{p_1 p_2 p_3 \in \mathcal{A} \\ x^{\kappa} < p_3 < p_2 < p_1 < p_2 p_3 }} 2 - \sum_{\substack{p_1 \cdots p_7 \in \mathcal{A} \\ x^{\kappa} < p_7 < \cdots < p_1 }} 20 + \sum_{\substack{p_1 \cdots p_5 \in \mathcal{A} \\ x^{\kappa} < p_5 < \cdots < p_1 }} \mu^2(n) \left(\sum_{\substack{d \mid n \\ d < \sqrt{n} }} \mu(d) \right) \psi\left(n, x^{\kappa} \right) + O\left(\frac{y}{q (\log x)^{4}}\right) \\
=&\ \sum_{p \in \mathcal{A}} 1 - S_{2521} - S_{2522} + \sum_{\substack{n = p_1 \cdots p_5 \in \mathcal{A} \\ x^{\kappa} < p_5 < \cdots < p_1 }} \mu^2(n) \left(\sum_{\substack{d \mid n \\ d < \sqrt{n} }} \mu(d) \right) \psi\left(n, x^{\kappa} \right) + O\left(\frac{y}{q (\log x)^{4}}\right).
\end{align}
Similarly, we can decompose $S_{252B}$ as
\begin{align}
\nonumber S_{252B} =&\ \sum_{n \in \mathcal{B}} \left(\sum_{\substack{d \mid n \\ d < \sqrt{x} }} \mu(d) \right) \psi\left(n, x^{\kappa} \right) \\
\nonumber =&\ \sum_{n \in \mathcal{B}} \mu^2(n) \left(\sum_{\substack{d \mid n \\ d < \sqrt{n} }} \mu(d) \right) \psi\left(n, x^{\kappa} \right) + O\left(\frac{y}{(\log x)^{4}}\right) \\
\nonumber =&\ \sum_{p \in \mathcal{B}} 1 - \sum_{\substack{p_1 p_2 p_3 \in \mathcal{B} \\ x^{\kappa} < p_3 < p_2 < p_1 < p_2 p_3 }} 2 - \sum_{\substack{p_1 \cdots p_7 \in \mathcal{B} \\ x^{\kappa} < p_7 < \cdots < p_1 }} 20 + \sum_{\substack{p_1 \cdots p_5 \in \mathcal{B} \\ x^{\kappa} < p_5 < \cdots < p_1 }} \mu^2(n) \left(\sum_{\substack{d \mid n \\ d < \sqrt{n} }} \mu(d) \right) \psi\left(n, x^{\kappa} \right) + O\left(\frac{y}{(\log x)^{4}}\right) \\
=&\ \sum_{p \in \mathcal{B}} 1 - S_{2523} - S_{2524} + \sum_{\substack{n = p_1 \cdots p_5 \in \mathcal{B} \\ x^{\kappa} < p_5 < \cdots < p_1 }} \mu^2(n) \left(\sum_{\substack{d \mid n \\ d < \sqrt{n} }} \mu(d) \right) \psi\left(n, x^{\kappa} \right) + O\left(\frac{y}{(\log x)^{4}}\right).
\end{align}

Now, if we have
\begin{equation}
\sum_{n \in \mathcal{A}} \left(\sum_{\substack{d \mid n \\ d < \sqrt{x} }} \mu(d) \right) \psi\left(n, x^{\kappa} \right) = \frac{1}{\varphi(q)} \sum_{n \in \mathcal{B}} \left(\sum_{\substack{d \mid n \\ d < \sqrt{x} }} \mu(d) \right) \psi\left(n, x^{\kappa} \right) + O\left(\frac{y}{q (\log x)^{3}}\right)
\end{equation}
for almost all $q \sim Q$, we can put (32) and (33) together. For the sum
$$
\sum_{\substack{n = p_1 \cdots p_5 \in \mathcal{A} \\ x^{\kappa} < p_5 < \cdots < p_1 }} \mu^2(n) \left(\sum_{\substack{d \mid n \\ d < \sqrt{n} }} \mu(d) \right) \psi\left(n, x^{\kappa} \right),
$$
we can use Lemma~\ref{l226} to get
\begin{equation}
\sum_{\substack{n = p_1 \cdots p_5 \in \mathcal{A} \\ x^{\kappa} < p_5 < \cdots < p_1 }} \mu^2(n) \left(\sum_{\substack{d \mid n \\ d < \sqrt{n} }} \mu(d) \right) \psi\left(n, x^{\kappa} \right) \leqslant \sum_{\substack{n = p_1 \cdots p_5 \in \mathcal{A} \\ x^{\kappa} < p_5 < \cdots < p_1 }} \mu^2(n) \left(\sum_{\substack{d \mid n \\ \mu^2(d) = 1 \\ \Omega(d) = 3 \\ \sqrt{n} < d < \sqrt{n P^{-}(d)} }} 2 \right) \psi\left(n, x^{\kappa} \right).
\end{equation}
Similar to (29)--(31), we have
\begin{align}
\nonumber &\ \sum_{\substack{n = p_1 \cdots p_5 \in \mathcal{A} \\ x^{\kappa} < p_5 < \cdots < p_1 }} \mu^2(n) \left(\sum_{\substack{d \mid n \\ d < \sqrt{n} }} \mu(d) \right) \psi\left(n, x^{\kappa} \right) \\
\nonumber \leqslant&\ \sum_{\substack{n = p_1 \cdots p_5 \in \mathcal{A} \\ x^{\kappa} < p_5 < \cdots < p_1 }} \mu^2(n) \left(\sum_{\substack{d \mid n \\ \mu^2(d) = 1 \\ \Omega(d) = 3 \\ \sqrt{n} < d < \sqrt{n P^{-}(d)} }} 2 \right) \psi\left(n, x^{\kappa} \right) \\
\nonumber \leqslant&\ \sum_{n \in \mathcal{A}} \left(\sum_{\substack{d \mid n \\ \mu^2(d) = 1 \\ \Omega(d) = 3 \\ \sqrt{x} < d < \sqrt{x P^{-}(d)} }} 2 \right) \psi\left(n, x^{\kappa} \right) + O\left(\sum_{n \in \mathcal{A}} \sum_{\substack{d \mid n \\ d^2 \in \mathcal{B} }} 1 \right) \\
=&\ \sum_{n \in \mathcal{A}} \left(\sum_{\substack{d \mid n \\ \mu^2(d) = 1 \\ \Omega(d) = 3 \\ \sqrt{x} < d < \sqrt{x P^{-}(d)} }} 2 \right) \psi\left(n, x^{\kappa} \right) + O\left(\frac{y}{q (\log x)^{4}}\right).
\end{align}

Now we consider the above main term with condition $p_1 \cdots p_5 \in \mathcal{A}$ replaced by $p_1 \cdots p_5 \in \mathcal{B}$. Revisiting the above process (36) but with opposite direction and by Lemma~\ref{l225}, we get
\begin{align}
\nonumber &\ \sum_{n \in \mathcal{B}} \left(\sum_{\substack{d \mid n \\ \mu^2(d) = 1 \\ \Omega(d) = 3 \\ \sqrt{x} < d < \sqrt{x P^{-}(d)} }} 2 \right) \psi\left(n, x^{\kappa} \right) \\
\nonumber =&\ \sum_{n \in \mathcal{B}} \left(\sum_{\substack{d \mid n \\ \mu^2(d) = 1 \\ \Omega(d) = 3 \\ \sqrt{n} < d < \sqrt{n P^{-}(d)} }} 2 \right) \psi\left(n, x^{\kappa} \right) + O\left(\frac{y}{(\log x)^{3}}\right) \\
\nonumber =&\ \sum_{n \in \mathcal{B}} \mu^2(n) \left(\sum_{\substack{d \mid n \\ \mu^2(d) = 1 \\ \Omega(d) = 3 \\ \sqrt{n} < d < \sqrt{n P^{-}(d)} }} 2 \right) \psi\left(n, x^{\kappa} \right) + O\left(\frac{y}{(\log x)^{3}}\right) \\
\nonumber \leqslant&\ \sum_{\substack{p_1 p_2 p_3 p_4 \in \mathcal{B} \\ x^{\kappa} < p_4 < p_3 < p_2 < p_1 \\ p_1 < p_2 p_3 p_4 \\ p_2 p_3 < p_1 }} 2 + \sum_{\substack{p_1 \cdots p_5 \in \mathcal{B} \\ x^{\kappa} < p_5 < p_4 < p_3 < p_2 < p_1 \\ p_2 p_3 < p_1 p_5 }} 2 + \sum_{\substack{p_1 \cdots p_6 \in \mathcal{B} \\ x^{\kappa} < p_6 < p_5 < p_4 < p_3 < p_2 < p_1 }} 20 \\
\nonumber &+ \sum_{\substack{p_1 \cdots p_5 \in \mathcal{B} \\ x^{\kappa} < p_5 < \cdots < p_1 }} \mu^2(n) \left(\sum_{\substack{d \mid n \\ d < \sqrt{n} }} \mu(d) \right) \psi\left(n, x^{\kappa} \right) + O\left(\frac{y}{(\log x)^{3}}\right) \\
=&\ S_{2525} + S_{2526} + S_{2527} + \sum_{\substack{n = p_1 \cdots p_5 \in \mathcal{B} \\ x^{\kappa} < p_5 < \cdots < p_1 }} \mu^2(n) \left(\sum_{\substack{d \mid n \\ d < \sqrt{n} }} \mu(d) \right) \psi\left(n, x^{\kappa} \right) + O\left(\frac{y}{(\log x)^{3}}\right).
\end{align}

Now we suppose that (34) holds true. Suppose also that we have
\begin{equation}
\sum_{n \in \mathcal{A}} \left(\sum_{\substack{d \mid n \\ \mu^2(d) = 1 \\ \Omega(d) = 3 \\ \sqrt{x} < d < \sqrt{x P^{-}(d)} }} 2 \right) \psi\left(n, x^{\kappa} \right) = \frac{1}{\varphi(q)} \sum_{n \in \mathcal{B}} \left(\sum_{\substack{d \mid n \\ \mu^2(d) = 1 \\ \Omega(d) = 3 \\ \sqrt{x} < d < \sqrt{x P^{-}(d)} }} 2 \right) \psi\left(n, x^{\kappa} \right) + O\left(\frac{y}{q (\log x)^{3}}\right)
\end{equation}
for almost all $q \sim Q$. By (36)--(38), we get
\begin{align}
\nonumber &\ \sum_{\substack{n = p_1 \cdots p_5 \in \mathcal{A} \\ x^{\kappa} < p_5 < \cdots < p_1 }} \mu^2(n) \left(\sum_{\substack{d \mid n \\ d < \sqrt{n} }} \mu(d) \right) \psi\left(n, x^{\kappa} \right) \\
\leqslant&\ \frac{1}{\varphi(q)} \sum_{\substack{n = p_1 \cdots p_5 \in \mathcal{B} \\ x^{\kappa} < p_5 < \cdots < p_1 }} \mu^2(n) \left(\sum_{\substack{d \mid n \\ d < \sqrt{n} }} \mu(d) \right) \psi\left(n, x^{\kappa} \right) + \frac{1}{\varphi(q)} \left(S_{2525} + S_{2526} + S_{2527}\right) + O\left(\frac{y}{q (\log x)^{3}}\right)
\end{align}
and thus
\begin{align}
\nonumber &\ \sum_{p \in \mathcal{A}} 1 - S_{2521} - S_{2522} + \sum_{\substack{n = p_1 \cdots p_5 \in \mathcal{A} \\ x^{\kappa} < p_5 < \cdots < p_1 }} \mu^2(n) \left(\sum_{\substack{d \mid n \\ d < \sqrt{n} }} \mu(d) \right) \psi\left(n, x^{\kappa} \right) \\
\nonumber \leqslant&\ \frac{1}{\varphi(q)} \sum_{\substack{n = p_1 \cdots p_5 \in \mathcal{B} \\ x^{\kappa} < p_5 < \cdots < p_1 }} \mu^2(n) \left(\sum_{\substack{d \mid n \\ d < \sqrt{n} }} \mu(d) \right) \psi\left(n, x^{\kappa} \right) \\
&+ \sum_{p \in \mathcal{A}} 1 - S_{2521} - S_{2522} + \frac{1}{\varphi(q)} \left(S_{2525} + S_{2526} + S_{2527} \right) + O\left(\frac{y}{q (\log x)^{3}}\right)
\end{align}
for almost all $q \sim Q$. By (32)--(34), we also have
\begin{align}
\nonumber &\ \sum_{p \in \mathcal{A}} 1 - S_{2521} - S_{2522} + \sum_{\substack{n = p_1 \cdots p_5 \in \mathcal{A} \\ x^{\kappa} < p_5 < \cdots < p_1 }} \mu^2(n) \left(\sum_{\substack{d \mid n \\ d < \sqrt{n} }} \mu(d) \right) \psi\left(n, x^{\kappa} \right) \\
=&\ \frac{1}{\varphi(q)} \sum_{\substack{n = p_1 \cdots p_5 \in \mathcal{B} \\ x^{\kappa} < p_5 < \cdots < p_1 }} \mu^2(n) \left(\sum_{\substack{d \mid n \\ d < \sqrt{n} }} \mu(d) \right) \psi\left(n, x^{\kappa} \right) + \frac{1}{\varphi(q)} \left(\sum_{p \in \mathcal{B}} 1 - S_{2523} - S_{2524} \right) + O\left(\frac{y}{q (\log x)^{3}}\right)
\end{align}
for almost all $q \sim Q$. Combining (40) and (41), for almost all $q \sim Q$ we get
\begin{align}
\nonumber \sum_{p \in \mathcal{A}} 1 \geqslant&\ S_{2521} + S_{2522} + \frac{1}{\varphi(q)} \left(\sum_{p \in \mathcal{B}} 1 - S_{2523} - S_{2524} - S_{2525} - S_{2526} - S_{2527} \right) + O\left(\frac{y}{q (\log x)^{3}}\right) \\
\nonumber =&\ \frac{1}{\varphi(q)} \sum_{p \in \mathcal{B}} 1 + \left(S_{2521} - \frac{1}{\varphi(q)} S_{2523} \right) + \left(S_{2522} - \frac{1}{\varphi(q)} S_{2524} \right) \\
& - \frac{1}{\varphi(q)} \left(S_{2525} + S_{2526} + S_{2527} \right) + O\left(\frac{y}{q (\log x)^{3}}\right).
\end{align}
Note that by Prime Number Theorem in short intervals and partial summation, we can calculate the loss from $S_{2525}$, $S_{2526}$ and $S_{2527}$:
\begin{align}
S_{2525} =&\ 2 (1+o(1)) \frac{y}{\log x} \left( \int_{(t_1, t_2, t_3) \in U_{2525}} \frac{1}{t_1 t_2 t_3 (1 - t_1 - t_2 - t_3)} d t_3 d t_2 d t_1 \right), \\
S_{2526} =&\ 2 (1+o(1)) \frac{y}{\log x} \left( \int_{(t_1, t_2, t_3, t_4) \in U_{2526}} \frac{1}{t_1 t_2 t_3 t_4 (1 - t_1 - t_2 - t_3 - t_4)} d t_4 d t_3 d t_2 d t_1 \right), \\
S_{2527} =&\ 20 (1+o(1)) \frac{y}{\log x} \left( \int_{(t_1, t_2, t_3, t_4, t_5) \in U_{2527}} \frac{1}{t_1 t_2 t_3 t_4 t_5 (1 - t_1 - t_2 - t_3 - t_4 - t_5)} d t_5 d t_4 d t_3 d t_2 d t_1 \right),
\end{align}
where
\begin{align}
\nonumber U_{2525}(\boldsymbol{\alpha}_{3}) :=&\ \left\{ \kappa < \alpha_3 < \alpha_2 < \alpha_1,\ \alpha_1 + \alpha_2 + \alpha_3 > \frac{1}{2},\ 2 \alpha_1 + 2 \alpha_2 + \alpha_3 < 1 \right\}, \\
\nonumber U_{2526}(\boldsymbol{\alpha}_{4}) :=&\ \left\{ \kappa < \alpha_4 < \alpha_3 < \alpha_2 < \alpha_1,\ 2 \alpha_1 + 2 \alpha_2 + \alpha_3 < 1 \right\}, \\
\nonumber U_{2527}(\boldsymbol{\alpha}_{5}) :=&\ \left\{ \kappa < \alpha_5 < \alpha_4 < \alpha_3 < \alpha_2 < \alpha_1,\ 2 \alpha_1 + \alpha_2 + \alpha_3 + \alpha_4 + \alpha_5 < 1 \right\}.
\end{align}

Now we only need to give lower bounds for
\begin{equation}
\sum_{\substack{p_1 p_2 p_3 \in \mathcal{A} \\ x^{\kappa} < p_3 < p_2 < p_1 < p_2 p_3 }} 2 - \frac{1}{\varphi(q)} \sum_{\substack{p_1 p_2 p_3 \in \mathcal{B} \\ x^{\kappa} < p_3 < p_2 < p_1 < p_2 p_3 }} 2
\end{equation}
and
\begin{equation}
\sum_{\substack{p_1 \cdots p_7 \in \mathcal{A} \\ x^{\kappa} < p_7 < \cdots < p_1 }} 20 - \frac{1}{\varphi(q)} \sum_{\substack{p_1 \cdots p_7 \in \mathcal{B} \\ x^{\kappa} < p_7 < \cdots < p_1 }} 20.
\end{equation}
Here, Lemmas~\ref{l25}--\ref{l28} and \ref{l216} (see the proof of [\cite{MaynardLargeModuliI}, Lemma 8.13]) are applicable for parts of (46) and (47). After removing cross conditions between prime variables (see [\cite{Mikawa}, Section 9]), we can use those lemmas to show parts of (46) and (47) have a negligible contribution for almost all $q \sim Q$. We discard the remaining parts of $S_{2521}$ and $S_{2522}$ since they are positive. Note that when $\theta < \frac{29}{56} - \varepsilon$, we have $\kappa \geqslant \frac{1}{7}$ and (47) equals zero. The loss from (46) and (47) is
\begin{align}
&\ 2 (1+o(1)) \frac{y}{\log x} \left( \int_{(t_1, t_2) \in U_{2523}} \frac{1}{t_1 t_2 (1 - t_1 - t_2)} d t_2 d t_1 \right) \\
+&\ 20 (1+o(1)) \frac{y}{\log x} \left( \int_{(t_1, t_2, t_3, t_4, t_5, t_6) \in U_{2524}} \frac{1}{t_1 t_2 t_3 t_4 t_5 t_6 (1 - t_1 - t_2 - t_3 - t_4 - t_5 - t_6)} d t_6 d t_5 d t_4 d t_3 d t_2 d t_1 \right),
\end{align}
where
\begin{align}
\nonumber U_{2523}(\boldsymbol{\alpha}_{2}) :=&\ \left\{ \kappa < \alpha_2 < \alpha_1,\ \alpha_1 + \alpha_2 > \frac{1}{2},\ 2 \alpha_1 + \alpha_2 < 1,\ \boldsymbol{\alpha}_2 \notin \boldsymbol{G}_2 \right\}, \\
\nonumber U_{2524}(\boldsymbol{\alpha}_{6}) :=&\ \left\{ \kappa < \alpha_6 < \alpha_5 < \alpha_4 < \alpha_3 < \alpha_2 < \alpha_1,\ 2 \alpha_1 + \alpha_2 + \alpha_3 + \alpha_4 + \alpha_5 + \alpha_6 < 1,\ \boldsymbol{\alpha}_6 \notin \boldsymbol{G}_6 \right\}.
\end{align}

Finally, we can prove lower bounds for $C_0^{*}(\theta)$ for almost all $q \sim Q$ in the range $\theta \in \left[\frac{1}{2},\ \frac{17}{32} - \varepsilon \right]$ under our assumptions (34) and (38) by subtracting the values of the 5 integrals in (43)--(45) and (48)--(49) from $1$ and summing over $\asymp (\log x)^5$ subintervals. The proof of (34) and (38) for the case $\theta = \frac{17}{32} - \varepsilon$ was given in [\cite{Mikawa}, Propositions 1 and 2]. In the proof,
\begin{equation}
\sum_{n \in \mathcal{A}} \left(\sum_{\substack{d \mid n \\ d < \sqrt{x} }} \mu(d) \right) \psi\left(n, x^{\kappa} \right) \quad \text{and} \quad \sum_{n \in \mathcal{A}} \left(\sum_{\substack{d \mid n \\ \mu^2(d) = 1 \\ \Omega(d) = 3 \\ \sqrt{x} < d < \sqrt{x P^{-}(d)} }} 2 \right) \psi\left(n, x^{\kappa} \right)
\end{equation}
was decomposed into Type-I and Type-II sums, where the coefficients are convolutions of the form, say,
$$
\psi * \psi * \psi * \Psi * \Psi
$$
and similar convolutions. After that, Lemma~\ref{l29} and Lemma~\ref{l23} are applied to estimate the Type-I case and the Type-II case. Note that Lemma~\ref{l29} does not require any non-fixed coefficient to satisfy \textbf{Condition A} (the Siegel--Walfisz condition). This fact is important in this part of the next section.

For $\theta < \frac{17}{32} - \varepsilon$, similar decompositions of (28) are still valid, and we can prove (34) and (38) using the same method (the readers can follow the proof in \cite{Mikawa} and check all the conditions we need, especially some upper bounds). Note that (34) and (38) are only applicable for Theorem~\ref{t22} but not Theorem~\ref{t21}. After decompositions, we get
\begin{center}
\begin{tabular}{|c|c|c|c|}
\hline \boldmath{$\theta$} & \boldmath{$C_0^{*}(\theta)$} & \boldmath{$\theta$} & \boldmath{$C_0^{*}(\theta)$} \\
\hline $0.501$ & $0.8025$ & $0.517$ & $0.6033$ \\
\hline $0.502$ & $0.8003$ & $0.518$ & $0.5851$ \\
\hline $0.503$ & $0.7981$ & $0.519$ & $0.5689$ \\
\hline $0.504$ & $0.7859$ & $0.520$ & $0.5487$ \\
\hline $0.505$ & $0.7699$ & $0.521$ & $0.5323$ \\
\hline $0.506$ & $0.7579$ & $0.522$ & $0.5139$ \\
\hline $0.507$ & $0.7439$ & $0.523$ & $0.4919$ \\
\hline $0.508$ & $0.7339$ & $0.524$ & $0.4699$ \\
\hline $0.509$ & $0.7199$ & $0.525$ & $0.4499$ \\
\hline $0.510$ & $0.7079$ & $0.526$ & $0.4299$ \\
\hline $0.511$ & $0.6919$ & $0.527$ & $0.4079$ \\
\hline $0.512$ & $0.6838$ & $0.528$ & $0.3839$ \\
\hline $0.513$ & $0.6658$ & $0.529$ & $0.3579$ \\
\hline $0.514$ & $0.6517$ & $0.530$ & $0.3339$ \\
\hline $0.515$ & $0.6355$ & $0.531$ & $0.3099$ \\
\hline $0.516$ & $0.6215$ &  & \\
\hline
\end{tabular} \\
\textbf{Table 2.6: Lower Bounds for }\boldmath{$C_0^{*}(\theta)$} \textbf{(Second Method, }\boldmath{$\frac{1}{2} < \theta \leqslant \frac{17}{32} - \varepsilon$}\textbf{)}
\end{center}
Note that the lower bound does not become trivial even when $\theta = \frac{17}{32} - \varepsilon$ using this method. However, this method collapses when $\theta \geqslant \frac{17}{32}$ and we cannot get any nontrivial result using this method. There is no ``grey area'' between a lower bound $>0.3$ and no result at all.
\begin{center}
\begin{tabular}{|c|c|c|c|c|c|}
\hline \boldmath{$\theta$} & \boldmath{$C_0^{*}(\theta)$} \textbf{ (Method 1)} & \boldmath{$C_0^{*}(\theta)$} \textbf{ (Method 2)} & \boldmath{$\theta$} & \boldmath{$C_0^{*}(\theta)$} \textbf{ (Method 1)} & \boldmath{$C_0^{*}(\theta)$} \textbf{ (Method 2)} \\
\hline $0.501$ & \boldmath{$0.8636$} & $0.8025$ & $0.517$ & $0.4455$ & \boldmath{$0.6033$} \\
\hline $0.502$ & \boldmath{$0.8455$} & $0.8003$ & $0.518$ & $0.4086$ & \boldmath{$0.5851$} \\
\hline $0.503$ & \boldmath{$0.8298$} & $0.7981$ & $0.519$ & $0.3744$ & \boldmath{$0.5689$} \\
\hline $0.504$ & \boldmath{$0.8135$} & $0.7859$ & $0.520$ & $0.3363$ & \boldmath{$0.5487$} \\
\hline $0.505$ & \boldmath{$0.7972$} & $0.7699$ & $0.521$ & $0.2973$ & \boldmath{$0.5323$} \\
\hline $0.506$ & $0.7299$ & \boldmath{$0.7579$} & $0.522$ & $0.2546$ & \boldmath{$0.5139$} \\
\hline $0.507$ & $0.7082$ & \boldmath{$0.7439$} & $0.523$ & $0.2148$ & \boldmath{$0.4919$} \\
\hline $0.508$ & $0.6894$ & \boldmath{$0.7339$} & $0.524$ & $0.1706$ & \boldmath{$0.4699$} \\
\hline $0.509$ & $0.6606$ & \boldmath{$0.7199$} & $0.525$ & $0.1218$ & \boldmath{$0.4499$} \\
\hline $0.510$ & $0.6369$ & \boldmath{$0.7079$} & $0.526$ & $0.0778$ & \boldmath{$0.4299$} \\
\hline $0.511$ & $0.6128$ & \boldmath{$0.6919$} & $0.527$ & $0.0267$ & \boldmath{$0.4079$} \\
\hline $0.512$ & $0.5907$ & \boldmath{$0.6838$} & $0.528$ & $0$ & \boldmath{$0.3839$} \\
\hline $0.513$ & $0.5567$ & \boldmath{$0.6658$} & $0.529$ & $0$ & \boldmath{$0.3579$} \\
\hline $0.514$ & $0.5332$ & \boldmath{$0.6517$} & $0.530$ & $0$ & \boldmath{$0.3339$} \\
\hline $0.515$ & $0.5040$ & \boldmath{$0.6355$} & $0.531$ & $0$ & \boldmath{$0.3099$} \\
\hline $0.516$ & $0.4762$ & \boldmath{$0.6215$} & & & \\
\hline
\end{tabular} \\
\textbf{Table 2.7: A Comparison of Two Methods on the Lower Bounds for }\boldmath{$C_0^{*}(\theta)$} \textbf{(}\boldmath{$\frac{1}{2} < \theta \leqslant \frac{17}{32} - \varepsilon$}\textbf{)}
\end{center}

In the end of this section, we mention an improvement over Lemma~\ref{l225} on the case $\Omega(n) = 6$. In this case, Lemma~\ref{l225} gives that
\begin{equation}
\sum_{\substack{\Omega(n) = 6 \\ \mu^2(n) = 1 \\ d \mid n \\ \mu^2(d) = 1 \\ \Omega(d) = 3 \\ \sqrt{n} < d < \sqrt{n P^{-}(d)} }} 1 \leqslant 10.
\end{equation}
We shall prove that the upper bound $10$ can be reduced under some conditions on the prime factors of $d$. Let $n = p_1 p_2 p_3 p_4 p_5 p_6$ with $p_1 > p_2 > p_3 > p_4 > p_5 > p_6$. The upper bound $10$ can be easily obtained: Let $d = p_i p_j p_k$ with $1 \leqslant i < j < k \leqslant 6$, we have $\binom{6}{3} = 20$ choices for $(i, j, k)$. Since only one of choices $d$ and $\frac{n}{d}$ is larger than $\sqrt{n}$, the number of $d$ counted is no more than $\frac{1}{2}\binom{6}{3} = 10$.

However, this upper bound ignores another restriction $d < \sqrt{n P^{-}(d)} = \sqrt{n p_k}$ in the sum. Taking this condition into our consideration, we can show that some of the $10$ possible combinations of $d$ are not acceptable when $\boldsymbol{\alpha}_6$ lies in some special regions. We give a table in the end of this section to reveal one possible $\boldsymbol{\alpha}_6$ for each $(i, j, k)$ such that $d = p_i p_j p_k$ will not be counted in the sum. Trivially, we know that the following 5 choices of $(i, j, k)$ are impossible for $d$ since $p_i p_j p_k < \frac{n}{p_i p_j p_k}$:
$$
(2, 4, 6),\ (2, 5, 6),\ (3, 4, 6),\ (3, 5, 6),\ (4, 5, 6).
$$
Thus, we only need to consider the remaining 15 choices.
\begin{center}
\begin{tabular}{|c|c|}
\hline $\boldsymbol{(i, j, k)}$ & $\boldsymbol{\alpha}_6$ \\
\hline $(1, 2, 3)$ & $(0.295, 0.143, 0.142, 0.141, 0.140, 0.139)$ \\
\hline $(1, 2, 4)$ & $(0.295, 0.143, 0.142, 0.141, 0.140, 0.139)$ \\
\hline $(1, 2, 5)$ & $(0.295, 0.143, 0.142, 0.141, 0.140, 0.139)$ \\
\hline $(1, 2, 6)$ & $(0.295, 0.143, 0.142, 0.141, 0.140, 0.139)$ \\
\hline $(1, 3, 4)$ & $(0.295, 0.143, 0.142, 0.141, 0.140, 0.139)$ \\
\hline $(1, 3, 5)$ & $(0.295, 0.143, 0.142, 0.141, 0.140, 0.139)$ \\
\hline $(1, 3, 6)$ & $(0.295, 0.143, 0.142, 0.141, 0.140, 0.139)$ \\
\hline $(1, 4, 5)$ & $(0.295, 0.143, 0.142, 0.141, 0.140, 0.139)$ \\
\hline $(1, 4, 6)$ & $(0.295, 0.143, 0.142, 0.141, 0.140, 0.139)$ \\
\hline $(1, 5, 6)$ & $(0.295, 0.143, 0.142, 0.141, 0.140, 0.139)$ \\
\hline $(2, 3, 4)$ & Do Not Exist \\
\hline $(2, 3, 5)$ & Do Not Exist \\
\hline $(2, 3, 6)$ & Do Not Exist \\
\hline $(2, 4, 5)$ & Do Not Exist \\
\hline $(3, 4, 5)$ & Do Not Exist \\
\hline
\end{tabular} \\
\textbf{Table 2.8: Examples of }\boldmath{$d = p_i p_j p_k$}\textbf{ such that }\boldmath{$d > \sqrt{n P^{-}(d)}$}
\end{center}
By this table, we know that $d$ will always be counted in the sum if $i \neq 1$. For a comparison, we also give a table that shows one possible $\boldsymbol{\alpha}_6$ for each $(i, j, k)$ such that $d = p_i p_j p_k$ will be counted in the sum.
\begin{center}
\begin{tabular}{|c|c|}
\hline $\boldsymbol{(i, j, k)}$ & $\boldsymbol{\alpha}_6$ \\
\hline $(1, 2, 3)$ & $(0.279, 0.147, 0.146, 0.145, 0.143, 0.140)$ \\
\hline $(1, 2, 4)$ & $(0.280, 0.147, 0.146, 0.145, 0.143, 0.139)$ \\
\hline $(1, 2, 5)$ & $(0.283, 0.146, 0.145, 0.143, 0.142, 0.141)$ \\
\hline $(1, 2, 6)$ & $(0.284, 0.146, 0.145, 0.143, 0.142, 0.140)$ \\
\hline $(1, 3, 4)$ & $(0.287, 0.151, 0.142, 0.141, 0.140, 0.139)$ \\
\hline $(1, 3, 5)$ & $(0.287, 0.151, 0.142, 0.141, 0.140, 0.139)$ \\
\hline $(1, 3, 6)$ & $(0.287, 0.151, 0.142, 0.141, 0.140, 0.139)$ \\
\hline $(1, 4, 5)$ & $(0.288, 0.150, 0.142, 0.141, 0.140, 0.139)$ \\
\hline $(1, 4, 6)$ & $(0.288, 0.150, 0.142, 0.141, 0.140, 0.139)$ \\
\hline $(1, 5, 6)$ & $(0.289, 0.149, 0.142, 0.141, 0.140, 0.139)$ \\
\hline $(2, 3, 4)$ & $(0.220, 0.217, 0.143, 0.141, 0.140, 0.139)$ \\
\hline $(2, 3, 5)$ & $(0.219, 0.217, 0.144, 0.141, 0.140, 0.139)$ \\
\hline $(2, 3, 6)$ & $(0.217, 0.216, 0.147, 0.141, 0.140, 0.139)$ \\
\hline $(2, 4, 5)$ & $(0.214, 0.213, 0.146, 0.145, 0.143, 0.139)$ \\
\hline $(3, 4, 5)$ & $(0.190, 0.170, 0.169, 0.168, 0.164, 0.139)$ \\
\hline
\end{tabular} \\
\textbf{Table 2.9: Examples of }\boldmath{$d = p_i p_j p_k$}\textbf{ such that }\boldmath{$\sqrt{n} < d < \sqrt{n P^{-}(d)}$}
\end{center}

\section{$2$-factored Moduli, 1}
In this section we focus on the first $2$-factored case, where the moduli $q = q_1 q_2$ with $q_1 \sim Q_1 = x^{\theta_1}$ and $q_2 \sim Q_2 = x^{\theta_2}$. By the definitions of the sieved sets $\mathcal{A}^{q}$, $\mathcal{B}^{q}$ and the sieve function $S\left(\mathcal{C}, z\right)$, using Prime Number Theorem, we have
\begin{equation}
\pi(2 x; q_1 q_2, a) - \pi(x; q_1 q_2, a) = \sum_{p \in \mathcal{A}^{q_1 q_2}} 1 = S\left(\mathcal{A}^{q_1 q_2}, (2x)^{\frac{1}{2}} \right) \quad \text{and} \quad S\left(\mathcal{B}^{q_1 q_2}, (2x)^{\frac{1}{2}} \right) = (1+o(1)) \frac{x}{\log x}.
\end{equation}
Our aim is again to show that the sparser set $\mathcal{A}^{q_1 q_2}$ contains the expected proportion of primes compared to the larger set $\mathcal{B}^{q_1 q_2}$, which requires us to decompose $S\left(\mathcal{A}^{q_1 q_2}, (2x)^{\frac{1}{2}} \right)$ and prove ``asymptotic formulas'' of the form
\begin{equation}
\sum_{\substack{q_1 \sim Q_1 \\ q_2 \sim Q_2 \\ (q_1 q_2, a) = 1}} \left| S\left(\mathcal{A}^{q_1 q_2}, z \right) - \frac{1}{\varphi(q_1 q_2)} S\left(\mathcal{B}^{q_1 q_2}, z \right) \right| \ll \frac{x}{(\log x)^A}.
\end{equation}
for some parts of it, and drop the remaining parts to construct a suitable majorant or minorant. For the majorant case we can only drop negative parts, while for the minorant case we can only drop positive parts. After the final decompositions, we can get the following result with some $0 < C_0(\theta_1, \theta_2) \leqslant 1$ and $C_1(\theta_1, \theta_2) \geqslant 1$:
\begin{theorem}\label{t31}
There exist functions $\rho_0$ and $\rho_1$ which satisfies the following properties:

(Majorant / Minorant). $\rho_0(n)$ is a minorant for the prime indicator function $\mathbbm{1}_{p}(n)$, and $\rho_1(n)$ is a majorant for the prime indicator function $\mathbbm{1}_{p}(n)$. That is, we have
$$
\rho_0(n) \leqslant \mathbbm{1}_{p}(n) \leqslant \rho_1(n).
$$

(Upper and Lower bounds). We have
$$
\sum_{n \leqslant x} \rho_0(n) \geqslant (1+o(1))\frac{C_0(\theta_1, \theta_2) x}{\log x} \quad \text{and} \quad \sum_{n \leqslant x} \rho_1(n) \leqslant (1+o(1))\frac{C_1(\theta_1, \theta_2) x}{\log x}
$$
for two functions $C_0(\theta_1, \theta_2)$ and $C_1(\theta_1, \theta_2)$ satisfy $0 < C_0(\theta_1, \theta_2) \leqslant 1$ and $C_1(\theta_1, \theta_2) \geqslant 1$.

(Distributions in Arithmetic Progressions). For any $a \in \mathbb{Z} \backslash \{0\}$ and any $A>0$, we have
$$
\sum_{\substack{q_1 \sim Q_1 \\ q_2 \sim Q_2 \\ (q_1 q_2, a) = 1}} \left| \sum_{\substack{n \leqslant x \\ n \equiv a (\bmod q_1 q_2)}} \rho_j(n) - \frac{1}{\varphi(q_1 q_2)} \sum_{\substack{n \leqslant x \\ (n, q_1 q_2) = 1}} \rho_j(n) \right| \ll \frac{x}{(\log x)^A}
$$
for $j = 0, 1$.
\end{theorem}

The ``asymptotic formula'' (53) yields the following asymptotic formula for almost all $q \sim Q$ such that $q_1 \mid q$ and $q_1 \sim Q_1$:
\begin{equation}
S\left(\mathcal{A}^{q_1 q_2}, z \right) = (1+o(1)) \frac{1}{\varphi(q_1 q_2)} S\left(\mathcal{B}^{q_1 q_2}, z \right).
\end{equation}
By a similar decomposing process, one can deduce the following result using asymptotic formulas of the form (54):
\begin{theorem}\label{t32}
For almost all  $q \sim Q$ such that $q_1 \mid q$ and $q_1 \sim Q_1$, we have
$$
\frac{C_0^{*}(\theta_1, \theta_2) x}{\varphi(q) \log x} \leqslant \pi(x; q, a) \leqslant \frac{C_1^{*}(\theta_1, \theta_2) x}{\varphi(q) \log x}
$$
for two functions $C_0^{*}(\theta_1, \theta_2)$ and $C_1^{*}(\theta_1, \theta_2)$ satisfy $0 < C_0(\theta_1, \theta_2) \leqslant C_0^{*}(\theta_1, \theta_2) \leqslant 1$ and $1 \leqslant C_1^{*}(\theta_1, \theta_2) \leqslant C_1(\theta_1, \theta_2)$.
\end{theorem}

In order to give asymptotic formulas (53) and (54) for sieve functions $S\left(\mathcal{A}^{q_1 q_2}, z \right)$, we need results of the form
\begin{equation}
\sum_{\substack{q_1 \sim Q_1 \\ q_2 \sim Q_2 \\ (q_1 q_2, a) = 1}} \left| \sum_{\substack{n \sim x \\ n \equiv a (\bmod q_1 q_2)}} f(n) - \frac{1}{\varphi(q_1 q_2)} \sum_{\substack{n \sim x \\ (n, q_1 q_2) = 1}} f(n) \right| \ll \frac{x}{(\log x)^A}.
\end{equation}
As in Section 2, we may want the coefficients to satisfy \textbf{Conditions A and B}.

\subsection{Preliminary Lemmas}
Before constructing the majorant and minorant, we need estimate results of the form (55). Note that the results from Section 2 are still applicable in the final decomposition, and the results here are still useful in the next two sections.

\subsubsection{Type-II estimate}
All of the following three lemmas come from \cite{MaynardLargeModuliI}. The first one helps us to increase the value of $\kappa$ defined in Section 2.
\begin{lemma}\label{l33} ([\cite{MaynardLargeModuliI}, Proposition 12.1]).
Let $M_1 M_2 \asymp x$. Let $a_{1, m_1}$ and $a_{2, m_2}$ be divisor-bounded complex sequences. Suppose that $a_{2, m_2}$ satisfies \textbf{Conditions A and B}. If we have
$$
Q_2 x^{\varepsilon} < M_2,\ Q_1^2 Q_2^3 M_2^6 < x^{2 - 15 \varepsilon},\ Q_1^3 Q_2^3 M_2^3 < x^{2 - 15 \varepsilon},\ Q_1^4 Q_2^3 M_2^3 < x^{\frac{5}{2} - 15 \varepsilon},\ Q_1^2 Q_2^2 < M_2 x^{1 - 4 \varepsilon},
$$
then
$$
\sum_{\substack{q_1 \sim Q_1 \\ q_2 \sim Q_2 \\ (q_1 q_2, a) = 1}} \left| \sum_{\substack{m_1 \sim M_1 \\ m_2 \sim M_2 \\ m_1 m_2 \equiv a (\bmod q_1 q_2)}} a_{1, m_1} a_{2, m_2} - \frac{1}{\varphi(q_1 q_2)} \sum_{\substack{m_1 \sim M_1 \\ m_2 \sim M_2 \\ (m_1 m_2, q_1 q_2) = 1}} a_{1, m_1} a_{2, m_2} \right| \ll \frac{x}{(\log x)^A}.
$$
\end{lemma}

The next one provides a ``small factor'' Type-II estimate that can be used to replace $\boldsymbol{U}_{j}$ with $\boldsymbol{S}_{j}$ in various conditions.
\begin{lemma}\label{l34} ([\cite{MaynardLargeModuliI}, Proposition 12.2]).
Let $M_1 M_2 \asymp x$. Let $a_{1, m_1}$ and $a_{2, m_2}$ be divisor-bounded complex sequences. Suppose that $a_{2, m_2}$ satisfies \textbf{Conditions A and B}. If we have
$$
Q_1 Q_2^2 < M_2 x^{1 - 7 \varepsilon},\ Q_1^8 Q_2^7 M_2^6 < x^{4 - 13 \varepsilon},
$$
then
$$
\sum_{\substack{q_1 \sim Q_1 \\ q_2 \sim Q_2 \\ (q_1 q_2, a) = 1}} \left| \sum_{\substack{m_1 \sim M_1 \\ m_2 \sim M_2 \\ m_1 m_2 \equiv a (\bmod q_1 q_2)}} a_{1, m_1} a_{2, m_2} - \frac{1}{\varphi(q_1 q_2)} \sum_{\substack{m_1 \sim M_1 \\ m_2 \sim M_2 \\ (m_1 m_2, q_1 q_2) = 1}} a_{1, m_1} a_{2, m_2} \right| \ll \frac{x}{(\log x)^A}.
$$
\end{lemma}

The last one is very important since it handles a Type-II sum with a product of variables lies in $\left[1 - \theta,\ \theta \right]$ in some cases, which is crucial in obtaining nontrivial lower bounds for $C_0(\theta_1, \theta_2)$.
\begin{lemma}\label{l35} ([\cite{MaynardLargeModuliI}, Proposition 8.2]).
Let $M_1 M_2 \asymp x$. Let $a_{1, m_1}$ and $a_{2, m_2}$ be divisor-bounded complex sequences. Suppose that $a_{2, m_2}$ satisfies \textbf{Conditions A and B}. If we have
$$
Q_1^7 Q_2^{12} < x^{4 - 19 \varepsilon},\ Q_1 x^{\varepsilon} < M_2 < Q_1^{-1} x^{1 - 6 \varepsilon},
$$
then
$$
\sum_{\substack{q_1 \sim Q_1 \\ q_2 \sim Q_2 \\ (q_1 q_2, a) = 1}} \left| \sum_{\substack{m_1 \sim M_1 \\ m_2 \sim M_2 \\ m_1 m_2 \equiv a (\bmod q_1 q_2)}} a_{1, m_1} a_{2, m_2} - \frac{1}{\varphi(q_1 q_2)} \sum_{\substack{m_1 \sim M_1 \\ m_2 \sim M_2 \\ (m_1 m_2, q_1 q_2) = 1}} a_{1, m_1} a_{2, m_2} \right| \ll \frac{x}{(\log x)^A}.
$$
\end{lemma}

Other results from \cite{MaynardLargeModuliIII}, such as [\cite{MaynardLargeModuliIII}, Proposition 5.1] and [\cite{MaynardLargeModuliIII}, Proposition 5.3], can also be used here; however, we decide not to use them since they are only applicable for $Q_1 Q_2 < x^{0.501}$.

\subsubsection{Type-I$_3$ estimate}
Now we provide some estimates for the triple divisor function proved in \cite{MaynardLargeModuliI} and \cite{Lichtman2}, which will be useful when dealing with sieve functions that count products of three large variables. These can be seen as variants of the three-dimensional Harman's sieve in Section 2. Note that we will use them in both cases $C_1(\theta_1, \theta_2)$ and $C_1^{*}(\theta_1, \theta_2)$ here, which is different from the case in Section 2 where we can only use the three-dimensional Harman's sieve in the case $C_1^{*}(\theta)$.
\begin{lemma}\label{l36} ([\cite{MaynardLargeModuliI}, Lemma 20.7]).
Let $x^{2 \varepsilon} \leqslant M_3 \leqslant M_2 \leqslant M_1$, $x^{\varepsilon} \leqslant M_0$ and $M_0 M_1 M_2 M_3 \asymp x$. Let $a_{m_0}$ be a divisor-bounded complex sequence, $z = \exp\left(\log x (\log \log x)^{-3} \right)$. Let $\mathbf{M}_1$, $\mathbf{M}_2$, $\mathbf{M}_3$ be intervals such that $\mathbf{M}_i \subseteq [M_i, 2 M_i]$. If we have
$$
Q_1 Q_2 \leqslant x^{1 - \varepsilon},\ \frac{M_0 Q_1^{\frac{5}{2}} Q_2^3}{x^{1 - 15 \varepsilon}} \leqslant M_1 \leqslant \frac{x^{2 - \varepsilon}}{Q_1^3 Q_2^2 M_0},
$$
then
$$
\sum_{\substack{q_1 \sim Q_1 \\ q_2 \sim Q_2 \\ (q_1 q_2, a) = 1}} \left| \sum_{\substack{m_0 \sim M_0 \\ m_1 \in \mathbf{M}_1 \\ m_2 \in \mathbf{M}_2 \\ m_3 \in \mathbf{M}_3 \\ \left(m_1 m_2 m_3, P(z)\right) = 1 \\ m_0 m_1 m_2 m_3 \equiv a (\bmod q_1 q_2)}} a_{m_0} - \frac{1}{\varphi(q_1 q_2)} \sum_{\substack{m_0 \sim M_0 \\ m_1 \in \mathbf{M}_1 \\ m_2 \in \mathbf{M}_2 \\ m_3 \in \mathbf{M}_3 \\ \left(m_1 m_2 m_3, P(z)\right) = 1 \\ (m_0 m_1 m_2 m_3, q_1 q_2) = 1}} a_{m_0} \right| \ll \frac{x}{(\log x)^A}.
$$
\end{lemma}

\begin{lemma}\label{l37} ([\cite{MaynardLargeModuliI}, Proposition 11.1]).
Let $x^{\varepsilon} \leqslant M_3 \leqslant M_2 \leqslant M_1 \leqslant x^{\frac{3}{7} + \varepsilon}$, $x^{\varepsilon} \leqslant M_0$ and $M_0 M_1 M_2 M_3 \asymp x$. Let $a_{m_0}$ be a divisor-bounded complex sequence, $z = \exp\left(\log x (\log \log x)^{-3} \right)$. Let $\mathbf{M}_1$, $\mathbf{M}_2$, $\mathbf{M}_3$ be intervals such that $\mathbf{M}_i \subseteq [M_i, 2 M_i]$. If we have
$$
Q_1^7 Q_2^9 < x^4,\ Q_1^9 Q_2 < x^{\frac{32}{7}},\ Q_1^{\frac{15}{8}} Q_2^{\frac{15}{8}} M_0 < x^{1 - 20 \varepsilon},
$$
then
$$
\sum_{\substack{q_1 \sim Q_1 \\ q_2 \sim Q_2 \\ (q_1 q_2, a) = 1}} \left| \sum_{\substack{m_0 \sim M_0 \\ m_1 \in \mathbf{M}_1 \\ m_2 \in \mathbf{M}_2 \\ m_3 \in \mathbf{M}_3 \\ \left(m_1 m_2 m_3, P(z)\right) = 1 \\ m_0 m_1 m_2 m_3 \equiv a (\bmod q_1 q_2)}} a_{m_0} - \frac{1}{\varphi(q_1 q_2)} \sum_{\substack{m_0 \sim M_0 \\ m_1 \in \mathbf{M}_1 \\ m_2 \in \mathbf{M}_2 \\ m_3 \in \mathbf{M}_3 \\ \left(m_1 m_2 m_3, P(z)\right) = 1 \\ (m_0 m_1 m_2 m_3, q_1 q_2) = 1}} a_{m_0} \right| \ll \frac{x}{(\log x)^A}.
$$
\end{lemma}

\begin{lemma}\label{l38} ([\cite{Lichtman2}, Proposition 12.2]).
Let $x^{\varepsilon} \leqslant M_3 \leqslant M_2 \leqslant M_1 \leqslant x^{\frac{3}{7} + \varepsilon}$, $M_0 = x^{\varepsilon}$ and $M_0 M_1 M_2 M_3 \asymp x$. Let $a_{m_0}$ be a divisor-bounded complex sequence, $z = \exp\left(\log x (\log \log x)^{-3} \right)$. Let $\mathbf{M}_1$, $\mathbf{M}_2$, $\mathbf{M}_3$ be intervals such that $\mathbf{M}_i \subseteq [M_i, 2 M_i]$. If we have
$$
Q_1^3 Q_2^2 < x^{\frac{11}{7} - 30 \varepsilon},\ Q_1^{11} Q_2^{12} < x^{6 - 30 \varepsilon},\ Q_1 Q_2 < x^{\frac{8}{15} - 30 \varepsilon},
$$
then
$$
\sum_{\substack{q_1 \sim Q_1 \\ q_2 \sim Q_2 \\ (q_1 q_2, a) = 1}} \left| \sum_{\substack{m_0 \sim M_0 \\ m_1 \in \mathbf{M}_1 \\ m_2 \in \mathbf{M}_2 \\ m_3 \in \mathbf{M}_3 \\ \left(m_1 m_2 m_3, P(z)\right) = 1 \\ m_0 m_1 m_2 m_3 \equiv a (\bmod q_1 q_2)}} a_{m_0} - \frac{1}{\varphi(q_1 q_2)} \sum_{\substack{m_0 \sim M_0 \\ m_1 \in \mathbf{M}_1 \\ m_2 \in \mathbf{M}_2 \\ m_3 \in \mathbf{M}_3 \\ \left(m_1 m_2 m_3, P(z)\right) = 1 \\ (m_0 m_1 m_2 m_3, q_1 q_2) = 1}} a_{m_0} \right| \ll \frac{x}{(\log x)^A}.
$$
\end{lemma}

\subsection{Sieve Asymptotic Formulas}
In this section, many asymptotic formulas used in the decompositions will be adopted from Section 2. We shall also use the following powerful lemma, which comes from \cite{MaynardLargeModuliI} and gives asymptotic formulas of the form (53) for all sums that count numbers with $4$ or more prime factors, all larger than $x^{\frac{1}{7}}$.
\begin{lemma}\label{l39} ([\cite{MaynardLargeModuliI}, Proposition 7.3]).
Let $j \geqslant 4$, $P_1 P_2 \cdots P_j \asymp x$ and $P_1 \geqslant P_2 \geqslant \cdots \geqslant P_j \geqslant x^{\frac{1}{7} + 10 \varepsilon}$. Suppose that
$$
2 \theta_1 + \theta_2 < 1 - 100 \varepsilon \quad \text{and} \quad 7 \theta_1 + 12 \theta_2 < 4 - 100 \varepsilon.
$$
Then
$$
\sum_{\substack{p_1, \ldots, p_j \\ p_i \sim P_i,\ 1 \leqslant i \leqslant j \\ p_1 \cdots p_j \equiv a (\bmod q_1 q_2) }} 1
$$
has an asymptotic formula of the form (53).
\end{lemma}
We can use Lemma~\ref{l37} and Lemma~\ref{l38} to construct a new three-dimensional Harman's sieve similar to Lemma~\ref{l220} and Lemma~\ref{l221} but applicable for Theorem~\ref{t31}, and we will discuss the construction later.

\subsection{Upper Bounds}
We shall construct the majorant $\rho_1(n)$ and prove upper bounds for $C_1(\theta_1, \theta_2)$ and $C_1^{*}(\theta_1, \theta_2)$ in this subsection. Before our final decompositions, we first mention some existing results of $C_1(\theta_1, \theta_2)$ and $C_1^{*}(\theta_1, \theta_2)$.
\begin{theorem}\label{t310}
The functions $C_1(\theta_1, \theta_2)$ and $C_1^{*}(\theta_1, \theta_2)$ satisfy the following conditions:

(1). $C_1(\theta_1, \theta_2) = C_1(\theta_2, \theta_1)$, $C_1^{*}(\theta_1, \theta_2) = C_1^{*}(\theta_2, \theta_1)$;

(2). $C_1(\theta_1, \theta_2) = C_1^{*}(\theta_1, \theta_2) = 1$ for all $\theta_1, \theta_2$ satisfy $\theta_1 + \theta_2 \leqslant 0.5 - \varepsilon$;

(3). $C_1(\theta_1, \theta_2) = C_1^{*}(\theta_1, \theta_2) = 1$ for all $\theta_1, \theta_2$ satisfy $2 \theta_1 + \theta_2 \leqslant 1 - \varepsilon$, $7 \theta_1 + 12 \theta_2 \leqslant 4 - \varepsilon$ and $19 \theta_1 + 20 \theta_2 \leqslant 10 - \varepsilon$;

(4). $C_1(\theta_1, \theta_2) \leqslant C_1(\theta_1 + \theta_2)$, $C_1^{*}(\theta_1, \theta_2) \leqslant C_1^{*}(\theta_1 + \theta_2)$ for $0.5 \leqslant \theta_1 + \theta_2 \leqslant 1$;

(5). $C_1^{*}(\theta_1, \theta_2) \leqslant C_1(\theta_1, \theta_2) \leqslant 1 + \varepsilon$ for all $\theta_1, \theta_2$ satisfy $\theta_1 + \theta_2 = 0.5$;

(6). $C_1^{*}(\theta_1, \theta_2) \leqslant C_1(\theta_1, \theta_2) \leqslant 1 + \varepsilon$ for all $\theta_1, \theta_2$ satisfy $\theta_1 \leqslant 0.5$ and $\theta_2 = \min\left(1 - 2 \theta_1, \frac{4 - 7 \theta_1}{12}, \frac{10 - 19 \theta_1}{20}\right)$.
\end{theorem}
\begin{proof}
Statement (1) is obvious. Statements (2)--(3) follow easily from the Bombieri--Vinogradov Theorem and [\cite{MaynardLargeModuliI}, Theorem 1.1]. Statement (4) holds trivially by the work done in Section 2. When there are no new arithmetic information inputs outside of those in Section 2, we use $C_1(\theta_1 + \theta_2)$ and $C_1^{*}(\theta_1 + \theta_2)$ as upper bounds for $C_1(\theta_1, \theta_2)$ and $C_1^{*}(\theta_1, \theta_2)$ respectively. Statements (5)--(6) hold from our final decompositions in this Subsection and arguments in \cite{676}, with a ``loss'' of size $O(\varepsilon)$ and a fact that a ``three-dimensional Harman's sieve''(which will be discussed later) is applicable when $\theta_1 + \theta_2 = 0.5$ (using Lemma~\ref{l38}) or $(\theta_1, \theta_2)$ lies in the boundary of the region defined by [\cite{MaynardLargeModuliI}, Theorem 1.1].
\end{proof}

From here to the end of this section, we assume that $\theta_1 \geqslant \theta_2$ to simplify the conditions. We also write $\theta = \theta_1 + \theta_2$. Before performing our final decompositions, we define several regions of the pair $(\theta_1, \theta_2)$ based on various arithmetic information inputs.
\begin{align}
\nonumber \boldsymbol{U} =&\ \left\{ (\theta_1, \theta_2) : 0 \leqslant \theta_2 \leqslant \theta_1 < 1,\ \theta_1 + \theta_2 < 1 \right\}, \\
\nonumber \boldsymbol{I} =&\ \left\{ (\theta_1, \theta_2) : (\theta_1, \theta_2) \in \boldsymbol{U};\ \theta_1 + \theta_2 \leqslant \frac{1}{2} - \varepsilon \right. \\
\nonumber & \left. \qquad \qquad \quad \text{ or } 2 \theta_1 + \theta_2 \leqslant 1 - \varepsilon,\ 7 \theta_1 + 12 \theta_2 \leqslant 4 - \varepsilon,\ 19 \theta_1 + 20 \theta_2 \leqslant 10 - \varepsilon \right\}, \\
\nonumber \boldsymbol{T}_{1} =&\ \left\{ (\theta_1, \theta_2) : (\theta_1, \theta_2) \in \boldsymbol{U};\ 7 \theta_1 + 9 \theta_2 < 4,\ 9 \theta_1 + \theta_2 < \frac{32}{7} \right\}, \\
\nonumber \boldsymbol{T}_{2} =&\ \left\{ (\theta_1, \theta_2) : (\theta_1, \theta_2) \in \boldsymbol{U} \backslash \boldsymbol{T}_{1};\ 3 \theta_1 + 2 \theta_2 < \frac{11}{7} - \varepsilon,\ 11 \theta_1 + 12 \theta_2 < 6 - \varepsilon,\ \theta_1 + \theta_2 < \frac{8}{15} - \varepsilon \right\}, \\
\nonumber \boldsymbol{T} =&\ \left\{ (\theta_1, \theta_2) : (\theta_1, \theta_2) \in \boldsymbol{U};\ (\theta_1, \theta_2) \in \boldsymbol{T}_{1} \cup \boldsymbol{T}_{2} \right\}, \\
\nonumber \boldsymbol{A} =&\ \left\{ (\theta_1, \theta_2) : (\theta_1, \theta_2) \in \boldsymbol{U} \backslash \boldsymbol{I};\ \frac{5}{14} < \theta_1 \leqslant \frac{2}{5},\ \frac{1}{2}(1 - 2 \theta_1) < \theta_2 < \frac{1}{9}(2 - 2 \theta_1) \right. \\
\nonumber & \qquad \qquad \quad \text{ or } \frac{2}{5} < \theta_1 \leqslant \frac{4}{9},\ \frac{1}{2}(1 - 2 \theta_1) < \theta_2 < \frac{1}{6}(2 - 3 \theta_1) \\
\nonumber & \qquad \qquad \quad \text{ or } \frac{4}{9} < \theta_1 \leqslant \frac{1}{2},\ \frac{1}{2}(1 - 2 \theta_1) < \theta_2 < \frac{1}{9}(5 - 9 \theta_1) \\
\nonumber & \left. \qquad \qquad \quad \text{ or } \frac{1}{2} < \theta_1 < \frac{11}{20},\ 0 < \theta_2 < \frac{1}{18}(11 - 20 \theta_1) \right\}, \\
\nonumber \boldsymbol{B} =&\ \left\{ (\theta_1, \theta_2) : (\theta_1, \theta_2) \in \boldsymbol{U} \backslash \boldsymbol{I};\ \frac{1}{4} < \theta_1 \leqslant \frac{10}{33},\ \frac{1}{2}(1 - 2 \theta_1) < \theta_2 < \theta_1 \right. \\
\nonumber & \qquad \qquad \quad \text{ or } \frac{10}{33} < \theta_1 \leqslant \frac{1}{2},\ \frac{1}{2}(1 - 2 \theta_1) < \theta_2 < \frac{1}{19}(10 - 14 \theta_1) \\
\nonumber & \left. \qquad \qquad \quad \text{ or } \frac{1}{2} < \theta_1 < \frac{5}{7},\ 0 < \theta_2 < \frac{1}{19}(10 - 14 \theta_1) \right\}, \\
\nonumber \boldsymbol{C} =&\ \left\{ (\theta_1, \theta_2) : (\theta_1, \theta_2) \in \boldsymbol{U} \backslash \boldsymbol{I};\ \frac{2}{5} < \theta_1 < \frac{1}{2},\ \frac{1}{2}(1 - 2 \theta_1) < \theta_2 < \frac{1}{12}(4 - 7 \theta_1) \right\}.
\end{align}
Here, $\boldsymbol{U}$ denote all possible pairs $(\theta_1, \theta_2)$ in our problem, $\boldsymbol{I}$ denote the region that $C_1(\theta_1, \theta_2) = 1$ follows by the Bombieri--Vinogradov Theorem or [\cite{MaynardLargeModuliI}, Theorem 1.1], and $\boldsymbol{T}$ denote a new ``three-dimensional Harman's sieve'' region corresponding to Lemma~\ref{l37} and Lemma~\ref{l38}. Region $\boldsymbol{A}$ corresponds to Lemma~\ref{l33}, region $\boldsymbol{B}$ corresponds to Lemma~\ref{l34}, and region $\boldsymbol{C}$ corresponds to Lemma~\ref{l35}. Region $\boldsymbol{B}$ covers both region $\boldsymbol{A}$ and region $\boldsymbol{C}$.

Before our discussions on each region, we first give a result of the three-dimensional Harman's sieve. We shall implicitly use this result in many decompositions below. Assume that $\theta \leqslant \frac{17}{32} - \varepsilon$. Put
$$
\mathcal{F}^{q} =\left\{m_1 m_2 m_3 : m_1 m_2 m_3 \sim x,\ m_1 m_2 m_3 \equiv a (\bmod q),\ x^{\varepsilon} \leqslant m_3 \leqslant m_2 \leqslant m_1 \leqslant x^{\frac{3}{7} + \varepsilon} \right\}.
$$
We want to give an asymptotic formula of the form (53) for
\begin{equation}
S\left(\mathcal{F}^{q_1 q_2}, x^{\kappa} \right)
\end{equation}
with $\kappa = \frac{5 - 8 \theta}{6} - \varepsilon$ or $\frac{2 - 3 \theta}{3} - \varepsilon$ or some other values. In order to give an asymptotic formula of the form (53) for (56), we follow the steps in [\cite{MaynardLargeModuliI}, Chapter 11]. We apply [\cite{MaynardLargeModuliI}, Lemma 10.2] 3 times with $z_1 = x^{\kappa}$ and $z_2 = \exp\left(\log x (\log \log x)^{-3} \right)$ on each variable $m_1$, $m_2$ and $m_3$, and we use Lemma~\ref{l211}, Lemma~\ref{l33} and Lemma~\ref{l34} to give asymptotic formulas of the form (53) for the corresponding Type-II sums (see $\Sigma_{2}$ in the Remark of Lemma~\ref{l213} and the sums $V^{(1)}(q_1 q_2)$, $V^{(2)}(q_1 q_2)$ and $V^{(3)}(q_1 q_2)$ in \cite{MaynardLargeModuliI}). After applying [\cite{MaynardLargeModuliI}, Lemma 10.2] 3 times, we use Lemma~\ref{l37} or Lemma~\ref{l38} to give an asymptotic formula of the form (53) for the remaining Type-I$_3$ sum (see [\cite{MaynardLargeModuliI}, (11.4)]). 

When the available Type-II range starts from $\left[2 \theta - 1 + \varepsilon,\ \cdots \right)$, we need $2 \theta - 1 + \varepsilon < 1 - \frac{15}{8} \theta - 20 \varepsilon$, or $\theta < \frac{16}{31} - 10 \varepsilon \approx 0.5161$ to give an asymptotic formula of the form (53) for (56). In this case we take $M_0 = x^{2 \theta - 1 + \varepsilon}$ in Lemma~\ref{l37} and $y = x^{2 \theta - 1 + \varepsilon}$ in [\cite{MaynardLargeModuliI}, Lemma 10.2]. Since $3 (2 \theta - 1 + \varepsilon) < \kappa$ when $\theta < \frac{16}{31} - 10 \varepsilon$, we can give an asymptotic formula of the form (53) for (56) if $(\theta_1, \theta_2) \in \boldsymbol{T}_{1}$. When the available Type-II range starts from $\left[2 \theta_1 + \theta_2 - 1 + \varepsilon,\ \cdots \right)$, we need the condition $44 \theta_1 + 26 \theta_2 < 23$ and the details are similar.

When the available Type-II range starts from $\left[\varepsilon,\ \cdots \right)$, we can use both Lemma~\ref{l37} and Lemma~\ref{l38} since we can now take $M_0 = y = x^{\varepsilon}$. For all $\theta < \frac{17}{32}$, the Type-II range $\left[\varepsilon,\ \kappa \right]$ is sufficient to give an asymptotic formula of the form (53) for (56) if $(\theta_1, \theta_2) \in \boldsymbol{T}$.

When applying this three-dimensional Harman's sieve, one can follow the process in \cite{676}: Suppose that $\kappa \geqslant \frac{1}{7} + \varepsilon$. For three ``large'' variables $p_1 p_2 m_3 \sim x$ such that $x^{\kappa} \leqslant p_1, p_2, m_3 \leqslant x^{\frac{3}{7} + \varepsilon}$, we have
\begin{align}
\nonumber \sum_{\kappa \leqslant \alpha_1, \alpha_2 \leqslant \frac{3}{7} + \varepsilon} S\left(\mathcal{A}^{q_1 q_2}_{p_1 p_2}, x^{\kappa}\right) =&\ \sum_{\substack{x^{\kappa} \leqslant p_2 < p_1 < x^{\frac{3}{7} + \varepsilon} \\ p_1 p_2 m_3 \in \mathcal{A}^{q_1 q_2} \\ \left(p_1 p_2 m_3, P(x^{\kappa})\right) = 1 }} 1 \\
\nonumber =&\ \sum_{\substack{x^{\kappa} \leqslant m_2 < m_1 < x^{\frac{3}{7} + \varepsilon} \\ m_1 m_2 m_3 \in \mathcal{A}^{q_1 q_2} \\ \left(m_1 m_2 m_3, P(x^{\kappa})\right) = 1 }} 1 - \sum_{\substack{x^{\kappa} \leqslant m_2 < m_1 < x^{\frac{3}{7} + \varepsilon} \\ m_1 m_2 m_3 \in \mathcal{A}^{q_1 q_2} \\ \left(m_1 m_2 m_3, P(x^{\kappa})\right) = 1 \\ \Omega(m_1 m_2) \geqslant 3 }} 1 \\
=&\ S\left(\mathcal{F}^{q_1 q_2}, x^{\kappa} \right) - \sum_{\substack{x^{\kappa} \leqslant m_2 < m_1 < x^{\frac{3}{7} + \varepsilon} \\ m_1 m_2 m_3 \in \mathcal{A}^{q_1 q_2} \\ \left(m_1 m_2 m_3, P(x^{\kappa})\right) = 1 \\ \Omega(m_1 m_2) \geqslant 3 }} 1.
\end{align}
We can give an asymptotic formula of the form (53) for the first sum on the right-hand side of (57). For the second sum, since $\kappa \geqslant \frac{1}{7} + \varepsilon$ and $m_1, m_2 \leqslant x^{\frac{3}{7} + \varepsilon}$, we know that $\Omega(m_1), \Omega(m_2) \leqslant 2$. Thus, the loss integrals come from this sum is similar to $I_1$ ($\Omega(m_1 m_2) = 3$) and $I_2$ ($\Omega(m_1 m_2) = 4$) in Lemma~\ref{l220} with modified integration regions. We note that in many applications of this device, we can use it on the whole of the two-dimensional region $C$ defined in Section 2 since it only requires the variables are smaller than $x^{\frac{3}{7} + \varepsilon}$, which holds naturally since we have $\alpha_2 < \frac{1}{3}$, $\alpha_1 < \frac{3}{7} + \varepsilon$ and $\alpha_1 + \alpha_2 > \frac{4}{7} - \varepsilon$ (note that $\left[\theta + \varepsilon,\ \frac{4}{7} - \varepsilon \right]$ is a Type-II range and $\alpha_1 + \alpha_2 < \theta + \varepsilon \Rightarrow \boldsymbol{\alpha}_2 \in A$) when $\boldsymbol{\alpha}_2 \in C$, $\theta \leqslant \frac{17}{32} - \varepsilon$. In addition, the above process is still applicable in some cases with $\kappa < \frac{1}{7} + \varepsilon$ in Section 4, and we shall discuss them in the next section. Sometimes we do not have enough Type-II information to apply this device, and we can go back to Lemma~\ref{l220} if we are in the case $C_1^{*}(\theta_1, \theta_2)$. Note that Lemma~\ref{l221} is not applicable in this and later Sections, since $q = q_1 q_2$ cannot be prime.

In the final decompositions in this Subsection, we need to consider both $C_1(\theta_1, \theta_2)$ and $C_1^{*}(\theta_1, \theta_2)$. In the case $C_1(\theta_1, \theta_2)$, we still need to discard the one-dimensional sum that summing over $p_1 \in \left(1 - \theta - \varepsilon,\ \frac{1}{2} \right)$ as in Section 2. Fortunately, Lemma~\ref{l35} gives an asymptotic formula of the form (53) for that sum when $p_1 \in \left[\theta_1 + \varepsilon,\ 1 - \theta_1 - \varepsilon \right]$ under the condition $7 \theta_1 + 12 \theta_2 < 4 - 19 \varepsilon$, and we only need to discard the remaining parts in this case. Note that when we have $2 \theta_1 + \theta_2 \leqslant 1 - 2 \varepsilon$, we can give an asymptotic formula of the form (53) for the whole sum that summing over $p_1 \in \left(1 - \theta - \varepsilon,\ \frac{1}{2} \right)$ since $\theta_1 + \varepsilon \leqslant 1 - \theta - \varepsilon$. In the case $C_1^{*}(\theta_1, \theta_2)$, we can use Lemma~\ref{l218} to handle this sum when $\theta < \frac{17}{32}$, and Lemma~\ref{l35} also helps us give asymptotic formulas of the form (54) for parts of the three-dimensional sum in Lemma~\ref{l218} under the condition $7 \theta_1 + 12 \theta_2 < 4 - 19 \varepsilon$. We shall not explicitly state the differences between cases $C_1(\theta_1, \theta_2)$ and $C_1^{*}(\theta_1, \theta_2)$ in the final decompositions, and the readers can easily find the differences from here. Note that their differences only come from two-dimensional and three-dimensional sieves.

Now we assume that $(\theta_1, \theta_2) \in \boldsymbol{A}$. We divide $\boldsymbol{A}$ into 17 subregions:
$$
\boldsymbol{A} = \boldsymbol{A}_{01} \cup \boldsymbol{A}_{02} \cup \boldsymbol{A}_{03} \cup \boldsymbol{A}_{04} \cup \boldsymbol{A}_{05} \cup \boldsymbol{A}_{06} \cup \boldsymbol{A}_{07} \cup \boldsymbol{A}_{08} \cup \boldsymbol{A}_{09} \cup \boldsymbol{A}_{10} \cup \boldsymbol{A}_{11} \cup \boldsymbol{A}_{12} \cup \boldsymbol{A}_{13} \cup \boldsymbol{A}_{14} \cup \boldsymbol{A}_{15} \cup \boldsymbol{A}_{16} \cup \boldsymbol{A}_{17},
$$
where
\begin{align}
\nonumber \boldsymbol{A}_{01} =&\ \left\{ (\theta_1, \theta_2) : \frac{5}{14} < \theta_1 \leqslant \frac{2}{5},\ \frac{1}{2}(1 - 2 \theta_1) < \theta_2 < \frac{1}{9}(2 - 2 \theta_1) \right\}, \\
\nonumber \boldsymbol{A}_{02} =&\ \left\{ (\theta_1, \theta_2) : \frac{2}{5} < \theta_1 \leqslant \frac{4}{9},\ \frac{1}{12}(4 - 7 \theta_1) < \theta_2 < \frac{1}{3}(2 - 4 \theta_1) \right\}, \\
\nonumber \boldsymbol{A}_{03} =&\ \left\{ (\theta_1, \theta_2) : \frac{2}{5} < \theta_1 \leqslant \frac{4}{9},\ \frac{1}{3}(2 - 4 \theta_1) \leqslant \theta_2 < \frac{1}{6}(2 - 3 \theta_1),\ \theta_2 < \frac{1}{13}(10 - 20 \theta_1) \right\}, \\
\nonumber \boldsymbol{A}_{04} =&\ \left\{ (\theta_1, \theta_2) : \frac{2}{5} < \theta_1 \leqslant \frac{4}{9},\ \frac{1}{13}(10 - 20 \theta_1) \leqslant \theta_2 < \frac{1}{6}(2 - 3 \theta_1),\ \theta_1 + \theta_2 < \frac{11}{20} \right\}, \\
\nonumber \boldsymbol{A}_{05} =&\ \left\{ (\theta_1, \theta_2) : \frac{2}{5} < \theta_1 \leqslant \frac{4}{9},\ \frac{11}{20} - \theta_1 \leqslant \theta_2 < \frac{1}{6}(2 - 3 \theta_1) \right\}, \\
\nonumber \boldsymbol{A}_{06} =&\ \left\{ (\theta_1, \theta_2) : \frac{4}{9} < \theta_1 \leqslant \frac{5}{11},\ \frac{1}{12}(4 - 7 \theta_1) < \theta_2 < \frac{1}{13}(10 - 20 \theta_1) \right\}, \\
\nonumber \boldsymbol{A}_{07} =&\ \left\{ (\theta_1, \theta_2) : \frac{4}{9} < \theta_1 \leqslant \frac{5}{11},\ \frac{1}{13}(10 - 20 \theta_1) \leqslant \theta_2 < 1 - 2 \theta_1,\ \theta_1 + \theta_2 < \frac{11}{20} \right\}, \\
\nonumber \boldsymbol{A}_{08} =&\ \left\{ (\theta_1, \theta_2) : \frac{4}{9} < \theta_1 \leqslant \frac{5}{11},\ \frac{11}{20} - \theta_1 \leqslant \theta_2 < 1 - 2 \theta_1 \right\}, \\
\nonumber \boldsymbol{A}_{09} =&\ \left\{ (\theta_1, \theta_2) : \frac{5}{11} < \theta_1 < \frac{1}{2},\ \frac{1}{20}(10 - 19 \theta_1) < \theta_2 < \frac{1}{13}(10 - 20 \theta_1) \right\}, \\
\nonumber \boldsymbol{A}_{10} =&\ \left\{ (\theta_1, \theta_2) : \frac{5}{11} < \theta_1 < \frac{1}{2},\ \frac{1}{20}(10 - 19 \theta_1) < \theta_2 < 1 - 2 \theta_1,\ \frac{1}{13}(10 - 20 \theta_1) \leqslant \theta_2 \right\}, \\
\nonumber \boldsymbol{A}_{11} =&\ \left\{ (\theta_1, \theta_2) : \frac{4}{9} < \theta_1 < \frac{1}{2},\ 1 - 2 \theta_1 \leqslant \theta_2 < \frac{1}{9}(5 - 9 \theta_1),\ \theta_1 + \theta_2 < \frac{11}{20},\ \theta_2 < \frac{1}{14}(10 - 19 \theta_1) \right\}, \\
\nonumber \boldsymbol{A}_{12} =&\ \left\{ (\theta_1, \theta_2) : \frac{4}{9} < \theta_1 < \frac{1}{2},\ 1 - 2 \theta_1 \leqslant \theta_2 < \frac{1}{9}(5 - 9 \theta_1),\ \theta_1 + \theta_2 < \frac{11}{20},\ \theta_2 \geqslant \frac{1}{14}(10 - 19 \theta_1) \right\}, \\
\nonumber \boldsymbol{A}_{13} =&\ \left\{ (\theta_1, \theta_2) : \frac{4}{9} < \theta_1 < \frac{1}{2},\ 1 - 2 \theta_1 \leqslant \theta_2 < \frac{1}{9}(5 - 9 \theta_1),\ \theta_1 + \theta_2 \geqslant \frac{11}{20},\ \theta_2 < \frac{1}{14}(10 - 19 \theta_1) \right\}, \\
\nonumber \boldsymbol{A}_{14} =&\ \left\{ (\theta_1, \theta_2) : \frac{4}{9} < \theta_1 < \frac{1}{2},\ 1 - 2 \theta_1 \leqslant \theta_2 < \frac{1}{9}(5 - 9 \theta_1),\ \theta_1 + \theta_2 \geqslant \frac{11}{20},\ \theta_2 \geqslant \frac{1}{14}(10 - 19 \theta_1) \right\}, \\
\nonumber \boldsymbol{A}_{15} =&\ \left\{ (\theta_1, \theta_2) : \frac{1}{2} \leqslant \theta_1 < \frac{10}{19},\ 0 < \theta_2 < \frac{1}{14}(10 - 19 \theta_1) \right\}, \\
\nonumber \boldsymbol{A}_{16} =&\ \left\{ (\theta_1, \theta_2) : \frac{1}{2} \leqslant \theta_1 < \frac{11}{20},\ \frac{1}{14}(10 - 19 \theta_1) \leqslant \theta_2 < \frac{11}{20} - \theta_1 \right\}, \\
\nonumber \boldsymbol{A}_{17} =&\ \left\{ (\theta_1, \theta_2) : \frac{1}{2} \leqslant \theta_1 < \frac{11}{20},\ \frac{11}{20} - \theta_1 \leqslant \theta_2 < \frac{1}{18}(11 - 20 \theta_1) \right\}.
\end{align}
Note that we have a Type-II range $\left[\frac{3}{7} + \varepsilon,\ 1 - \theta - \varepsilon \right]$ for all $\theta \leqslant \frac{127}{224} - \varepsilon$ by Lemma~\ref{l216}. For $\theta \leqslant \frac{45}{89} - \varepsilon$ we also have a Type-II range $\left((\log x)^{\varepsilon - 1},\ \kappa \right]$ by Lemma~\ref{l24}, and we shall not repeatedly state these two ranges in the decompositions for the sake of simplicity.

\subsubsection{$\boldsymbol{A}_{01}$}
For $(\theta_1, \theta_2) \in \boldsymbol{A}_{01}$ we have 3 available Type-II information ranges:
\begin{equation}
\left[2 \theta_1 + 2 \theta_2 - 1 + \varepsilon,\ \frac{1}{6}(5 - 8 \theta_1 - 8 \theta_2) - \varepsilon \right], \quad \left[\theta_2 + \varepsilon,\ \frac{1}{6}(2 - 2 \theta_1 - 3 \theta_2) - \varepsilon \right] \quad \text{and} \quad \left[\varepsilon,\ \frac{1}{6}(4 - 7 \theta_1 - 8 \theta_2) - \varepsilon \right].
\end{equation}
The first range comes from Lemma~\ref{l23}, the second comes from Lemma~\ref{l33}, and the third comes from Lemma~\ref{l34}. We divide $\boldsymbol{A}_{01}$ into 4 subregions based on the overlapping conditions of these ranges.
$$
\boldsymbol{A}_{01} = \boldsymbol{A}_{0101} \cup \boldsymbol{A}_{0102} \cup \boldsymbol{A}_{0103} \cup \boldsymbol{A}_{0104},
$$
where
\begin{align}
\nonumber \boldsymbol{A}_{0101} =&\ \left\{ (\theta_1, \theta_2) : \frac{5}{14} < \theta_1 \leqslant \frac{50}{131},\ \frac{1}{2}(1 - 2 \theta_1) < \theta_2 < \frac{1}{9}(2 - 2 \theta_1) \right. \\
\nonumber & \left. \qquad \qquad \quad \text{ or } \frac{50}{131} < \theta_1 \leqslant \frac{2}{5},\ \frac{1}{2}(1 - 2 \theta_1) < \theta_2 < \frac{1}{20}(10 - 19 \theta_1) \right\}, \\
\nonumber \boldsymbol{A}_{0102} =&\ \left\{ (\theta_1, \theta_2) : \frac{50}{131} < \theta_1 \leqslant \frac{17}{44},\ \frac{1}{20}(10 - 19 \theta_1) < \theta_2 < \frac{1}{9}(2 - 2 \theta_1) \right. \\
\nonumber & \left. \qquad \qquad \quad \text{ or } \frac{17}{44} < \theta_1 \leqslant \frac{2}{5},\ \frac{1}{20}(10 - 19 \theta_1) < \theta_2 < \frac{1}{5}(3 - 6 \theta_1) \right\}, \\
\nonumber \boldsymbol{A}_{0103} =&\ \left\{ (\theta_1, \theta_2) : \frac{17}{44} < \theta_1 \leqslant \frac{2}{5},\ \frac{1}{5}(3 - 6 \theta_1) < \theta_2 < \frac{1}{14}(5 - 8 \theta_1) \right\}, \\
\nonumber \boldsymbol{A}_{0104} =&\ \left\{ (\theta_1, \theta_2) : \frac{17}{44} < \theta_1 \leqslant \frac{2}{5},\ \frac{1}{14}(5 - 8 \theta_1) < \theta_2 < \frac{1}{9}(2 - 2 \theta_1) \right\}.
\end{align}
Note that we have $\theta < \frac{7}{13}$ for $(\theta_1, \theta_2) \in \boldsymbol{A}_{01}$, and $\theta < \frac{17}{32}$ for $(\theta_1, \theta_2) \in \boldsymbol{A}_{0101} \cup \boldsymbol{A}_{0102} \cup \boldsymbol{A}_{0103}$.

In $\boldsymbol{A}_{0101}$ we have 
$$
\frac{1}{6}(4 - 7 \theta_1 - 8 \theta_2) > 2 \theta_1 + 2 \theta_2 - 1 \quad \text{and} \quad \frac{1}{6}(5 - 8 \theta_1 - 8 \theta_2) > \frac{1}{6}(2 - 2 \theta_1 - 3 \theta_2).
$$
Hence, the Type-II range for $\boldsymbol{A}_{0101}$ is
\begin{equation}
\left[\varepsilon,\ \frac{1}{6}(5 - 8 \theta_1 - 8 \theta_2) - \varepsilon \right].
\end{equation}
Since the Type-II range starts from $\left[\varepsilon,\ \cdots \right)$, we want to prove that for $(\theta_1, \theta_2) \in \boldsymbol{A}_{0101}$,
\begin{equation}
\sum_{\substack{\alpha_j < \cdots < \alpha_1 \\ (\alpha_1, \ldots, \alpha_j, \varepsilon) \in \boldsymbol{S}_{j+1} }} S\left(\mathcal{A}^{q_1 q_2}_{p_1 \cdots p_j}, x^{\kappa}\right)
\end{equation}
has an asymptotic formula of the form (53). Using Buchstab's identity repeatedly as in the proof of Lemma~\ref{l213}, we have
\begin{align}
\nonumber \sum_{\substack{\alpha_j < \cdots < \alpha_1 \\ (\alpha_1, \ldots, \alpha_j, \varepsilon) \in \boldsymbol{S}_{j+1} }} S\left(\mathcal{A}^{q_1 q_2}_{p_1 \cdots p_j}, x^{\kappa}\right) =&\ \sum_{\substack{\alpha_j < \cdots < \alpha_1 \\ (\alpha_1, \ldots, \alpha_j, \varepsilon) \in \boldsymbol{S}_{j+1} }} S\left(\mathcal{A}^{q_1 q_2}_{p_1 \cdots p_j}, x^{2 \varepsilon}\right) \\
\nonumber & - \sum_{\substack{\alpha_j < \cdots < \alpha_1 \\ (\alpha_1, \ldots, \alpha_j, \varepsilon) \in \boldsymbol{S}_{j+1} \\ 2 \varepsilon \leqslant \alpha_{j+1} < \kappa }} S\left(\mathcal{A}^{q_1 q_2}_{p_1 \cdots p_{j+1}}, p_{j+1}\right) \\
\nonumber =&\ \sum_{\substack{\alpha_j < \cdots < \alpha_1 \\ (\alpha_1, \ldots, \alpha_j, \varepsilon) \in \boldsymbol{S}_{j+1} }} S\left(\mathcal{A}^{q_1 q_2}_{p_1 \cdots p_j}, \exp\left(\log x (\log \log x)^{-3} \right)\right) \\
\nonumber & - \sum_{\substack{\alpha_j < \cdots < \alpha_1 \\ (\alpha_1, \ldots, \alpha_j, \varepsilon) \in \boldsymbol{S}_{j+1} \\ \varepsilon \leqslant \alpha_{j+1} < 2 \varepsilon }} S\left(\mathcal{A}^{q_1 q_2}_{p_1 \cdots p_{j+1}}, p_{j+1}\right) \\
\nonumber & + \sum_{k \geqslant j+1} (-1)^{k-j} \sum_{\substack{\alpha_j < \cdots < \alpha_1 \\ (\alpha_1, \ldots, \alpha_j, \varepsilon) \in \boldsymbol{S}_{j+1} \\ \alpha_k < \cdots < \alpha_{j+1} \\ \alpha_{j+1} < 2 \varepsilon \\ \alpha_{j+1} + \cdots + \alpha_{k} < \varepsilon }} S\left(\mathcal{A}^{q_1 q_2}_{p_1 \cdots p_{k}}, \exp\left(\log x (\log \log x)^{-3} \right)\right) \\
\nonumber & + \sum_{k \geqslant j+2} (-1)^{k-j} \sum_{\substack{\alpha_j < \cdots < \alpha_1 \\ (\alpha_1, \ldots, \alpha_j, \varepsilon) \in \boldsymbol{S}_{j+1} \\ \alpha_k < \cdots < \alpha_{j+1} \\ \alpha_{j+1} < 2 \varepsilon \\ \alpha_{j+1} + \cdots + \alpha_{k-1} < \varepsilon \leqslant \alpha_{j+1} + \cdots + \alpha_{k} }} S\left(\mathcal{A}^{q_1 q_2}_{p_1 \cdots p_{k}}, p_{k}\right) \\
\nonumber & - \sum_{\substack{\alpha_j < \cdots < \alpha_1 \\ (\alpha_1, \ldots, \alpha_j, \varepsilon) \in \boldsymbol{S}_{j+1} \\ 2 \varepsilon \leqslant \alpha_{j+1} < \kappa }} S\left(\mathcal{A}^{q_1 q_2}_{p_1 \cdots p_{j+1}}, p_{j+1}\right) \\
=&\ \Sigma_{33011} - \Sigma_{33012} + \Sigma_{33013} + \Sigma_{33014} - \Sigma_{33015}.
\end{align}
We can give asymptotic formulas of the form (53) for $\Sigma_{33012}$ and $\Sigma_{33015}$ by Lemma~\ref{l211} and our Type-II range (59). We can give asymptotic formulas of the form (53) for $\Sigma_{33011}$ and $\Sigma_{33013}$ by the condition $(\alpha_1, \ldots, \alpha_j, \varepsilon) \in \boldsymbol{S}_{j+1}$ and an application of Lemma~\ref{l29}. Note that $\alpha_{j+1} + \cdots + \alpha_{k} \in \left[\varepsilon,\ 3 \varepsilon \right)$, we can give an asymptotic formula of the form (53) for $\Sigma_{33014}$ by Lemma~\ref{l211} and our Type-II range (59). Hence we know that the sum (60) has an asymptotic formula of the form (53).

Now, the decompositions in this case are very similar to the \textit{Case 1} in Subsection 2.4, except that we cannot use role-reversals when $\theta \geqslant \frac{45}{89}$. We split the range of $\theta_1 + \theta_2$ as in \cite{LRB679}, and replace the conditions $\boldsymbol{\alpha}_{j} \in \boldsymbol{U}_{j}$ with $(\alpha_1, \ldots, \alpha_j, \varepsilon) \in \boldsymbol{S}_{j+1}$.

In $\boldsymbol{A}_{0102}$ we have 
$$
\frac{1}{6}(4 - 7 \theta_1 - 8 \theta_2) < 2 \theta_1 + 2 \theta_2 - 1 \quad \text{and} \quad \frac{1}{6}(5 - 8 \theta_1 - 8 \theta_2) > \frac{1}{6}(2 - 2 \theta_1 - 3 \theta_2).
$$
Hence, the Type-II range for $\boldsymbol{A}_{0102}$ is
\begin{equation}
\left[\varepsilon,\ \frac{1}{6}(4 - 7 \theta_1 - 8 \theta_2) - \varepsilon \right] \cup \left[2 \theta_1 + 2 \theta_2 - 1 + \varepsilon,\ \frac{1}{6}(5 - 8 \theta_1 - 8 \theta_2) - \varepsilon \right].
\end{equation}
This case is similar to the \textit{Case 2} in Subsection 2.4. We ignore the first range $\left[\varepsilon,\ \frac{1}{6}(4 - 7 \theta_1 - 8 \theta_2) - \varepsilon \right]$ and use only the range comes from Lemma~\ref{l23}. The decompositions in this case are very similar to parts of the work done in \cite{LRB679}, where $\theta_1 + \theta_2 = \theta$ lies in some subranges.

In $\boldsymbol{A}_{0103}$ we have 
$$
\frac{1}{6}(4 - 7 \theta_1 - 8 \theta_2) < 2 \theta_1 + 2 \theta_2 - 1 \quad \text{and} \quad \theta_2 < \frac{1}{6}(5 - 8 \theta_1 - 8 \theta_2) < \frac{1}{6}(2 - 2 \theta_1 - 3 \theta_2).
$$
Hence, the Type-II range for $\boldsymbol{A}_{0103}$ is
\begin{equation}
\left[\varepsilon,\ \frac{1}{6}(4 - 7 \theta_1 - 8 \theta_2) - \varepsilon \right] \cup \left[2 \theta_1 + 2 \theta_2 - 1 + \varepsilon,\ \frac{1}{6}(2 - 2 \theta_1 - 3 \theta_2) - \varepsilon \right].
\end{equation}
This case is almost the same as the above case. The only difference here is that we replace the value of $\kappa = \frac{5 - 8 \theta}{6} - \varepsilon$ with $\frac{2 - 2 \theta_1 - 3 \theta_2}{6} - \varepsilon$. Since we have $\frac{5 - 8 \theta}{6} < \frac{2 - 2 \theta_1 - 3 \theta_2}{6}$ in this case, we can simply use Buchstab's identity to get
\begin{equation}
\sum_{\boldsymbol{\alpha}_{j} } S\left(\mathcal{A}^{q_1 q_2}_{p_1 \cdots p_j}, x^{\frac{2 - 2 \theta_1 - 3 \theta_2}{6} - \varepsilon}\right) = \sum_{\boldsymbol{\alpha}_{j} } S\left(\mathcal{A}^{q_1 q_2}_{p_1 \cdots p_j}, x^{\kappa}\right) - \sum_{\substack{ \boldsymbol{\alpha}_{j} \\ \kappa \leqslant \alpha_{j+1} < \frac{2 - 2 \theta_1 - 3 \theta_2}{6} - \varepsilon }} S\left(\mathcal{A}^{q_1 q_2}_{p_1 \cdots p_j p_{j+1}}, p_{j+1}\right).
\end{equation}
We can give an asymptotic formula of the form (53) for the last sum in (64) by our Type-II information in this case. By applying (64), we know that for $(\theta_1, \theta_2) \in \boldsymbol{A}_{0103}$, any $\boldsymbol{\alpha}_{j}$ such that
$$
\sum_{\boldsymbol{\alpha}_{j} } S\left(\mathcal{A}^{q_1 q_2}_{p_1 \cdots p_j}, x^{\kappa}\right)
$$
has an asymptotic formula of the form (53) also yields an asymptotic formula of the form (53) for the sum
$$
\sum_{\boldsymbol{\alpha}_{j} } S\left(\mathcal{A}^{q_1 q_2}_{p_1 \cdots p_j}, x^{\frac{2 - 2 \theta_1 - 3 \theta_2}{6} - \varepsilon}\right).
$$
Hence, we can change our ``starting point'' from $\kappa$ to $\frac{2 - 2 \theta_1 - 3 \theta_2}{6} - \varepsilon$. In the applications of two-dimensional and three-dimensional sieves, the details are similar. Now the decompositions can be easily performed as in the case $\boldsymbol{A}_{0102}$.

In $\boldsymbol{A}_{0104}$ we have 
$$
\frac{1}{6}(4 - 7 \theta_1 - 8 \theta_2) < 2 \theta_1 + 2 \theta_2 - 1 \quad \text{and} \quad \frac{1}{6}(5 - 8 \theta_1 - 8 \theta_2) < \theta_2.
$$
Hence, the Type-II range for $\boldsymbol{A}_{0104}$ is
\begin{equation}
\left[\varepsilon,\ \frac{1}{6}(4 - 7 \theta_1 - 8 \theta_2) - \varepsilon \right] \cup \left[2 \theta_1 + 2 \theta_2 - 1 + \varepsilon,\ \frac{1}{6}(5 - 8 \theta_1 - 8 \theta_2) - \varepsilon \right] \cup \left[\theta_2 + \varepsilon,\ \frac{1}{6}(2 - 2 \theta_1 - 3 \theta_2) - \varepsilon \right].
\end{equation}
Assume that $\theta < \frac{17}{32}$. We use the middle Type-II range $\left[2 \theta_1 + 2 \theta_2 - 1 + \varepsilon,\ \frac{1}{6}(5 - 8 \theta_1 - 8 \theta_2) - \varepsilon \right]$ to give a ``starting point'' $\kappa = \frac{5 - 8 \theta}{6} - \varepsilon$, and the third Type-II range $\left[\theta_2 + \varepsilon,\ \frac{1}{6}(2 - 2 \theta_1 - 3 \theta_2) - \varepsilon \right]$ is used to subtract the contributions of those sums with products of variables lie in this range. 
The decompositions are also similar to which in the case $\boldsymbol{A}_{0102}$.

When $\theta \geqslant \frac{17}{32}$, we follow the decompositions in \cite{LRB679} and Subsection 2.4, and use our Type-II ranges to give asymptotic formulas of the form (53) for sums with products of variables lie in the range.

\subsubsection{$\boldsymbol{A}_{02}$}
For $(\theta_1, \theta_2) \in \boldsymbol{A}_{02}$ we have 3 available Type-II information ranges:
\begin{equation}
\left[2 \theta_1 + 2 \theta_2 - 1 + \varepsilon,\ \frac{1}{6}(5 - 8 \theta_1 - 8 \theta_2) - \varepsilon \right], \quad \left[\theta_2 + \varepsilon,\ \frac{1}{6}(2 - 2 \theta_1 - 3 \theta_2) - \varepsilon \right] \quad \text{and} \quad \left[\varepsilon,\ \frac{1}{6}(4 - 7 \theta_1 - 8 \theta_2) - \varepsilon \right].
\end{equation}
The first range comes from Lemma~\ref{l23}, the second comes from Lemma~\ref{l33}, and the third comes from Lemma~\ref{l34}. We divide $\boldsymbol{A}_{02}$ into 4 subregions based on the overlapping conditions of these ranges.
$$
\boldsymbol{A}_{02} = \boldsymbol{A}_{0201} \cup \boldsymbol{A}_{0202} \cup \boldsymbol{A}_{0203} \cup \boldsymbol{A}_{0204},
$$
where
\begin{align}
\nonumber \boldsymbol{A}_{0201} =&\ \left\{ (\theta_1, \theta_2) : \frac{2}{5} < \theta_1 < \frac{16}{37},\ \frac{1}{12}(4 - 7 \theta_1) < \theta_2 < \frac{1}{5}(3 - 6 \theta_1) \right\}, \\
\nonumber \boldsymbol{A}_{0202} =&\ \left\{ (\theta_1, \theta_2) : \frac{2}{5} < \theta_1 \leqslant \frac{16}{37},\ \frac{1}{5}(3 - 6 \theta_1) \leqslant \theta_2 < \frac{1}{20}(10 - 19 \theta_1) \right. \\
\nonumber & \qquad \qquad \quad \text{ or } \frac{16}{37} < \theta_1 \leqslant \frac{10}{23},\ \frac{1}{12}(4 - 7 \theta_1) < \theta_2 < \frac{1}{20}(10 - 19 \theta_1) \\
\nonumber & \left. \qquad \qquad \quad \text{ or } \frac{10}{23} < \theta_1 < \frac{4}{9},\ \frac{1}{12}(4 - 7 \theta_1) < \theta_2 < \frac{1}{3}(2 - 4 \theta_1) \right\}, \\
\nonumber \boldsymbol{A}_{0203} =&\ \left\{ (\theta_1, \theta_2) : \frac{2}{5} < \theta_1 \leqslant \frac{13}{32},\ \frac{1}{20}(10 - 19 \theta_1) \leqslant \theta_2 < \frac{1}{14}(5 - 8 \theta_1) \right. \\
\nonumber & \left. \qquad \qquad \quad \text{ or } \frac{13}{32} < \theta_1 < \frac{10}{23},\ \frac{1}{20}(10 - 19 \theta_1) \leqslant \theta_2 < \frac{1}{3}(2 - 4 \theta_1) \right\}, \\
\nonumber \boldsymbol{A}_{0204} =&\ \left\{ (\theta_1, \theta_2) : \frac{2}{5} < \theta_1 < \frac{13}{32},\ \frac{1}{14}(5 - 8 \theta_1) \leqslant \theta_2 < \frac{1}{3}(2 - 4 \theta_1) \right\}.
\end{align}
Note that we have $\theta < \frac{7}{13}$ for $(\theta_1, \theta_2) \in \boldsymbol{A}_{02}$, and $\theta < \frac{17}{32}$ for $(\theta_1, \theta_2) \in \boldsymbol{A}_{0201} \cup \boldsymbol{A}_{0202} \cup \boldsymbol{A}_{0203}$.

In $\boldsymbol{A}_{0201}$ we have 
$$
\frac{1}{6}(4 - 7 \theta_1 - 8 \theta_2) > 2 \theta_1 + 2 \theta_2 - 1, \quad \theta_2 < \frac{1}{6}(5 - 8 \theta_1 - 8 \theta_2) \quad \text{and} \quad \frac{1}{6}(5 - 8 \theta_1 - 8 \theta_2) > \frac{1}{6}(2 - 2 \theta_1 - 3 \theta_2).
$$
Hence, the Type-II range for $\boldsymbol{A}_{0201}$ is
\begin{equation}
\left[\varepsilon,\ \frac{1}{6}(5 - 8 \theta_1 - 8 \theta_2) - \varepsilon \right].
\end{equation}
The decompositions are similar to which in the case $\boldsymbol{A}_{0101}$.

In $\boldsymbol{A}_{0202}$ we have 
$$
\frac{1}{6}(4 - 7 \theta_1 - 8 \theta_2) > 2 \theta_1 + 2 \theta_2 - 1 \quad \text{and} \quad \theta_2 < \frac{1}{6}(5 - 8 \theta_1 - 8 \theta_2) < \frac{1}{6}(2 - 2 \theta_1 - 3 \theta_2).
$$
Hence, the Type-II range for $\boldsymbol{A}_{0202}$ is
\begin{equation}
\left[\varepsilon,\ \frac{1}{6}(2 - 2 \theta_1 - 3 \theta_2) - \varepsilon \right].
\end{equation}
Using Buchstab's identity as in (64), one can replace the value of $\kappa = \frac{5 - 8 \theta}{6} - \varepsilon$ with $\frac{2 - 2 \theta_1 - 3 \theta_2}{6} - \varepsilon$. Again, the decompositions are similar to which in the case $\boldsymbol{A}_{0101}$.

In $\boldsymbol{A}_{0203}$ we have 
$$
\frac{1}{6}(4 - 7 \theta_1 - 8 \theta_2) < 2 \theta_1 + 2 \theta_2 - 1 \quad \text{and} \quad \theta_2 < \frac{1}{6}(5 - 8 \theta_1 - 8 \theta_2) < \frac{1}{6}(2 - 2 \theta_1 - 3 \theta_2).
$$
Hence, the Type-II range for $\boldsymbol{A}_{0203}$ is
\begin{equation}
\left[\varepsilon,\ \frac{1}{6}(4 - 7 \theta_1 - 8 \theta_2) - \varepsilon \right] \cup \left[2 \theta_1 + 2 \theta_2 - 1 + \varepsilon,\ \frac{1}{6}(2 - 2 \theta_1 - 3 \theta_2) - \varepsilon \right].
\end{equation}
The decompositions are similar to which in the case $\boldsymbol{A}_{0103}$.

In $\boldsymbol{A}_{0204}$ we have 
$$
\frac{1}{6}(4 - 7 \theta_1 - 8 \theta_2) < 2 \theta_1 + 2 \theta_2 - 1 \quad \text{and} \quad \frac{1}{6}(5 - 8 \theta_1 - 8 \theta_2) < \theta_2.
$$
Hence, the Type-II range for $\boldsymbol{A}_{0204}$ is
\begin{equation}
\left[\varepsilon,\ \frac{1}{6}(4 - 7 \theta_1 - 8 \theta_2) - \varepsilon \right] \cup \left[2 \theta_1 + 2 \theta_2 - 1 + \varepsilon,\ \frac{1}{6}(5 - 8 \theta_1 - 8 \theta_2) - \varepsilon \right] \cup \left[\theta_2 + \varepsilon,\ \frac{1}{6}(2 - 2 \theta_1 - 3 \theta_2) - \varepsilon \right].
\end{equation}
The decompositions are similar to which in the case $\boldsymbol{A}_{0104}$.

\subsubsection{$\boldsymbol{A}_{03}$}
For $(\theta_1, \theta_2) \in \boldsymbol{A}_{03}$ we have 3 available Type-II information ranges:
\begin{equation}
\left[2 \theta_1 + 2 \theta_2 - 1 + \varepsilon,\ \frac{1}{6}(5 - 8 \theta_1 - 8 \theta_2) - \varepsilon \right], \quad \left[\theta_2 + \varepsilon,\ \frac{1}{3}(2 - 3 \theta_1 - 3 \theta_2) - \varepsilon \right] \quad \text{and} \quad \left[\varepsilon,\ \frac{1}{6}(4 - 7 \theta_1 - 8 \theta_2) - \varepsilon \right].
\end{equation}
The first range comes from Lemma~\ref{l23}, the second comes from Lemma~\ref{l33}, and the third comes from Lemma~\ref{l34}. We divide $\boldsymbol{A}_{03}$ into 3 subregions based on the overlapping conditions of these ranges.
$$
\boldsymbol{A}_{03} = \boldsymbol{A}_{0301} \cup \boldsymbol{A}_{0302} \cup \boldsymbol{A}_{0303},
$$
where
\begin{align}
\nonumber \boldsymbol{A}_{0301} =&\ \left\{ (\theta_1, \theta_2) : \frac{10}{23} < \theta_1 \leqslant \frac{4}{9},\ \frac{1}{3}(2 - 4 \theta_1) \leqslant \theta_2 < \frac{1}{20}(10 - 19 \theta_1) \right\}, \\
\nonumber \boldsymbol{A}_{0302} =&\ \left\{ (\theta_1, \theta_2) : \frac{13}{32} < \theta_1 \leqslant \frac{75}{176},\ \frac{1}{3}(2 - 4 \theta_1) \leqslant \theta_2 < \frac{1}{14}(5 - 8 \theta_1) \right. \\
\nonumber & \qquad \qquad \quad \text{ or } \frac{75}{176} < \theta_1 < \frac{10}{23},\ \frac{1}{3}(2 - 4 \theta_1) < \theta_2 < \frac{1}{13}(10 - 20 \theta_1) \\
\nonumber & \left. \qquad \qquad \quad \text{ or } \frac{10}{23} \leqslant \theta_1 \leqslant \frac{4}{9},\ \frac{1}{20}(10 - 19 \theta_1) < \theta_2 < \frac{1}{13}(10 - 20 \theta_1) \right\}, \\
\nonumber \boldsymbol{A}_{0303} =&\ \left\{ (\theta_1, \theta_2) : \frac{2}{5} < \theta_1 < \frac{13}{32},\ \frac{1}{3}(2 - 4 \theta_1) \leqslant \theta_2 < \frac{1}{6}(2 - 3 \theta_1) \right. \\
\nonumber & \qquad \qquad \quad \text{ or } \frac{13}{32} \leqslant \theta_1 \leqslant \frac{34}{81},\ \frac{1}{14}(5 - 8 \theta_1) < \theta_2 < \frac{1}{6}(2 - 3 \theta_1) \\
\nonumber & \left. \qquad \qquad \quad \text{ or } \frac{34}{81} < \theta_1 < \frac{75}{176},\ \frac{1}{14}(5 - 8 \theta_1) < \theta_2 < \frac{1}{13}(10 - 20 \theta_1) \right\}.
\end{align}
Note that we have $\theta < \frac{17}{32}$ for $(\theta_1, \theta_2) \in \boldsymbol{A}_{0301}$.

In $\boldsymbol{A}_{0301}$ we have 
$$
\frac{1}{6}(4 - 7 \theta_1 - 8 \theta_2) > 2 \theta_1 + 2 \theta_2 - 1 \quad \text{and} \quad \frac{1}{6}(5 - 8 \theta_1 - 8 \theta_2) > \theta_2.
$$
Hence, the Type-II range for $\boldsymbol{A}_{0301}$ is
\begin{equation}
\left[\varepsilon,\ \frac{1}{3}(2 - 3 \theta_1 - 3 \theta_2) - \varepsilon \right].
\end{equation}
Similar to (64), one can replace the value of $\kappa = \frac{5 - 8 \theta}{6} - \varepsilon$ with $\frac{2 - 3 \theta}{3} - \varepsilon$ by applying Buchstab's identity. The decompositions are similar to which in the case $\boldsymbol{A}_{0101}$.

In $\boldsymbol{A}_{0302}$ we have 
$$
\frac{1}{6}(4 - 7 \theta_1 - 8 \theta_2) < 2 \theta_1 + 2 \theta_2 - 1 \quad \text{and} \quad \frac{1}{6}(5 - 8 \theta_1 - 8 \theta_2) > \theta_2.
$$
Hence, the Type-II range for $\boldsymbol{A}_{0302}$ is
\begin{equation}
\left[\varepsilon,\ \frac{1}{6}(4 - 7 \theta_1 - 8 \theta_2) - \varepsilon \right] \cup \left[2 \theta_1 + 2 \theta_2 - 1 + \varepsilon,\ \frac{1}{3}(2 - 3 \theta_1 - 3 \theta_2) - \varepsilon \right].
\end{equation}
For $\theta < \frac{17}{32}$, we replace $\kappa$ with $\frac{2 - 3 \theta}{3} - \varepsilon$, and the decompositions are similar to which in the case $\boldsymbol{A}_{0103}$. For $\theta \geqslant \frac{17}{32}$, we follow the decompositions in \cite{LRB679} and Subsection 2.4, and we use our Type-II range to give extra asymptotic formulas of the form (53) for sums with products of variables lie in the range $\left(\kappa,\ \frac{2 - 3 \theta}{3} - \varepsilon \right]$.

In $\boldsymbol{A}_{0303}$ we have 
$$
\frac{1}{6}(4 - 7 \theta_1 - 8 \theta_2) < 2 \theta_1 + 2 \theta_2 - 1 \quad \text{and} \quad \frac{1}{6}(5 - 8 \theta_1 - 8 \theta_2) < \theta_2.
$$
Hence, the Type-II range for $\boldsymbol{A}_{0303}$ is
\begin{equation}
\left[\varepsilon,\ \frac{1}{6}(4 - 7 \theta_1 - 8 \theta_2) - \varepsilon \right] \cup \left[2 \theta_1 + 2 \theta_2 - 1 + \varepsilon,\ \frac{1}{6}(5 - 8 \theta_1 - 8 \theta_2) - \varepsilon \right] \cup \left[\theta_2 + \varepsilon,\ \frac{1}{3}(2 - 3 \theta_1 - 3 \theta_2) - \varepsilon \right].
\end{equation}
The decompositions are similar to which in the case $\boldsymbol{A}_{0302}$, but now we cannot replace $\kappa$ with $\frac{2 - 3 \theta}{3} - \varepsilon$ when $\theta < \frac{17}{32}$, and we can only give asymptotic formulas of the form (53) for sums with products of variables lie in the range $\left[\theta_2 + \varepsilon,\ \frac{1}{3}(2 - 3 \theta_1 - 3 \theta_2) - \varepsilon \right]$.

\subsubsection{$\boldsymbol{A}_{04}$}
For $(\theta_1, \theta_2) \in \boldsymbol{A}_{04}$ we have 4 available Type-II information ranges:
\begin{equation}
\left[2 \theta_1 + 2 \theta_2 - 1 + \varepsilon,\ \frac{1}{6}(5 - 8 \theta_1 - 8 \theta_2) - \varepsilon \right], \quad \left[\theta_2 + \varepsilon,\ \frac{1}{3}(2 - 3 \theta_1 - 3 \theta_2) - \varepsilon \right],
\end{equation}
\begin{equation}
\left[\varepsilon,\ \frac{1}{6}(4 - 8 \theta_1 - 7 \theta_2) - \varepsilon \right] \quad \text{and} \quad \left[2 \theta_1 + \theta_2 - 1 + \varepsilon,\ \frac{1}{6}(4 - 7 \theta_1 - 8 \theta_2) - \varepsilon \right].
\end{equation}
The first range comes from Lemma~\ref{l23}, the second comes from Lemma~\ref{l33}, and the third and fourth come from Lemma~\ref{l34}. Since we have
$$
\frac{1}{6}(4 - 8 \theta_1 - 7 \theta_2) < 0 \quad \text{and} \quad 2 \theta_1 + \theta_2 - 1 < 0
$$
in this region, (75)--(76) are equivalent to 3 Type-II information ranges
\begin{equation}
\left[2 \theta_1 + 2 \theta_2 - 1 + \varepsilon,\ \frac{1}{6}(5 - 8 \theta_1 - 8 \theta_2) - \varepsilon \right], \quad \left[\theta_2 + \varepsilon,\ \frac{1}{3}(2 - 3 \theta_1 - 3 \theta_2) - \varepsilon \right] \quad \text{and} \quad \left[\varepsilon,\ \frac{1}{6}(4 - 7 \theta_1 - 8 \theta_2) - \varepsilon \right].
\end{equation}
We divide $\boldsymbol{A}_{04}$ into 2 subregions based on the overlapping conditions of these ranges.
$$
\boldsymbol{A}_{04} = \boldsymbol{A}_{0401} \cup \boldsymbol{A}_{0402},
$$
where
\begin{align}
\nonumber \boldsymbol{A}_{0401} =&\ \left\{ (\theta_1, \theta_2) : \frac{75}{176} < \theta_1 \leqslant \frac{4}{9},\ \frac{1}{13}(10 - 20 \theta_1) \leqslant \theta_2 < \frac{1}{14}(5 - 8 \theta_1) \right\}, \\
\nonumber \boldsymbol{A}_{0402} =&\ \left\{ (\theta_1, \theta_2) : \frac{34}{81} < \theta_1 < \frac{75}{176},\ \frac{1}{13}(10 - 20 \theta_1) \leqslant \theta_2 < \frac{1}{6}(2 - 3 \theta_1) \right. \\
\nonumber & \qquad \qquad \quad \text{ or } \frac{75}{176} \leqslant \theta_1 \leqslant \frac{13}{30},\ \frac{1}{14}(5 - 8 \theta_1) < \theta_2 < \frac{1}{6}(2 - 3 \theta_1) \\
\nonumber & \left. \qquad \qquad \quad \text{ or } \frac{13}{30} < \theta_1 \leqslant \frac{4}{9},\ \frac{1}{14}(5 - 8 \theta_1) < \theta_2 < \frac{1}{20}(11 - 20 \theta_1) \right\}.
\end{align}
Note that we have $\theta > \frac{7}{13}$ for $(\theta_1, \theta_2) \in \boldsymbol{A}_{0402}$.

In $\boldsymbol{A}_{0401}$ we have 
$$
\frac{1}{6}(5 - 8 \theta_1 - 8 \theta_2) > \theta_2.
$$
Hence, the Type-II range for $\boldsymbol{A}_{0401}$ is
\begin{equation}
\left[\varepsilon,\ \frac{1}{6}(4 - 7 \theta_1 - 8 \theta_2) - \varepsilon \right] \cup \left[2 \theta_1 + 2 \theta_2 - 1 + \varepsilon,\ \frac{1}{3}(2 - 3 \theta_1 - 3 \theta_2) - \varepsilon \right].
\end{equation}
The decompositions are similar to which in the case $\boldsymbol{A}_{0302}$.

In $\boldsymbol{A}_{0402}$ we have 
$$
\frac{1}{6}(5 - 8 \theta_1 - 8 \theta_2) < \theta_2.
$$
Hence, the Type-II range for $\boldsymbol{A}_{0402}$ is
\begin{equation}
\left[\varepsilon,\ \frac{1}{6}(4 - 7 \theta_1 - 8 \theta_2) - \varepsilon \right] \cup \left[2 \theta_1 + 2 \theta_2 - 1 + \varepsilon,\ \frac{1}{6}(5 - 8 \theta_1 - 8 \theta_2) - \varepsilon \right] \cup \left[\theta_2 + \varepsilon,\ \frac{1}{3}(2 - 3 \theta_1 - 3 \theta_2) - \varepsilon \right].
\end{equation}
The decompositions are similar to which in the case $\boldsymbol{A}_{0303}$.

\subsubsection{$\boldsymbol{A}_{05}$}
For $(\theta_1, \theta_2) \in \boldsymbol{A}_{05}$ we have 3 available Type-II information ranges:
\begin{equation}
\left[\theta_2 + \varepsilon,\ \frac{1}{3}(2 - 3 \theta_1 - 3 \theta_2) - \varepsilon \right], \quad \left[\varepsilon,\ \frac{1}{6}(4 - 8 \theta_1 - 7 \theta_2) - \varepsilon \right] \quad \text{and} \quad \left[2 \theta_1 + \theta_2 - 1 + \varepsilon,\ \frac{1}{6}(4 - 7 \theta_1 - 8 \theta_2) - \varepsilon \right].
\end{equation}
The first range comes from Lemma~\ref{l33}, and the second and third come from Lemma~\ref{l34}. Note that Lemma~\ref{l23} becomes trivial in this region since $\theta_1 + \theta_2 \geqslant \frac{11}{20}$. Since we have
$$
\frac{1}{6}(4 - 8 \theta_1 - 7 \theta_2) < 0 \quad \text{and} \quad 2 \theta_1 + \theta_2 - 1 < 0
$$
in this region, (80) is equivalent to 2 Type-II information ranges
\begin{equation}
\left[\theta_2 + \varepsilon,\ \frac{1}{3}(2 - 3 \theta_1 - 3 \theta_2) - \varepsilon \right] \quad \text{and} \quad \left[\varepsilon,\ \frac{1}{6}(4 - 7 \theta_1 - 8 \theta_2) - \varepsilon \right].
\end{equation}
Hence, the Type-II range for $\boldsymbol{A}_{05}$ is
\begin{equation}
\left[\varepsilon,\ \frac{1}{6}(4 - 7 \theta_1 - 8 \theta_2) - \varepsilon \right] \cup \left[\theta_2 + \varepsilon,\ \frac{1}{3}(2 - 3 \theta_1 - 3 \theta_2) - \varepsilon \right].
\end{equation}
We follow the decomposing process in \cite{LRB679}, and use our Type-II information to give asymptotic formulas of the form (53) for sums with products of variables lie in the range $\left[\theta_2 + \varepsilon,\ \frac{1}{3}(2 - 3 \theta_1 - 3 \theta_2) - \varepsilon \right]$.

\subsubsection{$\boldsymbol{A}_{06}$}
For $(\theta_1, \theta_2) \in \boldsymbol{A}_{06}$ we have 3 available Type-II information ranges:
\begin{equation}
\left[2 \theta_1 + 2 \theta_2 - 1 + \varepsilon,\ \frac{1}{6}(5 - 8 \theta_1 - 8 \theta_2) - \varepsilon \right], \quad \left[\theta_2 + \varepsilon,\ \frac{1}{3}(2 - 3 \theta_1 - 3 \theta_2) - \varepsilon \right] \quad \text{and} \quad \left[\varepsilon,\ \frac{1}{6}(4 - 7 \theta_1 - 8 \theta_2) - \varepsilon \right].
\end{equation}
The first range comes from Lemma~\ref{l23}, the second comes from Lemma~\ref{l33}, and the third comes from Lemma~\ref{l34}. We divide $\boldsymbol{A}_{06}$ into 2 subregions based on the overlapping conditions of these ranges.
$$
\boldsymbol{A}_{06} = \boldsymbol{A}_{0601} \cup \boldsymbol{A}_{0602},
$$
where
\begin{align}
\nonumber \boldsymbol{A}_{0601} =&\ \left\{ (\theta_1, \theta_2) : \frac{4}{9} < \theta_1 < \frac{5}{11},\ \frac{1}{12}(4 - 7 \theta_1) < \theta_2 < \frac{1}{20}(10 - 19 \theta_1) \right\}, \\
\nonumber \boldsymbol{A}_{0602} =&\ \left\{ (\theta_1, \theta_2) : \frac{4}{9} < \theta_1 \leqslant \frac{5}{11},\ \frac{1}{20}(10 - 19 \theta_1) \leqslant \theta_2 < \frac{1}{13}(10 - 20 \theta_1) \right\}.
\end{align}

In $\boldsymbol{A}_{0601}$ we have 
$$
\frac{1}{6}(4 - 7 \theta_1 - 8 \theta_2) > 2 \theta_1 + 2 \theta_2 - 1.
$$
Hence, the Type-II range for $\boldsymbol{A}_{0601}$ is
\begin{equation}
\left[\varepsilon,\ \frac{1}{3}(2 - 3 \theta_1 - 3 \theta_2) - \varepsilon \right].
\end{equation}
The decompositions are similar to which in the case $\boldsymbol{A}_{0301}$.

In $\boldsymbol{A}_{0602}$ we have 
$$
\frac{1}{6}(4 - 7 \theta_1 - 8 \theta_2) < 2 \theta_1 + 2 \theta_2 - 1.
$$
Hence, the Type-II range for $\boldsymbol{A}_{0602}$ is
\begin{equation}
\left[\varepsilon,\ \frac{1}{6}(4 - 7 \theta_1 - 8 \theta_2) - \varepsilon \right] \cup \left[2 \theta_1 + 2 \theta_2 - 1 + \varepsilon,\ \frac{1}{3}(2 - 3 \theta_1 - 3 \theta_2) - \varepsilon \right].
\end{equation}
The decompositions are similar to which in the case $\boldsymbol{A}_{0302}$.

\subsubsection{$\boldsymbol{A}_{07}$}
For $(\theta_1, \theta_2) \in \boldsymbol{A}_{07}$ we have 4 available Type-II information ranges:
\begin{equation}
\left[2 \theta_1 + 2 \theta_2 - 1 + \varepsilon,\ \frac{1}{6}(5 - 8 \theta_1 - 8 \theta_2) - \varepsilon \right], \quad \left[\theta_2 + \varepsilon,\ \frac{1}{3}(2 - 3 \theta_1 - 3 \theta_2) - \varepsilon \right],
\end{equation}
\begin{equation}
\left[\varepsilon,\ \frac{1}{6}(4 - 8 \theta_1 - 7 \theta_2) - \varepsilon \right] \quad \text{and} \quad \left[2 \theta_1 + \theta_2 - 1 + \varepsilon,\ \frac{1}{6}(4 - 7 \theta_1 - 8 \theta_2) - \varepsilon \right].
\end{equation}
The first range comes from Lemma~\ref{l23}, the second comes from Lemma~\ref{l33}, and the third and fourth come from Lemma~\ref{l34}. Since we have
$$
\frac{1}{6}(4 - 8 \theta_1 - 7 \theta_2) < 0 \quad \text{and} \quad 2 \theta_1 + \theta_2 - 1 < 0
$$
in this region, (86)--(87) are equivalent to 3 Type-II information ranges
\begin{equation}
\left[2 \theta_1 + 2 \theta_2 - 1 + \varepsilon,\ \frac{1}{6}(5 - 8 \theta_1 - 8 \theta_2) - \varepsilon \right], \quad \left[\theta_2 + \varepsilon,\ \frac{1}{3}(2 - 3 \theta_1 - 3 \theta_2) - \varepsilon \right] \quad \text{and} \quad \left[\varepsilon,\ \frac{1}{6}(4 - 7 \theta_1 - 8 \theta_2) - \varepsilon \right].
\end{equation}
We divide $\boldsymbol{A}_{07}$ into 2 subregions based on the overlapping conditions of these ranges.
$$
\boldsymbol{A}_{07} = \boldsymbol{A}_{0701} \cup \boldsymbol{A}_{0702},
$$
where
\begin{align}
\nonumber \boldsymbol{A}_{0701} =&\ \left\{ (\theta_1, \theta_2) : \frac{4}{9} < \theta_1 \leqslant \frac{9}{20},\ \frac{1}{13}(10 - 20 \theta_1) \leqslant \theta_2 < \frac{1}{14}(5 - 8 \theta_1) \right. \\
\nonumber & \left. \qquad \qquad \quad \text{ or } \frac{9}{20} < \theta_1 \leqslant \frac{5}{11},\ \frac{1}{13}(10 - 20 \theta_1) \leqslant \theta_2 < 1 - 2 \theta_1 \right\}, \\
\nonumber \boldsymbol{A}_{0702} =&\ \left\{ (\theta_1, \theta_2) : \frac{4}{9} < \theta_1 < \frac{9}{20},\ \frac{1}{14}(5 - 8 \theta_1) \leqslant \theta_2 < \frac{1}{20}(11 - 20 \theta_1) \right\}.
\end{align}

In $\boldsymbol{A}_{0701}$ we have 
$$
\frac{1}{6}(5 - 8 \theta_1 - 8 \theta_2) > \theta_2.
$$
Hence, the Type-II range for $\boldsymbol{A}_{0701}$ is
\begin{equation}
\left[\varepsilon,\ \frac{1}{6}(4 - 7 \theta_1 - 8 \theta_2) - \varepsilon \right] \cup \left[2 \theta_1 + 2 \theta_2 - 1 + \varepsilon,\ \frac{1}{3}(2 - 3 \theta_1 - 3 \theta_2) - \varepsilon \right].
\end{equation}
The decompositions are similar to which in the case $\boldsymbol{A}_{0302}$.

In $\boldsymbol{A}_{0702}$ we have 
$$
\frac{1}{6}(5 - 8 \theta_1 - 8 \theta_2) < \theta_2.
$$
Hence, the Type-II range for $\boldsymbol{A}_{0702}$ is
\begin{equation}
\left[\varepsilon,\ \frac{1}{6}(4 - 7 \theta_1 - 8 \theta_2) - \varepsilon \right] \cup \left[2 \theta_1 + 2 \theta_2 - 1 + \varepsilon,\ \frac{1}{6}(5 - 8 \theta_1 - 8 \theta_2) - \varepsilon \right] \cup \left[\theta_2 + \varepsilon,\ \frac{1}{3}(2 - 3 \theta_1 - 3 \theta_2) - \varepsilon \right].
\end{equation}
The decompositions are similar to which in the case $\boldsymbol{A}_{0303}$.

\subsubsection{$\boldsymbol{A}_{08}$}
For $(\theta_1, \theta_2) \in \boldsymbol{A}_{08}$ we have 3 available Type-II information ranges:
\begin{equation}
\left[\theta_2 + \varepsilon,\ \frac{1}{3}(2 - 3 \theta_1 - 3 \theta_2) - \varepsilon \right], \quad \left[\varepsilon,\ \frac{1}{6}(4 - 8 \theta_1 - 7 \theta_2) - \varepsilon \right] \quad \text{and} \quad \left[2 \theta_1 + \theta_2 - 1 + \varepsilon,\ \frac{1}{6}(4 - 7 \theta_1 - 8 \theta_2) - \varepsilon \right].
\end{equation}
The first range comes from Lemma~\ref{l33}, and the second and third come from Lemma~\ref{l34}. Note that Lemma~\ref{l23} becomes trivial in this region since $\theta_1 + \theta_2 \geqslant \frac{11}{20}$. Since we have
$$
\frac{1}{6}(4 - 8 \theta_1 - 7 \theta_2) < 0 \quad \text{and} \quad 2 \theta_1 + \theta_2 - 1 < 0
$$
in this region, (91) is equivalent to 2 Type-II information ranges
\begin{equation}
\left[\theta_2 + \varepsilon,\ \frac{1}{3}(2 - 3 \theta_1 - 3 \theta_2) - \varepsilon \right] \quad \text{and} \quad \left[\varepsilon,\ \frac{1}{6}(4 - 7 \theta_1 - 8 \theta_2) - \varepsilon \right].
\end{equation}
Hence, the Type-II range for $\boldsymbol{A}_{08}$ is
\begin{equation}
\left[\varepsilon,\ \frac{1}{6}(4 - 7 \theta_1 - 8 \theta_2) - \varepsilon \right] \cup \left[\theta_2 + \varepsilon,\ \frac{1}{3}(2 - 3 \theta_1 - 3 \theta_2) - \varepsilon \right].
\end{equation}
The decompositions are similar to which in the case $\boldsymbol{A}_{05}$.

\subsubsection{$\boldsymbol{A}_{09}$}
For $(\theta_1, \theta_2) \in \boldsymbol{A}_{09}$ we have 3 available Type-II information ranges:
\begin{equation}
\left[2 \theta_1 + 2 \theta_2 - 1 + \varepsilon,\ \frac{1}{6}(5 - 8 \theta_1 - 8 \theta_2) - \varepsilon \right], \quad \left[\theta_2 + \varepsilon,\ \frac{1}{3}(2 - 3 \theta_1 - 3 \theta_2) - \varepsilon \right] \quad \text{and} \quad \left[\varepsilon,\ \frac{1}{6}(4 - 7 \theta_1 - 8 \theta_2) - \varepsilon \right].
\end{equation}
The first range comes from Lemma~\ref{l23}, the second comes from Lemma~\ref{l33}, and the third comes from Lemma~\ref{l34}. Note that in $\boldsymbol{A}_{09}$ we have
$$
2 \theta_1 + 2 \theta_2 - 1 < \theta_2 < \frac{1}{6}(5 - 8 \theta_1 - 8 \theta_2) < \frac{1}{3}(2 - 3 \theta_1 - 3 \theta_2).
$$
Hence, the Type-II range for $\boldsymbol{A}_{09}$ is
\begin{equation}
\left[\varepsilon,\ \frac{1}{6}(4 - 7 \theta_1 - 8 \theta_2) - \varepsilon \right] \cup \left[2 \theta_1 + 2 \theta_2 - 1 + \varepsilon,\ \frac{1}{3}(2 - 3 \theta_1 - 3 \theta_2) - \varepsilon \right].
\end{equation}
We shall discuss the decompositions on this interval in next subsubsection, together with the case $(\theta_1, \theta_2) \in \boldsymbol{A}_{10}$.

\subsubsection{$\boldsymbol{A}_{10}$}
For $(\theta_1, \theta_2) \in \boldsymbol{A}_{10}$ we have 4 available Type-II information ranges:
\begin{equation}
\left[2 \theta_1 + 2 \theta_2 - 1 + \varepsilon,\ \frac{1}{6}(5 - 8 \theta_1 - 8 \theta_2) - \varepsilon \right], \quad \left[\theta_2 + \varepsilon,\ \frac{1}{3}(2 - 3 \theta_1 - 3 \theta_2) - \varepsilon \right],
\end{equation}
\begin{equation}
\left[\varepsilon,\ \frac{1}{6}(4 - 8 \theta_1 - 7 \theta_2) - \varepsilon \right] \quad \text{and} \quad \left[2 \theta_1 + \theta_2 - 1 + \varepsilon,\ \frac{1}{6}(4 - 7 \theta_1 - 8 \theta_2) - \varepsilon \right].
\end{equation}
The first range comes from Lemma~\ref{l23}, the second comes from Lemma~\ref{l33}, and the third and fourth come from Lemma~\ref{l34}. Since we have
$$
\frac{1}{6}(4 - 8 \theta_1 - 7 \theta_2) < 0 \quad \text{and} \quad 2 \theta_1 + \theta_2 - 1 < 0
$$
in this region, (96)--(97) are equivalent to 3 Type-II information ranges
\begin{equation}
\left[2 \theta_1 + 2 \theta_2 - 1 + \varepsilon,\ \frac{1}{6}(5 - 8 \theta_1 - 8 \theta_2) - \varepsilon \right], \quad \left[\theta_2 + \varepsilon,\ \frac{1}{3}(2 - 3 \theta_1 - 3 \theta_2) - \varepsilon \right] \quad \text{and} \quad \left[\varepsilon,\ \frac{1}{6}(4 - 7 \theta_1 - 8 \theta_2) - \varepsilon \right].
\end{equation}
Note that in $\boldsymbol{A}_{10}$ we have
$$
2 \theta_1 + 2 \theta_2 - 1 < \theta_2 < \frac{1}{6}(5 - 8 \theta_1 - 8 \theta_2) < \frac{1}{3}(2 - 3 \theta_1 - 3 \theta_2).
$$
Hence, the Type-II range for $\boldsymbol{A}_{10}$ is
\begin{equation}
\left[\varepsilon,\ \frac{1}{6}(4 - 7 \theta_1 - 8 \theta_2) - \varepsilon \right] \cup \left[2 \theta_1 + 2 \theta_2 - 1 + \varepsilon,\ \frac{1}{3}(2 - 3 \theta_1 - 3 \theta_2) - \varepsilon \right].
\end{equation}

In $\boldsymbol{A}_{09} \cup \boldsymbol{A}_{10}$, we can replace the value of $\kappa = \frac{5 - 8 \theta}{6} - \varepsilon$ with $\frac{2 - 3 \theta}{3} - \varepsilon$, and thus $\kappa \geqslant \frac{1}{7} + \varepsilon$ when $\theta \leqslant \frac{11}{21} - \varepsilon$.

Comparing to previous cases, another region $\boldsymbol{C}$ covers $\boldsymbol{A}_{09} \cup \boldsymbol{A}_{10}$ when $\theta_2 < \frac{4 - 7 \theta_1 - \varepsilon}{12}$. Thus, in $\boldsymbol{A}_{09} \cup \boldsymbol{A}_{10}$ we can use an extra Type-II range
\begin{equation}
\left[\theta_1 + \varepsilon,\ 1 - \theta_1 - \varepsilon \right]
\end{equation}
comes from Lemma~\ref{l35} on those parts covered by $\boldsymbol{C}$.

Another important tool we can use in $\boldsymbol{A}_{09} \cup \boldsymbol{A}_{10}$ is Lemma~\ref{l39}, which gives asymptotic formulas of the form (53) for any sum that counts numbers with more than $4$ prime factors and all prime factors are larger than $x^{\frac{1}{7}}$. This lemma is very useful on estimating high-dimensional sums with $\boldsymbol{\alpha}_2 \in A \cup B$. Suppose that $7 \theta_1 + 12 \theta_2 \leqslant 4 - \varepsilon$ and $\theta \leqslant \frac{11}{21} - \varepsilon$. After using Buchstab's identity twice on the sum
\begin{equation}
\sum_{\kappa \leqslant \alpha_1 \leqslant \frac{3}{7} + \varepsilon} S\left(\mathcal{A}^{q_1 q_2}_{p_1}, p_1 \right),
\end{equation}
we get
\begin{equation}
\sum_{\substack{\kappa \leqslant \alpha_1 \leqslant \frac{3}{7} + \varepsilon \\ \kappa \leqslant \alpha_2 < \min\left(\alpha_1, \frac{1}{2}(1 - \alpha_1)\right) \\ \kappa \leqslant \alpha_3 < \min\left(\alpha_2, \frac{1}{2}(1 - \alpha_1 - \alpha_2)\right) }} S\left(\mathcal{A}^{q_1 q_2}_{p_1 p_2 p_3}, p_3 \right).
\end{equation}
We know that sum (102) only counts numbers with $4$ or more prime factors. Since $\kappa \geqslant \frac{1}{7} + \varepsilon$, we can use Lemma~\ref{l39} to give an asymptotic formula of the form (53).

\subsubsection{$\boldsymbol{A}_{11}$}
For $(\theta_1, \theta_2) \in \boldsymbol{A}_{11}$ we have 3 available Type-II information ranges:
\begin{equation}
\left[2 \theta - 1 + \varepsilon,\ \frac{1}{6}(5 - 8 \theta) - \varepsilon \right], \ \left[2 \theta - 1 + \varepsilon,\ \frac{1}{3}(2 - 3 \theta) - \varepsilon \right] \quad \text{and} \quad  \left[2 \theta_1 + \theta_2 - 1 + \varepsilon,\ \frac{1}{6}(4 - 7 \theta_1 - 8 \theta_2) - \varepsilon \right].
\end{equation}
The first range comes from Lemma~\ref{l23}, the second comes from Lemma~\ref{l33}, and the third comes from Lemma~\ref{l34}. We divide $\boldsymbol{A}_{11}$ into 2 subregions based on the overlapping conditions of these ranges.
$$
\boldsymbol{A}_{11} = \boldsymbol{A}_{1101} \cup \boldsymbol{A}_{1102},
$$
where
\begin{align}
\nonumber \boldsymbol{A}_{1101} =&\ \left\{ (\theta_1, \theta_2) : \frac{10}{21} < \theta_1 < \frac{1}{2},\ 1 - 2 \theta_1 \leqslant \theta_2 < \frac{1}{20}(10 - 19 \theta_1) \right\}, \\
\nonumber \boldsymbol{A}_{1102} =&\ \left\{ (\theta_1, \theta_2) : \frac{9}{20} < \theta_1 \leqslant \frac{23}{50},\ 1 - 2 \theta_1 \leqslant \theta_2 < \frac{1}{20}(11 - 20 \theta_1) \right. \\
\nonumber & \qquad \qquad \quad \text{ or } \frac{23}{50} < \theta_1 < \frac{10}{21},\ 1 - 2 \theta_1 \leqslant \theta_2 < \frac{1}{14}(10 - 19 \theta_1) \\
\nonumber & \left. \qquad \qquad \quad \text{ or } \frac{10}{21} \leqslant \theta_1 < \frac{1}{2},\ \frac{1}{20}(10 - 19 \theta_1) < \theta_2 < \frac{1}{14}(10 - 19 \theta_1) \right\}.
\end{align}
Note that we have $\frac{1}{6}(5 - 8 \theta_1 - 8 \theta_2) < \frac{1}{3}(2 - 3 \theta_1 - 3 \theta_2)$ for $(\theta_1, \theta_2) \in \boldsymbol{A}_{11}$ and $\theta < \frac{10}{19} < \frac{17}{32}$ for $(\theta_1, \theta_2) \in \boldsymbol{A}_{1101}$.

In $\boldsymbol{A}_{1101}$ we have 
$$
\frac{1}{6}(4 - 7 \theta_1 - 8 \theta_2) > 2 \theta_1 + 2 \theta_2 - 1.
$$
Hence, the Type-II range for $\boldsymbol{A}_{1101}$ is
\begin{equation}
\left[2 \theta_1 + \theta_2 - 1 + \varepsilon,\ \frac{1}{3}(2 - 3 \theta_1 - 3 \theta_2) - \varepsilon \right].
\end{equation}
The decompositions are similar to which in the case $\boldsymbol{A}_{10}$. A modification on the decompositions in this case is that we want to relax the condition $(\alpha_1, \ldots, \alpha_j, 2\theta-1+\varepsilon) \in \boldsymbol{S}_{j+1}$ in the definition of $\boldsymbol{U}_{j}$. Since the ``left endpoint'' of our Type-II interval now becomes $2 \theta_1 + \theta_2 - 1 + \varepsilon$, we want to prove a variant of Lemma~\ref{l213}, where $(\alpha_1, \ldots, \alpha_j, 2 \theta - 1 + \varepsilon) \in \boldsymbol{S}_{j+1}$ was replaced by a new condition $(\alpha_1, \ldots, \alpha_j, 2 \theta_1 + \theta_2 - 1 + \varepsilon) \in \boldsymbol{S}_{j+1}$ when $(\theta_1, \theta_2) \in \boldsymbol{A}_{1101}$. Note that we have done similar things in case $\boldsymbol{A}_{01}$, where that condition was replaced by $(\alpha_1, \ldots, \alpha_j, \varepsilon) \in \boldsymbol{S}_{j+1}$. Define
\begin{equation}
\boldsymbol{U}_{j}^{\prime \prime}(\theta_1, \theta_2) = \left\{\boldsymbol{\alpha}_{j}: \boldsymbol{\alpha}_{j} \in \boldsymbol{A}_{j},\ (\alpha_1, \ldots, \alpha_j, 2 \theta_1 + \theta_2 - 1 + \varepsilon) \in \boldsymbol{S}_{j+1} \right\}.
\end{equation}
Now we prove that
\begin{equation}
\sum_{\boldsymbol{\alpha}_{j} \in \boldsymbol{U}_{j}^{\prime \prime}} S\left(\mathcal{A}^{q_1 q_2}_{p_1 \cdots p_j}, x^{\frac{2 - 3 \theta}{3} - \varepsilon}\right)
\end{equation}
has an asymptotic formula of the form (53) when $(\theta_1, \theta_2) \in \boldsymbol{A}_{1101}$. Using Buchstab's identity repeatedly as in the proof of Lemma~\ref{l213}, we have
\begin{align}
\nonumber \sum_{\boldsymbol{\alpha}_{j} \in \boldsymbol{U}_{j}^{\prime \prime}} S\left(\mathcal{A}^{q_1 q_2}_{p_1 \cdots p_j}, x^{\frac{2 - 3 \theta}{3} - \varepsilon}\right) =&\ \sum_{\boldsymbol{\alpha}_{j} \in \boldsymbol{U}_{j}^{\prime \prime}} S\left(\mathcal{A}^{q_1 q_2}_{p_1 \cdots p_j}, \exp\left(\log x (\log \log x)^{-3} \right)\right) \\
\nonumber & - \sum_{\substack{\boldsymbol{\alpha}_{j} \in \boldsymbol{U}_{j}^{\prime \prime} \\ \boldsymbol{\alpha}_{j+1} \in \boldsymbol{A}_{j+1} \\ 2 \theta_1 + \theta_2 - 1 + \varepsilon \leqslant \alpha_{j+1} < \frac{2 - 3 \theta}{3} - \varepsilon }} S\left(\mathcal{A}^{q_1 q_2}_{p_1 \cdots p_{j+1}}, p_{j+1}\right) \\
\nonumber & + \sum_{k \geqslant j+1} (-1)^{k-j} \sum_{\substack{\boldsymbol{\alpha}_{j} \in \boldsymbol{U}_{j}^{\prime \prime} \\ \boldsymbol{\alpha}_{k} \in \boldsymbol{A}_{k} \\ \alpha_{j+1} < \frac{2 - 3 \theta}{3} - \varepsilon \\ \alpha_{j+1} + \cdots + \alpha_{k} < 2 \theta_1 + \theta_2 - 1 + \varepsilon }} S\left(\mathcal{A}^{q_1 q_2}_{p_1 \cdots p_{k}}, \exp\left(\log x (\log \log x)^{-3} \right)\right) \\
\nonumber & + \sum_{k \geqslant j+2} (-1)^{k-j} \sum_{\substack{\boldsymbol{\alpha}_{j} \in \boldsymbol{U}_{j}^{\prime \prime} \\ \boldsymbol{\alpha}_{k} \in \boldsymbol{A}_{k} \\ \alpha_{j+1} < \frac{2 - 3 \theta}{3} - \varepsilon \\ \alpha_{j+1} + \cdots + \alpha_{k-1} < 2 \theta_1 + \theta_2 - 1 + \varepsilon \leqslant \alpha_{j+1} + \cdots + \alpha_{k} }} S\left(\mathcal{A}^{q_1 q_2}_{p_1 \cdots p_{k}}, p_{k}\right) \\
=&\ S^{q_1 q_2}_{33111} - S^{q_1 q_2}_{33112} + S^{q_1 q_2}_{33113} + S^{q_1 q_2}_{33114}.
\end{align}
We can give an asymptotic formula of the form (53) for $S^{q_1 q_2}_{33112}$ by Lemma~\ref{l211} and our Type-II range (104). We can give asymptotic formulas of the form (53) for $S^{q_1 q_2}_{33111}$ and $S^{q_1 q_2}_{33113}$ by the condition $(\alpha_1, \ldots, \alpha_j, 2 \theta_1 + \theta_2 - 1 + \varepsilon) \in \boldsymbol{S}_{j+1}$ and an application of Lemma~\ref{l29}. In order to show that the sum (106) has an asymptotic formula of the form (53), we only need to give an asymptotic formula of the form (53) for $S^{q_1 q_2}_{33114}$. For this sum, we need to prove a variant of [\cite{676}, Lemma 11]: Let $\theta < \frac{17}{32}$ and $\boldsymbol{\alpha}_{k} \in \boldsymbol{A}_{k}$. Suppose that
$$
\alpha_{j+1} + \cdots + \alpha_{k-1} < 2 \theta_1 + \theta_2 - 1 + \varepsilon \leqslant \alpha_{j+1} + \cdots + \alpha_{k}
$$
and
$$
\alpha_{j+1} < \frac{2 - 3 \theta}{3} - \varepsilon
$$
for some $j$ ($0 \leqslant j \leqslant k - 1$), then $\boldsymbol{\alpha}_{k} \in \boldsymbol{G}_{k}$.

When $\alpha_{k} < \frac{5 - 9 \theta_1 - 6 \theta_2}{3} - 2 \varepsilon$, then
$$
2 \theta_1 + \theta_2 - 1 + \varepsilon \leqslant \alpha_{j+1} + \cdots + \alpha_{k} < (2 \theta_1 + \theta_2 - 1 + \varepsilon) + \frac{5 - 9 \theta_1 - 6 \theta_2}{3} - 2 \varepsilon = \frac{2 - 3 \theta}{3} - \varepsilon
$$
and $\boldsymbol{\alpha}_{k} \in \boldsymbol{G}_{k}$.

Suppose that $\alpha_{k} \geqslant \frac{5 - 9 \theta_1 - 6 \theta_2}{3} - 2 \varepsilon$. Since $\boldsymbol{\alpha}_{k} \in \boldsymbol{A}_{k}$, we have $\alpha_{k} < \alpha_{j+1} < \frac{2 - 3 \theta}{3} - \varepsilon$. Now we only need to prove that
$$
\frac{5 - 9 \theta_1 - 6 \theta_2}{3} - 2 \varepsilon \geqslant 2 \theta_1 + \theta_2 - 1 + \varepsilon,
$$
or
$$
15 \theta_1 + 9 \theta_2 \leqslant 8 - 9 \varepsilon
$$
when $(\theta_1, \theta_2) \in \boldsymbol{A}_{1101}$. A simple verification then completes the proof, and we know that the sum (106) has an asymptotic formula of the form (53).

Now, we can replace all $\boldsymbol{U}_{j}$ with $\boldsymbol{U}_{j}^{\prime \prime}$ in the decompositions. Since the whole $\boldsymbol{A}_{1101}$ is covered by $\boldsymbol{C}$, we can use the Type-II range $\left[\theta_1 + \varepsilon,\ 1 - \theta_1 - \varepsilon \right]$ from Lemma~\ref{l35} to give extra asymptotic formulas of the form (53). Since we have $(\theta_1, \theta_2) \in \boldsymbol{T}_{1}$ and $44 \theta_1 + 26 \theta_2 < 23$ for $(\theta_1, \theta_2) \in \boldsymbol{A}_{1101}$, the new three-dimensional Harman's sieve is also applicable here.

One difference between the decompositions for $(\theta_1, \theta_2) \in \boldsymbol{A}_{1101}$ and $(\theta_1, \theta_2) \in \boldsymbol{A}_{10}$ is that we can use Lemma~\ref{l39} to give lots of asymptotic formulas of the form (53) when $(\theta_1, \theta_2) \in \boldsymbol{A}_{10}$, but we cannot use Lemma~\ref{l39} when $(\theta_1, \theta_2) \in \boldsymbol{A}_{1101}$ since $2 \theta_1 + \theta_2 \geqslant 1$.

If we can use Lemma~\ref{l39}, then we need to have
\begin{equation}
2 \theta_1 + \theta_2 \leqslant 1 - 2 \varepsilon
\end{equation}
and
\begin{equation}
7 \theta_1 + 12 \theta_2 \leqslant 4 - \varepsilon.
\end{equation}
When $\theta \geqslant \frac{1}{2}$, (109) also yields $(\theta_1, \theta_2) \in \boldsymbol{C}$, and we can use the extra Type-II information given by Lemma~\ref{l35}. Note that (108) means that $\theta_1 + \varepsilon \leqslant 1 - \theta - \varepsilon$, the new Type-II range $\left(\theta_1 + \varepsilon,\ 1 - \theta_1 - \varepsilon \right)$ covers the whole interval $\left(1 - \theta - \varepsilon,\ \frac{1}{2}\right)$. If we want to use the new three-dimensional Harman's sieve with $\theta \geqslant \frac{16}{31}$, then the ``left endpoint'' of our Type-II range cannot be $2 \theta - 1 + \varepsilon$, which means that we need Lemma~\ref{l34} to make our Type-II range starts from $2 \theta_1 + \theta_2 - 1 + \varepsilon$. The requirement of that is
$$
\frac{1}{6}(4 - 7 \theta_1 - 8 \theta_2) - \varepsilon \geqslant 2 \theta - 1 + \varepsilon,
$$
or
\begin{equation}
19 \theta_1 + 20 \theta_2 \leqslant 10 - 12 \varepsilon.
\end{equation}
The readers can find that (108)--(110) are the conditions (1.3)--(1.5) in [\cite{MaynardLargeModuliI}, Theorem 1.1]. In fact, the above process gives a proof of [\cite{MaynardLargeModuliI}, Theorem 1.1]: Since $\theta \leqslant \frac{11}{21} - \varepsilon$ is satisfied automatically by (108)--(110), we have $\frac{2 - 3 \theta}{3} - \varepsilon \geqslant \frac{1}{7}$ and the sums after using the new three-dimensional Harman's sieve will have more than $4$ prime factors by earlier discussions. Lemma~\ref{l35} gives an asymptotic formula of the form (53) for the case $1 - \theta - \varepsilon < \alpha_1 < \frac{1}{2}$, and the sums
$$
S\left(\mathcal{A}^{q_1 q_2}, x^{\frac{2 - 3 \theta}{3} - \varepsilon} \right), \quad \sum_{\alpha_1 < \frac{3}{7} + \varepsilon} S\left(\mathcal{A}^{q_1 q_2}_{p_1}, x^{\frac{2 - 3 \theta}{3} - \varepsilon} \right), \quad \sum_{\boldsymbol{\alpha}_2 \in A \cup B} S\left(\mathcal{A}^{q_1 q_2}_{p_1 p_2}, x^{\frac{2 - 3 \theta}{3} - \varepsilon} \right)
$$
after straightforward decompositions have asymptotic formulas of the form (53) by Lemma~\ref{l212}, Lemma~\ref{l217} and our Type-II range $\left[\varepsilon,\ \frac{2 - 3 \theta}{3} - \varepsilon \right]$ under (108)--(110). The remaining sums (without the new three-dimensional Harman's sieve) can be dealt with Lemma~\ref{l39} since they only count numbers with $4$ or more prime factors.

In $\boldsymbol{A}_{1102}$ we have 
$$
\frac{1}{6}(4 - 7 \theta_1 - 8 \theta_2) < 2 \theta_1 + 2 \theta_2 - 1.
$$
Hence, the Type-II range for $\boldsymbol{A}_{1102}$ is
\begin{equation}
\left[2 \theta_1 + \theta_2 - 1 + \varepsilon,\ \frac{1}{6}(4 - 7 \theta_1 - 8 \theta_2) - \varepsilon \right] \cup \left[2 \theta_1 + 2 \theta_2 - 1 + \varepsilon,\ \frac{1}{3}(2 - 3 \theta_1 - 3 \theta_2) - \varepsilon \right].
\end{equation}
The decompositions are similar to which in the case $\boldsymbol{A}_{10}$. We can still use the Type-II range $\left[\theta_1 + \varepsilon,\ 1 - \theta_1 - \varepsilon \right]$ from Lemma~\ref{l35} to give extra asymptotic formulas of the form (53) for parts of $\boldsymbol{A}_{1102}$.

\subsubsection{$\boldsymbol{A}_{12}$}
For $(\theta_1, \theta_2) \in \boldsymbol{A}_{12}$ we have 2 available Type-II information ranges:
\begin{equation}
\left[2 \theta_1 + 2 \theta_2 - 1 + \varepsilon,\ \frac{1}{6}(5 - 8 \theta_1 - 8 \theta_2) - \varepsilon \right] \quad \text{and} \quad \left[2 \theta_1 + 2 \theta_2 - 1 + \varepsilon,\ \frac{1}{3}(2 - 3 \theta_1 - 3 \theta_2) - \varepsilon \right].
\end{equation}
The first range comes from Lemma~\ref{l23}, and the second comes from Lemma~\ref{l33}. Note that in $\boldsymbol{A}_{12}$ we have 
$$
\frac{1}{6}(5 - 8 \theta_1 - 8 \theta_2) < \frac{1}{3}(2 - 3 \theta_1 - 3 \theta_2).
$$
Hence, the Type-II range for $\boldsymbol{A}_{12}$ is
\begin{equation}
\left[2 \theta_1 + 2 \theta_2 - 1 + \varepsilon,\ \frac{1}{3}(2 - 3 \theta_1 - 3 \theta_2) - \varepsilon \right].
\end{equation}
The decompositions are similar to which in the case $\boldsymbol{A}_{10}$. We can still use the Type-II range $\left[\theta_1 + \varepsilon,\ 1 - \theta_1 - \varepsilon \right]$ from Lemma~\ref{l35} to give extra asymptotic formulas of the form (53) for parts of $\boldsymbol{A}_{12}$.

\subsubsection{$\boldsymbol{A}_{13}$}
For $(\theta_1, \theta_2) \in \boldsymbol{A}_{13}$ we have 2 available Type-II information ranges:
\begin{equation}
\left[2 \theta_1 + 2 \theta_2 - 1 + \varepsilon,\ \frac{1}{3}(2 - 3 \theta_1 - 3 \theta_2) - \varepsilon \right] \quad \text{and} \quad \left[2 \theta_1 + \theta_2 - 1 + \varepsilon,\ \frac{1}{6}(4 - 7 \theta_1 - 8 \theta_2) - \varepsilon \right].
\end{equation}
The first range comes from Lemma~\ref{l33}, and the second comes from Lemma~\ref{l34}. Note that Lemma~\ref{l23} becomes trivial in this region since $\theta_1 + \theta_2 \geqslant \frac{11}{20}$. Hence, the Type-II range for $\boldsymbol{A}_{13}$ is
\begin{equation}
\left[2 \theta_1 + \theta_2 - 1 + \varepsilon,\ \frac{1}{6}(4 - 7 \theta_1 - 8 \theta_2) - \varepsilon \right] \cup \left[2 \theta_1 + 2 \theta_2 - 1 + \varepsilon,\ \frac{1}{3}(2 - 3 \theta_1 - 3 \theta_2) - \varepsilon \right].
\end{equation}
The decompositions are similar to which in the case $\boldsymbol{A}_{05}$.

\subsubsection{$\boldsymbol{A}_{14}$}
For $(\theta_1, \theta_2) \in \boldsymbol{A}_{14}$, we only have 1 available Type-II information range comes from Lemma~\ref{l33}:
\begin{equation}
\left[2 \theta_1 + 2 \theta_2 - 1 + \varepsilon,\ \frac{1}{3}(2 - 3 \theta_1 - 3 \theta_2) - \varepsilon \right].
\end{equation}
Again, the decompositions are similar to which in the case $\boldsymbol{A}_{05}$.

\subsubsection{$\boldsymbol{A}_{15}$}
For $(\theta_1, \theta_2) \in \boldsymbol{A}_{15}$ we have 3 available Type-II information ranges:
\begin{equation}
\left[2 \theta - 1 + \varepsilon,\ \frac{1}{6}(5 - 8 \theta) - \varepsilon \right], \ \left[2 \theta - 1 + \varepsilon,\ \frac{1}{6}(5 - 8 \theta_1 - 6 \theta_2) - \varepsilon \right] \quad \text{and} \quad \left[2 \theta_1 + \theta_2 - 1 + \varepsilon,\ \frac{1}{6}(4 - 7 \theta_1 - 8 \theta_2) - \varepsilon \right].
\end{equation}
The first range comes from Lemma~\ref{l23}, the second comes from Lemma~\ref{l33}, and the third comes from Lemma~\ref{l34}. We divide $\boldsymbol{A}_{15}$ into 2 subregions based on the overlapping conditions of these ranges.
$$
\boldsymbol{A}_{15} = \boldsymbol{A}_{1501} \cup \boldsymbol{A}_{1502},
$$
where
\begin{align}
\nonumber \boldsymbol{A}_{1501} =&\ \left\{ (\theta_1, \theta_2) : \frac{1}{2} \leqslant \theta_1 < \frac{10}{19},\ 0 < \theta_2 < \frac{1}{20}(10 - 19 \theta_1) \right\}, \\
\nonumber \boldsymbol{A}_{1502} =&\ \left\{ (\theta_1, \theta_2) : \frac{1}{2} \leqslant \theta_1 < \frac{10}{19},\ \frac{1}{20}(10 - 19 \theta_1) \leqslant \theta_2 < \frac{1}{14}(10 - 19 \theta_1) \right\}.
\end{align}
Note that we have $\theta < \frac{10}{19} < \frac{17}{32}$ for $(\theta_1, \theta_2) \in \boldsymbol{A}_{1501}$.

In $\boldsymbol{A}_{1501}$ we have 
$$
\frac{1}{6}(4 - 7 \theta_1 - 8 \theta_2) > 2 \theta_1 + 2 \theta_2 - 1.
$$
Hence, the Type-II range for $\boldsymbol{A}_{1501}$ is
\begin{equation}
\left[2 \theta_1 + \theta_2 - 1 + \varepsilon,\ \frac{1}{6}(5 - 8 \theta_1 - 6 \theta_2) - \varepsilon \right].
\end{equation}
The decompositions are similar to which in the case $\boldsymbol{A}_{1101}$. We want to replace $\boldsymbol{U}_{j}$ with $\boldsymbol{U}_{j}^{\prime \prime}$ in the decompositions as in $\boldsymbol{A}_{1101}$. Thus, we need to prove that
\begin{equation}
\sum_{\boldsymbol{\alpha}_{j} \in \boldsymbol{U}_{j}^{\prime \prime}} S\left(\mathcal{A}^{q_1 q_2}_{p_1 \cdots p_j}, x^{\frac{5 - 8 \theta_1 - 6 \theta_2}{6} - \varepsilon}\right)
\end{equation}
has an asymptotic formula of the form (53) when $(\theta_1, \theta_2) \in \boldsymbol{A}_{1501}$. Using Buchstab's identity repeatedly as in the proof of Lemma~\ref{l213}, we have
\begin{align}
\nonumber \sum_{\boldsymbol{\alpha}_{j} \in \boldsymbol{U}_{j}^{\prime \prime}} S\left(\mathcal{A}^{q_1 q_2}_{p_1 \cdots p_j}, x^{\frac{5 - 8 \theta_1 - 6 \theta_2}{6} - \varepsilon}\right) =&\ \sum_{\boldsymbol{\alpha}_{j} \in \boldsymbol{U}_{j}^{\prime \prime}} S\left(\mathcal{A}^{q_1 q_2}_{p_1 \cdots p_j}, \exp\left(\log x (\log \log x)^{-3} \right)\right) \\
\nonumber & - \sum_{\substack{\boldsymbol{\alpha}_{j} \in \boldsymbol{U}_{j}^{\prime \prime} \\ \boldsymbol{\alpha}_{j+1} \in \boldsymbol{A}_{j+1} \\ 2 \theta_1 + \theta_2 - 1 + \varepsilon \leqslant \alpha_{j+1} < \frac{5 - 8 \theta_1 - 6 \theta_2}{6} - \varepsilon }} S\left(\mathcal{A}^{q_1 q_2}_{p_1 \cdots p_{j+1}}, p_{j+1}\right) \\
\nonumber & + \sum_{k \geqslant j+1} (-1)^{k-j} \sum_{\substack{\boldsymbol{\alpha}_{j} \in \boldsymbol{U}_{j}^{\prime \prime} \\ \boldsymbol{\alpha}_{k} \in \boldsymbol{A}_{k} \\ \alpha_{j+1} < \frac{5 - 8 \theta_1 - 6 \theta_2}{6} - \varepsilon \\ \alpha_{j+1} + \cdots + \alpha_{k} < 2 \theta_1 + \theta_2 - 1 + \varepsilon }} S\left(\mathcal{A}^{q_1 q_2}_{p_1 \cdots p_{k}}, \exp\left(\log x (\log \log x)^{-3} \right)\right) \\
\nonumber & + \sum_{k \geqslant j+2} (-1)^{k-j} \sum_{\substack{\boldsymbol{\alpha}_{j} \in \boldsymbol{U}_{j}^{\prime \prime} \\ \boldsymbol{\alpha}_{k} \in \boldsymbol{A}_{k} \\ \alpha_{j+1} < \frac{5 - 8 \theta_1 - 6 \theta_2}{6} - \varepsilon \\ \alpha_{j+1} + \cdots + \alpha_{k-1} < 2 \theta_1 + \theta_2 - 1 + \varepsilon \leqslant \alpha_{j+1} + \cdots + \alpha_{k} }} S\left(\mathcal{A}^{q_1 q_2}_{p_1 \cdots p_{k}}, p_{k}\right) \\
=&\ S^{q_1 q_2}_{33151} - S^{q_1 q_2}_{33152} + S^{q_1 q_2}_{33153} + S^{q_1 q_2}_{33154}.
\end{align}
We can give an asymptotic formula of the form (53) for $S^{q_1 q_2}_{33152}$ by Lemma~\ref{l211} and our Type-II range (129). We can give asymptotic formulas of the form (53) for $S^{q_1 q_2}_{33151}$ and $S^{q_1 q_2}_{33153}$ by the condition $(\alpha_1, \ldots, \alpha_j, 2 \theta_1 + \theta_2 - 1 + \varepsilon) \in \boldsymbol{S}_{j+1}$ and an application of Lemma~\ref{l29}. Now we only need to show that $S^{q_1 q_2}_{33154}$ has an asymptotic formula of the form (53). Let $\theta < \frac{17}{32}$ and $\boldsymbol{\alpha}_{k} \in \boldsymbol{A}_{k}$. Suppose that
$$
\alpha_{j+1} + \cdots + \alpha_{k-1} < 2 \theta_1 + \theta_2 - 1 + \varepsilon \leqslant \alpha_{j+1} + \cdots + \alpha_{k}
$$
and
$$
\alpha_{j+1} < \frac{5 - 8 \theta_1 - 6 \theta_2}{6} - \varepsilon
$$
for some $j$ ($0 \leqslant j \leqslant k - 1$), then $\boldsymbol{\alpha}_{k} \in \boldsymbol{G}_{k}$.

When $\alpha_{k} < \frac{11 - 20 \theta_1 - 12 \theta_2}{6} - 2 \varepsilon$, then
$$
2 \theta_1 + \theta_2 - 1 + \varepsilon \leqslant \alpha_{j+1} + \cdots + \alpha_{k} < (2 \theta_1 + \theta_2 - 1 + \varepsilon) + \frac{11 - 20 \theta_1 - 12 \theta_2}{6} - 2 \varepsilon = \frac{5 - 8 \theta_1 - 6 \theta_2}{6} - \varepsilon
$$
and $\boldsymbol{\alpha}_{k} \in \boldsymbol{G}_{k}$.

Suppose that $\alpha_{k} \geqslant \frac{11 - 20 \theta_1 - 12 \theta_2}{6} - 2 \varepsilon$. Since $\boldsymbol{\alpha}_{k} \in \boldsymbol{A}_{k}$, we have $\alpha_{k} < \alpha_{j+1} < \frac{5 - 8 \theta_1 - 6 \theta_2}{6} - \varepsilon$. Now we only need to prove that
$$
\frac{11 - 20 \theta_1 - 12 \theta_2}{6} - 2 \varepsilon \geqslant 2 \theta_1 + \theta_2 - 1 + \varepsilon,
$$
or
$$
32 \theta_1 + 18 \theta_2 \leqslant 17 - 18 \varepsilon
$$
when $(\theta_1, \theta_2) \in \boldsymbol{A}_{1501}$. A simple verification then completes the proof, and we know that the sum (119) has an asymptotic formula of the form (53).

Since $\theta_1 \geqslant \frac{1}{2}$, Lemma~\ref{l35} is not applicable here. However, we can still use the new three-dimensional Harman's sieve for parts of $\boldsymbol{A}_{1501}$ covered by $\boldsymbol{T}_{1}$, since we have $44 \theta_1 + 26 \theta_2 < 23$ for $(\theta_1, \theta_2) \in \boldsymbol{A}_{1501} \cap \boldsymbol{T}_{1}$.

In $\boldsymbol{A}_{1502}$ we have 
$$
\frac{1}{6}(4 - 7 \theta_1 - 8 \theta_2) < 2 \theta_1 + 2 \theta_2 - 1.
$$
Hence, the Type-II range for $\boldsymbol{A}_{1502}$ is
\begin{equation}
\left[2 \theta_1 + \theta_2 - 1 + \varepsilon,\ \frac{1}{6}(4 - 7 \theta_1 - 8 \theta_2) - \varepsilon \right] \cup \left[2 \theta_1 + 2 \theta_2 - 1 + \varepsilon,\ \frac{1}{6}(5 - 8 \theta_1 - 6 \theta_2) - \varepsilon \right].
\end{equation}
The decompositions are similar to which in the case $\boldsymbol{A}_{0302}$.

\subsubsection{$\boldsymbol{A}_{16}$}
For $(\theta_1, \theta_2) \in \boldsymbol{A}_{16}$ we have 2 available Type-II information ranges:
\begin{equation}
\left[2 \theta_1 + 2 \theta_2 - 1 + \varepsilon,\ \frac{1}{6}(5 - 8 \theta_1 - 8 \theta_2) - \varepsilon \right] \quad \text{and} \quad \left[2 \theta_1 + 2 \theta_2 - 1 + \varepsilon,\ \frac{1}{6}(5 - 8 \theta_1 - 6 \theta_2) - \varepsilon \right].
\end{equation}
The first range comes from Lemma~\ref{l23}, and the second comes from Lemma~\ref{l33}. Since $\theta_2 > 0$, the Type-II range for $\boldsymbol{A}_{16}$ is
\begin{equation}
\left[2 \theta_1 + 2 \theta_2 - 1 + \varepsilon,\ \frac{1}{6}(5 - 8 \theta_1 - 6 \theta_2) - \varepsilon \right].
\end{equation}
The decompositions are similar to which in the case $\boldsymbol{A}_{1502}$.

\subsubsection{$\boldsymbol{A}_{17}$}
For $(\theta_1, \theta_2) \in \boldsymbol{A}_{17}$, we only have 1 available Type-II information range comes from Lemma~\ref{l33}:
\begin{equation}
\left[2 \theta_1 + 2 \theta_2 - 1 + \varepsilon,\ \frac{1}{6}(5 - 8 \theta_1 - 6 \theta_2) - \varepsilon \right].
\end{equation}
The decompositions are similar to which in the case $\boldsymbol{A}_{05}$.

Next we assume that $(\theta_1, \theta_2) \in \boldsymbol{B} \backslash \boldsymbol{A}$. From the condition $Q_1^8 Q_2^7 M_2^6 < x^{4 - 13 \varepsilon}$ in Lemma~\ref{l34}, we can find that the Type-II ranges come from Lemma~\ref{l34} are of the form
\begin{equation}
\left(\cdots, \frac{1}{6}(4 - 7 \theta_1 - 8 \theta_2) - \varepsilon \right] \quad \text{and} \quad \left(\cdots, \frac{1}{6}(4 - 8 \theta_1 - 7 \theta_2) - \varepsilon \right],
\end{equation}
where the ``left endpoint'' may be $\varepsilon$ or $2 \theta_1 + \theta_2 - 1 + \varepsilon$. Since Lemma~\ref{l33} and Lemma~\ref{l35} are not applicable in this case (we have $\boldsymbol{C} \subset \boldsymbol{A}$), the only situation that Lemma~\ref{l34} brings improvements is when we have
\begin{equation}
\frac{1}{6}(4 - 7 \theta_1 - 8 \theta_2) > 2 \theta - 1 \quad \text{or} \quad \frac{1}{6}(4 - 8 \theta_1 - 7 \theta_2) > 2 \theta - 1.
\end{equation}
In this case, we can make the ``left endpoint'' of our Type-II range smaller, hence to relax the condition $(\alpha_1, \ldots, \alpha_j, 2\theta-1+\varepsilon) \in \boldsymbol{S}_{j+1}$ in the definition of $\boldsymbol{U}_{j}$. We may replace $\boldsymbol{\alpha}_{j} \in \boldsymbol{U}_{j}$ with $(\alpha_1, \ldots, \alpha_j, \varepsilon) \in \boldsymbol{S}_{j+1}$ or $\boldsymbol{\alpha}_{j} \in \boldsymbol{U}_{j}^{\prime \prime}$ when the ``left endpoint'' becomes $\varepsilon$ or $2 \theta_1 + \theta_2 - 1 + \varepsilon$. After checking the conditions, we find that a subregion of $\boldsymbol{B} \backslash \boldsymbol{A}$ that satisfies (126) is
\begin{align}
\nonumber \boldsymbol{B}_{01} =&\ \left\{ (\theta_1, \theta_2) : \frac{1}{4} < \theta_1 < \frac{10}{39},\ \frac{1}{2}(1 - 2 \theta_1) < \theta_2 \leqslant \theta_1 \right. \\
\nonumber & \qquad \qquad \quad \text{ or } \frac{10}{39} \leqslant \theta_1 \leqslant \frac{5}{14},\ \frac{1}{2}(1 - 2 \theta_1) < \theta_2 < \frac{1}{20}(10 - 19 \theta_1) \\
\nonumber & \left. \qquad \qquad \quad \text{ or } \frac{5}{14} < \theta_1 < \frac{50}{131},\ \frac{1}{9}(2 - 2 \theta_1) < \theta_2 < \frac{1}{20}(10 - 19 \theta_1) \right\}.
\end{align}
In this region, we have a Type-II range
\begin{equation}
\left[\varepsilon,\ \frac{1}{6}(5 - 8 \theta_1 - 8 \theta_2) - \varepsilon \right].
\end{equation}
The decompositions are similar to which in the case $\boldsymbol{A}_{0101}$. The decompositions in remaining parts of $\boldsymbol{B} \backslash \boldsymbol{A}$ stay the same as in \cite{LRB679} and Subsection 2.4.

Working on each case above carefully, we can get the following upper bounds for $C_1(\theta_1, \theta_2)$ and $C_1^{*}(\theta_1, \theta_2)$ ($0.5 < \theta \leqslant 0.55$):
\begin{center}
\begin{tabular}{|>{\centering\arraybackslash}p{0.6cm}|>{\centering\arraybackslash}p{0.6cm}|>{\centering\arraybackslash}p{0.6cm}|>{\centering\arraybackslash}p{0.6cm}|>{\centering\arraybackslash}p{0.6cm}|>{\centering\arraybackslash}p{0.6cm}|>{\centering\arraybackslash}p{0.6cm}|>{\centering\arraybackslash}p{0.6cm}|>{\centering\arraybackslash}p{0.6cm}|>{\centering\arraybackslash}p{0.6cm}|>{\centering\arraybackslash}p{0.6cm}|>{\centering\arraybackslash}p{0.6cm}|>{\centering\arraybackslash}p{0.6cm}|>{\centering\arraybackslash}p{0.6cm}|>{\centering\arraybackslash}p{0.6cm}|}
\hline \boldmath{$0.25$} & \tiny $1.0486$ & \tiny $1.7394$ & \tiny $1.8338$ & \tiny $1.9993$ & \tiny $2.3996$ & --- & --- & --- & --- & --- & --- & --- & --- & --- \\
\hline \boldmath{$0.24$} & $1 + \varepsilon$ & \tiny $1.0486$ & \tiny $1.7394$ & \tiny $1.8338$ & \tiny $1.9993$ & \tiny $2.3996$ & --- & --- & --- & --- & --- & --- & --- & --- \\
\hline \boldmath{$0.23$} & $1$ & $1 + \varepsilon$ & \tiny $1.0486$ & \tiny $1.7394$ & \tiny $1.8338$ & \tiny $1.9993$ & \tiny $2.3996$ & --- & --- & --- & --- & --- & --- & --- \\
\hline \boldmath{$0.22$} & $1$ & $1$ & $1 + \varepsilon$ & \tiny $1.0486$ & \tiny $1.7394$ & \tiny $1.8338$ & \tiny $1.9993$ & \tiny $2.3996$ & --- & --- & --- & --- & --- & --- \\
\hline \boldmath{$0.21$} & $1$ & $1$ & $1$ & $1 + \varepsilon$ & \tiny $1.0486$ & \tiny $1.7394$ & \tiny $1.8338$ & \tiny $1.9993$ & \tiny $2.3996$ & --- & --- & --- & --- & --- \\
\hline \boldmath{$0.20$} & $1$ & $1$ & $1$ & $1$ & $1 + \varepsilon$ & \tiny $1.0486$ & \tiny $1.7394$ & \tiny $1.8338$ & \tiny $1.9993$ & \tiny $2.3996$ & --- & --- & --- & --- \\
\hline \boldmath{$0.19$} & $1$ & $1$ & $1$ & $1$ & $1$ & $1 + \varepsilon$ & \tiny $1.0486$ & \tiny $1.7394$ & \tiny $1.8338$ & \tiny $1.9993$ & \tiny $2.3996$ & --- & --- & --- \\
\hline \boldmath{$0.18$} & $1$ & $1$ & $1$ & $1$ & $1$ & $1$ & $1 + \varepsilon$ & \tiny $1.0486$ & \tiny $1.7394$ & \tiny $1.8338$ & \tiny $1.9993$ & \tiny $2.3996$ & --- & --- \\
\hline \boldmath{$0.17$} & $1$ & $1$ & $1$ & $1$ & $1$ & $1$ & $1$ & $1 + \varepsilon$ & \tiny $1.0486$ & \tiny $1.7394$ & \tiny $1.8338$ & \tiny $1.9993$ & \tiny $2.3996$ & --- \\
\hline \boldmath{$0.16$} & $1$ & $1$ & $1$ & $1$ & $1$ & $1$ & $1$ & $1$ & $1 + \varepsilon$ & \tiny $1.0486$ & \tiny $1.7394$ & \tiny $1.8338$ & \tiny $1.9993$ & \tiny $2.3996$ \\
\hline \boldmath{$0.15$} & $1$ & $1$ & $1$ & $1$ & $1$ & $1$ & $1$ & $1$ & $1$ & $1 + \varepsilon$ & \tiny $1.0486$ & \tiny $1.7394$ & \tiny $1.8338$ & \tiny $1.9993$ \\
\hline \boldmath{$0.14$} & $1$ & $1$ & $1$ & $1$ & $1$ & $1$ & $1$ & $1$ & $1$ & $1$ & $1 + \varepsilon$ & \tiny $1.0486$ & \tiny $1.7394$ & \tiny $1.8338$ \\
\hline \boldmath{$0.13$} & $1$ & $1$ & $1$ & $1$ & $1$ & $1$ & $1$ & $1$ & $1$ & $1$ & $1$ & $1 + \varepsilon$ & \tiny $1.0486$ & \tiny $1.7394$ \\
\hline \boldmath{$0.12$} & $1$ & $1$ & $1$ & $1$ & $1$ & $1$ & $1$ & $1$ & $1$ & $1$ & $1$ & $1$ & $1 + \varepsilon$ & \tiny $1.0486$ \\
\hline \small \boldmath{$\theta_2 \backslash \theta_1$} & \boldmath{$0.26$} & \boldmath{$0.27$} & \boldmath{$0.28$} & \boldmath{$0.29$} & \boldmath{$0.30$} & \boldmath{$0.31$} & \boldmath{$0.32$} & \boldmath{$0.33$} & \boldmath{$0.34$} & \boldmath{$0.35$} & \boldmath{$0.36$} & \boldmath{$0.37$} & \boldmath{$0.38$} & \boldmath{$0.39$} \\
\hline
\end{tabular} \\
\textbf{Table 3.1: Upper Bounds for }\boldmath{$C_1(\theta_1, \theta_2)$} \textbf{(}\boldmath{$0.5 < \theta \leqslant 0.55$}\textbf{) 1/2}
\end{center}
\begin{center}
\begin{tabular}{|>{\centering\arraybackslash}p{0.6cm}|>{\centering\arraybackslash}p{0.6cm}|>{\centering\arraybackslash}p{0.6cm}|>{\centering\arraybackslash}p{0.6cm}|>{\centering\arraybackslash}p{0.6cm}|>{\centering\arraybackslash}p{0.6cm}|>{\centering\arraybackslash}p{0.6cm}|>{\centering\arraybackslash}p{0.6cm}|>{\centering\arraybackslash}p{0.6cm}|>{\centering\arraybackslash}p{0.6cm}|>{\centering\arraybackslash}p{0.6cm}|>{\centering\arraybackslash}p{0.6cm}|>{\centering\arraybackslash}p{0.6cm}|>{\centering\arraybackslash}p{0.6cm}|>{\centering\arraybackslash}p{0.6cm}|>{\centering\arraybackslash}p{0.6cm}|}
\hline \boldmath{$0.15$} & \tiny $2.3996$ & --- & --- & --- & --- & --- & --- & --- & --- & --- & --- & --- & --- & --- & --- \\
\hline \boldmath{$0.14$} & \tiny $1.9993$ & \tiny $2.3996$ & --- & --- & --- & --- & --- & --- & --- & --- & --- & --- & --- & --- & --- \\
\hline \boldmath{$0.13$} & \tiny $1.8250$ & \tiny $1.9993$ & \tiny $2.3996$ & --- & --- & --- & --- & --- & --- & --- & --- & --- & --- & --- & --- \\
\hline \boldmath{$0.12$} & \tiny $1.7394$ & \tiny $1.8182$ & \tiny $1.9851$ & \tiny $2.3996$ & --- & --- & --- & --- & --- & --- & --- & --- & --- & --- & --- \\
\hline \boldmath{$0.11$} & \tiny $1.0486$ & \tiny $1.1067$ & \tiny $1.8182$ & \tiny $1.9486$ & \tiny $2.3852$ & --- & --- & --- & --- & --- & --- & --- & --- & --- & --- \\
\hline \boldmath{$0.10$} & $1 + \varepsilon$ & \tiny $1.0486$ & \tiny $1.1111$ & \tiny $1.8182$ & \tiny $1.9486$ & \tiny $2.3511$ & --- & --- & --- & --- & --- & --- & --- & --- & --- \\
\hline \boldmath{$0.09$} & $1$ & $1$ & \tiny $1.0486$ & \tiny $1.1111$ & \tiny $1.8182$ & \tiny $1.9486$ & \tiny $2.3511$ & --- & --- & --- & --- & --- & --- & --- & --- \\
\hline \boldmath{$0.08$} & $1$ & $1$ & $1$ & $1$ & \tiny $1.1099$ & \tiny $1.8182$ & \tiny $1.9486$ & \tiny $2.3511$ & --- & --- & --- & --- & --- & --- & --- \\
\hline \boldmath{$0.07$} & $1$ & $1$ & $1$ & $1$ & $1$ & $1$ & \tiny $1.8182$ & \tiny $1.9486$ & \tiny $2.3511$ & --- & --- & --- & --- & --- & --- \\
\hline \boldmath{$0.06$} & $1$ & $1$ & $1$ & $1$ & $1$ & $1$ & $1$ & \tiny $1.8182$ & \tiny $1.9486$ & \tiny $2.3511$ & --- & --- & --- & --- & --- \\
\hline \boldmath{$0.05$} & $1$ & $1$ & $1$ & $1$ & $1$ & $1$ & $1$ & $1$ & \tiny $1.7215$ & \tiny $1.9486$ & \tiny $2.3511$ & --- & --- & --- & --- \\
\hline \boldmath{$0.04$} & $1$ & $1$ & $1$ & $1$ & $1$ & $1$ & $1$ & $1$ & $1 + \varepsilon$ & \tiny $1.7676$ & \tiny $1.9486$ & \tiny $2.3549$ & --- & --- & --- \\
\hline \boldmath{$0.03$} & $1$ & $1$ & $1$ & $1$ & $1$ & $1$ & $1$ & $1$ & $1$ & \tiny $1.0503$ & \tiny $1.8182$ & \tiny $1.9615$ & \tiny $2.3648$ & --- & --- \\
\hline \boldmath{$0.02$} & $1$ & $1$ & $1$ & $1$ & $1$ & $1$ & $1$ & $1$ & $1$ & $1 + \varepsilon$ & \tiny $1.1099$ & \tiny $1.8212$ & \tiny $1.9765$ & \tiny $2.3928$ & --- \\
\hline \boldmath{$0.01$} & $1$ & $1$ & $1$ & $1$ & $1$ & $1$ & $1$ & $1$ & $1$ & $1$ & \tiny $1.0482$ & \tiny $1.7415$ & \tiny $1.8271$ & \tiny $1.9893$ & \tiny $2.3935$ \\
\hline \small \boldmath{$\theta_2 \backslash \theta_1$} & \boldmath{$0.40$} & \boldmath{$0.41$} & \boldmath{$0.42$} & \boldmath{$0.43$} & \boldmath{$0.44$} & \boldmath{$0.45$} & \boldmath{$0.46$} & \boldmath{$0.47$} & \boldmath{$0.48$} & \boldmath{$0.49$} & \boldmath{$0.50$} & \boldmath{$0.51$} & \boldmath{$0.52$} & \boldmath{$0.53$} & \boldmath{$0.54$} \\
\hline
\end{tabular} \\
\textbf{Table 3.2: Upper Bounds for }\boldmath{$C_1(\theta_1, \theta_2)$} \textbf{(}\boldmath{$0.5 < \theta \leqslant 0.55$}\textbf{) 2/2}
\end{center}
\begin{center}
\begin{tabular}{|>{\centering\arraybackslash}p{0.6cm}|>{\centering\arraybackslash}p{0.6cm}|>{\centering\arraybackslash}p{0.6cm}|>{\centering\arraybackslash}p{0.6cm}|>{\centering\arraybackslash}p{0.6cm}|>{\centering\arraybackslash}p{0.6cm}|>{\centering\arraybackslash}p{0.6cm}|>{\centering\arraybackslash}p{0.6cm}|>{\centering\arraybackslash}p{0.6cm}|>{\centering\arraybackslash}p{0.6cm}|>{\centering\arraybackslash}p{0.6cm}|>{\centering\arraybackslash}p{0.6cm}|>{\centering\arraybackslash}p{0.6cm}|>{\centering\arraybackslash}p{0.6cm}|>{\centering\arraybackslash}p{0.6cm}|}
\hline \boldmath{$0.25$} & \tiny $1.0150$ & \tiny $1.7028$ & \tiny $1.8169$ & \tiny $1.9993$ & \tiny $2.3996$ & --- & --- & --- & --- & --- & --- & --- & --- & --- \\
\hline \boldmath{$0.24$} & $1 + \varepsilon$ & \tiny $1.0150$ & \tiny $1.7028$ & \tiny $1.8169$ & \tiny $1.9993$ & \tiny $2.3996$ & --- & --- & --- & --- & --- & --- & --- & --- \\
\hline \boldmath{$0.23$} & $1$ & $1 + \varepsilon$ & \tiny $1.0150$ & \tiny $1.7028$ & \tiny $1.8169$ & \tiny $1.9993$ & \tiny $2.3996$ & --- & --- & --- & --- & --- & --- & --- \\
\hline \boldmath{$0.22$} & $1$ & $1$ & $1 + \varepsilon$ & \tiny $1.0150$ & \tiny $1.7028$ & \tiny $1.8169$ & \tiny $1.9993$ & \tiny $2.3996$ & --- & --- & --- & --- & --- & --- \\
\hline \boldmath{$0.21$} & $1$ & $1$ & $1$ & $1 + \varepsilon$ & \tiny $1.0150$ & \tiny $1.7028$ & \tiny $1.8169$ & \tiny $1.9993$ & \tiny $2.3996$ & --- & --- & --- & --- & --- \\
\hline \boldmath{$0.20$} & $1$ & $1$ & $1$ & $1$ & $1 + \varepsilon$ & \tiny $1.0150$ & \tiny $1.7028$ & \tiny $1.8169$ & \tiny $1.9993$ & \tiny $2.3996$ & --- & --- & --- & --- \\
\hline \boldmath{$0.19$} & $1$ & $1$ & $1$ & $1$ & $1$ & $1 + \varepsilon$ & \tiny $1.0150$ & \tiny $1.7028$ & \tiny $1.8169$ & \tiny $1.9993$ & \tiny $2.3996$ & --- & --- & --- \\
\hline \boldmath{$0.18$} & $1$ & $1$ & $1$ & $1$ & $1$ & $1$ & $1 + \varepsilon$ & \tiny $1.0150$ & \tiny $1.7028$ & \tiny $1.8169$ & \tiny $1.9993$ & \tiny $2.3996$ & --- & --- \\
\hline \boldmath{$0.17$} & $1$ & $1$ & $1$ & $1$ & $1$ & $1$ & $1$ & $1 + \varepsilon$ & \tiny $1.0150$ & \tiny $1.7028$ & \tiny $1.8169$ & \tiny $1.9993$ & \tiny $2.3996$ & --- \\
\hline \boldmath{$0.16$} & $1$ & $1$ & $1$ & $1$ & $1$ & $1$ & $1$ & $1$ & $1 + \varepsilon$ & \tiny $1.0150$ & \tiny $1.7028$ & \tiny $1.8169$ & \tiny $1.9993$ & \tiny $2.3996$ \\
\hline \boldmath{$0.15$} & $1$ & $1$ & $1$ & $1$ & $1$ & $1$ & $1$ & $1$ & $1$ & $1 + \varepsilon$ & \tiny $1.0150$ & \tiny $1.7028$ & \tiny $1.8169$ & \tiny $1.9993$ \\
\hline \boldmath{$0.14$} & $1$ & $1$ & $1$ & $1$ & $1$ & $1$ & $1$ & $1$ & $1$ & $1$ & $1 + \varepsilon$ & \tiny $1.0150$ & \tiny $1.7028$ & \tiny $1.8169$ \\
\hline \boldmath{$0.13$} & $1$ & $1$ & $1$ & $1$ & $1$ & $1$ & $1$ & $1$ & $1$ & $1$ & $1$ & $1 + \varepsilon$ & \tiny $1.0150$ & \tiny $1.7028$ \\
\hline \boldmath{$0.12$} & $1$ & $1$ & $1$ & $1$ & $1$ & $1$ & $1$ & $1$ & $1$ & $1$ & $1$ & $1$ & $1 + \varepsilon$ & \tiny $1.0150$ \\
\hline \small \boldmath{$\theta_2 \backslash \theta_1$} & \boldmath{$0.26$} & \boldmath{$0.27$} & \boldmath{$0.28$} & \boldmath{$0.29$} & \boldmath{$0.30$} & \boldmath{$0.31$} & \boldmath{$0.32$} & \boldmath{$0.33$} & \boldmath{$0.34$} & \boldmath{$0.35$} & \boldmath{$0.36$} & \boldmath{$0.37$} & \boldmath{$0.38$} & \boldmath{$0.39$} \\
\hline
\end{tabular} \\
\textbf{Table 3.3: Upper Bounds for }\boldmath{$C_1^{*}(\theta_1, \theta_2)$} \textbf{(}\boldmath{$0.5 < \theta \leqslant 0.55$}\textbf{) 1/2}
\end{center}
\begin{center}
\begin{tabular}{|>{\centering\arraybackslash}p{0.6cm}|>{\centering\arraybackslash}p{0.6cm}|>{\centering\arraybackslash}p{0.6cm}|>{\centering\arraybackslash}p{0.6cm}|>{\centering\arraybackslash}p{0.6cm}|>{\centering\arraybackslash}p{0.6cm}|>{\centering\arraybackslash}p{0.6cm}|>{\centering\arraybackslash}p{0.6cm}|>{\centering\arraybackslash}p{0.6cm}|>{\centering\arraybackslash}p{0.6cm}|>{\centering\arraybackslash}p{0.6cm}|>{\centering\arraybackslash}p{0.6cm}|>{\centering\arraybackslash}p{0.6cm}|>{\centering\arraybackslash}p{0.6cm}|>{\centering\arraybackslash}p{0.6cm}|>{\centering\arraybackslash}p{0.6cm}|}
\hline \boldmath{$0.15$} & \tiny $2.3996$ & --- & --- & --- & --- & --- & --- & --- & --- & --- & --- & --- & --- & --- & --- \\
\hline \boldmath{$0.14$} & \tiny $1.9993$ & \tiny $2.3996$ & --- & --- & --- & --- & --- & --- & --- & --- & --- & --- & --- & --- & --- \\
\hline \boldmath{$0.13$} & \tiny $1.8169$ & \tiny $1.9993$ & \tiny $2.3996$ & --- & --- & --- & --- & --- & --- & --- & --- & --- & --- & --- & --- \\
\hline \boldmath{$0.12$} & \tiny $1.7028$ & \tiny $1.8120$ & \tiny $1.9851$ & \tiny $2.3996$ & --- & --- & --- & --- & --- & --- & --- & --- & --- & --- & --- \\
\hline \boldmath{$0.11$} & \tiny $1.0150$ & \tiny $1.0678$ & \tiny $1.8120$ & \tiny $1.9486$ & \tiny $2.3852$ & --- & --- & --- & --- & --- & --- & --- & --- & --- & --- \\
\hline \boldmath{$0.10$} & $1 + \varepsilon$ & \tiny $1.0150$ & \tiny $1.0687$ & \tiny $1.8120$ & \tiny $1.9486$ & \tiny $2.3511$ & --- & --- & --- & --- & --- & --- & --- & --- & --- \\
\hline \boldmath{$0.09$} & $1$ & $1$ & \tiny $1.0150$ & \tiny $1.0670$ & \tiny $1.8120$ & \tiny $1.9486$ & \tiny $2.3511$ & --- & --- & --- & --- & --- & --- & --- & --- \\
\hline \boldmath{$0.08$} & $1$ & $1$ & $1$ & $1$ & \tiny $1.0636$ & \tiny $1.8120$ & \tiny $1.9486$ & \tiny $2.3511$ & --- & --- & --- & --- & --- & --- & --- \\
\hline \boldmath{$0.07$} & $1$ & $1$ & $1$ & $1$ & $1$ & $1$ & \tiny $1.8120$ & \tiny $1.9486$ & \tiny $2.3511$ & --- & --- & --- & --- & --- & --- \\
\hline \boldmath{$0.06$} & $1$ & $1$ & $1$ & $1$ & $1$ & $1$ & $1$ & \tiny $1.8120$ & \tiny $1.9486$ & \tiny $2.3511$ & --- & --- & --- & --- & --- \\
\hline \boldmath{$0.05$} & $1$ & $1$ & $1$ & $1$ & $1$ & $1$ & $1$ & $1$ & \tiny $1.6992$ & \tiny $1.9486$ & \tiny $2.3511$ & --- & --- & --- & --- \\
\hline \boldmath{$0.04$} & $1$ & $1$ & $1$ & $1$ & $1$ & $1$ & $1$ & $1$ & $1 + \varepsilon$ & \tiny $1.7406$ & \tiny $1.9486$ & \tiny $2.3549$ & --- & --- & --- \\
\hline \boldmath{$0.03$} & $1$ & $1$ & $1$ & $1$ & $1$ & $1$ & $1$ & $1$ & $1$ & \tiny $1.0187$ & \tiny $1.8120$ & \tiny $1.9615$ & \tiny $2.3648$ & --- & --- \\
\hline \boldmath{$0.02$} & $1$ & $1$ & $1$ & $1$ & $1$ & $1$ & $1$ & $1$ & $1$ & $1 + \varepsilon$ & \tiny $1.0636$ & \tiny $1.8169$ & \tiny $1.9765$ & \tiny $2.3928$ & --- \\
\hline \boldmath{$0.01$} & $1$ & $1$ & $1$ & $1$ & $1$ & $1$ & $1$ & $1$ & $1$ & $1$ & \tiny $1.0150$ & \tiny $1.6991$ & \tiny $1.8169$ & \tiny $1.9893$ & \tiny $2.3935$ \\
\hline \small \boldmath{$\theta_2 \backslash \theta_1$} & \boldmath{$0.40$} & \boldmath{$0.41$} & \boldmath{$0.42$} & \boldmath{$0.43$} & \boldmath{$0.44$} & \boldmath{$0.45$} & \boldmath{$0.46$} & \boldmath{$0.47$} & \boldmath{$0.48$} & \boldmath{$0.49$} & \boldmath{$0.50$} & \boldmath{$0.51$} & \boldmath{$0.52$} & \boldmath{$0.53$} & \boldmath{$0.54$} \\
\hline
\end{tabular} \\
\textbf{Table 3.4: Upper Bounds for }\boldmath{$C_1^{*}(\theta_1, \theta_2)$} \textbf{(}\boldmath{$0.5 < \theta \leqslant 0.55$}\textbf{) 2/2}
\end{center}

\subsection{Lower Bounds}
We shall construct the minorant $\rho_0(n)$ and prove lower bounds for $C_0(\theta_1, \theta_2)$ and $C_0^{*}(\theta_1, \theta_2)$ in this subsection. Before our final decompositions, we first mention some existing results of $C_0(\theta_1, \theta_2)$ and $C_0^{*}(\theta_1, \theta_2)$.
\begin{theorem}\label{t311}
The functions $C_0(\theta_1, \theta_2)$ and $C_0^{*}(\theta_1, \theta_2)$ satisfy the following conditions:

(1). $C_0(\theta_1, \theta_2) = C_0(\theta_2, \theta_1)$, $C_0^{*}(\theta_1, \theta_2) = C_0^{*}(\theta_2, \theta_1)$;

(2). $C_0(\theta_1, \theta_2) = C_0^{*}(\theta_1, \theta_2) = 1$ for all $\theta_1, \theta_2$ satisfy $\theta_1 + \theta_2 \leqslant 0.5 - \varepsilon$; 

(3). $C_0(\theta_1, \theta_2) = C_0^{*}(\theta_1, \theta_2) = 1$ for all $\theta_1, \theta_2$ satisfy $2 \theta_1 + \theta_2 \leqslant 1 - \varepsilon$, $7 \theta_1 + 12 \theta_2 \leqslant 4 - \varepsilon$ and $19 \theta_1 + 20 \theta_2 \leqslant 10 - \varepsilon$;

(4). $C_0(\theta_1, \theta_2) \geqslant C_0(\theta_1 + \theta_2)$, $C_0^{*}(\theta_1, \theta_2) \geqslant C_0^{*}(\theta_1 + \theta_2)$ for $0.5 \leqslant \theta_1 + \theta_2 \leqslant 1$;

(5). $C_0^{*}(\theta_1, \theta_2) \geqslant 1 - \varepsilon$ for all $\theta_1, \theta_2$ satisfy $\theta_1 + \theta_2 = 0.5$.
\end{theorem}
\begin{proof}
Statement (1) is obvious. Statements (2)--(3) follow easily from the Bombieri--Vinogradov Theorem and [\cite{MaynardLargeModuliI}, Theorem 1.1]. Statement (4) holds trivially by the work done in Section 2. When there are no new arithmetic information inputs outside of those in Section 2, we use $C_0(\theta_1 + \theta_2)$ and $C_0^{*}(\theta_1 + \theta_2)$ as lower bounds for $C_0(\theta_1, \theta_2)$ and $C_0^{*}(\theta_1, \theta_2)$ respectively. Statement (5) follows from the Statement (2) of Theorem~\ref{t223} together with Statement (4).
\end{proof}

We first consider the case $C_0(\theta_1, \theta_2)$. By the discussions in Subsection 2.5, we need to give an asymptotic formula of the form (53) for the sum
\begin{equation}
\sum_{1 - \theta - \varepsilon < \alpha_1 < \frac{1}{2}} S\left(\mathcal{A}^{q_1 q_2}_{p_1}, p_1 \right).
\end{equation}
The only result that gives arithmetic information for this sum is Lemma~\ref{l35}, hence we need (108) and (109) to give an asymptotic formula of the form (53) for the sum (128). Assume that (108) and (109) hold for $(\theta_1, \theta_2)$. If (110) also holds for $(\theta_1, \theta_2)$, then we have $C_0(\theta_1, \theta_2) = 1$ by the Statement (3) of Theorem~\ref{t311}. If we have $19 \theta_1 + 20 \theta_2 > 10 - \varepsilon$, however, we can obtain nontrivial lower bounds for $C_0(\theta_1, \theta_2)$ less than 1. By the discussions in Subsection 2.5, we cannot get any nontrivial lower bound for $C_0(\theta_1, \theta_2)$ if either (108) or (109) is not fulfilled.

Now we focus on the case $C_0^{*}(\theta_1, \theta_2)$. From here to the end of this section, we assume that $\theta \leqslant \frac{17}{32} - \varepsilon$. We still use two different methods to give lower bounds for $C_0^{*}(\theta_1, \theta_2)$: The first is Harman's sieve, and the second is due to Mikawa \cite{Mikawa}.

\subsubsection{First Method}
The first method is to use Harman's sieve as in \cite{676} and Subsubsection 2.5.1. Again, we can only discard positive terms that do not have asymptotic formulas of the form (54) in this case. The main steps remain the same as in Subsubsection 2.5.1, but now we can use the new Type-II information corresponding to different $(\theta_1, \theta_2)$ given in Subsection 3.3. Modifications need to do in the lower bound case are similar to those in the upper bound case; However, we do not need to consider the validity of three-dimensional sieves (Lemma~\ref{l220} and the new three-dimensional Harman's sieve given in the upper bound Subsections) since they only give upper bounds for positive terms. Lemma~\ref{l39} is still applicable for parts of $\boldsymbol{A}_{09}$ and $\boldsymbol{A}_{10}$, and we can apply Lemma~\ref{l35} in parts of $\boldsymbol{A}_{09}$--$\boldsymbol{A}_{12}$. Working on each region and subregion carefully, we can obtain the following lower bounds for $C_0^{*}(\theta_1, \theta_2)$ ($0.5 < \theta \leqslant 0.53$):
\begin{center}
\begin{tabular}{|>{\centering\arraybackslash}p{0.6cm}|>{\centering\arraybackslash}p{0.6cm}|>{\centering\arraybackslash}p{0.6cm}|>{\centering\arraybackslash}p{0.6cm}|>{\centering\arraybackslash}p{0.6cm}|>{\centering\arraybackslash}p{0.6cm}|>{\centering\arraybackslash}p{0.6cm}|>{\centering\arraybackslash}p{0.6cm}|>{\centering\arraybackslash}p{0.6cm}|>{\centering\arraybackslash}p{0.6cm}|>{\centering\arraybackslash}p{0.6cm}|>{\centering\arraybackslash}p{0.6cm}|>{\centering\arraybackslash}p{0.6cm}|>{\centering\arraybackslash}p{0.6cm}|>{\centering\arraybackslash}p{0.6cm}|}
\hline \boldmath{$0.25$} & \tiny $0.6369$ & \tiny $0.3363$ & --- & --- & --- & --- & --- & --- & --- & --- & --- & --- & --- & --- \\
\hline \boldmath{$0.24$} & $1 - \varepsilon$ & \tiny $0.6369$ & \tiny $0.3363$ & --- & --- & --- & --- & --- & --- & --- & --- & --- & --- & --- \\
\hline \boldmath{$0.23$} & $1$ & $1 - \varepsilon$ & \tiny $0.6369$ & \tiny $0.3363$ & --- & --- & --- & --- & --- & --- & --- & --- & --- & --- \\
\hline \boldmath{$0.22$} & $1$ & $1$ & $1 - \varepsilon$ & \tiny $0.6369$ & \tiny $0.3363$ & --- & --- & --- & --- & --- & --- & --- & --- & --- \\
\hline \boldmath{$0.21$} & $1$ & $1$ & $1$ & $1 - \varepsilon$ & \tiny $0.6369$ & \tiny $0.3363$ & --- & --- & --- & --- & --- & --- & --- & --- \\
\hline \boldmath{$0.20$} & $1$ & $1$ & $1$ & $1$ & $1 - \varepsilon$ & \tiny $0.6369$ & \tiny $0.3363$ & --- & --- & --- & --- & --- & --- & --- \\
\hline \boldmath{$0.19$} & $1$ & $1$ & $1$ & $1$ & $1$ & $1 - \varepsilon$ & \tiny $0.6369$ & \tiny $0.3363$ & --- & --- & --- & --- & --- & --- \\
\hline \boldmath{$0.18$} & $1$ & $1$ & $1$ & $1$ & $1$ & $1$ & $1 - \varepsilon$ & \tiny $0.6369$ & \tiny $0.3363$ & --- & --- & --- & --- & --- \\
\hline \boldmath{$0.17$} & $1$ & $1$ & $1$ & $1$ & $1$ & $1$ & $1$ & $1 - \varepsilon$ & \tiny $0.6369$ & \tiny $0.3363$ & --- & --- & --- & --- \\
\hline \boldmath{$0.16$} & $1$ & $1$ & $1$ & $1$ & $1$ & $1$ & $1$ & $1$ & $1 - \varepsilon$ & \tiny $0.6369$ & \tiny $0.3363$ & --- & --- & --- \\
\hline \boldmath{$0.15$} & $1$ & $1$ & $1$ & $1$ & $1$ & $1$ & $1$ & $1$ & $1$ & $1 - \varepsilon$ & \tiny $0.6369$ & \tiny $0.3363$ & --- & --- \\
\hline \boldmath{$0.14$} & $1$ & $1$ & $1$ & $1$ & $1$ & $1$ & $1$ & $1$ & $1$ & $1$ & $1 - \varepsilon$ & \tiny $0.6369$ & \tiny $0.3363$ & --- \\
\hline \boldmath{$0.13$} & $1$ & $1$ & $1$ & $1$ & $1$ & $1$ & $1$ & $1$ & $1$ & $1$ & $1$ & $1 - \varepsilon$ & \tiny $0.6369$ & \tiny $0.3363$ \\
\hline \boldmath{$0.12$} & $1$ & $1$ & $1$ & $1$ & $1$ & $1$ & $1$ & $1$ & $1$ & $1$ & $1$ & $1$ & $1 - \varepsilon$ & \tiny $0.6369$ \\
\hline \small \boldmath{$\theta_2 \backslash \theta_1$} & \boldmath{$0.26$} & \boldmath{$0.27$} & \boldmath{$0.28$} & \boldmath{$0.29$} & \boldmath{$0.30$} & \boldmath{$0.31$} & \boldmath{$0.32$} & \boldmath{$0.33$} & \boldmath{$0.34$} & \boldmath{$0.35$} & \boldmath{$0.36$} & \boldmath{$0.37$} & \boldmath{$0.38$} & \boldmath{$0.39$} \\
\hline
\end{tabular} \\
\textbf{Table 3.5: Lower Bounds for }\boldmath{$C_0^{*}(\theta_1, \theta_2)$} \textbf{(First Method, }\boldmath{$0.5 < \theta \leqslant 0.53$}\textbf{) 1/2}
\end{center}
\begin{center}
\begin{tabular}{|>{\centering\arraybackslash}p{0.6cm}|>{\centering\arraybackslash}p{0.6cm}|>{\centering\arraybackslash}p{0.6cm}|>{\centering\arraybackslash}p{0.6cm}|>{\centering\arraybackslash}p{0.6cm}|>{\centering\arraybackslash}p{0.6cm}|>{\centering\arraybackslash}p{0.6cm}|>{\centering\arraybackslash}p{0.6cm}|>{\centering\arraybackslash}p{0.6cm}|>{\centering\arraybackslash}p{0.6cm}|>{\centering\arraybackslash}p{0.6cm}|>{\centering\arraybackslash}p{0.6cm}|>{\centering\arraybackslash}p{0.6cm}|>{\centering\arraybackslash}p{0.6cm}|}
\hline \boldmath{$0.13$} & \tiny $-0.12$ & --- & --- & --- & --- & --- & --- & --- & --- & --- & --- & --- & --- \\
\hline \boldmath{$0.12$} & \tiny $0.3363$ & \tiny $0.0055$ & --- & --- & --- & --- & --- & --- & --- & --- & --- & --- & --- \\
\hline \boldmath{$0.11$} & \tiny $0.6369$ & \tiny $0.3506$ & \tiny $0.0055$ & --- & --- & --- & --- & --- & --- & --- & --- & --- & --- \\
\hline \boldmath{$0.10$} & $1 - \varepsilon$ & \tiny $0.6369$ & \tiny $0.3616$ & \tiny $0.0055$ & --- & --- & --- & --- & --- & --- & --- & --- & --- \\
\hline \boldmath{$0.09$} & $1$ & $1$ & \tiny $0.6369$ & \tiny $0.3786$ & \tiny $0.0055$ & --- & --- & --- & --- & --- & --- & --- & --- \\
\hline \boldmath{$0.08$} & $1$ & $1$ & $1$ & $1$ & \tiny $0.3828$ & \tiny $0.0055$ & --- & --- & --- & --- & --- & --- & --- \\
\hline \boldmath{$0.07$} & $1$ & $1$ & $1$ & $1$ & $1$ & $1$ & \tiny $0.0055$ & --- & --- & --- & --- & --- & --- \\
\hline \boldmath{$0.06$} & $1$ & $1$ & $1$ & $1$ & $1$ & $1$ & $1$ & \tiny $0.0055$ & --- & --- & --- & --- & --- \\
\hline \boldmath{$0.05$} & $1$ & $1$ & $1$ & $1$ & $1$ & $1$ & $1$ & $1$ & \tiny $0.5212$ & --- & --- & --- & --- \\
\hline \boldmath{$0.04$} & $1$ & $1$ & $1$ & $1$ & $1$ & $1$ & $1$ & $1$ & \tiny $0.7961$ & \tiny $0.2650$ & --- & --- & --- \\
\hline \boldmath{$0.03$} & $1$ & $1$ & $1$ & $1$ & $1$ & $1$ & $1$ & $1$ & $1$ & \tiny $0.5903$ & \tiny $0.0055$ & --- & --- \\
\hline \boldmath{$0.02$} & $1$ & $1$ & $1$ & $1$ & $1$ & $1$ & $1$ & $1$ & $1$ & \tiny $0.8186$ & \tiny $0.3828$ & \tiny $-0.04$ & --- \\
\hline \boldmath{$0.01$} & $1$ & $1$ & $1$ & $1$ & $1$ & $1$ & $1$ & $1$ & $1$ & $1$ & \tiny $0.6530$ & \tiny $0.3616$ & \tiny $-0.09$ \\
\hline \small \boldmath{$\theta_2 \backslash \theta_1$} & \boldmath{$0.40$} & \boldmath{$0.41$} & \boldmath{$0.42$} & \boldmath{$0.43$} & \boldmath{$0.44$} & \boldmath{$0.45$} & \boldmath{$0.46$} & \boldmath{$0.47$} & \boldmath{$0.48$} & \boldmath{$0.49$} & \boldmath{$0.50$} & \boldmath{$0.51$} & \boldmath{$0.52$} \\
\hline
\end{tabular} \\
\textbf{Table 3.6: Lower Bounds for }\boldmath{$C_0^{*}(\theta_1, \theta_2)$} \textbf{(First Method, }\boldmath{$0.5 < \theta \leqslant 0.53$}\textbf{) 2/2}
\end{center}

\subsubsection{Second Method}
The second method is to use Mikawa's modified sieve developed in \cite{Mikawa}, and the whole process is discussed in Section 2. When doing the decomposing process in \cite{Mikawa}, we need the Type-II range $[r_0,\ r_1]$ satisfies $r_1 > 2 r_0$. In the case $q \sim Q$ in Section 2, the Type-II range is fixed on $\left[2 \theta - 1 + \varepsilon,\ \frac{5 - 8 \theta}{6} - \varepsilon \right]$ (when $\theta < \frac{45}{89}$, the Type-II range $\left((\log x)^{\varepsilon - 1},\ 2 \theta - 1 + \varepsilon \right)$ would not bring useful improvements here), and that is why this method is not applicable when $\theta \geqslant \frac{17}{32}$. In the first $2$-factored moduli case, we can enlarge $r_1$ to values like $\frac{2 - 3 \theta}{3} - \varepsilon$ and $\frac{5 - 8 \theta_1 - 6 \theta_2}{6} - \varepsilon$ for some special $(\theta_1, \theta_2)$. For example, when $\kappa$ is replaced by $\frac{2 - 3 \theta}{3} - \varepsilon$, then we only need
$$
\frac{2 - 3 \theta}{3} - \varepsilon > 2 (2 \theta - 1 + \varepsilon), \quad \text{or} \quad \theta < \frac{8}{15} - 9 \varepsilon
$$
for the Type-II requirements. However, the method in \cite{Mikawa} still becomes invalid when $\theta > \frac{17}{32}$ even in this $2$-factored case. The real problem is not on the Type-II requirements but on the Type-I requirements. In the estimate of $S_{I}^{\prime \prime}$ in \cite{Mikawa} (see [\cite{Mikawa}, Page 148]), we need the following three conditions corresponding to Lemma~\ref{l29}:
$$
\frac{1}{2} - (2 \theta - 1) < 1 - \theta, \quad \frac{1}{2} + 15 (2 \theta - 1) < 2 - \theta \quad \text{and} \quad \frac{1}{2} + 7 (2 \theta - 1) < 2 - 2 \theta.
$$
The last condition above is equivalent to $\theta < \frac{17}{32}$. Replacing $2 \theta - 1 + \varepsilon$ with smaller $2 \theta_1 + \theta_2 - 1 + \varepsilon$ or $\varepsilon$ is meaningless here since these replacements also need the condition $\theta < \frac{17}{32}$ or even stronger $\theta < \frac{10}{19}$. Another idea is to use Lemmas~\ref{l25}--\ref{l28} to cover the region that Lemma~\ref{l29} cannot cover; Indeed Lemma~\ref{l26} covers the remaining region for $\theta < \frac{17}{30}$. However, the coefficients in the Type-I sum $S_{I}^{\prime \prime}$ in \cite{Mikawa} are convolutions involving the $\Psi$ function, which means that they do not satisfy \textbf{Condition B} (No small prime factors), and thus Lemmas~\ref{l25}--\ref{l28} are not applicable here. If we can prove Lemma~\ref{l26} without the \textbf{Condition B} on $a_{1, m_1}$ and $a_{3, m_3}$, then we can enlarge the applicable $\theta$ to some value larger than $\frac{17}{32}$.

In the first $2$-factored moduli case, we fail to extend the ``Mikawa applicable range'' of $\theta$ to $\geqslant \frac{17}{32}$. However, in this case we can still make some minor improvements over those results in Section 2. In many subregions we can use the Type-II range $\left(\frac{5 - 8 \theta}{6} - \varepsilon,\ \frac{2 - 3 \theta}{3} - \varepsilon \right]$ or similar ranges to give more asymptotic formulas of the form (54) for (46) and (47), and we can also use Lemma~\ref{l35} if $(\theta_1, \theta_2) \in \boldsymbol{C}$. Working on each region and subregion carefully, we can obtain the following lower bounds for $C_0^{*}(\theta_1, \theta_2)$ ($0.5 < \theta \leqslant 0.53$):
\begin{center}
\begin{tabular}{|>{\centering\arraybackslash}p{0.6cm}|>{\centering\arraybackslash}p{0.6cm}|>{\centering\arraybackslash}p{0.6cm}|>{\centering\arraybackslash}p{0.6cm}|>{\centering\arraybackslash}p{0.6cm}|>{\centering\arraybackslash}p{0.6cm}|>{\centering\arraybackslash}p{0.6cm}|>{\centering\arraybackslash}p{0.6cm}|>{\centering\arraybackslash}p{0.6cm}|>{\centering\arraybackslash}p{0.6cm}|>{\centering\arraybackslash}p{0.6cm}|>{\centering\arraybackslash}p{0.6cm}|>{\centering\arraybackslash}p{0.6cm}|>{\centering\arraybackslash}p{0.6cm}|}
\hline \boldmath{$0.13$} & \tiny $0.3477$ & --- & --- & --- & --- & --- & --- & --- & --- & --- & --- & --- & --- \\
\hline \boldmath{$0.12$} & \tiny $0.5487$ & \tiny $0.3625$ & --- & --- & --- & --- & --- & --- & --- & --- & --- & --- & --- \\
\hline \boldmath{$0.11$} & \tiny $0.7079$ & \tiny $0.5519$ & \tiny $0.3625$ & --- & --- & --- & --- & --- & --- & --- & --- & --- & --- \\
\hline \boldmath{$0.10$} & $1 - \varepsilon$ & \tiny $0.7079$ & \tiny $0.5551$ & \tiny $0.3625$ & --- & --- & --- & --- & --- & --- & --- & --- & --- \\
\hline \boldmath{$0.09$} & $1$ & $1$ & \tiny $0.7079$ & \tiny $0.5569$ & \tiny $0.3625$ & --- & --- & --- & --- & --- & --- & --- & --- \\
\hline \boldmath{$0.08$} & $1$ & $1$ & $1$ & $1$ & \tiny $0.5633$ & \tiny $0.3625$ & --- & --- & --- & --- & --- & --- & --- \\
\hline \boldmath{$0.07$} & $1$ & $1$ & $1$ & $1$ & $1$ & $1$ & \tiny $0.3625$ & --- & --- & --- & --- & --- & --- \\
\hline \boldmath{$0.06$} & $1$ & $1$ & $1$ & $1$ & $1$ & $1$ & $1$ & \tiny $0.3625$ & --- & --- & --- & --- & --- \\
\hline \boldmath{$0.05$} & $1$ & $1$ & $1$ & $1$ & $1$ & $1$ & $1$ & $1$ & \tiny $0.5237$ & --- & --- & --- & --- \\
\hline \boldmath{$0.04$} & $1$ & $1$ & $1$ & $1$ & $1$ & $1$ & $1$ & $1$ & \tiny $0.7051$ & \tiny $0.4419$ & --- & --- & --- \\
\hline \boldmath{$0.03$} & $1$ & $1$ & $1$ & $1$ & $1$ & $1$ & $1$ & $1$ & $1$ & \tiny $0.6339$ & \tiny $0.3625$ & --- & --- \\
\hline \boldmath{$0.02$} & $1$ & $1$ & $1$ & $1$ & $1$ & $1$ & $1$ & $1$ & $1$ & \tiny $0.7741$ & \tiny $0.5634$ & \tiny $0.3509$ & --- \\
\hline \boldmath{$0.01$} & $1$ & $1$ & $1$ & $1$ & $1$ & $1$ & $1$ & $1$ & $1$ & $1$ & \tiny $0.7097$ & \tiny $0.5552$ & \tiny $0.3449$ \\
\hline \small \boldmath{$\theta_2 \backslash \theta_1$} & \boldmath{$0.40$} & \boldmath{$0.41$} & \boldmath{$0.42$} & \boldmath{$0.43$} & \boldmath{$0.44$} & \boldmath{$0.45$} & \boldmath{$0.46$} & \boldmath{$0.47$} & \boldmath{$0.48$} & \boldmath{$0.49$} & \boldmath{$0.50$} & \boldmath{$0.51$} & \boldmath{$0.52$} \\
\hline
\end{tabular} \\
\textbf{Table 3.7: Lower Bounds for }\boldmath{$C_0^{*}(\theta_1, \theta_2)$} \textbf{(Second Method, }\boldmath{$0.5 < \theta \leqslant 0.53$}\textbf{)}
\end{center}
\begin{center}
\begin{tabular}{|>{\centering\arraybackslash}p{0.7cm}|>{\centering\arraybackslash}p{0.7cm}|>{\centering\arraybackslash}p{0.7cm}|>{\centering\arraybackslash}p{0.7cm}|>{\centering\arraybackslash}p{0.7cm}|>{\centering\arraybackslash}p{0.7cm}|>{\centering\arraybackslash}p{0.7cm}|>{\centering\arraybackslash}p{0.7cm}|>{\centering\arraybackslash}p{0.7cm}|>{\centering\arraybackslash}p{0.7cm}|>{\centering\arraybackslash}p{0.7cm}|>{\centering\arraybackslash}p{0.7cm}|>{\centering\arraybackslash}p{0.7cm}|>{\centering\arraybackslash}p{0.7cm}|}
\hline \makecell{ \\ \boldmath{$0.13$}} & \tiny \makecell{$-0.12$ \\ \boldmath{$0.3477$}} & \makecell{ \\ ---} & \makecell{ \\ ---} & \makecell{ \\ ---} & \makecell{ \\ ---} & \makecell{ \\ ---} & \makecell{ \\ ---} & \makecell{ \\ ---} & \makecell{ \\ ---} & \makecell{ \\ ---} & \makecell{ \\ ---} & \makecell{ \\ ---} & \makecell{ \\ ---} \\
\hline \makecell{ \\ \boldmath{$0.12$}} & \tiny \makecell{$0.3363$ \\ \boldmath{$0.5487$}} & \tiny \makecell{$0.0055$ \\ \boldmath{$0.3625$}} & \makecell{ \\ ---} & \makecell{ \\ ---} & \makecell{ \\ ---} & \makecell{ \\ ---} & \makecell{ \\ ---} & \makecell{ \\ ---} & \makecell{ \\ ---} & \makecell{ \\ ---} & \makecell{ \\ ---} & \makecell{ \\ ---} & \makecell{ \\ ---} \\
\hline \makecell{ \\ \boldmath{$0.11$}} & \tiny \makecell{$0.6369$ \\ \boldmath{$0.7079$}} & \tiny \makecell{$0.3506$ \\ \boldmath{$0.5519$}} & \tiny \makecell{$0.0055$ \\ \boldmath{$0.3625$}} & \makecell{ \\ ---} & \makecell{ \\ ---} & \makecell{ \\ ---} & \makecell{ \\ ---} & \makecell{ \\ ---} & \makecell{ \\ ---} & \makecell{ \\ ---} & \makecell{ \\ ---} & \makecell{ \\ ---} & \makecell{ \\ ---} \\
\hline \makecell{ \\ \boldmath{$0.10$}} & \makecell{ \\ $1 - \varepsilon$} & \tiny \makecell{$0.6369$ \\ \boldmath{$0.7079$}} & \tiny \makecell{$0.3616$ \\ \boldmath{$0.5551$}} & \tiny \makecell{$0.0055$ \\ \boldmath{$0.3625$}} & \makecell{ \\ ---} & \makecell{ \\ ---} & \makecell{ \\ ---} & \makecell{ \\ ---} & \makecell{ \\ ---} & \makecell{ \\ ---} & \makecell{ \\ ---} & \makecell{ \\ ---} & \makecell{ \\ ---} \\
\hline \makecell{ \\ \boldmath{$0.09$}} & \makecell{ \\ $1$} & \makecell{ \\ $1$} & \tiny \makecell{$0.6369$ \\ \boldmath{$0.7079$}} & \tiny \makecell{$0.3786$ \\ \boldmath{$0.5569$}} & \tiny \makecell{$0.0055$ \\ \boldmath{$0.3625$}} & \makecell{ \\ ---} & \makecell{ \\ ---} & \makecell{ \\ ---} & \makecell{ \\ ---} & \makecell{ \\ ---} & \makecell{ \\ ---} & \makecell{ \\ ---} & \makecell{ \\ ---} \\
\hline \makecell{ \\ \boldmath{$0.08$}} & \makecell{ \\ $1$} & \makecell{ \\ $1$} & \makecell{ \\ $1$} & \makecell{ \\ $1$} & \tiny \makecell{$0.3828$ \\ \boldmath{$0.5633$}} & \tiny \makecell{$0.0055$ \\ \boldmath{$0.3625$}} & \makecell{ \\ ---} & \makecell{ \\ ---} & \makecell{ \\ ---} & \makecell{ \\ ---} & \makecell{ \\ ---} & \makecell{ \\ ---} & \makecell{ \\ ---} \\
\hline \makecell{ \\ \boldmath{$0.07$}} & \makecell{ \\ $1$} & \makecell{ \\ $1$} & \makecell{ \\ $1$} & \makecell{ \\ $1$} & \makecell{ \\ $1$} & \makecell{ \\ $1$} & \tiny \makecell{$0.0055$ \\ \boldmath{$0.3625$}} & \makecell{ \\ ---} & \makecell{ \\ ---} & \makecell{ \\ ---} & \makecell{ \\ ---} & \makecell{ \\ ---} & \makecell{ \\ ---} \\
\hline \makecell{ \\ \boldmath{$0.06$}} & \makecell{ \\ $1$} & \makecell{ \\ $1$} & \makecell{ \\ $1$} & \makecell{ \\ $1$} & \makecell{ \\ $1$} & \makecell{ \\ $1$} & \makecell{ \\ $1$} & \tiny \makecell{$0.0055$ \\ \boldmath{$0.3625$}} & \makecell{ \\ ---} & \makecell{ \\ ---} & \makecell{ \\ ---} & \makecell{ \\ ---} & \makecell{ \\ ---} \\
\hline \makecell{ \\ \boldmath{$0.05$}} & \makecell{ \\ $1$} & \makecell{ \\ $1$} & \makecell{ \\ $1$} & \makecell{ \\ $1$} & \makecell{ \\ $1$} & \makecell{ \\ $1$} & \makecell{ \\ $1$} & \makecell{ \\ $1$} & \tiny \makecell{$0.5212$ \\ \boldmath{$0.5237$}} & \makecell{ \\ ---} & \makecell{ \\ ---} & \makecell{ \\ ---} & \makecell{ \\ ---} \\
\hline \makecell{ \\ \boldmath{$0.04$}} & \makecell{ \\ $1$} & \makecell{ \\ $1$} & \makecell{ \\ $1$} & \makecell{ \\ $1$} & \makecell{ \\ $1$} & \makecell{ \\ $1$} & \makecell{ \\ $1$} & \makecell{ \\ $1$} & \tiny \makecell{\boldmath{$0.7961$} \\ $0.7051$} & \tiny \makecell{$0.2650$ \\ \boldmath{$0.4419$}} & \makecell{ \\ ---} & \makecell{ \\ ---} & \makecell{ \\ ---} \\
\hline \makecell{ \\ \boldmath{$0.03$}} & \makecell{ \\ $1$} & \makecell{ \\ $1$} & \makecell{ \\ $1$} & \makecell{ \\ $1$} & \makecell{ \\ $1$} & \makecell{ \\ $1$} & \makecell{ \\ $1$} & \makecell{ \\ $1$} & \makecell{ \\ $1$} & \tiny \makecell{$0.5903$ \\ \boldmath{$0.6339$}} & \tiny \makecell{$0.0055$ \\ \boldmath{$0.3625$}} & \makecell{ \\ ---} & \makecell{ \\ ---} \\
\hline \makecell{ \\ \boldmath{$0.02$}} & \makecell{ \\ $1$} & \makecell{ \\ $1$} & \makecell{ \\ $1$} & \makecell{ \\ $1$} & \makecell{ \\ $1$} & \makecell{ \\ $1$} & \makecell{ \\ $1$} & \makecell{ \\ $1$} & \makecell{ \\ $1$} & \tiny \makecell{\boldmath{$0.8186$} \\ $0.7741$} & \tiny \makecell{$0.3828$ \\ \boldmath{$0.5634$}} & \tiny \makecell{$-0.04$ \\ \boldmath{$0.3509$}} & \makecell{ \\ ---} \\
\hline \makecell{ \\ \boldmath{$0.01$}} & \makecell{ \\ $1$} & \makecell{ \\ $1$} & \makecell{ \\ $1$} & \makecell{ \\ $1$} & \makecell{ \\ $1$} & \makecell{ \\ $1$} & \makecell{ \\ $1$} & \makecell{ \\ $1$} & \makecell{ \\ $1$} & \makecell{ \\ $1$} & \tiny \makecell{$0.6530$ \\ \boldmath{$0.7097$}} & \tiny \makecell{$0.3616$ \\ \boldmath{$0.5552$}} & \tiny \makecell{$-0.09$ \\ \boldmath{$0.3449$}} \\
\hline \small \boldmath{$\theta_2 \backslash \theta_1$} & \boldmath{$0.40$} & \boldmath{$0.41$} & \boldmath{$0.42$} & \boldmath{$0.43$} & \boldmath{$0.44$} & \boldmath{$0.45$} & \boldmath{$0.46$} & \boldmath{$0.47$} & \boldmath{$0.48$} & \boldmath{$0.49$} & \boldmath{$0.50$} & \boldmath{$0.51$} & \boldmath{$0.52$} \\
\hline
\end{tabular} \\
\textbf{Table 3.8: A Comparison of Two Methods on the Lower Bounds for} \boldmath{$C_0^{*}(\theta_1, \theta_2)$} \textbf{(}\boldmath{$0.5 < \theta \leqslant 0.53$}\textbf{)}
\end{center}

\section{$2$-factored Moduli, 2}
In this section we focus on the second $2$-factored case with bilinear weights, where the absolute values in the first case are replaced by divisor-bounded coefficients. We also call this case ``the bilinear case''. The initial setups on the sieves are similar to the first case. We want to get the following result with some $0 < C_0^{\prime}(\theta_1, \theta_2) \leqslant 1$ and $C_1^{\prime}(\theta_1, \theta_2) \geqslant 1$:
\begin{theorem}\label{t41}
There exist functions $\rho_0$ and $\rho_1$ which satisfies the following properties:

(Majorant / Minorant). $\rho_0(n)$ is a minorant for the prime indicator function $\mathbbm{1}_{p}(n)$, and $\rho_1(n)$ is a majorant for the prime indicator function $\mathbbm{1}_{p}(n)$. That is, we have
$$
\rho_0(n) \leqslant \mathbbm{1}_{p}(n) \leqslant \rho_1(n).
$$

(Upper and Lower bounds). We have
$$
\sum_{n \leqslant x} \rho_0(n) \geqslant (1+o(1))\frac{C_0^{\prime}(\theta_1, \theta_2) x}{\log x} \quad \text{and} \quad \sum_{n \leqslant x} \rho_1(n) \leqslant (1+o(1))\frac{C_1^{\prime}(\theta_1, \theta_2) x}{\log x}
$$
for two functions $C_0^{\prime}(\theta_1, \theta_2)$ and $C_1^{\prime}(\theta_1, \theta_2)$ satisfy $0 < C_0^{\prime}(\theta_1, \theta_2) \leqslant 1$ and $C_1^{\prime}(\theta_1, \theta_2) \geqslant 1$.

(Distributions in Arithmetic Progressions). Let $\lambda_{1, q_1}$ and $\lambda_{2, q_2}$ be divisor-bounded complex sequences. For any $a \in \mathbb{Z} \backslash \{0\}$ and any $A>0$, we have
$$
\sum_{\substack{q_1 \sim Q_1 \\ q_2 \sim Q_2 \\ (q_1 q_2, a) = 1}} \lambda_{1, q_1} \lambda_{2, q_2} \left( \sum_{\substack{n \leqslant x \\ n \equiv a (\bmod q_1 q_2)}} \rho_j(n) - \frac{1}{\varphi(q_1 q_2)} \sum_{\substack{n \leqslant x \\ (n, q_1 q_2) = 1}} \rho_j(n) \right) \ll \frac{x}{(\log x)^A}
$$
for $j = 0, 1$.
\end{theorem}

In order to prove Theorem~\ref{t41} with suitable $C_0^{\prime}(\theta_1, \theta_2)$ and $C_1^{\prime}(\theta_1, \theta_2)$, we need results of the form
\begin{equation}
\sum_{\substack{q_1 \sim Q_1 \\ q_2 \sim Q_2 \\ (q_1 q_2, a) = 1}} \lambda_{1, q_1} \lambda_{2, q_2} \left( \sum_{\substack{n \sim x \\ n \equiv a (\bmod q_1 q_2)}} f(n) - \frac{1}{\varphi(q_1 q_2)} \sum_{\substack{n \sim x \\ (n, q_1 q_2) = 1}} f(n) \right) \ll \frac{x}{(\log x)^A}.
\end{equation}
Again, we may want the coefficients to satisfy \textbf{Conditions A and B} mentioned in Section 2.

\subsection{Preliminary Lemmas}
Before constructing the majorant and minorant, we need estimate results of the form (129). The results from Section 2 and Section 3 are still applicable in the final decomposition (except for the results in Subsection 2.3), and the results here are still useful in the later Section 5. 

\subsubsection{Type-II estimate}
The first lemma comes from \cite{FouvryA2}. Note that Case (1) of this lemma can be deduced easily by Lemma~\ref{l34}.
\begin{lemma}\label{l42} ([\cite{FouvryA2}, Théorème]).
Let $M_1 M_2 \asymp x$ and $M_2 \geqslant x^{\varepsilon}$. Let $a_{1, m_1}$, $a_{2, m_2}$, $\lambda_{1, q_1}$ and $\lambda_{2, q_2}$ be divisor-bounded complex sequences. Suppose that $a_{2, m_2}$ satisfies \textbf{Conditions A and B}. If any of the following conditions
\begin{align}\label{l42}
\nonumber (1). \qquad&\ Q_1 Q_2^2 \leqslant M_2 x^{1 - \varepsilon},\ Q_1^8 Q_2^7 M_2^6 \leqslant x^{4 - \varepsilon}; \\
\nonumber (2). \qquad&\ Q_1 Q_2^2 \leqslant x^{1 - \varepsilon},\ Q_1^8 Q_2^7 M_2^5 \leqslant x^{4 - \varepsilon}; \\
\nonumber (3). \qquad&\ Q_1 Q_2^2 \leqslant M_2 x^{1 - \varepsilon},\ Q_1^5 Q_2^2 M_2^6 \leqslant x^{2 - \varepsilon},\ Q_1^4 Q_2^3 M_2^3 \leqslant x^{2 - \varepsilon},\ Q_1^9 Q_2^8 M_2^6 \leqslant x^{5 - \varepsilon}; \\
\nonumber (4). \qquad&\ Q_1 Q_2^2 \leqslant x^{1 - \varepsilon},\ Q_1^5 Q_2^2 M_2^5 \leqslant x^{2 - \varepsilon},\ Q_1^8 Q_2^6 M_2^5 \leqslant x^{4 - \varepsilon}; \\
\nonumber (5). \qquad&\ Q_1 Q_2^2 M_2 \leqslant x^{1 - \varepsilon},\ Q_1^5 Q_2^2 M_2^4 \leqslant x^{2 - \varepsilon},\ Q_1^8 Q_2^6 M_2^5 \leqslant x^{4 - \varepsilon}
\end{align}
holds, then
$$
\sum_{\substack{q_1 \sim Q_1 \\ q_2 \sim Q_2 \\ (q_1 q_2, a) = 1}} \lambda_{1, q_1} \lambda_{2, q_2} \left( \sum_{\substack{m_1 \sim M_1 \\ m_2 \sim M_2 \\ m_1 m_2 \equiv a (\bmod q_1 q_2)}} a_{1, m_1} a_{2, m_2} - \frac{1}{\varphi(q_1 q_2)} \sum_{\substack{m_1 \sim M_1 \\ m_2 \sim M_2 \\ (m_1 m_2, q_1 q_2) = 1}} a_{1, m_1} a_{2, m_2} \right) \ll \frac{x}{(\log x)^A}.
$$
\end{lemma}
The next lemma comes from \cite{IwaniecPomykala}, and it was used in \cite{Lichtman2}.
\begin{lemma}\label{l43} ([\cite{IwaniecPomykala}, Proposition]).
Let $M_1 M_2 \asymp x$. Let $a_{1, m_1}$, $a_{2, m_2}$, $\lambda_{1, q_1}$ and $\lambda_{2, q_2}$ be divisor-bounded complex sequences. Suppose that $a_{2, m_2}$ satisfies \textbf{Conditions A and B}. If we have
$$
Q_1 x^{\varepsilon} < M_2,\ Q_1 Q_2^2 < M_2 x^{1 - \varepsilon},\ Q_1^4 Q_2^8 M_2^6 < x^{5 - \varepsilon},
$$
then
$$
\sum_{\substack{q_1 \sim Q_1 \\ q_2 \sim Q_2 \\ (q_1 q_2, a) = 1}} \lambda_{1, q_1} \lambda_{2, q_2} \left( \sum_{\substack{m_1 \sim M_1 \\ m_2 \sim M_2 \\ m_1 m_2 \equiv a (\bmod q_1 q_2)}} a_{1, m_1} a_{2, m_2} - \frac{1}{\varphi(q_1 q_2)} \sum_{\substack{m_1 \sim M_1 \\ m_2 \sim M_2 \\ (m_1 m_2, q_1 q_2) = 1}} a_{1, m_1} a_{2, m_2} \right) \ll \frac{x}{(\log x)^A}.
$$
\end{lemma}

The next two lemmas come from \cite{BFI} and were used in the proof of [\cite{BFI}, Theorem 8], [\cite{FouvryA2}, Corollaire 5] and [\cite{677}, Theorem 3].
\begin{lemma}\label{l44} ([\cite{BFI}, Theorem 1]).
Let $M_1 M_2 \asymp x$ and $\min(M_1, M_2) > x^{\varepsilon}$. Let $a_{1, m_1}$, $a_{2, m_2}$, $\lambda_{1, q_1}$ and $\lambda_{2, q_2}$ be divisor-bounded complex sequences. Suppose that $a_{2, m_2}$ satisfies \textbf{Conditions A and B}. If we have
$$
Q_2 x^{\varepsilon} \ll M_2 \ll x^{-\varepsilon} \min\left(x^{\frac{1}{2}} Q_1^{-\frac{1}{2}}, x^2 Q_1^{-5} Q_2^{-1}, x Q_1^{-2} Q_2^{-\frac{1}{2}} \right),
$$
then
$$
\sum_{\substack{q_2 \sim Q_2 \\ m_1 \sim M_1 \\ (q_2, a m_1) = 1}} \left( \sum_{\substack{q_1 \sim Q_1 \\ (q_1, a m_1) = 1}} \lambda_{1, q_1} \left( \sum_{\substack{m_2 \sim M_2 \\ m_1 m_2 \equiv a (\bmod q_1 q_2)}} a_{2, m_2} - \frac{1}{\varphi(q_1 q_2)} \sum_{\substack{m_2 \sim M_2 \\ (m_2, q_1 q_2) = 1}} a_{2, m_2} \right) \right)^2 \ll \|a_2\|^2 x Q_2^{-1} (\log x)^{-A},
$$
where $\|a_2\| = \left(\sum_{m_2} |a_{2, m_2}|^2 \right)^{\frac{1}{2}}$ is the $l^2$ norm.
\end{lemma}
\begin{lemma}\label{l45} ([\cite{BFI}, Theorem 2]).
Let $M_1 M_2 \asymp x$ and $\min(M_1, M_2) > x^{\varepsilon}$. Let $a_{1, m_1}$, $a_{2, m_2}$, $\lambda_{1, q_1}$ and $\lambda_{2, q_2}$ be divisor-bounded complex sequences. Suppose that $a_{2, m_2}$ satisfies \textbf{Conditions A and B}. Suppose also that
$$
M_2^{1 - \varepsilon} \sum_{m_2} |a_{2, m_2}|^4 \ll \left(\sum_{m_2} |a_{2, m_2}|^2 \right)^2.
$$
If we have
$$
Q_2 x^{\varepsilon} \ll M_2 \ll x^{-\varepsilon} \min\left(x^{\frac{1}{2}} Q_1^{-1} Q_2^{\frac{1}{2}}, x^{\frac{2}{5}} Q_1^{-\frac{2}{5}}, x^{\frac{1}{2}} Q_1^{-\frac{3}{4}} \right),
$$
then
$$
\sum_{\substack{q_2 \sim Q_2 \\ m_1 \sim M_1 \\ (q_2, a m_1) = 1}} \left( \sum_{\substack{q_1 \sim Q_1 \\ (q_1, a m_1) = 1}} \lambda_{1, q_1} \left( \sum_{\substack{m_2 \sim M_2 \\ m_1 m_2 \equiv a (\bmod q_1 q_2)}} a_{2, m_2} - \frac{1}{\varphi(q_1 q_2)} \sum_{\substack{m_2 \sim M_2 \\ (m_2, q_1 q_2) = 1}} a_{2, m_2} \right) \right)^2 \ll \|a_2\|^2 x Q_2^{-1} (\log x)^{-A}.
$$
\end{lemma}
\begin{remark*}
Before applying Lemma~\ref{l44} and Lemma~\ref{l45}, one need to use Cauchy's inequality as the following:
\begin{align}
\nonumber &\ \left( \sum_{\substack{q_1 \sim Q_1 \\ q_2 \sim Q_2 \\ (q_1 q_2, a) = 1}} \lambda_{1, q_1} \lambda_{2, q_2} \left( \sum_{\substack{m_1 \sim M_1 \\ m_2 \sim M_2 \\ m_1 m_2 \equiv a (\bmod q_1 q_2)}} a_{1, m_1} a_{2, m_2} - \frac{1}{\varphi(q_1 q_2)} \sum_{\substack{m_1 \sim M_1 \\ m_2 \sim M_2 \\ (m_1 m_2, q_1 q_2) = 1}} a_{1, m_1} a_{2, m_2} \right) \right)^2 \\
\nonumber \leqslant&\ \|\lambda_2\|^2 \|a_1\|^2 \sum_{\substack{q_2 \sim Q_2 \\ m_1 \sim M_1 \\ (q_2, a m_1) = 1}} \left( \sum_{\substack{q_1 \sim Q_1 \\ (q_1, a m_1) = 1}} \lambda_{1, q_1} \left( \sum_{\substack{m_2 \sim M_2 \\ m_1 m_2 \equiv a (\bmod q_1 q_2)}} a_{2, m_2} - \frac{1}{\varphi(q_1 q_2)} \sum_{\substack{m_2 \sim M_2 \\ (m_2, q_1 q_2) = 1}} a_{2, m_2} \right) \right)^2.
\end{align}
\end{remark*}

\subsubsection{Type-I estimate}
The next lemma comes from \cite{BFI}, and it was used in \cite{677}.
\begin{lemma}\label{l46} ([\cite{BFI}, Theorems 5 and 5$^{*}$]).
Let $M_1 M_2 \asymp x$ and $z \ll \exp\left(\log x (\log \log x)^{-1} \right)$. Let $a_{1, m_1}$, $a_{2, m_2}$, $\lambda_{1, q_1}$ and $\lambda_{2, q_2}$ be divisor-bounded complex sequences. Suppose that
$$
a_{1, m_1} = \mathbbm{1}_{m_1 \in \mathbf{M}} \qquad \text{or} \qquad a_{1, m_1} = \mathbbm{1}_{\substack{m_1 \in \mathbf{M} \\ \left(m_1, P(z)\right) = 1}}
$$
for some interval $\mathbf{M} \subseteq [M_1, 2 M_1]$. If we have
$$
Q_1 M_2 < x^{1 - \varepsilon}, \quad Q_1 Q_2^4 M_2 < x^{2 - \varepsilon}, \quad Q_1 Q_2^2 M_2^2 < x^{2 - 2 \varepsilon}, \quad Q_1^3 Q_2^4 M_2 < x^{3 - \varepsilon},
$$
then
$$
\sum_{\substack{q_1 \sim Q_1 \\ q_2 \sim Q_2 \\ (q_1 q_2, a) = 1}} \lambda_{1, q_1} \lambda_{2, q_2} \left( \sum_{\substack{m_1 \in \mathbf{M} \\ m_2 \sim M_2 \\ m_1 m_2 \equiv a (\bmod q_1 q_2)}} a_{1, m_1} a_{2, m_2} - \frac{1}{\varphi(q_1 q_2)} \sum_{\substack{m_1 \in \mathbf{M} \\ m_2 \sim M_2 \\ (m_1 m_2, q_1 q_2) = 1}} a_{1, m_1} a_{2, m_2} \right) \ll \|a_2\| M_1^{\frac{1}{2}} x^{\frac{1}{2} - \varepsilon}.
$$
\end{lemma}

\subsection{Upper Bounds}
We shall construct the majorant $\rho_1(n)$ and prove upper bounds for $C_1^{\prime}(\theta_1, \theta_2)$ in this subsection. Before constructing, we first mention some existing results of $C_1^{\prime}(\theta_1, \theta_2)$.
\begin{theorem}\label{t47}
The function $C_1^{\prime}(\theta_1, \theta_2)$ satisfies the following conditions:

(1). $C_1^{\prime}(\theta_1, \theta_2) = C_1^{\prime}(\theta_2, \theta_1)$;

(2). $C_1^{\prime}(\theta_1, \theta_2) = 1$ for all $\theta_1, \theta_2$ satisfy $\theta_1 + \theta_2 \leqslant 0.5 - \varepsilon$; 

(3). $C_1^{\prime}(\theta_1, \theta_2) = 1$ for all $\theta_1, \theta_2$ satisfy $2 \theta_1 + \theta_2 \leqslant 1 - \varepsilon$, $7 \theta_1 + 12 \theta_2 \leqslant 4 - \varepsilon$ and $19 \theta_1 + 20 \theta_2 \leqslant 10 - \varepsilon$;

(4). $C_1^{\prime}(\theta_1, \theta_2) = 1$ for all $\theta_1, \theta_2$ satisfy $\theta_1 < \frac{1}{3}$, $\theta_2 < \frac{1}{5}$ and $\theta_1 + \theta_2 < \frac{29}{56}$;

(5). $C_1^{\prime}(\theta_1, \theta_2) = 1$ for all $\theta_1, \theta_2$ satisfy $\theta_1 + 3 \theta_2 < 1$, $\theta_1 + \theta_2 < \frac{29}{56}$ and $\theta_2 < \max\left(\frac{1 - 2 \theta_1}{2}, \frac{2 - 2 \theta_1}{5}\right)$;

(6). $C_1^{\prime}(\theta_1, \theta_2) = 1$ for all $\theta_1, \theta_2$ satisfy $\theta_1 + 3 \theta_2 < 1$, $\theta_1 + \theta_2 < \frac{29}{56}$, $4 \theta_1 + \theta_2 < \frac{403}{266}$ and $\frac{7}{4} \theta_1 + \theta_2 < \frac{403}{532}$;

(7). $C_1^{\prime}(\theta_1, \theta_2) \leqslant C_1(\theta_1, \theta_2) \leqslant C_1(\theta_1 + \theta_2)$ for $0.5 \leqslant \theta_1 + \theta_2 \leqslant 1$;

(8). $C_1^{\prime}(\theta_1, \theta_2) \leqslant 1 + \varepsilon$ for all $\theta_1, \theta_2$ satisfy $\theta_1 + \theta_2 = 0.5$;

(9). $C_1^{\prime}(\theta_1, \theta_2) \leqslant 1 + \varepsilon$ for all $\theta_1, \theta_2$ satisfy $\theta_1 \leqslant 0.5$ and $\theta_2 = \min\left(1 - 2 \theta_1, \frac{4 - 7 \theta_1}{12}, \frac{10 - 19 \theta_1}{20}\right)$.
\end{theorem}
\begin{proof}
Statement (1) is obvious. Statements (2)--(6) follow easily from the Bombieri--Vinogradov Theorem, [\cite{MaynardLargeModuliI}, Theorem 1.1], [\cite{677}, Theorem 3] and [\cite{FouvryA2}, Page 621 and Corollaire 5]. Statement (7) holds trivially by the work done in Section 2 and Section 3. When there are no new arithmetic information inputs outside of those in previous sections, we use $C_1(\theta_1, \theta_2)$ as an upper bound for $C_1^{\prime}(\theta_1, \theta_2)$. Statements (8)--(9) follows from Statement (7) together with Statements (5)--(6) of Theorem~\ref{t310}.
\end{proof}

From here to the end of this section, we assume that $\theta_1 \geqslant \theta_2$ to simplify the conditions. We also write $\theta = \theta_1 + \theta_2$. Before performing our final decompositions, we define several regions of the pair $(\theta_1, \theta_2)$ based on various arithmetic information inputs. We also use many other regions defined before, and the readers can find the definitions of them in previous sections. Moreover, the region $\boldsymbol{S}$ defined in Section 2 can be enlarged in this section (and also in Section 5) using Lemma~\ref{l46} above.
\begin{align}
\nonumber \boldsymbol{S} =&\ \left\{ (s, t): s \leqslant 1 - \theta - \varepsilon,\ s + 2t \leqslant 2 - 2 \theta - \varepsilon,\ s + 4t \leqslant 2 - \theta - \varepsilon \right. \\
\nonumber & \left. \qquad \qquad \text{ or } s + t < \min\left(1 - \theta_1, 2 - \theta_1 - 4 \theta_2, 1 - \frac{1}{2} \theta_1 - \theta_2, 3 - 3 \theta_1 - 4 \theta_2 \right) - \varepsilon \right\}, \\
\nonumber \boldsymbol{J} =&\ \left\{ (\theta_1, \theta_2) : (\theta_1, \theta_2) \in \boldsymbol{U};\ \theta_1 + \theta_2 \leqslant \frac{1}{2} - \varepsilon \right. \\
\nonumber & \qquad \qquad \quad \text{ or } 2 \theta_1 + \theta_2 \leqslant 1 - \varepsilon,\ 7 \theta_1 + 12 \theta_2 \leqslant 4 - \varepsilon,\ 19 \theta_1 + 20 \theta_2 \leqslant 10 - \varepsilon \\
\nonumber & \qquad \qquad \quad \text{ or } \theta_1 < \frac{1}{3},\ \theta_2 < \frac{1}{5},\ \theta_1 + \theta_2 < \frac{29}{56} \\
\nonumber & \qquad \qquad \quad \text{ or } \theta_1 + 3 \theta_2 < 1,\ \theta_1 + \theta_2 < \frac{29}{56},\ \theta_2 < \max\left(\frac{1 - 2 \theta_1}{2}, \frac{2 - 2 \theta_1}{5}\right) \\
\nonumber & \left. \qquad \qquad \quad \text{ or } \theta_1 + 3 \theta_2 < 1,\ \theta_1 + \theta_2 < \frac{29}{56},\ 4 \theta_1 + \theta_2 < \frac{403}{266},\ \frac{7}{4} \theta_1 + \theta_2 < \frac{403}{532} \right\}, \\
\nonumber \boldsymbol{E} =&\ \left\{ (\theta_1, \theta_2) : (\theta_1, \theta_2) \in \boldsymbol{U} \backslash \boldsymbol{J};\ \frac{1}{4} < \theta_1 \leqslant \frac{2}{7},\ \frac{1}{2}(1 - 2 \theta_1) < \theta_2 \leqslant \theta_1 \right. \\
\nonumber & \qquad \qquad \quad \text{ or } \frac{2}{7} < \theta_1 \leqslant \frac{2}{5},\ \frac{1}{2}(1 - 2 \theta_1) < \theta_2 < \frac{1}{4}(2 - 3 \theta_1) \\
\nonumber & \qquad \qquad \quad \text{ or } \frac{2}{5} < \theta_1 < \frac{1}{2},\ \frac{1}{2}(1 - 2 \theta_1) < \theta_2 < \frac{1}{15}(11 - 20 \theta_1) \\
\nonumber & \left. \qquad \qquad \quad \text{ or } \frac{1}{2} \leqslant \theta_1 < \frac{11}{20},\ 0 < \theta_2 < \frac{1}{15}(11 - 20 \theta_1) \right\}, \\
\nonumber \boldsymbol{F}_1 =&\ \left\{ (\theta_1, \theta_2) : (\theta_1, \theta_2) \in \boldsymbol{U} \backslash \boldsymbol{J};\ \frac{1}{4} < \theta_1 \leqslant \frac{5}{18},\ \frac{1}{2}(1 - 2 \theta_1) < \theta_2 \leqslant \theta_1 \right. \\
\nonumber & \qquad \qquad \quad \text{ or } \frac{5}{18} < \theta_1 < \frac{1}{2},\ \frac{1}{2}(1 - 2 \theta_1) < \theta_2 < \frac{1}{10}(5 - 8 \theta_1) \\
\nonumber & \left. \qquad \qquad \quad \text{ or } \frac{1}{2} \leqslant \theta_1 < \frac{11}{20},\ 0 < \theta_2 < \frac{1}{10}(11 - 20 \theta_1) \right\}, \\
\nonumber \boldsymbol{F}_2 =&\ \left\{ (\theta_1, \theta_2) : (\theta_1, \theta_2) \in \boldsymbol{U} \backslash \boldsymbol{J};\ \frac{1}{4} < \theta_1 < \frac{2}{7},\ \frac{1}{2}(1 - 2 \theta_1) < \theta_2 \leqslant \theta_1 \right. \\
\nonumber & \qquad \qquad \quad \text{ or } \frac{2}{7} \leqslant \theta_1 \leqslant \frac{2}{5},\ \frac{1}{2}(1 - 2 \theta_1) < \theta_2 < \frac{1}{4}(2 - 3 \theta_1) \\
\nonumber & \left. \qquad \qquad \quad \text{ or } \frac{2}{5} < \theta_1 < \frac{1}{2},\ \frac{1}{2}(1 - 2 \theta_1) < \theta_2 < 1 - 2 \theta_1 \right\}.
\end{align}
Here, $\boldsymbol{J}$ denote the region that $C_1^{\prime}(\theta_1, \theta_2) = 1$ follows by the Bombieri--Vinogradov Theorem or theorems mentioned in the proof of Statements (3)--(6) of Theorem~\ref{t47}. Region $\boldsymbol{E}$ corresponds to Lemma~\ref{l42}, region $\boldsymbol{F}_1$ corresponds to Lemma~\ref{l43}, and region $\boldsymbol{F}_2$ corresponds to Lemma~\ref{l44} and Lemma~\ref{l45}. The results in Section 3 will also be applied when $(\theta_1, \theta_2)$ is in $\boldsymbol{A}$, $\boldsymbol{B}$, $\boldsymbol{C}$ or $\boldsymbol{T}$. We note that the three-dimensional Harman's sieve corresponds to region $\boldsymbol{T}_{2}$ (see Lemma~\ref{l38} and the discussions in Subsection 3.3) will be used a lot since we have many Type-II ranges start from $\varepsilon$, compare to the cases in Section 3. Again, we shall implicitly use the new  three-dimensional Harman's sieve in many decompositions below. Sometimes we use it even when $\kappa < \frac{1}{7} + \varepsilon$; in this case, we only need to make some modifications to the ``loss integrals'': Suppose that $\frac{1}{8} < \kappa < \frac{1}{7} + \varepsilon$, and we apply a three-dimensional Harman's sieve on a region $R$. Since $\alpha_1 < \frac{1}{2} = \frac{4}{8}$ and $\alpha_2 < \frac{1}{3} < \frac{3}{8}$, we have $\Omega(m_1) \leqslant 3$ and $\Omega(m_2) \leqslant 2$, and $\Omega(m_1 m_2)$ can be $3$, $4$ or $5$. Now we have 3 cases based on different values of $\Omega(m_1 m_2)$:

(1). $\Omega(m_1) = 1,\ \Omega(m_2) = 2$ or $\Omega(m_1) = 2,\ \Omega(m_2) = 1$;

(2). $\Omega(m_1) = 2,\ \Omega(m_2) = 2$ or $\Omega(m_1) = 3,\ \Omega(m_2) = 1$;

(3). $\Omega(m_1) = 3,\ \Omega(m_2) = 2$.

\noindent These 3 cases correspond to 3 ``loss integrals'':
\begin{equation}
\frac{1}{\kappa} \int_{\substack{t_1, t_2, t_3 \geqslant \kappa \\ (t_1, t_2 + t_3) \in R,\ t_2 \geqslant t_3 \\ \text{or } (t_1 + t_2, t_3) \in R,\ t_1 \geqslant t_2 \\ (t_1, t_2, t_3) \notin \boldsymbol{G}_{3} }} \frac{\omega\left(\frac{1 - t_1 - t_2 - t_3}{\kappa}\right)}{t_1 t_2 t_3} d t_3 d t_2 d t_1,
\end{equation}
\begin{equation}
\frac{1}{\kappa} \int_{\substack{t_1, t_2, t_3, t_4 \geqslant \kappa \\ (t_1 + t_2, t_3 + t_4) \in R,\ t_1 \geqslant t_2,\ t_3 \geqslant t_4 \\ \text{or } (t_1 + t_2 + t_3, t_4) \in R,\ t_1 \geqslant t_2 \geqslant t_3 \\ (t_1, t_2, t_3, t_4) \notin \boldsymbol{G}_{4} }} \frac{\omega\left(\frac{1 - t_1 - t_2 - t_3 - t_4}{\kappa}\right)}{t_1 t_2 t_3 t_4} d t_4 d t_3 d t_2 d t_1,
\end{equation}
and
\begin{equation}
\frac{1}{\kappa} \int_{\substack{t_1, t_2, t_3, t_4, t_5 \geqslant \kappa \\ (t_1 + t_2 + t_3, t_4 + t_5) \in R \\ t_1 \geqslant t_2 \geqslant t_3,\ t_4 \geqslant t_5 \\ (t_1, t_2, t_3, t_4, t_5) \notin \boldsymbol{G}_{5} }} \frac{\omega\left(\frac{1 - t_1 - t_2 - t_3 - t_4 - t_5}{\kappa}\right)}{t_1 t_2 t_3 t_4 t_5} d t_5 d t_4 d t_3 d t_2 d t_1.
\end{equation}
One can compare them with the ``loss integrals'' in Lemma~\ref{l220} and Lemma~\ref{l221}. When $\frac{11}{21} - \varepsilon < \theta \leqslant \frac{17}{32} - \varepsilon$, the new three-dimensional Harman's sieve should be applied on both $B$ and $C$ (see Lemma~\ref{l217}) if we use it, since we can only perform straightforward decompositions on $A$ where $\alpha_1 + \alpha_2 < \theta$. Since we have $\alpha_2 < \frac{1}{3}$, $\alpha_1 < \frac{3}{7} + \varepsilon$ and $\alpha_1 + \alpha_2 > \frac{4}{7} - \varepsilon$ for $\boldsymbol{\alpha}_{2} \in B \cup C$, the new three-dimensional Harman's sieve is available here. Now we want to prove that $\boldsymbol{\alpha}_{2} \in \boldsymbol{U}_{2}$ for $\boldsymbol{\alpha}_{2} \in A$. The proof is trivial when $\alpha_1 < \tau = \frac{2}{7}$. When $\frac{2}{7} \leqslant \alpha_1 \leqslant \frac{3}{7} + \varepsilon$, we need to show that $\boldsymbol{\alpha}_{2} \in \boldsymbol{U}_{2}$ or $(\alpha_1, \alpha_2, 2 \theta - 1 + \varepsilon) \in \boldsymbol{S}_{3}$. When $\theta < \frac{48}{91} - \varepsilon$, we have
\begin{equation}
\alpha_1 \leqslant \frac{3}{7} + \varepsilon < 1 - \theta - \varepsilon,
\end{equation}
\begin{equation}
\alpha_1 + 2 (\alpha_2 + 2 \theta - 1 + \varepsilon) < 6 \theta - \alpha_1 - 2 + 2 \varepsilon \leqslant 6 \theta - \frac{16}{7} + 2 \varepsilon < 2 - 2 \theta - \varepsilon,
\end{equation}
and
\begin{equation}
\alpha_1 + 4 (\alpha_2 + 2 \theta - 1 + \varepsilon) < 12 \theta - 3 \alpha_1 - 4 + 4 \varepsilon \leqslant 12 \theta - \frac{34}{7} + 4 \varepsilon < 2 - \theta - \varepsilon.
\end{equation}
When $\theta \geqslant \frac{48}{91} - \varepsilon$, however, we need an extra condition $\alpha_1 < 2 - 3 \theta - 3 \varepsilon$ to get
\begin{equation}
\alpha_1 + 2 \theta - 1 + \varepsilon < 1 - \theta - \varepsilon,
\end{equation}
\begin{equation}
\alpha_1 + 2 \theta - 1 + \varepsilon + 2 \alpha_2 < 4 \theta - \alpha_1 - 1 + \varepsilon 
\leqslant 4 \theta - \frac{9}{7} + \varepsilon < 2 - 2 \theta - \varepsilon,
\end{equation}
and
\begin{equation}
\alpha_1 + 2 \theta - 1 + \varepsilon + 4 \alpha_2 < 6 \theta - 3 \alpha_1 - 1 + \varepsilon 
\leqslant 6 \theta - \frac{13}{7} + \varepsilon < 2 - \theta - \varepsilon.
\end{equation}
Since we have $\alpha_2 > \frac{5 - 8 \theta}{6} - \varepsilon$ and $\theta - \frac{5 - 8 \theta}{6} + 2 \varepsilon < 2 - 3 \theta - 3 \varepsilon$ for all $\theta < \frac{17}{32} - \varepsilon$, we know that the condition $\alpha_1 < 2 - 3 \theta - 3 \varepsilon$ is satisfied automatically when $\boldsymbol{\alpha}_{2} \in A$. Of course, we often decide to use it instead of discarding the region $\frac{2}{7} \leqslant \alpha_1 \leqslant \frac{3}{7} + \varepsilon$ since the higher dimensional loss comes from integrals (130)--(132) is usually smaller than the one-dimensional loss.

Now we assume that $(\theta_1, \theta_2) \in \boldsymbol{E}$. We divide $\boldsymbol{E}$ into 11 subregions:
$$
\boldsymbol{E} = \boldsymbol{E}_{01} \cup \boldsymbol{E}_{02} \cup \boldsymbol{E}_{03} \cup \boldsymbol{E}_{04} \cup \boldsymbol{E}_{05} \cup \boldsymbol{E}_{06} \cup \boldsymbol{E}_{07} \cup \boldsymbol{E}_{08} \cup \boldsymbol{E}_{09} \cup \boldsymbol{E}_{10} \cup \boldsymbol{E}_{11},
$$
where
\begin{align}
\nonumber \boldsymbol{E}_{01} =&\ \left\{ (\theta_1, \theta_2) : \frac{1}{4} < \theta_1 \leqslant \frac{2}{7},\ \frac{1}{2}(1 - 2 \theta_1) < \theta_2 \leqslant \theta_1 \right\}, \\
\nonumber \boldsymbol{E}_{02} =&\ \left\{ (\theta_1, \theta_2) : \frac{2}{7} < \theta_1 \leqslant \frac{5}{14},\ \frac{1}{2}(1 - 2 \theta_1) < \theta_2 < \frac{1}{7}(6 - 14 \theta_1) \right\}, \\
\nonumber \boldsymbol{E}_{03} =&\ \left\{ (\theta_1, \theta_2) : \frac{2}{7} < \theta_1 \leqslant \frac{5}{14},\ \frac{1}{7}(6 - 14 \theta_1) < \theta_2 < \frac{1}{4}(2 - 3 \theta_1) \right\}, \\
\nonumber \boldsymbol{E}_{04} =&\ \left\{ (\theta_1, \theta_2) : \frac{5}{14} < \theta_1 \leqslant \frac{3}{8},\ \frac{1}{2}(1 - 2 \theta_1) < \theta_2 < \frac{1}{3}(-1 + 4 \theta_1) \right\}, \\
\nonumber \boldsymbol{E}_{05} =&\ \left\{ (\theta_1, \theta_2) : \frac{5}{14} < \theta_1 \leqslant \frac{3}{8},\ \frac{1}{3}(-1 + 4 \theta_1) < \theta_2 < \frac{1}{4}(2 - 3 \theta_1) \right\}, \\
\nonumber \boldsymbol{E}_{06} =&\ \left\{ (\theta_1, \theta_2) : \frac{3}{8} < \theta_1 \leqslant \frac{2}{5},\ \frac{1}{2}(1 - 2 \theta_1) < \theta_2 < \frac{1}{3}(2 - 4 \theta_1) \right\}, \\
\nonumber \boldsymbol{E}_{07} =&\ \left\{ (\theta_1, \theta_2) : \frac{3}{8} < \theta_1 \leqslant \frac{2}{5},\ \frac{1}{3}(2 - 4 \theta_1) < \theta_2 < \frac{1}{4}(2 - 3 \theta_1) \right\}, \\
\nonumber \boldsymbol{E}_{08} =&\ \left\{ (\theta_1, \theta_2) : \frac{2}{5} < \theta_1 < \frac{1}{2},\ \frac{1}{2}(1 - 2 \theta_1) < \theta_2 < \frac{1}{3}(2 - 4 \theta_1) \right\}, \\
\nonumber \boldsymbol{E}_{09} =&\ \left\{ (\theta_1, \theta_2) : \frac{2}{5} < \theta_1 < \frac{1}{2},\ \frac{1}{3}(2 - 4 \theta_1) < \theta_2 < 1 - 2 \theta_1 \right\}, \\
\nonumber \boldsymbol{E}_{10} =&\ \left\{ (\theta_1, \theta_2) : \frac{2}{5} < \theta_1 < \frac{1}{2},\ 1 - 2 \theta_1 < \theta_2 < \frac{1}{15}(11 - 20 \theta_1) \right\}, \\
\nonumber \boldsymbol{E}_{11} =&\ \left\{ (\theta_1, \theta_2) : \frac{1}{2} \leqslant \theta_1 < \frac{11}{20},\ 0 < \theta_2 < \frac{1}{15}(11 - 20 \theta_1) \right\}.
\end{align}
Note that we have a Type-II range $\left[\frac{3}{7} + \varepsilon,\ 1 - \theta - \varepsilon \right]$ for all $\theta \leqslant \frac{127}{224} - \varepsilon$ by Lemma~\ref{l216}. For $\theta \leqslant \frac{45}{89} - \varepsilon$ we also have a Type-II range $\left((\log x)^{\varepsilon - 1},\ \kappa \right]$ by Lemma~\ref{l24}, and we shall not repeatedly state these two ranges in the decompositions for the sake of simplicity.

\subsubsection{$\boldsymbol{E}_{01}$}
We divide $\boldsymbol{E}_{01}$ into 4 subregions:
$$
\boldsymbol{E}_{01} = \boldsymbol{E}_{0101} \cup \boldsymbol{E}_{0102} \cup \boldsymbol{E}_{0103} \cup \boldsymbol{E}_{0104},
$$
where
\begin{align}
\nonumber \boldsymbol{E}_{0101} =&\ \left\{ (\theta_1, \theta_2) : \frac{1}{4} < \theta_1 \leqslant \frac{2}{7},\ \frac{1}{2}(1 - 2 \theta_1) < \theta_2 < \frac{1}{3}(1 - \theta_1) \right\}, \\
\nonumber \boldsymbol{E}_{0102} =&\ \left\{ (\theta_1, \theta_2) : \frac{1}{4} < \theta_1 < \frac{6}{23},\ \frac{1}{3}(1 - \theta_1) \leqslant \theta_2 \leqslant \theta_1 \right. \\
\nonumber & \left. \qquad \qquad \quad \text{ or } \frac{6}{23} \leqslant \theta_1 \leqslant \frac{2}{7},\ \frac{1}{3}(1 - \theta_1) \leqslant \theta_2 < \frac{1}{13}(6 - 10 \theta_1) \right\}, \\
\nonumber \boldsymbol{E}_{0103} =&\ \left\{ (\theta_1, \theta_2) : \frac{6}{23} < \theta_1 < \frac{11}{40},\ \frac{1}{13}(6 - 10 \theta_1) \leqslant \theta_2 \leqslant \theta_1 \right. \\
\nonumber & \left. \qquad \qquad \quad \text{ or } \frac{11}{40} \leqslant \theta_1 \leqslant \frac{2}{7},\ \frac{1}{13}(6 - 10 \theta_1) \leqslant \theta_2 < \frac{1}{20}(11 - 20 \theta_1) \right\}, \\
\nonumber \boldsymbol{E}_{0104} =&\ \left\{ (\theta_1, \theta_2) : \frac{11}{40} \leqslant \theta_1 \leqslant \frac{2}{7},\ \frac{1}{20}(11 - 20 \theta_1) \leqslant \theta_2 \leqslant \theta_1 \right\}.
\end{align}
Note that we have $\theta < \frac{11}{20}$ for $(\theta_1, \theta_2) \in \boldsymbol{E}_{0101} \cup \boldsymbol{E}_{0102} \cup \boldsymbol{E}_{0103}$.

The Type-II range for $\boldsymbol{E}_{0101}$ is
\begin{equation}
\left[\varepsilon,\ \frac{1}{6}(5 - 8 \theta_1 - 8 \theta_2) - \varepsilon \right] \cup \left[\theta_2 + \varepsilon,\ \frac{1}{2}(1 - \theta_2) - \varepsilon \right].
\end{equation}
The decompositions in this case will be discussed later together with the case $\boldsymbol{E}_{0202}$.

The Type-II range for $\boldsymbol{E}_{0102}$ is
\begin{equation}
\left[\varepsilon,\ \frac{1}{6}(5 - 8 \theta_1 - 8 \theta_2) - \varepsilon \right] \cup \left[\theta_2 + \varepsilon,\ \frac{1}{2}(2 - \theta_1 - 4 \theta_2) - \varepsilon \right].
\end{equation}
The decompositions in this case will be discussed later together with the case $\boldsymbol{E}_{0203}$.

The Type-II range for $\boldsymbol{E}_{0103}$ is
\begin{equation}
\left[\varepsilon,\ \frac{1}{4}(2 - 2 \theta_1 - 5 \theta_2) - \varepsilon \right] \cup \left[2 \theta_1 + 2 \theta_2 - 1 + \varepsilon,\ \frac{1}{6}(5 - 8 \theta_1 - 8 \theta_2) - \varepsilon \right] \cup \left[\theta_2 + \varepsilon,\ \frac{1}{2}(2 - \theta_1 - 4 \theta_2) - \varepsilon \right].
\end{equation}
The decompositions are similar to which in the case $\boldsymbol{A}_{0104}$ in Section 3. When $\theta < \frac{17}{32}$, we use the middle Type-II range $\left[2 \theta - 1 + \varepsilon,\ \frac{5 - 8 \theta}{6} - \varepsilon \right]$ to give a ``starting point'' $\kappa = \frac{5 - 8 \theta}{6} - \varepsilon$, and the third Type-II range $\left[\theta_2 + \varepsilon,\ \frac{2 - \theta_1 - 4 \theta_2}{2} - \varepsilon \right]$ is used to subtract the contributions of those sums with products of variables lie in this range. When $\theta \geqslant \frac{17}{32}$, both the second and the third Type-II ranges are used explicitly to discard suitable sums. In this case $\kappa$ is reduced to $\frac{5 - 8 \theta}{12} - 3 \varepsilon$ or $\frac{3 - 5 \theta}{7} - 2 \varepsilon$, and if $\kappa < \frac{2 - 2 \theta_1 - 5 \theta_2}{4} - \varepsilon$, we can replace $\kappa$ with a larger $\frac{2 - 2 \theta_1 - 5 \theta_2}{4} - \varepsilon$ using the first Type-II range and Buchstab's identity as in (64). In many cases below (especially cases with $\theta \geqslant \frac{11}{20}$), we can also compare the value of $\kappa$ with $\frac{2 - 2 \theta_1 - 5 \theta_2}{4} - \varepsilon$ (or maybe other ``end point values''), and we shall not state similar process again.

The Type-II range for $\boldsymbol{E}_{0104}$ is
\begin{equation}
\left[\varepsilon,\ \frac{1}{4}(2 - 2 \theta_1 - 5 \theta_2) - \varepsilon \right] \cup \left[\theta_2 + \varepsilon,\ \frac{1}{2}(2 - \theta_1 - 4 \theta_2) - \varepsilon \right].
\end{equation}
The decompositions are similar to which in the case $\boldsymbol{A}_{05}$ in Section 3.

\subsubsection{$\boldsymbol{E}_{02}$}
We divide $\boldsymbol{E}_{02}$ into 15 subregions:
$$
\boldsymbol{E}_{02} = \boldsymbol{E}_{0201} \cup \boldsymbol{E}_{0202} \cup \boldsymbol{E}_{0203} \cup \boldsymbol{E}_{0204} \cup \boldsymbol{E}_{0205} \cup \boldsymbol{E}_{0206} \cup \boldsymbol{E}_{0207} \cup \boldsymbol{E}_{0208} \cup \boldsymbol{E}_{0209} \cup \boldsymbol{E}_{0210} \cup \boldsymbol{E}_{0211} \cup \boldsymbol{E}_{0212} \cup \boldsymbol{E}_{0213} \cup \boldsymbol{E}_{0214} \cup \boldsymbol{E}_{0215},
$$
where
\begin{align}
\nonumber \boldsymbol{E}_{0201} =&\ \left\{ (\theta_1, \theta_2) : \frac{1}{3} < \theta_1 \leqslant \frac{7}{20},\ \frac{1}{2}(1 - 2 \theta_1) < \theta_2 < \frac{1}{14}(5 - 8 \theta_1) \right. \\
\nonumber & \left. \qquad \qquad \quad \text{ or } \frac{7}{20} < \theta_1 < \frac{5}{14},\ \frac{1}{2}(1 - 2 \theta_1) < \theta_2 < \frac{1}{7}(6 - 14 \theta_1) \right\}, \\
\nonumber \boldsymbol{E}_{0202} =&\ \left\{ (\theta_1, \theta_2) : \frac{2}{7} < \theta_1 < \frac{5}{17},\ \frac{1}{2}(1 - 2 \theta_1) < \theta_2 < \frac{1}{3}(1 - \theta_1) \right. \\
\nonumber & \left. \qquad \qquad \quad \text{ or } \frac{5}{17} \leqslant \theta_1 < \frac{3}{10},\ \frac{1}{2}(1 - 2 \theta_1) < \theta_2 < 2 - 6 \theta_1 \right\}, \\
\nonumber \boldsymbol{E}_{0203} =&\ \left\{ (\theta_1, \theta_2) : \frac{2}{7} < \theta_1 < \frac{5}{17},\ \frac{1}{3}(1 - \theta_1) \leqslant \theta_2 < \frac{1}{13}(6 - 10 \theta_1) \right\}, \\
\nonumber \boldsymbol{E}_{0204} =&\ \left\{ (\theta_1, \theta_2) : \frac{2}{7} < \theta_1 \leqslant \frac{29}{100},\ \frac{1}{13}(6 - 10 \theta_1) \leqslant \theta_2 < \frac{1}{20}(11 - 20 \theta_1) \right. \\
\nonumber & \left. \qquad \qquad \quad \text{ or } \frac{29}{100} < \theta_1 < \frac{5}{17},\ \frac{1}{13}(6 - 10 \theta_1) \leqslant \theta_2 < 2 - 6 \theta_1 \right\}, \\
\nonumber \boldsymbol{E}_{0205} =&\ \left\{ (\theta_1, \theta_2) : \frac{2}{7} < \theta_1 < \frac{29}{100},\ \frac{1}{20}(11 - 20 \theta_1) \leqslant \theta_2 < 2 - 6 \theta_1 \right\}, \\
\nonumber \boldsymbol{E}_{0206} =&\ \left\{ (\theta_1, \theta_2) : \frac{5}{17} \leqslant \theta_1 \leqslant \frac{79}{266},\ 2 - 6 \theta_1 \leqslant \theta_2 < \frac{1}{13}(6 - 10 \theta_1) \right. \\
\nonumber & \left. \qquad \qquad \quad \text{ or } \frac{79}{266} < \theta_1 < \frac{3}{10},\ 2 - 6 \theta_1 \leqslant \theta_2 < \frac{1}{2}(7 - 22 \theta_1) \right\}, \\
\nonumber \boldsymbol{E}_{0207} =&\ \left\{ (\theta_1, \theta_2) : \frac{5}{17} \leqslant \theta_1 \leqslant \frac{19}{64},\ \frac{1}{13}(6 - 10 \theta_1) \leqslant \theta_2 < \frac{1}{3}(1 - \theta_1) \right. \\
\nonumber & \left. \qquad \qquad \quad \text{ or } \frac{19}{64} < \theta_1 < \frac{79}{266},\ \frac{1}{13}(6 - 10 \theta_1) \leqslant \theta_2 < \frac{1}{2}(7 - 22 \theta_1) \right\}, \\
\nonumber \boldsymbol{E}_{0208} =&\ \left\{ (\theta_1, \theta_2) : \frac{29}{100} < \theta_1 \leqslant \frac{19}{65},\ 2 - 6 \theta_1 \leqslant \theta_2 < \frac{1}{20}(11 - 20 \theta_1) \right. \\
\nonumber & \qquad \qquad \quad \text{ or } \frac{19}{65} < \theta_1 < \frac{5}{17},\ 2 - 6 \theta_1 \leqslant \theta_2 \leqslant \frac{1}{4}(6 - 17 \theta_1) \\
\nonumber & \qquad \qquad \quad \text{ or } \frac{5}{17} \leqslant \theta_1 \leqslant \frac{8}{27},\ \frac{1}{3}(1 - \theta_1) \leqslant \theta_2 \leqslant \frac{1}{4}(6 - 17 \theta_1) \\
\nonumber & \left. \qquad \qquad \quad \text{ or } \frac{8}{27} < \theta_1 \leqslant \frac{19}{64},\ \frac{1}{3}(1 - \theta_1) \leqslant \theta_2 \leqslant \frac{1}{2}(7 - 22 \theta_1) \right\}, \\
\nonumber \boldsymbol{E}_{0209} =&\ \left\{ (\theta_1, \theta_2) : \frac{2}{7} < \theta_1 \leqslant \frac{29}{100},\ 2 - 6 \theta_1 \leqslant \theta_2 < \frac{1}{4}(6 - 17 \theta_1) \right. \\
\nonumber & \left. \qquad \qquad \quad \text{ or } \frac{29}{100} < \theta_1 \leqslant \frac{19}{65},\ \frac{1}{20}(11 - 20 \theta_1) \leqslant \theta_2 \leqslant \frac{1}{4}(6 - 17 \theta_1) \right\}, \\
\nonumber \boldsymbol{E}_{0210} =&\ \left\{ (\theta_1, \theta_2) : \frac{79}{266} < \theta_1 < \frac{3}{10},\ \frac{1}{2}(7 - 22 \theta_1) \leqslant \theta_2 < \frac{1}{13}(6 - 10 \theta_1) \right. \\
\nonumber & \qquad \qquad \quad \text{ or } \frac{3}{10} \leqslant \theta_1 \leqslant \frac{9}{28},\ \frac{1}{2}(1 - 2 \theta_1) < \theta_2 < \frac{1}{13}(6 - 10 \theta_1) \\
\nonumber & \qquad \qquad \quad \text{ or } \frac{9}{28} < \theta_1 \leqslant \frac{1}{3},\ \frac{1}{2}(1 - 2 \theta_1) < \theta_2 < \frac{1}{7}(6 - 14 \theta_1) \\
\nonumber & \left. \qquad \qquad \quad \text{ or } \frac{1}{3} < \theta_1 < \frac{7}{20},\ \frac{1}{14}(5 - 8 \theta_1) < \theta_2 < \frac{1}{7}(6 - 14 \theta_1) \right\}, \\
\nonumber \boldsymbol{E}_{0211} =&\ \left\{ (\theta_1, \theta_2) : \frac{19}{64} < \theta_1 \leqslant \frac{79}{266},\ \frac{1}{2}(7 - 22 \theta_1) \leqslant \theta_2 < \frac{1}{3}(1 - \theta_1) \right. \\
\nonumber & \qquad \qquad \quad \text{ or } \frac{79}{266} < \theta_1 < \frac{4}{13},\ \frac{1}{13}(6 - 10 \theta_1) \leqslant \theta_2 < \frac{1}{3}(1 - \theta_1) \\
\nonumber & \qquad \qquad \quad \text{ or } \frac{4}{13} \leqslant \theta_1 < \frac{20}{63},\ \frac{1}{13}(6 - 10 \theta_1) \leqslant \theta_2 \leqslant \frac{1}{8}(4 - 7 \theta_1) \\
\nonumber & \left. \qquad \qquad \quad \text{ or } \frac{20}{63} \leqslant \theta_1 < \frac{9}{28},\ \frac{1}{13}(6 - 10 \theta_1) \leqslant \theta_2 < \frac{1}{7}(6 - 14 \theta_1) \right\}, \\
\nonumber \boldsymbol{E}_{0212} =&\ \left\{ (\theta_1, \theta_2) : \frac{8}{27} < \theta_1 \leqslant \frac{19}{64},\ \frac{1}{2}(7 - 22 \theta_1) < \theta_2 \leqslant \frac{1}{8}(4 - 7 \theta_1) \right. \\
\nonumber & \left. \qquad \qquad \quad \text{ or } \frac{19}{64} < \theta_1 < \frac{4}{13},\ \frac{1}{3}(1 - \theta_1) \leqslant \theta_2 \leqslant \frac{1}{8}(4 - 7 \theta_1) \right\}, \\
\nonumber \boldsymbol{E}_{0213} =&\ \left\{ (\theta_1, \theta_2) : \frac{4}{13} < \theta_1 \leqslant \frac{11}{35},\ \frac{1}{8}(4 - 7 \theta_1) < \theta_2 < \frac{1}{3}(1 - \theta_1) \right. \\
\nonumber & \left. \qquad \qquad \quad \text{ or } \frac{11}{35} < \theta_1 < \frac{20}{63},\ \frac{1}{8}(4 - 7 \theta_1) < \theta_2 < \frac{1}{7}(6 - 14 \theta_1) \right\}, \\
\nonumber \boldsymbol{E}_{0214} =&\ \left\{ (\theta_1, \theta_2) : \frac{19}{65} < \theta_1 \leqslant \frac{8}{27},\ \frac{1}{4}(6 - 17 \theta_1) < \theta_2 < \frac{1}{20}(11 - 20 \theta_1) \right. \\
\nonumber & \qquad \qquad \quad \text{ or } \frac{8}{27} < \theta_1 \leqslant \frac{43}{140},\ \frac{1}{8}(4 - 7 \theta_1) < \theta_2 < \frac{1}{20}(11 - 20 \theta_1) \\
\nonumber & \qquad \qquad \quad \text{ or } \frac{43}{140} < \theta_1 \leqslant \frac{4}{13},\ \frac{1}{8}(4 - 7 \theta_1) < \theta_2 < \frac{1}{7}(6 - 14 \theta_1) \\
\nonumber & \left. \qquad \qquad \quad \text{ or } \frac{4}{13} < \theta_1 < \frac{11}{35},\ \frac{1}{3}(1 - \theta_1) \leqslant \theta_2 < \frac{1}{7}(6 - 14 \theta_1) \right\}, \\
\nonumber \boldsymbol{E}_{0215} =&\ \left\{ (\theta_1, \theta_2) : \frac{2}{7} < \theta_1 \leqslant \frac{19}{65},\ \frac{1}{4}(6 - 17 \theta_1) < \theta_2 < \frac{1}{7}(6 - 14 \theta_1) \right. \\
\nonumber & \left. \qquad \qquad \quad \text{ or } \frac{19}{65} < \theta_1 < \frac{43}{140},\ \frac{1}{20}(11 - 20 \theta_1) \leqslant \theta_2 < \frac{1}{7}(6 - 14 \theta_1) \right\}.
\end{align}
Note that we have $\theta < \frac{11}{20}$ for $(\theta_1, \theta_2) \in \boldsymbol{E}_{0201} \cup \boldsymbol{E}_{0202} \cup \boldsymbol{E}_{0203} \cup \boldsymbol{E}_{0204} \cup \boldsymbol{E}_{0206} \cup \boldsymbol{E}_{0207} \cup \boldsymbol{E}_{0208} \cup \boldsymbol{E}_{0210} \cup \boldsymbol{E}_{0211} \cup \boldsymbol{E}_{0212} \cup \boldsymbol{E}_{0213} \cup \boldsymbol{E}_{0214}$.

The Type-II range for $\boldsymbol{E}_{0201}$ is
\begin{equation}
\left[\varepsilon,\ \frac{1}{6}(5 - 8 \theta_1 - 4 \theta_2) - \varepsilon \right] \cup \left[\theta_1 + \varepsilon,\ \frac{1}{2}(1 - \theta_2) - \varepsilon \right].
\end{equation}
The decompositions in this case will be discussed later together with the case $\boldsymbol{E}_{0501}$.

The Type-II range for $\boldsymbol{E}_{0202}$ is
\begin{equation}
\left[\varepsilon,\ \frac{1}{6}(5 - 8 \theta_1 - 8 \theta_2) - \varepsilon \right] \cup \left[\theta_2 + \varepsilon,\ \frac{1}{2}(1 - \theta_2) - \varepsilon \right].
\end{equation}
In this case we discuss both $\boldsymbol{E}_{0101}$ and $\boldsymbol{E}_{0202}$, since the Type-II ranges for them are same:
\begin{align}
\nonumber \boldsymbol{\mathcal{Z}_1} = \boldsymbol{E}_{0101} \cup \boldsymbol{E}_{0202} =&\ \left\{ (\theta_1, \theta_2) : \frac{1}{4} < \theta_1 < \frac{5}{17},\ \frac{1}{2}(1 - 2 \theta_1) < \theta_2 < \frac{1}{3}(1 - \theta_1) \right. \\
\nonumber & \left. \qquad \qquad \quad \text{ or } \frac{5}{17} \leqslant \theta_1 < \frac{3}{10},\ \frac{1}{2}(1 - 2 \theta_1) < \theta_2 < 2 - 6 \theta_1 \right\}.
\end{align}

First, we shall prove the following theorem.
\begin{theorem}\label{case1}
Let $(\theta_1, \theta_2) \in \boldsymbol{\mathcal{Z}_1}$. Suppose that we have
$$
\theta_1 + \theta_2 \leqslant \frac{11}{21} - \varepsilon.
$$
Then we have
$$
C_1^{\prime}(\theta_1, \theta_2) \leqslant 1 + \int_{1 - \theta}^{\frac{1}{2}} \frac{\omega\left(\frac{1-t}{t} \right)}{t^2} d t + O(\varepsilon) = 1 + \log\left(\frac{\theta}{1 - \theta} \right) + O(\varepsilon).
$$
\end{theorem}
\begin{proof}
Since we have $(\theta_1, \theta_2) \in \boldsymbol{\mathcal{Z}_1}$ and $\theta_1 + \theta_2 \leqslant \frac{11}{21} - \varepsilon$, all of the following conditions hold true:
$$
\kappa = \frac{5 - 8 \theta}{6} - \varepsilon > \frac{17}{126} - \varepsilon > \frac{1}{8} + \varepsilon, \quad 11 \theta_1 + 12 \theta_2 < 6 - \varepsilon, \quad 8 \theta_1 + 11 \theta_2 < 5 - 8 \varepsilon,
$$
$$
3 \theta_1 + 2 \theta_2 < \frac{11}{7} - \varepsilon, \quad \theta_1 + 3 \theta_2 \leqslant 1, \quad \theta_2 < \frac{1 - \theta}{2} \leqslant \frac{1}{4} \quad \text{and} \quad \frac{1}{2}(1 - \theta_2) > \frac{3}{8}.
$$
We can simplify our Type-II range in this case to
\begin{equation}
\left[\varepsilon,\ \frac{17}{126} - \varepsilon \right] \cup \left[\frac{1}{4} + \varepsilon,\ \frac{3}{8} - \varepsilon \right] \cup \left[\frac{3}{7} + \varepsilon,\ 1 - \theta - \varepsilon \right].
\end{equation}

Now, we can decompose our $\mathbbm{1}_{n \sim x, n = p}(n) = \psi\left(n, (2 x)^{\frac{1}{2}}\right)$ in a way similar to the decompositions in [\cite{LRB679}, Sections 6.1--6.5]. By Buchstab's identity, we have
\begin{align}
\nonumber \psi\left(n, (2 x)^{\frac{1}{2}}\right) =&\ \psi\left(n, x^{\kappa}\right) - \sum_{\substack{n = p_1 \beta \\ \kappa \leqslant \alpha_1 < \frac{1}{2} }} \psi\left(\beta, p_1 \right) \\
\nonumber =&\ \psi\left(n, x^{\kappa}\right) - \sum_{\substack{n = p_1 \beta \\ \kappa \leqslant \alpha_1 < \frac{3}{7} + \varepsilon }} \psi\left(\beta, p_1 \right) - \sum_{\substack{n = p_1 \beta \\ \frac{3}{7} + \varepsilon \leqslant \alpha_1 \leqslant 1 - \theta - \varepsilon }} \psi\left(\beta, p_1 \right) - \sum_{\substack{n = p_1 \beta \\ 1 - \theta - \varepsilon < \alpha_1 < \frac{1}{2} }} \psi\left(\beta, p_1 \right) \\
\nonumber =&\ \psi\left(n, x^{\kappa}\right) - \sum_{\substack{n = p_1 \beta \\ \kappa \leqslant \alpha_1 < \frac{3}{7} + \varepsilon }} \psi\left(\beta, x^{\kappa} \right) + \sum_{\substack{n = p_1 p_2 \beta \\ \kappa \leqslant \alpha_1 < \frac{3}{7} + \varepsilon \\ \kappa \leqslant \alpha_2 < \min\left(\alpha_1, \frac{1}{2}(1 - \alpha_1) \right) \\ \boldsymbol{\alpha}_2 \in \boldsymbol{G}_2 }} \psi\left(\beta, p_2 \right) + \sum_{\substack{n = p_1 p_2 \beta \\ \kappa \leqslant \alpha_1 < \frac{3}{7} + \varepsilon \\ \kappa \leqslant \alpha_2 < \min\left(\alpha_1, \frac{1}{2}(1 - \alpha_1) \right) \\ \boldsymbol{\alpha}_2 \in A \cup B }} \psi\left(\beta, p_2 \right) \\
\nonumber & + \sum_{\substack{n = p_1 p_2 \beta \\ \kappa \leqslant \alpha_1 < \frac{3}{7} + \varepsilon \\ \kappa \leqslant \alpha_2 < \min\left(\alpha_1, \frac{1}{2}(1 - \alpha_1) \right) \\ \boldsymbol{\alpha}_2 \in C }} \psi\left(\beta, p_2 \right) - \sum_{\substack{n = p_1 \beta \\ \frac{3}{7} + \varepsilon \leqslant \alpha_1 \leqslant 1 - \theta - \varepsilon }} \psi\left(\beta, p_1 \right) - \sum_{\substack{n = p_1 \beta \\ 1 - \theta - \varepsilon < \alpha_1 < \frac{1}{2} }} \psi\left(\beta, p_1 \right) \\
\nonumber =&\ \psi\left(n, x^{\kappa}\right) - \sum_{\substack{n = p_1 \beta \\ \kappa \leqslant \alpha_1 < \frac{3}{7} + \varepsilon }} \psi\left(\beta, x^{\kappa} \right) + \sum_{\substack{n = p_1 p_2 \beta \\ \kappa \leqslant \alpha_1 < \frac{3}{7} + \varepsilon \\ \kappa \leqslant \alpha_2 < \min\left(\alpha_1, \frac{1}{2}(1 - \alpha_1) \right) \\ \boldsymbol{\alpha}_2 \in \boldsymbol{G}_2 }} \psi\left(\beta, p_2 \right) \\
\nonumber & + \sum_{\substack{n = p_1 p_2 \beta \\ \kappa \leqslant \alpha_1 < \frac{3}{7} + \varepsilon \\ \kappa \leqslant \alpha_2 < \min\left(\alpha_1, \frac{1}{2}(1 - \alpha_1) \right) \\ \boldsymbol{\alpha}_2 \in A \cup B }} \psi\left(\beta, x^{\kappa} \right) - \sum_{\substack{n = p_1 p_2 p_3 \beta \\ \kappa \leqslant \alpha_1 < \frac{3}{7} + \varepsilon \\ \kappa \leqslant \alpha_2 < \min\left(\alpha_1, \frac{1}{2}(1 - \alpha_1) \right) \\ \boldsymbol{\alpha}_2 \in A \cup B \cup C \\ \kappa \leqslant \alpha_3 < \min\left(\alpha_2, \frac{1}{2}(1 - \alpha_1 - \alpha_2) \right) }} \psi\left(\beta, p_3 \right) \\
\nonumber & + \sum_{\substack{n = p_1 p_2 \beta \\ \kappa \leqslant \alpha_1 < \frac{3}{7} + \varepsilon \\ \kappa \leqslant \alpha_2 < \min\left(\alpha_1, \frac{1}{2}(1 - \alpha_1) \right) \\ \boldsymbol{\alpha}_2 \in C }} \psi\left(\beta, x^{\kappa} \right) - \sum_{\substack{n = p_1 \beta \\ \frac{3}{7} + \varepsilon \leqslant \alpha_1 \leqslant 1 - \theta - \varepsilon }} \psi\left(\beta, p_1 \right) - \sum_{\substack{n = p_1 \beta \\ 1 - \theta - \varepsilon < \alpha_1 < \frac{1}{2} }} \psi\left(\beta, p_1 \right) \\
=&\ \Sigma_{42201} - \Sigma_{42202} + \Sigma_{42203} + \Sigma_{42204} - \Sigma_{42205} + \Sigma_{42206} - \Sigma_{42207} - \Sigma_{42208}.
\end{align}
By Lemma~\ref{l212}, (129) holds for $f(n) = \Sigma_{42201}$ and $f(n) = \Sigma_{42202}$. By Lemma~\ref{l211} and Lemma~\ref{l216}, (129) holds for $f(n) = \Sigma_{42203}$ and $f(n) = \Sigma_{42207}$. By Lemma~\ref{l217}, (129) holds for $f(n) = \Sigma_{42204}$. We discard the whole of $\Sigma_{42208}$, leading to the added integral in Theorem~\ref{case1}. For the remaining sums, $\Sigma_{42205}$ only counts numbers with $4$ or more prime factors.

For $\Sigma_{42206}$, since we have $3 \theta_1 + 2 \theta_2 < \frac{11}{7} - \varepsilon$, $11 \theta_1 + 12 \theta_2 < 6 - \varepsilon$, $\theta < \frac{11}{21} < \frac{8}{15} - \varepsilon$ and a Type-II range $\left[\varepsilon,\ \frac{5 - 8 \theta}{6} - \varepsilon \right]$, we can use Lemma~\ref{l38} and a three-dimensional Harman's sieve to get a ``loss term''
\begin{equation}
\Sigma_{42209} = \sum_{\substack{n = m_1 m_2 m_3 \\ \kappa \leqslant \alpha_1 < \frac{3}{7} + \varepsilon \\ \kappa \leqslant \alpha_2 < \min\left(\alpha_1, \frac{1}{2}(1 - \alpha_1) \right) \\ \boldsymbol{\alpha}_2 \in C \\ \Omega(m_1 m_2) \geqslant 3 }} \psi\left(m_1 m_2 m_3, x^{\kappa} \right).
\end{equation}
Since $\Omega(m_1 m_2 m_3) \geqslant \Omega(m_1 m_2) + 1 \geqslant 4$, $\Sigma_{42209}$ only counts numbers with $4$ or more prime factors.

Now, the proof of Theorem~\ref{case1} reduces to showing that (129) holds for $f(n) = $ sums that count numbers with $4$ or more prime factors. Since $\kappa > \frac{1}{8}$, we have $\Omega(n) \leqslant 7$. Assume that $\Omega(n) \geqslant 4$, and we only need to consider the following 7 cases:

\textbf{Case 1: $\Omega(n) = 7$.}
Suppose that $n = p_1 \cdots p_7$ and $\alpha_1 \geqslant \alpha_2 \geqslant \cdots \geqslant \alpha_7$. Now we have $\alpha_6 + \alpha_7 \leqslant \frac{2}{7}$. Since we have $\theta_2 \leqslant \frac{1}{4} < \frac{2}{7}$ and $\frac{1}{2}(1 - \theta_2) > \frac{3}{8} > \frac{2}{7}$, if $\alpha_6 + \alpha_7 \in \left[\theta_2 + \varepsilon,\ \frac{3}{8} - \varepsilon \right]$, then (129) holds for $f(n)$. Otherwise we have $\alpha_6 + \alpha_7 < \theta_2 + \varepsilon$. But since we have $8 \theta_1 + 11 \theta_2 < 5 - 8 \varepsilon$, we get
\begin{equation}
\alpha_6 + \alpha_7 < \theta_2 + \varepsilon < \frac{5 - 8 \theta}{3} - 2 \varepsilon < \alpha_6 + \alpha_7,
\end{equation}
making a contradiction. Hence (129) holds for $f(n)$ in \textbf{Case 1}.

\textbf{Case 2: $\Omega(n) = 6$.}
Suppose that $n = p_1 \cdots p_6$ and $\alpha_1 \geqslant \alpha_2 \geqslant \cdots \geqslant \alpha_6$. Now we have $\alpha_5 + \alpha_6 \leqslant \frac{1}{3}$. Since we have $\theta_2 \leqslant \frac{1}{4} < \frac{1}{3}$ and $\frac{1}{2}(1 - \theta_2) > \frac{3}{8} > \frac{1}{3}$, if $\alpha_5 + \alpha_6 \in \left[\theta_2 + \varepsilon,\ \frac{3}{8} - \varepsilon \right]$, then (129) holds for $f(n)$. Otherwise we have $\alpha_5 + \alpha_6 < \theta_2 + \varepsilon$. But since we have $8 \theta_1 + 11 \theta_2 < 5 - 8 \varepsilon$, we get
\begin{equation}
\alpha_5 + \alpha_6 < \theta_2 + \varepsilon < \frac{5 - 8 \theta}{3} - 2 \varepsilon < \alpha_5 + \alpha_6,
\end{equation}
making a contradiction. Hence (129) holds for $f(n)$ in \textbf{Case 2}.

\textbf{Case 3: $\Omega(n) = 4 \text{ or } 5$, no product of variables lies in $\left(1 - \theta - \varepsilon,\ \theta + \varepsilon \right)$.}
Suppose that $n = p_1 \cdots p_k$ and $\alpha_1 > \alpha_2 > \cdots > \alpha_k$, while the remaining part $\alpha_i = \alpha_j$ gives a ``loss'' of size $O(\varepsilon)$. In this case we can ``view'' $\left[\frac{3}{7} + \varepsilon,\ \frac{4}{7} - \varepsilon \right]$ as a ``fake'' Type-II range. Since we have $\theta < \frac{11}{21} < \frac{9}{17}$, the results in [\cite{MaynardLargeModuliI}, Chapter 9] can be applied here if any of the following conditions holds:

(1). $\alpha_k \geqslant \frac{1}{7} + \varepsilon$;

(2). $\Omega(n) = 5$, $\alpha_3 + \alpha_4 + \alpha_5 > \frac{4}{7} - \varepsilon$;

(3). $\Omega(n) = 4$, $\alpha_1 \leqslant \frac{3}{7} + \varepsilon$, $\alpha_1 + \alpha_4 > \frac{4}{7} - \varepsilon$.

By Condition (1), we can assume that $\frac{1}{8} < \alpha_k < \frac{1}{7} + \varepsilon$. We first consider the subcase $\Omega(n) = 5$. By Condition (2) and the ``fake'' Type-II range $\left[\frac{3}{7} + \varepsilon,\ \frac{4}{7} - \varepsilon \right]$, we can also assume that $\alpha_3 + \alpha_4 + \alpha_5 < \frac{3}{7} + \varepsilon$. This means that $\alpha_1 + \alpha_2 > \frac{4}{7} - \varepsilon$, and thus $\alpha_1 > \frac{2}{7} - \frac{1}{2} \varepsilon$. Since $\theta_2 < \frac{1}{4} < \frac{2}{7} - \frac{1}{2} \varepsilon < \frac{3}{8} < \frac{1}{2}(1 - \theta_2)$, we can assume that $\alpha_1 > \frac{3}{8} - \varepsilon$. Now we have $\frac{3}{7} + \varepsilon < \frac{1}{2} - \varepsilon = \frac{3}{8} - \varepsilon + \frac{1}{8} < \alpha_1 + \alpha_5$, and we can assume that $\alpha_1 + \alpha_5 > \frac{4}{7} - \varepsilon$. Since $\alpha_5 < \frac{1}{7} + \varepsilon$, we have $\alpha_1 > \frac{3}{7} - 2 \varepsilon$. If $\alpha_1 \geqslant \frac{3}{7} + \varepsilon$, we can assume that $\alpha_1 > \frac{4}{7} - \varepsilon$. But now we have $1 = \alpha_1 + \alpha_2 + \alpha_3 + \alpha_4 + \alpha_5 > \frac{4}{7} - \varepsilon + 4 \cdot \frac{1}{8} > 1$, making a contradiction. The remaining part is $\frac{3}{7} - 2 \varepsilon < \alpha_1 < \frac{3}{7} + \varepsilon$, leading to a ``loss'' of size $O(\varepsilon)$. Hence (129) holds for $f(n)$ with an error $O(\varepsilon)$ in \textbf{Case 3} if $\Omega(n) = 5$.

Next, we consider the second subcase $\Omega(n) = 4$. Suppose that $\alpha_1 + \alpha_4 \leqslant \frac{4}{7} - \varepsilon$. By Condition (3) and the ``fake'' Type-II range $\left[\frac{3}{7} + \varepsilon,\ \frac{4}{7} - \varepsilon \right]$, we can assume that $\alpha_1 + \alpha_4 < \frac{3}{7} + \varepsilon$. Since $\alpha_4 > \frac{1}{8}$, we have $\alpha_1 < \frac{17}{56} + \varepsilon = \frac{3}{7} + \varepsilon - \frac{1}{8}$. Since $\theta_2 < \frac{1}{4} < \frac{17}{56} + \varepsilon < \frac{3}{8} < \frac{1}{2}(1 - \theta_2)$, we can assume that $\alpha_1 \leqslant \frac{1}{4} + \varepsilon$. If $\alpha_1 \leqslant \frac{1}{4}$, then we have $1 = \alpha_1 + \alpha_2 + \alpha_3 + \alpha_4 < 4 \cdot \frac{1}{4} = 1$, making a contradiction. The remaining part is $\frac{1}{4} < \alpha_1 \leqslant \frac{1}{4} + \varepsilon$, leading to a ``loss'' of size $O(\varepsilon)$.

Now, by the ``fake'' Type-II range $\left[\frac{3}{7} + \varepsilon,\ \frac{4}{7} - \varepsilon \right]$, we only need to consider the subcase $\Omega(n) = 4$ with $\alpha_1 > \frac{4}{7} - \varepsilon$. This means that $\alpha_2 + \alpha_3 + \alpha_4 < \frac{3}{7} + \varepsilon$. Since $\alpha_2 > \alpha_3 > \alpha_4 > \frac{1}{8} + \varepsilon$, we have $\alpha_3 + \alpha_4 > \frac{1}{4} + \varepsilon > \theta_2 + \varepsilon$. Thus, we can assume that $\alpha_3 + \alpha_4 > \frac{3}{8} - \varepsilon$. But now we have
\begin{equation}
\frac{1}{2} < \frac{3}{2} \left(\frac{3}{8} - \varepsilon \right) < \frac{3}{2}(\alpha_3 + \alpha_4) < \alpha_2 + \alpha_3 + \alpha_4 < \frac{3}{7} + \varepsilon,
\end{equation}
making a contradiction. Hence (129) holds for $f(n)$ with an error $O(\varepsilon)$ in \textbf{Case 3} if $\Omega(n) = 4$.

\textbf{Case 4: $\Omega(n) = 5$, one variable lies in $\left(1 - \theta - \varepsilon,\ \theta + \varepsilon \right)$.}
Suppose that $n = p_1 \cdots p_5$ and $\alpha_1 \geqslant \alpha_2 \geqslant \cdots \geqslant \alpha_5$. Now we have $\alpha_1 > 1 - \theta - \varepsilon$ and $\alpha_2 \geqslant \cdots \geqslant \alpha_5 \geqslant \frac{5 - 8 \theta}{6} - \varepsilon$. But since $\theta < \frac{11}{21} < \frac{10}{19} - \varepsilon$, we have
\begin{equation}
1 - \theta - \varepsilon + 4 \left(\frac{5 - 8 \theta}{6} - \varepsilon \right) > 1,
\end{equation}
making a contradiction. Hence \textbf{Case 4} is empty.

\textbf{Case 5: $\Omega(n) = 5$, a product of two variables lies in $\left(1 - \theta - \varepsilon,\ \theta + \varepsilon \right)$.}
Suppose that $n = p_1 \cdots p_5$ and $\alpha_1 + \alpha_2 \in \left(1 - \theta - \varepsilon,\ \theta + \varepsilon \right)$ (of course, $\alpha_3 + \alpha_4 + \alpha_5 \in \left(1 - \theta - \varepsilon,\ \theta + \varepsilon \right)$ too). Since $\alpha_i + \alpha_j \geqslant \frac{5 - 8 \theta}{3} - 2 \varepsilon > \theta_2 + \varepsilon$ (see \textbf{Case 1}), we can assume that $\alpha_i + \alpha_j > \frac{1}{2}(1 - \theta_2) - \varepsilon > \frac{3}{8} - \varepsilon$ for all $i,j \in \left\{1,2,3,4,5\right\}$ (otherwise $\alpha_i + \alpha_j \in \left[\theta_2 + \varepsilon,\ \frac{1}{2}(1 - \theta_2) - \varepsilon \right]$, and (129) holds for $f(n)$). But now we have
\begin{equation}
\alpha_3 + \alpha_4 + \alpha_5 = \frac{(\alpha_3 + \alpha_4) + (\alpha_3 + \alpha_5) + (\alpha_4 + \alpha_5)}{2} > \frac{3}{2} \left(\frac{3}{8} - \varepsilon \right) > \frac{9}{16} - 2 \varepsilon > \frac{11}{21} + 2 \varepsilon > \theta + \varepsilon,
\end{equation}
making a contradiction. Hence (129) holds for $f(n)$ in \textbf{Case 5}.

\textbf{Case 6: $\Omega(n) = 4$, a product of two variables lies in $\left(1 - \theta - \varepsilon,\ \theta + \varepsilon \right)$.}
Suppose that $n = p_1 p_2 p_3 p_4$ and $\alpha_1 + \alpha_2 \in \left(1 - \theta - \varepsilon,\ \theta + \varepsilon \right)$ (of course, $\alpha_3 + \alpha_4 \in \left(1 - \theta - \varepsilon,\ \theta + \varepsilon \right)$ too). Without loss of generality, we further assume that $\alpha_1 > \alpha_2$ and $\alpha_3 > \alpha_4$, while the remaining part $\alpha_i = \alpha_j$ gives a ``loss'' of size $O(\varepsilon)$. Now we have $\alpha_2, \alpha_4 < \frac{\theta + \varepsilon}{2}$. Since $\frac{1}{2} < \theta \leqslant \frac{11}{21}$, we have $\frac{\theta + \varepsilon}{2} < \frac{11}{42} + \varepsilon < \frac{3}{8} - \varepsilon < \frac{1}{2}(1 - \theta_2) - \varepsilon$. Now we can assume that $\alpha_2, \alpha_4 < \theta_2 + \varepsilon$ (otherwise $\alpha_2 \text{ or } \alpha_4 \in \left[\theta_2 + \varepsilon,\ \frac{1}{2}(1 - \theta_2) - \varepsilon \right]$, and (129) holds for $f(n)$). Since we can assume that $\alpha_2 + \alpha_4 > \frac{3}{8} - \varepsilon$ (see \textbf{Case 5}), if $\alpha_1, \alpha_3 > \frac{3}{8} - \varepsilon$, then $1 = \alpha_1 + \alpha_3 + (\alpha_2 + \alpha_4) > \frac{9}{8} - 3 \varepsilon > 1$, making a contradiction. Now we assume that $\min(\alpha_1, \alpha_3) \leqslant \frac{3}{8} - \varepsilon < \frac{1}{2}(1 - \theta_2) - \varepsilon$. Without loss of generality, we further assume that $\alpha_1 > \alpha_3$, while the remaining part $\alpha_1 = \alpha_3$ gives a ``loss'' of size $O(\varepsilon)$. If $\alpha_3 \in \left[\theta_2 + \varepsilon,\ \frac{1}{2}(1 - \theta_2) - \varepsilon \right]$, then (129) holds for $f(n)$. Otherwise we have $\alpha_2, \alpha_3, \alpha_4 < \theta_2 + \varepsilon$. If $\alpha_3 < \theta_2 - 2 \varepsilon$, then we have $\alpha_3 + \alpha_4 \leqslant 2 \theta_2 - \varepsilon \leqslant 1 - \theta - \varepsilon$, making a contradiction with the assumption $\alpha_3 + \alpha_4 \in \left(1 - \theta - \varepsilon,\ \theta + \varepsilon \right)$. The remaining part is $\theta_2 - 2 \varepsilon < \alpha_3 \leqslant \theta_2 + \varepsilon$, leading to a ``loss'' of size $O(\varepsilon)$. Hence (129) holds for $f(n)$ with an error $O(\varepsilon)$ in \textbf{Case 6}.

\textbf{Case 7: $\Omega(n) = 4$, one variable lies in $\left(1 - \theta - \varepsilon,\ \theta + \varepsilon \right)$.}
Suppose that $n = p_0 p_1 p_2 p_3$ and $\alpha_0 \geqslant \alpha_1 \geqslant \alpha_2 \geqslant \alpha_3$. Now we have $\alpha_0 \in \left(1 - \theta - \varepsilon,\ \theta + \varepsilon \right)$ and $\alpha_1 + \alpha_2 + \alpha_3 \in \left(1 - \theta - \varepsilon,\ \theta + \varepsilon \right)$. We also have $\alpha_1 + \alpha_3 \leqslant \frac{1}{2}$ and $\alpha_2 \leqslant \frac{1}{3}$. Note that $\Sigma_{42209}$ only counts numbers with all prime factors smaller than $\max(m_1, m_2, m_3) \leqslant x^{\frac{3}{7} + \varepsilon}$, one can easily find that this type of $n$ will not be counted in $\Sigma_{42209}$ by a simple observation. Hence, this type of $n$ will only be counted in one part of $\Sigma_{42205}$:
\begin{equation}
\sum_{\substack{n = p_1 p_2 p_3 \beta \\ 1 - \theta - \varepsilon < \alpha_1 + \alpha_2 + \alpha_3 < \theta + \varepsilon }} \psi\left(\beta, p_3 \right),
\end{equation}
where $\beta$ in (153) is a large prime. Using Buchstab's identity, we have
\begin{equation}
\sum_{\substack{n = p_1 p_2 p_3 \beta \\ 1 - \theta - \varepsilon < \alpha_1 + \alpha_2 + \alpha_3 < \theta + \varepsilon }} \psi\left(\beta, p_3 \right) = \sum_{\substack{n = p_1 p_2 p_3 \beta \\ 1 - \theta - \varepsilon < \alpha_1 + \alpha_2 + \alpha_3 < \theta + \varepsilon }} \psi\left(\beta, x^{\kappa} \right) - \sum_{\substack{n = p_1 p_2 p_3 p_4 \beta \\ 1 - \theta - \varepsilon < \alpha_1 + \alpha_2 + \alpha_3 < \theta + \varepsilon \\ \kappa \leqslant \alpha_4 < \min\left(\alpha_3, \frac{1}{2}(1 - \alpha_1 - \alpha_2 - \alpha_3) \right) }} \psi\left(\beta, p_4 \right).
\end{equation}
By Lemma~\ref{l217}, (129) holds for the first sum in (154) since $(\alpha_1 + \alpha_2, \alpha_3) \in A$, and we only need to deal with the second sum in (154) that counts numbers with $5$ or more prime factors, with a product of $3$ variables lies in $\left(1 - \theta - \varepsilon,\ \theta + \varepsilon \right)$. Now we reduce this case to \textbf{Cases 1, 2 and 5}. Since the second sum in (154) is a positive sum that cannot be dropped directly, we need (129) holds for $f(n)$ in \textbf{Cases 1, 2 and 5} without an error $O(\varepsilon)$. By the discussions above, we know that (129) holds for $f(n)$ in \textbf{Cases 1, 2 and 5} without an error $O(\varepsilon)$, and thus (129) holds for $f(n)$ in \textbf{Case 7}.

Combining all the 7 cases above, the proof of Theorem~\ref{case1} is completed.
\end{proof}
\begin{remark*}
When $(\theta_1, \theta_2)$ lies in the boundary of this region:
$$
\theta_1 \in \left[\frac{1}{4},\ \frac{2}{7} - \varepsilon \right] \cup \left[\frac{2}{7} + \varepsilon,\ \frac{3}{10} - \varepsilon \right], \quad \theta_2 = \min\left(\frac{1 - \theta_1}{3}, \frac{11 - 21 \theta_1}{21}, 2 - 6 \theta_1 \right),
$$
we can use the exactly same decomposing process to prove that
$$
C_1^{\prime}(\theta_1, \theta_2) \leqslant 1 + \int_{1 - \theta}^{\frac{1}{2}} \frac{\omega\left(\frac{1-t}{t} \right)}{t^2} d t + O(\varepsilon) = 1 + \log\left(\frac{\theta}{1 - \theta} \right) + O(\varepsilon).
$$
For example, we have $C_1^{\prime}(0.28, 0.24) \leqslant 1 + \log\left(\frac{13}{12} \right) + O(\varepsilon) < 1.0801$. In order to prove this, we need to focus on the conditions needed in the proof of Theorem~\ref{case1}. When $\theta_1 + \theta_2 = \frac{11}{21}$, we move the part of $\Sigma_{42204}$ with $\boldsymbol{\alpha}_2 \in B$ into $\Sigma_{42206}$ (one can see the sums $\Sigma_{42904}$ and $\Sigma_{42906}$ in the proof of Theorem~\ref{case4} later). A simple verification yields $(\theta_1, \theta_2) \in \boldsymbol{T}$ except for $\theta_1 \in \left(\frac{2}{7} - \varepsilon, \frac{2}{7} + \varepsilon \right)$, and we can use the new three-dimensional Harman's sieve to handle $\Sigma_{42206}$. By the discussions above, we know that (129) holds for $f(n) = \Sigma_{42201}$, $f(n) = \Sigma_{42202}$, $f(n) = \Sigma_{42203}$, $f(n) = \Sigma_{42204}$ and $f(n) = \Sigma_{42207}$, and we only need to focus on the sums with $\Omega(n) \geqslant 4$. Since $8 \theta_1 + 11 \theta_2 < 5 - 8 \varepsilon$ still holds, we can use our Type-II range (145) to show that (129) holds for $f(n)$ in \textbf{Cases 1, 2, 4, 5 and 7} and for $f(n)$ in \textbf{Cases 3 and 6} with an error $O(\varepsilon)$.
\end{remark*}

For the remaining parts of $\boldsymbol{\mathcal{Z}_1}$, the decompositions are similar to which in the case $\boldsymbol{A}_{0101}$ in Section 3. We also use the Type-II range $\left[\theta_2 + \varepsilon,\ \frac{1}{2}(1 - \theta_2) - \varepsilon \right]$ to discard sums with products of variables lie in this range.

The Type-II range for $\boldsymbol{E}_{0203}$ is
\begin{equation}
\left[\varepsilon,\ \frac{1}{6}(5 - 8 \theta_1 - 8 \theta_2) - \varepsilon \right] \cup \left[\theta_2 + \varepsilon,\ \frac{1}{2}(2 - \theta_1 - 4 \theta_2) - \varepsilon \right].
\end{equation}
In this case we discuss both $\boldsymbol{E}_{0102}$ and $\boldsymbol{E}_{0203}$, since the Type-II ranges for them are same:
\begin{align}
\nonumber \boldsymbol{\mathcal{Z}_2} = \boldsymbol{E}_{0102} \cup \boldsymbol{E}_{0203} =&\ \left\{ (\theta_1, \theta_2) : \frac{1}{4} < \theta_1 < \frac{6}{23},\ \frac{1}{3}(1 - \theta_1) \leqslant \theta_2 \leqslant \theta_1 \right. \\
\nonumber & \left. \qquad \qquad \quad \text{ or } \frac{6}{23} \leqslant \theta_1 < \frac{5}{17},\ \frac{1}{3}(1 - \theta_1) \leqslant \theta_2 < \frac{1}{13}(6 - 10 \theta_1) \right\}.
\end{align}

First, we shall prove the following theorem.
\begin{theorem}\label{case2}
Let $(\theta_1, \theta_2) \in \boldsymbol{\mathcal{Z}_2}$. Suppose that we have
$$
\theta_2 \leqslant \frac{1}{4} - \varepsilon \quad \text{and} \quad 17 \theta_1 + 26 \theta_2 \leqslant 11 - 14 \varepsilon.
$$
Then we have
$$
C_1^{\prime}(\theta_1, \theta_2) \leqslant 1 + \int_{1 - \theta}^{\frac{1}{2}} \frac{\omega\left(\frac{1-t}{t} \right)}{t^2} d t + O(\varepsilon) = 1 + \log\left(\frac{\theta}{1 - \theta} \right) + O(\varepsilon).
$$
\end{theorem}
\begin{proof}
Since we have $(\theta_1, \theta_2) \in \boldsymbol{\mathcal{Z}_2}$, $\theta_2 \leqslant \frac{1}{4} - \varepsilon$ and $17 \theta_1 + 26 \theta_2 \leqslant 11 - 14 \varepsilon$, all of the following conditions hold true:
$$
\kappa = \frac{5 - 8 \theta}{6} - \varepsilon > \frac{17}{126} - \varepsilon > \frac{1}{8} + \varepsilon, \quad 11 \theta_1 + 12 \theta_2 < 6 - \varepsilon, \quad 8 \theta_1 + 11 \theta_2 < 5 - 8 \varepsilon,
$$
$$
3 \theta_1 + 2 \theta_2 < \frac{11}{7} - \varepsilon, \quad \frac{1}{2} \leqslant \theta_1 + \theta_2 \leqslant \frac{11}{21} - \varepsilon, \quad 7 \theta_1 + 16 \theta_2 < 6 - 8 \varepsilon \quad \text{and} \quad \frac{1}{2}(2 - \theta_1 - 4 \theta_2) > \frac{9}{25}.
$$
We can simplify our Type-II range in this case to
\begin{equation}
\left[\varepsilon,\ \frac{17}{126} - \varepsilon \right] \cup \left[\frac{1}{4} + \varepsilon,\ \frac{9}{25} - \varepsilon \right] \cup \left[\frac{3}{7} + \varepsilon,\ 1 - \theta - \varepsilon \right].
\end{equation}

Now, we can decompose our $\mathbbm{1}_{n \sim x, n = p}(n) = \psi\left(n, (2 x)^{\frac{1}{2}}\right)$ in a way similar to the decompositions in the proof of Theorem~\ref{case1}. By Buchstab's identity, we have
\begin{align}
\nonumber \psi\left(n, (2 x)^{\frac{1}{2}}\right) =&\ \psi\left(n, x^{\kappa}\right) - \sum_{\substack{n = p_1 \beta \\ \kappa \leqslant \alpha_1 < \frac{1}{2} }} \psi\left(\beta, p_1 \right) \\
\nonumber =&\ \psi\left(n, x^{\kappa}\right) - \sum_{\substack{n = p_1 \beta \\ \kappa \leqslant \alpha_1 < \frac{3}{7} + \varepsilon }} \psi\left(\beta, p_1 \right) - \sum_{\substack{n = p_1 \beta \\ \frac{3}{7} + \varepsilon \leqslant \alpha_1 \leqslant 1 - \theta - \varepsilon }} \psi\left(\beta, p_1 \right) - \sum_{\substack{n = p_1 \beta \\ 1 - \theta - \varepsilon < \alpha_1 < \frac{1}{2} }} \psi\left(\beta, p_1 \right) \\
\nonumber =&\ \psi\left(n, x^{\kappa}\right) - \sum_{\substack{n = p_1 \beta \\ \kappa \leqslant \alpha_1 < \frac{3}{7} + \varepsilon }} \psi\left(\beta, x^{\kappa} \right) + \sum_{\substack{n = p_1 p_2 \beta \\ \kappa \leqslant \alpha_1 < \frac{3}{7} + \varepsilon \\ \kappa \leqslant \alpha_2 < \min\left(\alpha_1, \frac{1}{2}(1 - \alpha_1) \right) \\ \boldsymbol{\alpha}_2 \in \boldsymbol{G}_2 }} \psi\left(\beta, p_2 \right) + \sum_{\substack{n = p_1 p_2 \beta \\ \kappa \leqslant \alpha_1 < \frac{3}{7} + \varepsilon \\ \kappa \leqslant \alpha_2 < \min\left(\alpha_1, \frac{1}{2}(1 - \alpha_1) \right) \\ \boldsymbol{\alpha}_2 \in A \cup B }} \psi\left(\beta, p_2 \right) \\
\nonumber & + \sum_{\substack{n = p_1 p_2 \beta \\ \kappa \leqslant \alpha_1 < \frac{3}{7} + \varepsilon \\ \kappa \leqslant \alpha_2 < \min\left(\alpha_1, \frac{1}{2}(1 - \alpha_1) \right) \\ \boldsymbol{\alpha}_2 \in C }} \psi\left(\beta, p_2 \right) - \sum_{\substack{n = p_1 \beta \\ \frac{3}{7} + \varepsilon \leqslant \alpha_1 \leqslant 1 - \theta - \varepsilon }} \psi\left(\beta, p_1 \right) - \sum_{\substack{n = p_1 \beta \\ 1 - \theta - \varepsilon < \alpha_1 < \frac{1}{2} }} \psi\left(\beta, p_1 \right) \\
\nonumber =&\ \psi\left(n, x^{\kappa}\right) - \sum_{\substack{n = p_1 \beta \\ \kappa \leqslant \alpha_1 < \frac{3}{7} + \varepsilon }} \psi\left(\beta, x^{\kappa} \right) + \sum_{\substack{n = p_1 p_2 \beta \\ \kappa \leqslant \alpha_1 < \frac{3}{7} + \varepsilon \\ \kappa \leqslant \alpha_2 < \min\left(\alpha_1, \frac{1}{2}(1 - \alpha_1) \right) \\ \boldsymbol{\alpha}_2 \in \boldsymbol{G}_2 }} \psi\left(\beta, p_2 \right) \\
\nonumber & + \sum_{\substack{n = p_1 p_2 \beta \\ \kappa \leqslant \alpha_1 < \frac{3}{7} + \varepsilon \\ \kappa \leqslant \alpha_2 < \min\left(\alpha_1, \frac{1}{2}(1 - \alpha_1) \right) \\ \boldsymbol{\alpha}_2 \in A \cup B }} \psi\left(\beta, x^{\kappa} \right) - \sum_{\substack{n = p_1 p_2 p_3 \beta \\ \kappa \leqslant \alpha_1 < \frac{3}{7} + \varepsilon \\ \kappa \leqslant \alpha_2 < \min\left(\alpha_1, \frac{1}{2}(1 - \alpha_1) \right) \\ \boldsymbol{\alpha}_2 \in A \cup B \cup C \\ \kappa \leqslant \alpha_3 < \min\left(\alpha_2, \frac{1}{2}(1 - \alpha_1 - \alpha_2) \right) }} \psi\left(\beta, p_3 \right) \\
\nonumber & + \sum_{\substack{n = p_1 p_2 \beta \\ \kappa \leqslant \alpha_1 < \frac{3}{7} + \varepsilon \\ \kappa \leqslant \alpha_2 < \min\left(\alpha_1, \frac{1}{2}(1 - \alpha_1) \right) \\ \boldsymbol{\alpha}_2 \in C }} \psi\left(\beta, x^{\kappa} \right) - \sum_{\substack{n = p_1 \beta \\ \frac{3}{7} + \varepsilon \leqslant \alpha_1 \leqslant 1 - \theta - \varepsilon }} \psi\left(\beta, p_1 \right) - \sum_{\substack{n = p_1 \beta \\ 1 - \theta - \varepsilon < \alpha_1 < \frac{1}{2} }} \psi\left(\beta, p_1 \right) \\
=&\ \Sigma_{42301} - \Sigma_{42302} + \Sigma_{42303} + \Sigma_{42304} - \Sigma_{42305} + \Sigma_{42306} - \Sigma_{42307} - \Sigma_{42308}.
\end{align}
By Lemma~\ref{l212}, (129) holds for $f(n) = \Sigma_{42301}$ and $f(n) = \Sigma_{42302}$. By Lemma~\ref{l211} and Lemma~\ref{l216}, (129) holds for $f(n) = \Sigma_{42303}$ and $f(n) = \Sigma_{42307}$. By Lemma~\ref{l217}, (129) holds for $f(n) = \Sigma_{42304}$. We discard the whole of $\Sigma_{42308}$, leading to the added integral in Theorem~\ref{case2}. For the remaining sums, $\Sigma_{42305}$ only counts numbers with $4$ or more prime factors.

For $\Sigma_{42306}$, since we have $3 \theta_1 + 2 \theta_2 < \frac{11}{7} - \varepsilon$, $11 \theta_1 + 12 \theta_2 < 6 - \varepsilon$, $\theta < \frac{11}{21} < \frac{8}{15} - \varepsilon$ and a Type-II range $\left[\varepsilon,\ \frac{5 - 8 \theta}{6} - \varepsilon \right]$, we can use Lemma~\ref{l38} and a three-dimensional Harman's sieve to get a ``loss term''
\begin{equation}
\Sigma_{42309} = \sum_{\substack{n = m_1 m_2 m_3 \\ \kappa \leqslant \alpha_1 < \frac{3}{7} + \varepsilon \\ \kappa \leqslant \alpha_2 < \min\left(\alpha_1, \frac{1}{2}(1 - \alpha_1) \right) \\ \boldsymbol{\alpha}_2 \in C \\ \Omega(m_1 m_2) \geqslant 3 }} \psi\left(m_1 m_2 m_3, x^{\kappa} \right).
\end{equation}
Since $\Omega(m_1 m_2 m_3) \geqslant \Omega(m_1 m_2) + 1 \geqslant 4$, $\Sigma_{42309}$ only counts numbers with $4$ or more prime factors.

Now, the proof of Theorem~\ref{case2} reduces to showing that (129) holds for $f(n) = $ sums that count numbers with $4$ or more prime factors. Since $\kappa > \frac{1}{8}$, we have $\Omega(n) \leqslant 7$. Assume that $\Omega(n) \geqslant 4$, and we only need to consider the following 7 cases:

\textbf{Case 1: $\Omega(n) = 7$.}
Suppose that $n = p_1 \cdots p_7$ and $\alpha_1 \geqslant \alpha_2 \geqslant \cdots \geqslant \alpha_7$. Now we have $\alpha_6 + \alpha_7 \leqslant \frac{2}{7}$. Since we have $\theta_2 \leqslant \frac{1}{4} < \frac{2}{7}$ and $\frac{1}{2}(2 - \theta_1 - 4 \theta_2) > \frac{9}{25} > \frac{2}{7}$, if $\alpha_6 + \alpha_7 \in \left[\theta_2 + \varepsilon,\ \frac{9}{25} - \varepsilon \right]$, then (129) holds for $f(n)$. Otherwise we have $\alpha_6 + \alpha_7 < \theta_2 + \varepsilon$. But since we have $8 \theta_1 + 11 \theta_2 < 5 - 8 \varepsilon$, we get
\begin{equation}
\alpha_6 + \alpha_7 < \theta_2 + \varepsilon < \frac{5 - 8 \theta}{3} - 2 \varepsilon < \alpha_6 + \alpha_7,
\end{equation}
making a contradiction. Hence (129) holds for $f(n)$ in \textbf{Case 1}.

\textbf{Case 2: $\Omega(n) = 6$.}
Suppose that $n = p_1 \cdots p_6$ and $\alpha_1 \geqslant \alpha_2 \geqslant \cdots \geqslant \alpha_6$. Now we have $\alpha_5 + \alpha_6 \leqslant \frac{1}{3}$. Since we have $\theta_2 \leqslant \frac{1}{4} < \frac{1}{3}$ and $\frac{1}{2}(2 - \theta_1 - 4 \theta_2) > \frac{9}{25} > \frac{1}{3}$, if $\alpha_5 + \alpha_6 \in \left[\theta_2 + \varepsilon,\ \frac{9}{25} - \varepsilon \right]$, then (129) holds for $f(n)$. Otherwise we have $\alpha_5 + \alpha_6 < \theta_2 + \varepsilon$. But since we have $8 \theta_1 + 11 \theta_2 < 5 - 8 \varepsilon$, we get
\begin{equation}
\alpha_5 + \alpha_6 < \theta_2 + \varepsilon < \frac{5 - 8 \theta}{3} - 2 \varepsilon < \alpha_5 + \alpha_6,
\end{equation}
making a contradiction. Hence (129) holds for $f(n)$ in \textbf{Case 2}.

\textbf{Case 3: $\Omega(n) = 4 \text{ or } 5$, no product of variables lies in $\left(1 - \theta - \varepsilon,\ \theta + \varepsilon \right)$.}
Suppose that $n = p_1 \cdots p_k$ and $\alpha_1 > \alpha_2 > \cdots > \alpha_k$, while the remaining part $\alpha_i = \alpha_j$ gives a ``loss'' of size $O(\varepsilon)$. In this case we can ``view'' $\left[\frac{3}{7} + \varepsilon,\ \frac{4}{7} - \varepsilon \right]$ as a ``fake'' Type-II range. Since we have $\theta < \frac{11}{21} < \frac{9}{17}$, the results in [\cite{MaynardLargeModuliI}, Chapter 9] can be applied here if any of the following conditions holds:

(1). $\alpha_k \geqslant \frac{1}{7} + \varepsilon$;

(2). $\Omega(n) = 5$, $\alpha_3 + \alpha_4 + \alpha_5 > \frac{4}{7} - \varepsilon$;

(3). $\Omega(n) = 4$, $\alpha_1 \leqslant \frac{3}{7} + \varepsilon$, $\alpha_1 + \alpha_4 > \frac{4}{7} - \varepsilon$.

By Condition (1), we can assume that $\frac{1}{8} < \alpha_k < \frac{1}{7} + \varepsilon$. We first consider the subcase $\Omega(n) = 5$. By Condition (2) and the ``fake'' Type-II range $\left[\frac{3}{7} + \varepsilon,\ \frac{4}{7} - \varepsilon \right]$, we can also assume that $\alpha_3 + \alpha_4 + \alpha_5 < \frac{3}{7} + \varepsilon$. This means that $\alpha_1 + \alpha_2 > \frac{4}{7} - \varepsilon$, and thus $\alpha_1 > \frac{2}{7} - \frac{1}{2} \varepsilon$. Since $\theta_2 < \frac{1}{4} < \frac{2}{7} - \frac{1}{2} \varepsilon < \frac{9}{25} < \frac{1}{2}(2 - \theta_1 - 4 \theta_2)$, we can assume that $\alpha_1 > \frac{9}{25} - \varepsilon$. Now we have $\frac{3}{7} + \varepsilon < \frac{97}{200} - \varepsilon = \frac{9}{25} - \varepsilon + \frac{1}{8} < \alpha_1 + \alpha_5$, and we can assume that $\alpha_1 + \alpha_5 > \frac{4}{7} - \varepsilon$. Since $\alpha_5 < \frac{1}{7} + \varepsilon$, we have $\alpha_1 > \frac{3}{7} - 2 \varepsilon$. If $\alpha_1 \geqslant \frac{3}{7} + \varepsilon$, we can assume that $\alpha_1 > \frac{4}{7} - \varepsilon$. But now we have $1 = \alpha_1 + \alpha_2 + \alpha_3 + \alpha_4 + \alpha_5 > \frac{4}{7} - \varepsilon + 4 \cdot \frac{1}{8} > 1$, making a contradiction. The remaining part is $\frac{3}{7} - 2 \varepsilon < \alpha_1 < \frac{3}{7} + \varepsilon$, leading to a ``loss'' of size $O(\varepsilon)$. Hence (129) holds for $f(n)$ with an error $O(\varepsilon)$ in \textbf{Case 3} if $\Omega(n) = 5$.

Next, we consider the second subcase $\Omega(n) = 4$. Suppose that $\alpha_1 + \alpha_4 \leqslant \frac{4}{7} - \varepsilon$. By Condition (3) and the ``fake'' Type-II range $\left[\frac{3}{7} + \varepsilon,\ \frac{4}{7} - \varepsilon \right]$, we can assume that $\alpha_1 + \alpha_4 < \frac{3}{7} + \varepsilon$. Since $\alpha_4 > \frac{1}{8}$, we have $\alpha_1 < \frac{17}{56} + \varepsilon = \frac{3}{7} + \varepsilon - \frac{1}{8}$. Since $\theta_2 < \frac{1}{4} < \frac{17}{56} + \varepsilon < \frac{9}{25} < \frac{1}{2}(2 - \theta_1 - 4 \theta_2)$, we can assume that $\alpha_1 \leqslant \frac{1}{4} + \varepsilon$. If $\alpha_1 \leqslant \frac{1}{4}$, then we have $1 = \alpha_1 + \alpha_2 + \alpha_3 + \alpha_4 < 4 \cdot \frac{1}{4} = 1$, making a contradiction. The remaining part is $\frac{1}{4} < \alpha_1 \leqslant \frac{1}{4} + \varepsilon$, leading to a ``loss'' of size $O(\varepsilon)$.

Now, by the ``fake'' Type-II range $\left[\frac{3}{7} + \varepsilon,\ \frac{4}{7} - \varepsilon \right]$, we only need to consider the subcase $\Omega(n) = 4$ with $\alpha_1 > \frac{4}{7} - \varepsilon$. This means that $\alpha_2 + \alpha_3 + \alpha_4 < \frac{3}{7} + \varepsilon$. Since $\alpha_2 > \alpha_3 > \alpha_4 > \frac{1}{8} + \varepsilon$, we have $\alpha_3 + \alpha_4 > \frac{1}{4} + \varepsilon > \theta_2 + \varepsilon$. Thus, we can assume that $\alpha_3 + \alpha_4 > \frac{9}{25} - \varepsilon$. But now we have
\begin{equation}
\frac{1}{2} < \frac{3}{2} \left(\frac{9}{25} - \varepsilon \right) < \frac{3}{2}(\alpha_3 + \alpha_4) < \alpha_2 + \alpha_3 + \alpha_4 < \frac{3}{7} + \varepsilon,
\end{equation}
making a contradiction. Hence (129) holds for $f(n)$ with an error $O(\varepsilon)$ in \textbf{Case 3} if $\Omega(n) = 4$.

\textbf{Case 4: $\Omega(n) = 5$, one variable lies in $\left(1 - \theta - \varepsilon,\ \theta + \varepsilon \right)$.}
Suppose that $n = p_1 \cdots p_5$ and $\alpha_1 \geqslant \alpha_2 \geqslant \cdots \geqslant \alpha_5$. Now we have $\alpha_1 > 1 - \theta - \varepsilon$ and $\alpha_2 \geqslant \cdots \geqslant \alpha_5 \geqslant \frac{5 - 8 \theta}{6} - \varepsilon$. But since $\theta < \frac{11}{21} < \frac{10}{19} - \varepsilon$, we have
\begin{equation}
1 - \theta - \varepsilon + 4 \left(\frac{5 - 8 \theta}{6} - \varepsilon \right) > 1,
\end{equation}
making a contradiction. Hence \textbf{Case 4} is empty.

\textbf{Case 5: $\Omega(n) = 5$, a product of two variables lies in $\left(1 - \theta - \varepsilon,\ \theta + \varepsilon \right)$.}
Suppose that $n = p_1 \cdots p_5$ and $\alpha_1 + \alpha_2 \in \left(1 - \theta - \varepsilon,\ \theta + \varepsilon \right)$ (of course, $\alpha_3 + \alpha_4 + \alpha_5 \in \left(1 - \theta - \varepsilon,\ \theta + \varepsilon \right)$ too). Since $\alpha_i + \alpha_j \geqslant \frac{5 - 8 \theta}{3} - 2 \varepsilon > \theta_2 + \varepsilon$ (see \textbf{Case 1}), we can assume that $\alpha_i + \alpha_j > \frac{1}{2}(2 - \theta_1 - 4 \theta_2) - \varepsilon$ for all $i,j \in \left\{1,2,3,4,5\right\}$ (otherwise $\alpha_i + \alpha_j \in \left[\theta_2 + \varepsilon,\ \frac{1}{2}(2 - \theta_1 - 4 \theta_2) - \varepsilon \right]$, and (129) holds for $f(n)$). Now we have 
\begin{equation}
\alpha_3 + \alpha_4 + \alpha_5 = \frac{(\alpha_3 + \alpha_4) + (\alpha_3 + \alpha_5) + (\alpha_4 + \alpha_5)}{2} > \frac{3}{4}(2 - \theta_1 - 4 \theta_2) - \varepsilon.
\end{equation}
Since $7 \theta_1 + 16 \theta_2 < 6 - 8 \varepsilon$, we have
\begin{equation}
\frac{3}{4}(2 - \theta_1 - 4 \theta_2) - \varepsilon > \theta + \varepsilon,
\end{equation}
making a contradiction with $\alpha_3 + \alpha_4 + \alpha_5 \in \left(1 - \theta - \varepsilon,\ \theta + \varepsilon \right)$. Hence (129) holds for $f(n)$ in \textbf{Case 5}.

\textbf{Case 6: $\Omega(n) = 4$, a product of two variables lies in $\left(1 - \theta - \varepsilon,\ \theta + \varepsilon \right)$.}
Suppose that $n = p_1 p_2 p_3 p_4$ and $\alpha_1 + \alpha_2 \in \left(1 - \theta - \varepsilon,\ \theta + \varepsilon \right)$ (of course, $\alpha_3 + \alpha_4 \in \left(1 - \theta - \varepsilon,\ \theta + \varepsilon \right)$ too). Without loss of generality, we further assume that $\max(\alpha_1, \alpha_2, \alpha_3, \alpha_4) = \alpha_1$, while the remaining part $\alpha_i = \alpha_j$ gives a ``loss'' of size $O(\varepsilon)$. Since $\alpha_1 > \frac{1}{4} \geqslant \theta_2 + \varepsilon$, we can assume that $\alpha_1 > \frac{1}{2}(2 - \theta_1 - 4 \theta_2) - \varepsilon$ (otherwise $\alpha_1 \in \left[\theta_2 + \varepsilon,\ \frac{1}{2}(2 - \theta_1 - 4 \theta_2) - \varepsilon \right]$, and (129) holds for $f(n)$). Now we have
\begin{equation}
\alpha_1 + \alpha_2 > \frac{2 - \theta_1 - 4 \theta_2}{2} + \frac{5 - 8 \theta}{6} - 2 \varepsilon = \frac{11 - 11 \theta_1 - 20 \theta_2}{6} - 2 \varepsilon.
\end{equation}
Since $17 \theta_1 + 26 \theta_2 < 11 - 14 \varepsilon$, we have
\begin{equation}
\frac{11 - 11 \theta_1 - 20 \theta_2}{6} - 2 \varepsilon > \theta + \varepsilon,
\end{equation}
making a contradiction with the assumption $\alpha_1 + \alpha_2 \in \left[1 - \theta,\ \theta \right]$. Hence (129) holds for $f(n)$ with an error $O(\varepsilon)$ in \textbf{Case 6}.

\textbf{Case 7: $\Omega(n) = 4$, one variable lies in $\left(1 - \theta - \varepsilon,\ \theta + \varepsilon \right)$.}
Suppose that $n = p_0 p_1 p_2 p_3$ and $\alpha_0 \geqslant \alpha_1 \geqslant \alpha_2 \geqslant \alpha_3$. Now we have $\alpha_0 \in \left(1 - \theta - \varepsilon,\ \theta + \varepsilon \right)$ and $\alpha_1 + \alpha_2 + \alpha_3 \in \left(1 - \theta - \varepsilon,\ \theta + \varepsilon \right)$. We also have $\alpha_1 + \alpha_3 \leqslant \frac{1}{2}$ and $\alpha_2 \leqslant \frac{1}{3}$. Note that $\Sigma_{42309}$ only counts numbers with all prime factors smaller than $\max(m_1, m_2, m_3) \leqslant x^{\frac{3}{7} + \varepsilon}$, one can easily find that this type of $n$ will not be counted in $\Sigma_{42309}$ by a simple observation. Hence, this type of $n$ will only be counted in one part of $\Sigma_{42305}$:
\begin{equation}
\sum_{\substack{n = p_1 p_2 p_3 \beta \\ 1 - \theta - \varepsilon < \alpha_1 + \alpha_2 + \alpha_3 < \theta + \varepsilon }} \psi\left(\beta, p_3 \right),
\end{equation}
where $\beta$ in (167) is a large prime. Using Buchstab's identity, we have
\begin{equation}
\sum_{\substack{n = p_1 p_2 p_3 \beta \\ 1 - \theta - \varepsilon < \alpha_1 + \alpha_2 + \alpha_3 < \theta + \varepsilon}} \psi\left(\beta, p_3 \right) = \sum_{\substack{n = p_1 p_2 p_3 \beta \\ 1 - \theta - \varepsilon < \alpha_1 + \alpha_2 + \alpha_3 < \theta + \varepsilon }} \psi\left(\beta, x^{\kappa} \right) - \sum_{\substack{n = p_1 p_2 p_3 p_4 \beta \\ 1 - \theta - \varepsilon < \alpha_1 + \alpha_2 + \alpha_3 < \theta + \varepsilon \\ \kappa \leqslant \alpha_4 < \min\left(\alpha_3, \frac{1}{2}(1 - \alpha_1 - \alpha_2 - \alpha_3) \right) }} \psi\left(\beta, p_4 \right).
\end{equation}
By Lemma~\ref{l217}, (129) holds for the first sum in (168) since $(\alpha_1 + \alpha_2, \alpha_3) \in A$, and we only need to deal with the second sum in (168) that counts numbers with $5$ or more prime factors, with a product of $3$ variables lies in $\left(1 - \theta - \varepsilon,\ \theta + \varepsilon \right)$. Now we reduce this case to \textbf{Cases 1, 2 and 5}. Since the second sum in (168) is a positive sum that cannot be dropped directly, we need (129) holds for $f(n)$ in \textbf{Cases 1, 2 and 5} without an error $O(\varepsilon)$. By the discussions above, we know that (129) holds for $f(n)$ in \textbf{Cases 1, 2 and 5} without an error $O(\varepsilon)$, and thus (129) holds for $f(n)$ in \textbf{Case 7}.

Combining all the 7 cases above, the proof of Theorem~\ref{case2} is completed.
\end{proof}
\begin{remark*}
When $(\theta_1, \theta_2)$ lies in the boundary of this region:
$$
\frac{1}{4} \leqslant \theta_1 < \frac{7}{25}, \quad \theta_2 = \min\left(\frac{1}{4}, \frac{11 - 17 \theta_1}{26}\right),
$$
we can use the exactly same decomposing process to prove that
$$
C_1^{\prime}(\theta_1, \theta_2) \leqslant 1 + \int_{1 - \theta}^{\frac{1}{2}} \frac{\omega\left(\frac{1-t}{t} \right)}{t^2} d t + O(\varepsilon) = 1 + \log\left(\frac{\theta}{1 - \theta} \right) + O(\varepsilon).
$$
For example, we have $C_1^{\prime}(0.26, 0.25) \leqslant 1 + \log\left(\frac{51}{49} \right) + O(\varepsilon) < 1.0401$. In order to prove this, we only need to focus on the conditions needed in the proof of Theorem~\ref{case2}. A simple verification yields $(\theta_1, \theta_2) \in \boldsymbol{T}$, and we can use the new three-dimensional Harman's sieve to handle $\Sigma_{42306}$. By the discussions above, we know that (129) holds for $f(n) = \Sigma_{42301}$, $f(n) = \Sigma_{42302}$, $f(n) = \Sigma_{42303}$, $f(n) = \Sigma_{42304}$ and $f(n) = \Sigma_{42307}$, and we only need to focus on the sums with $\Omega(n) \geqslant 4$. Since $8 \theta_1 + 11 \theta_2 < 5 - 8 \varepsilon$ and $7 \theta_1 + 16 \theta_2 < 6 - 8 \varepsilon$ still hold, we can use our Type-II range (156) to show that (129) holds for $f(n)$ in \textbf{Cases 1, 2, 4, 5 and 7} and for $f(n)$ in \textbf{Cases 3 and 6} with an error $O(\varepsilon)$.
\end{remark*}

For the remaining parts of $\boldsymbol{\mathcal{Z}_2}$, the decompositions are similar to which in the non-asymptotic parts in the case $\boldsymbol{\mathcal{Z}_1}$. We also use the Type-II range $\left[\theta_2 + \varepsilon,\ \frac{1}{2}(2 - \theta_1 - 4 \theta_2) - \varepsilon \right]$ to discard sums with products of variables lie in this range.

The Type-II range for $\boldsymbol{E}_{0204}$ is
\begin{equation}
\left[\varepsilon,\ \frac{1}{4}(2 - 2 \theta_1 - 5 \theta_2) - \varepsilon \right] \cup \left[2 \theta_1 + 2 \theta_2 - 1 + \varepsilon,\ \frac{1}{6}(5 - 8 \theta_1 - 8 \theta_2) - \varepsilon \right] \cup \left[\theta_2 + \varepsilon,\ \frac{1}{2}(2 - \theta_1 - 4 \theta_2) - \varepsilon \right].
\end{equation}
The decompositions are similar to which in the case $\boldsymbol{E}_{0103}$.

The Type-II range for $\boldsymbol{E}_{0205}$ is
\begin{equation}
\left[\varepsilon,\ \frac{1}{4}(2 - 2 \theta_1 - 5 \theta_2) - \varepsilon \right] \cup \left[\theta_2 + \varepsilon,\ \frac{1}{2}(2 - \theta_1 - 4 \theta_2) - \varepsilon \right].
\end{equation}
The decompositions are similar to which in the case $\boldsymbol{E}_{0104}$.

The Type-II range for $\boldsymbol{E}_{0206}$ is
\begin{equation}
\left[\varepsilon,\ \frac{1}{6}(5 - 8 \theta_1 - 8 \theta_2) - \varepsilon \right] \cup \left[\theta_2 + \varepsilon,\ 2 - 5 \theta_1 - \theta_2 - \varepsilon \right] \cup \left[\theta_1 + \varepsilon,\ \frac{1}{2}(1 - \theta_2) - \varepsilon \right].
\end{equation}
The decompositions are similar to which in the non-asymptotic parts in the case $\boldsymbol{\mathcal{Z}_1}$.

The Type-II range for $\boldsymbol{E}_{0207}$ is
\begin{equation}
\left[\varepsilon,\ \frac{1}{4}(2 - 2 \theta_1 - 5 \theta_2) - \varepsilon \right] \cup \left[2 \theta - 1 + \varepsilon,\ \frac{1}{6}(5 - 8 \theta) - \varepsilon \right] \cup \left[\theta_2 + \varepsilon,\ 2 - 5 \theta_1 - \theta_2 - \varepsilon \right] \cup \left[\theta_1 + \varepsilon,\ \frac{1}{2}(1 - \theta_2) - \varepsilon \right].
\end{equation}
The decompositions are similar to which in the case $\boldsymbol{E}_{0103}$.

The Type-II range for $\boldsymbol{E}_{0208}$ is
\begin{equation}
\left[\varepsilon,\ \frac{1}{4}(2 - 2 \theta_1 - 5 \theta_2) - \varepsilon \right] \cup \left[2 \theta - 1 + \varepsilon,\ \frac{1}{6}(5 - 8 \theta) - \varepsilon \right] \cup \left[\theta_2 + \varepsilon,\ 2 - 5 \theta_1 - \theta_2 - \varepsilon \right] \cup \left[\theta_1 + \varepsilon,\ \frac{1}{2}(2 - \theta_1 - 4 \theta_2) - \varepsilon \right].
\end{equation}
The decompositions are similar to which in the case $\boldsymbol{E}_{0103}$.

The Type-II range for $\boldsymbol{E}_{0209}$ is
\begin{equation}
\left[\varepsilon,\ \frac{1}{4}(2 - 2 \theta_1 - 5 \theta_2) - \varepsilon \right] \cup \left[\theta_2 + \varepsilon,\ 2 - 5 \theta_1 - \theta_2 - \varepsilon \right] \cup \left[\theta_1 + \varepsilon,\ \frac{1}{2}(2 - \theta_1 - 4 \theta_2) - \varepsilon \right].
\end{equation}
The decompositions are similar to which in the case $\boldsymbol{E}_{0104}$.

The Type-II range for $\boldsymbol{E}_{0210}$ is
\begin{equation}
\left[\varepsilon,\ \frac{1}{6}(5 - 8 \theta_1 - 8 \theta_2) - \varepsilon \right] \cup \left[\theta_2 + \varepsilon,\ \frac{1}{6}(5 - 8 \theta_1 - 4 \theta_2) - \varepsilon \right] \cup \left[\theta_1 + \varepsilon,\ \frac{1}{2}(1 - \theta_2) - \varepsilon \right].
\end{equation}
The decompositions are similar to which in the non-asymptotic parts in the case $\boldsymbol{\mathcal{Z}_1}$. We shall discuss this case in detail since it covers most parts of the region
\begin{equation}
\left\{ (\theta_1, \theta_2) : \theta_1 \leqslant \frac{1}{3} - \varepsilon,\ \theta_2 \leqslant \frac{1}{5},\ \theta \geqslant \frac{1}{2} \right\},
\end{equation}
where Baker and Harman [\cite{677}, Theorem 3] proved that (129) holds for $f(n) = \mathbbm{1}_{p}(n)$ in region (176) with $\theta < \frac{29}{56}$. The case $\boldsymbol{E}_{0210}$ covers the region (162) with $\theta < \frac{11}{21}$ (and a small part with $\theta \geqslant \frac{11}{21}$). From here to the end of the case $\boldsymbol{E}_{0210}$, we suppose that $\theta_1 \leqslant \frac{1}{3} - \varepsilon$, $\theta_2 \leqslant \frac{1}{5}$ and $\frac{1}{2} \leqslant \theta \leqslant \frac{11}{21} - \varepsilon$. Now we have
$$
\kappa = \frac{5 - 8 \theta}{6} - \varepsilon > \frac{17}{126} - \varepsilon > \frac{1}{8} + \varepsilon,
$$
$$
\frac{5 - 8 \theta_1 - 4 \theta_2}{6} - \varepsilon > \frac{5 - 8 \left(\frac{1}{3} - \varepsilon \right) - 4 \left(\frac{11}{21} - \varepsilon - \frac{1}{3} + \varepsilon \right)}{6} - \varepsilon = \frac{11}{42} + \frac{1}{3} \varepsilon > \frac{11}{42} > \frac{1}{4} + 2 \varepsilon
$$
and
$$
\frac{1 - \theta_2}{2} \geqslant \frac{1 - \frac{1}{5}}{2} = 0.4,
$$
and we can simplify our Type-II range in this case to
\begin{equation}
\left[\varepsilon,\ \frac{17}{126} - \varepsilon \right] \cup \left[\frac{1}{5} + \varepsilon,\ \frac{11}{42} \right] \cup \left[\frac{1}{3},\ 0.4 - \varepsilon \right] \cup \left[\frac{3}{7} + \varepsilon,\ 1 - \theta - \varepsilon \right].
\end{equation}
We also have
$$
11 \theta_1 + 12 \theta_2 < 6 - \varepsilon \quad \text{and} \quad 3 \theta_1 + 2 \theta_2 < \frac{11}{7} - \varepsilon,
$$
which are required when using Lemma~\ref{l38} together with a three-dimensional Harman's sieve.

Now, we can decompose our $\mathbbm{1}_{n \sim x, n = p}(n) = \psi\left(n, (2 x)^{\frac{1}{2}}\right)$ in a way similar to the decompositions in the proof of Theorem~\ref{case1}. By Buchstab's identity, we have
\begin{align}
\nonumber \psi\left(n, (2 x)^{\frac{1}{2}}\right) =&\ \psi\left(n, x^{\kappa}\right) - \sum_{\substack{n = p_1 \beta \\ \kappa \leqslant \alpha_1 < \frac{1}{2} }} \psi\left(\beta, p_1 \right) \\
\nonumber =&\ \psi\left(n, x^{\kappa}\right) - \sum_{\substack{n = p_1 \beta \\ \kappa \leqslant \alpha_1 < \frac{3}{7} + \varepsilon }} \psi\left(\beta, p_1 \right) - \sum_{\substack{n = p_1 \beta \\ \frac{3}{7} + \varepsilon \leqslant \alpha_1 \leqslant 1 - \theta - \varepsilon }} \psi\left(\beta, p_1 \right) - \sum_{\substack{n = p_1 \beta \\ 1 - \theta - \varepsilon < \alpha_1 < \frac{1}{2} }} \psi\left(\beta, p_1 \right) \\
\nonumber =&\ \psi\left(n, x^{\kappa}\right) - \sum_{\substack{n = p_1 \beta \\ \kappa \leqslant \alpha_1 < \frac{3}{7} + \varepsilon }} \psi\left(\beta, x^{\kappa} \right) + \sum_{\substack{n = p_1 p_2 \beta \\ \kappa \leqslant \alpha_1 < \frac{3}{7} + \varepsilon \\ \kappa \leqslant \alpha_2 < \min\left(\alpha_1, \frac{1}{2}(1 - \alpha_1) \right) \\ \boldsymbol{\alpha}_2 \in \boldsymbol{G}_2 }} \psi\left(\beta, p_2 \right) + \sum_{\substack{n = p_1 p_2 \beta \\ \kappa \leqslant \alpha_1 < \frac{3}{7} + \varepsilon \\ \kappa \leqslant \alpha_2 < \min\left(\alpha_1, \frac{1}{2}(1 - \alpha_1) \right) \\ \boldsymbol{\alpha}_2 \in A \cup B }} \psi\left(\beta, p_2 \right) \\
\nonumber & + \sum_{\substack{n = p_1 p_2 \beta \\ \kappa \leqslant \alpha_1 < \frac{3}{7} + \varepsilon \\ \kappa \leqslant \alpha_2 < \min\left(\alpha_1, \frac{1}{2}(1 - \alpha_1) \right) \\ \boldsymbol{\alpha}_2 \in C }} \psi\left(\beta, p_2 \right) - \sum_{\substack{n = p_1 \beta \\ \frac{3}{7} + \varepsilon \leqslant \alpha_1 \leqslant 1 - \theta - \varepsilon }} \psi\left(\beta, p_1 \right) - \sum_{\substack{n = p_1 \beta \\ 1 - \theta - \varepsilon < \alpha_1 < \frac{1}{2} }} \psi\left(\beta, p_1 \right) \\
\nonumber =&\ \psi\left(n, x^{\kappa}\right) - \sum_{\substack{n = p_1 \beta \\ \kappa \leqslant \alpha_1 < \frac{3}{7} + \varepsilon }} \psi\left(\beta, x^{\kappa} \right) + \sum_{\substack{n = p_1 p_2 \beta \\ \kappa \leqslant \alpha_1 < \frac{3}{7} + \varepsilon \\ \kappa \leqslant \alpha_2 < \min\left(\alpha_1, \frac{1}{2}(1 - \alpha_1) \right) \\ \boldsymbol{\alpha}_2 \in \boldsymbol{G}_2 }} \psi\left(\beta, p_2 \right) \\
\nonumber & + \sum_{\substack{n = p_1 p_2 \beta \\ \kappa \leqslant \alpha_1 < \frac{3}{7} + \varepsilon \\ \kappa \leqslant \alpha_2 < \min\left(\alpha_1, \frac{1}{2}(1 - \alpha_1) \right) \\ \boldsymbol{\alpha}_2 \in A \cup B }} \psi\left(\beta, x^{\kappa} \right) - \sum_{\substack{n = p_1 p_2 p_3 \beta \\ \kappa \leqslant \alpha_1 < \frac{3}{7} + \varepsilon \\ \kappa \leqslant \alpha_2 < \min\left(\alpha_1, \frac{1}{2}(1 - \alpha_1) \right) \\ \boldsymbol{\alpha}_2 \in A \cup B \cup C \\ \kappa \leqslant \alpha_3 < \min\left(\alpha_2, \frac{1}{2}(1 - \alpha_1 - \alpha_2) \right) }} \psi\left(\beta, p_3 \right) \\
\nonumber & + \sum_{\substack{n = p_1 p_2 \beta \\ \kappa \leqslant \alpha_1 < \frac{3}{7} + \varepsilon \\ \kappa \leqslant \alpha_2 < \min\left(\alpha_1, \frac{1}{2}(1 - \alpha_1) \right) \\ \boldsymbol{\alpha}_2 \in C }} \psi\left(\beta, x^{\kappa} \right) - \sum_{\substack{n = p_1 \beta \\ \frac{3}{7} + \varepsilon \leqslant \alpha_1 \leqslant 1 - \theta - \varepsilon }} \psi\left(\beta, p_1 \right) - \sum_{\substack{n = p_1 \beta \\ 1 - \theta - \varepsilon < \alpha_1 < \frac{1}{2} }} \psi\left(\beta, p_1 \right) \\
=&\ \Sigma_{42401} - \Sigma_{42402} + \Sigma_{42403} + \Sigma_{42404} - \Sigma_{42405} + \Sigma_{42406} - \Sigma_{42407} - \Sigma_{42408}.
\end{align}
By Lemma~\ref{l212}, (129) holds for $f(n) = \Sigma_{42401}$ and $f(n) = \Sigma_{42402}$. By Lemma~\ref{l211} and Lemma~\ref{l216}, (129) holds for $f(n) = \Sigma_{42403}$ and $f(n) = \Sigma_{42407}$. By Lemma~\ref{l217}, (129) holds for $f(n) = \Sigma_{42404}$. We discard the whole of $\Sigma_{42408}$. For the remaining sums, $\Sigma_{42405}$ only counts numbers with $4$ or more prime factors.

For $\Sigma_{42406}$, since we have $3 \theta_1 + 2 \theta_2 < \frac{11}{7} - \varepsilon$, $11 \theta_1 + 12 \theta_2 < 6 - \varepsilon$, $\theta < \frac{11}{21} < \frac{8}{15} - \varepsilon$ and a Type-II range $\left[\varepsilon,\ \frac{5 - 8 \theta}{6} - \varepsilon \right]$, we can use Lemma~\ref{l38} and a three-dimensional Harman's sieve to get a ``loss term''
\begin{equation}
\Sigma_{42409} = \sum_{\substack{n = m_1 m_2 m_3 \\ \kappa \leqslant \alpha_1 < \frac{3}{7} + \varepsilon \\ \kappa \leqslant \alpha_2 < \min\left(\alpha_1, \frac{1}{2}(1 - \alpha_1) \right) \\ \boldsymbol{\alpha}_2 \in C \\ \Omega(m_1 m_2) \geqslant 3 }} \psi\left(m_1 m_2 m_3, x^{\kappa} \right).
\end{equation}
Since $\Omega(m_1 m_2 m_3) \geqslant \Omega(m_1 m_2) + 1 \geqslant 4$, $\Sigma_{42409}$ only counts numbers with $4$ or more prime factors.

Now, we want to show that (129) holds for $f(n) = $ sums that count numbers with $4$ or more prime factors. Since $\kappa > \frac{1}{8}$, we have $\Omega(n) \leqslant 7$. Assume that $\Omega(n) \geqslant 4$, and we only need to consider 8 different cases as in the proof of Theorem~\ref{case1} and Theorem~\ref{case2}. We shall prove that (129) holds for $f(n)$ in 7 of the following 8 cases (sometimes with an error $O(\varepsilon)$), and the ``loss integrals'' for $C_1^{\prime}(\theta_1, \theta_2)$ will only come from that remaining case and the discarded $\Sigma_{42408}$.

\textbf{Case 1: $\Omega(n) = 7$.}
Suppose that $n = p_1 \cdots p_7$ and $\alpha_1 > \alpha_2 > \cdots > \alpha_7$, while the remaining part $\alpha_i = \alpha_j$ gives a ``loss'' of size $O(\varepsilon)$. Now we have $\alpha_4 + \alpha_5 + \alpha_6 + \alpha_7 < \frac{4}{7}$. If $\alpha_4 + \alpha_5 + \alpha_6 + \alpha_7 \in \left[\theta + \varepsilon,\ \frac{4}{7} - \varepsilon \right]$, then (129) holds for $f(n)$. Assume that $\alpha_4 + \alpha_5 + \alpha_6 + \alpha_7 < \theta + \varepsilon$. But since we have $\theta < \frac{11}{21} < \frac{10}{19} - \varepsilon$, we get
\begin{equation}
\alpha_4 + \alpha_5 + \alpha_6 + \alpha_7 < \theta + \varepsilon < \frac{10 - 16 \theta}{3} - 4 \varepsilon = 4 \left(\frac{5 - 8 \theta}{6} - \varepsilon \right) < \alpha_4 + \alpha_5 + \alpha_6 + \alpha_7,
\end{equation}
making a contradiction. Now we can assume that $\alpha_4 + \alpha_5 + \alpha_6 + \alpha_7 \in \left(\frac{4}{7} - \varepsilon,\ \frac{4}{7} \right)$, and thus $\alpha_1 + \alpha_2 + \alpha_3 \in \left(\frac{3}{7},\ \frac{3}{7} + \varepsilon \right)$. Now we have $\alpha_3 \leqslant \frac{1}{7} + \frac{1}{3} \varepsilon$ and $\alpha_4 \geqslant \frac{1}{7} - \frac{1}{4} \varepsilon$. Since $\alpha_3 > \alpha_4$, we have $\frac{1}{7} - \frac{1}{4} \varepsilon < \alpha_3 \leqslant \frac{1}{7} + \frac{1}{3} \varepsilon$, leading to a ``loss'' of size $O(\varepsilon)$. Hence (129) holds for $f(n)$ with an error $O(\varepsilon)$ in \textbf{Case 1}.

\textbf{Case 2: $\Omega(n) = 6$.}
Suppose that $n = p_1 \cdots p_6$ and $\alpha_1 > \alpha_2 > \cdots > \alpha_6$, while the remaining part $\alpha_i = \alpha_j$ gives a ``loss'' of size $O(\varepsilon)$. Now we have $\alpha_1 + \alpha_2 > \frac{1}{3}$. If $\alpha_1 + \alpha_2 \in \left[\frac{1}{3},\ 0.4 - \varepsilon \right]$, then (129) holds for $f(n)$. Otherwise we have $\alpha_1 + \alpha_2 > 0.4 - \varepsilon$. If $\alpha_1 + \alpha_2 \in \left[\frac{3}{7} + \varepsilon,\ 1 - \theta - \varepsilon \right]$, then (129) holds for $f(n)$. If $\alpha_1 + \alpha_2 > 1 - \theta - \varepsilon$, we have
\begin{equation}
\alpha_1 + \cdots + \alpha_6 > 1 - \theta - \varepsilon + 4 \left(\frac{5 - 8 \theta}{6} - \varepsilon \right) > 1,
\end{equation}
making a contradiction since $\theta < \frac{11}{21} < \frac{10}{19} - \varepsilon$ (see \textbf{Case 1}). Thus, we can assume that $\alpha_1 + \alpha_2 \in \left(0.4 - \varepsilon,\ \frac{3}{7} + \varepsilon \right)$.

Suppose that $\alpha_6 \geqslant \frac{1}{7} + \varepsilon$. Since $\alpha_1 > \alpha_2 > \cdots > \alpha_6 \geqslant \frac{1}{7} + \varepsilon$, we have $\alpha_1 + \alpha_2 + \alpha_3 > \frac{3}{7} + 3 \varepsilon$ and $\alpha_4 + \alpha_5 + \alpha_6 > \frac{3}{7} + 3 \varepsilon$. Hence $\alpha_1 + \alpha_2 + \alpha_3 \in \left[\frac{3}{7} + \varepsilon,\ \frac{4}{7} - \varepsilon \right]$. If $\alpha_1 + \alpha_2 + \alpha_3 \in \left[\frac{3}{7} + \varepsilon,\ 1 - \theta - \varepsilon \right] \cup \left[\theta + \varepsilon,\ \frac{4}{7} - \varepsilon \right]$, then (129) holds for $f(n)$. Otherwise we have $\alpha_1 + \alpha_2 + \alpha_3 \in \left(1 - \theta - \varepsilon,\ \theta + \varepsilon \right)$. But since $\alpha_1 + \alpha_2 > 0.4 - \varepsilon$, we have
\begin{equation}
\alpha_1 + \alpha_2 + \alpha_3 > 0.4 - \varepsilon + \frac{1}{7} > 0.54 > \theta + \varepsilon,
\end{equation}
making a contradiction.

Suppose that $\alpha_6 < \frac{1}{7} - 2 \varepsilon$. Now we have $\alpha_1 + \alpha_2 \in \left(0.4 - \varepsilon,\ \frac{3}{7} + \varepsilon \right)$ and $\frac{1}{8} < \alpha_6 < \frac{1}{7} - 2 \varepsilon$. Thus,
\begin{equation}
\alpha_1 + \alpha_2 + \alpha_6 \in \left(0.525 - \varepsilon,\ \frac{4}{7} - \varepsilon \right).
\end{equation}
Since $\theta + \varepsilon < 0.525 - \varepsilon$, we have (129) holds for $f(n)$.

The remaining part is $\frac{1}{7} - 2 \varepsilon \leqslant \alpha_6 < \frac{1}{7} + \varepsilon$, leading to a ``loss'' of size $O(\varepsilon)$. Hence (129) holds for $f(n)$ with an error $O(\varepsilon)$ in \textbf{Case 2}.

\textbf{Case 3: $\Omega(n) = 5$, a product of two variables lies in $\left(1 - \theta - \varepsilon,\ \theta + \varepsilon \right)$.}
Suppose that $n = p_1 \cdots p_5$ and $\alpha_1 + \alpha_2 \in \left(1 - \theta - \varepsilon,\ \theta + \varepsilon \right)$ (of course, $\alpha_3 + \alpha_4 + \alpha_5 \in \left(1 - \theta - \varepsilon,\ \theta + \varepsilon \right)$ too). Without loss of generality, we further assume that $\alpha_1 \geqslant \alpha_2$ and $\alpha_3 \geqslant \alpha_4 \geqslant \alpha_5$. Now we have $\alpha_1 > \frac{1 - \theta - \varepsilon}{2}$ and $\alpha_2 < \frac{\theta + \varepsilon}{2}$. Since $\frac{1}{5} + \varepsilon < \frac{1 - \theta - \varepsilon}{2} < \frac{\theta + \varepsilon}{2} \leqslant \frac{11}{42}$, we can assume that $\alpha_1 > \frac{11}{42} - \varepsilon$ and $\alpha_2 < \frac{1}{5} + \varepsilon$. Now we have $\alpha_1 > (1 - \theta) - \left(\frac{1}{5} + \varepsilon \right) = \frac{4}{5} - \theta - \varepsilon > \frac{4}{5} - \frac{11}{21} - 2 \varepsilon > 0.276$. If $\alpha_1 \geqslant \frac{1}{3}$, we can assume that $\alpha_1 > 0.4 - \varepsilon$ and $\alpha_1 + \alpha_2 > 0.4 - \varepsilon + \frac{1}{8} = 0.525 - \varepsilon > \theta + \varepsilon$, making a contradiction with the assumption $\alpha_1 + \alpha_2 \in \left[1 - \theta,\ \theta \right]$. Hence we also assume that $\alpha_1 < \frac{1}{3}$.

Since $\alpha_3 + \alpha_4 + \alpha_5 \in \left(1 - \theta - \varepsilon,\ \theta + \varepsilon \right)$, we also have $\alpha_3 + \alpha_4 \geqslant \frac{2}{3}(1 - \theta - \varepsilon) > \frac{20}{63} - \varepsilon$ and $\alpha_3 \geqslant \frac{10}{63} - \varepsilon$. If $\alpha_3 + \alpha_4 \geqslant \frac{1}{3}$, then we can assume that $\alpha_3 + \alpha_4 > 0.4 - \varepsilon$. Thus we have $\alpha_3 + \alpha_4 + \alpha_5 > 0.4 - \varepsilon + \frac{1}{8} = 0.525 - \varepsilon > \theta + \varepsilon$, contradicting with the assumption $\alpha_3 + \alpha_4 + \alpha_5 \in \left(1 - \theta - \varepsilon,\ \theta + \varepsilon \right)$. Hence we can also assume that $\alpha_3 + \alpha_4 \in \left(\frac{20}{63} - \varepsilon,\ \frac{1}{3} \right)$.

Since $\alpha_1 > 0.276$ and $\alpha_3 \geqslant \frac{10}{63} - \varepsilon$, we have $\alpha_1 + \alpha_3 > 0.276 + \frac{10}{63} - \varepsilon > \frac{3}{7} + \varepsilon$, and we can assume that $\alpha_1 + \alpha_3 > 1 - \theta - \varepsilon$. If $\alpha_1 + \alpha_3 \geqslant \theta + \varepsilon$, then we can assume that $\alpha_1 + \alpha_3 > \frac{4}{7} - \varepsilon$. Since $\alpha_1 < \frac{1}{3}$, we have $\alpha_3 > \frac{4}{7} - \varepsilon - \frac{1}{3} = \frac{5}{21} - \varepsilon > \frac{1}{5} + \varepsilon$, and we can assume that $\alpha_3 > \frac{11}{42} - \varepsilon$. But now we have $\alpha_3 + \alpha_4 + \alpha_5 > \frac{11}{42} - \varepsilon + 2 \left( \frac{17}{126} - \varepsilon \right) > 0.53 > \theta + \varepsilon$, making a contradiction with the assumption $\alpha_3 + \alpha_4 + \alpha_5 \in \left(1 - \theta - \varepsilon,\ \theta + \varepsilon \right)$. Thus, we assume that $\alpha_1 + \alpha_3 \in \left(1 - \theta - \varepsilon,\ \theta + \varepsilon \right)$ and $\alpha_2 + \alpha_4 + \alpha_5 \in \left(1 - \theta - \varepsilon,\ \theta + \varepsilon \right)$.

When $\alpha_3 \geqslant \alpha_4 \geqslant \alpha_5 \geqslant \alpha_2$, $\alpha_3 \geqslant \alpha_4 \geqslant \alpha_2 \geqslant \alpha_5$ or $\alpha_3 \geqslant \alpha_2 \geqslant \alpha_4 \geqslant \alpha_5$ holds, we know that
\begin{equation}
\alpha_3 + \alpha_4 \geqslant \frac{1}{2}(1 - \alpha_1) > \frac{1}{3},
\end{equation}
and we can assume that $\alpha_3 + \alpha_4 > 0.4 - \varepsilon$. But now we have $\alpha_3 + \alpha_4 + \alpha_5 > 0.4 - \varepsilon + \frac{1}{8} = 0.525 - \varepsilon > \theta + \varepsilon$, contradicting with the assumption $\alpha_3 + \alpha_4 + \alpha_5 \in \left(1 - \theta - \varepsilon,\ \theta + \varepsilon \right)$.

When $\alpha_2 \geqslant \alpha_3 \geqslant \alpha_4 \geqslant \alpha_5$ holds, we know that
\begin{equation}
\alpha_2 + \alpha_4 \geqslant \frac{1}{2}(1 - \alpha_1) > \frac{1}{3},
\end{equation}
and we can assume that $\alpha_2 + \alpha_4 > 0.4 - \varepsilon$. But now we have $\alpha_2 + \alpha_4 + \alpha_5 > 0.4 - \varepsilon + \frac{1}{8} = 0.525 - \varepsilon > \theta + \varepsilon$, contradicting with the assumption $\alpha_2 + \alpha_4 + \alpha_5 \in \left(1 - \theta - \varepsilon,\ \theta + \varepsilon \right)$.

Hence (129) holds for $f(n)$ in \textbf{Case 3}.

\textbf{Case 4: $\Omega(n) = 5$, one variable lies in $\left(1 - \theta - \varepsilon,\ \theta + \varepsilon \right)$.}
Suppose that $n = p_1 \cdots p_5$ and $\alpha_1 \geqslant \alpha_2 \geqslant \cdots \geqslant \alpha_5$. Now we have $\alpha_1 > 1 - \theta - \varepsilon$ and $\alpha_2 \geqslant \cdots \geqslant \alpha_5 \geqslant \frac{5 - 8 \theta}{6} - \varepsilon$. But since $\theta < \frac{11}{21} < \frac{10}{19} - \varepsilon$, we have
\begin{equation}
1 - \theta - \varepsilon + 4 \left(\frac{5 - 8 \theta}{6} - \varepsilon \right) > 1,
\end{equation}
making a contradiction. Hence \textbf{Case 4} is empty.

\textbf{Case 5: $\Omega(n) = 5$, no product of variables lies in $\left(1 - \theta - \varepsilon,\ \theta + \varepsilon \right)$.}
Suppose that $n = p_1 \cdots p_5$ and $\alpha_1 > \alpha_2 > \cdots > \alpha_5$, while the remaining part $\alpha_i = \alpha_j$ gives a ``loss'' of size $O(\varepsilon)$. In this case we can ``view'' $\left[\frac{3}{7} + \varepsilon,\ \frac{4}{7} - \varepsilon \right]$ as a ``fake'' Type-II range. Since we have $\theta < \frac{11}{21} < \frac{9}{17}$, the results in [\cite{MaynardLargeModuliI}, Chapter 9] can be applied here if any of the following conditions holds:

(1). $\alpha_5 \geqslant \frac{1}{7} + \varepsilon$;

(2). $\alpha_3 + \alpha_4 + \alpha_5 > \frac{4}{7} - \varepsilon$.

By Condition (1), we can assume that $\frac{1}{8} < \alpha_5 < \frac{1}{7} + \varepsilon$. By Condition (2) and the ``fake'' Type-II range $\left[\frac{3}{7} + \varepsilon,\ \frac{4}{7} - \varepsilon \right]$, we can also assume that $\alpha_3 + \alpha_4 + \alpha_5 < \frac{3}{7} + \varepsilon$. This means that $\alpha_1 + \alpha_2 > \frac{4}{7} - \varepsilon$ and $\alpha_1 > \frac{2}{7} - \frac{1}{2} \varepsilon$.

If we have $\alpha_1 + \alpha_5 \geqslant \frac{3}{7} + \varepsilon$, then we can assume that $\alpha_1 + \alpha_5 > \frac{4}{7} - \varepsilon$. Since $\alpha_5 < \frac{1}{7} + \varepsilon$, we have $\alpha_1 > \frac{3}{7} - 2 \varepsilon$. If $\alpha_1 \geqslant \frac{3}{7} + \varepsilon$, then we can assume that $\alpha_1 > \frac{4}{7} - \varepsilon$. Now we have $1 = \alpha_1 + \alpha_2 + \alpha_3 + \alpha_4 + \alpha_5 > \frac{4}{7} - \varepsilon + 4 \cdot \frac{1}{8} > 1$, making a contradiction. The remaining part is $\frac{3}{7} - 2 \varepsilon < \alpha_1 < \frac{3}{7} + \varepsilon$, leading to a ``loss'' of size $O(\varepsilon)$.

Hence we can assume that $\alpha_1 + \alpha_5 < \frac{3}{7} + \varepsilon$, and thus $\alpha_2 + \alpha_3 + \alpha_4 > \frac{4}{7} - \varepsilon$. Now we have $\alpha_2 + \alpha_3 > \frac{2}{3}(\alpha_2 + \alpha_3 + \alpha_4) > \frac{8}{21} - \varepsilon$. Since $\frac{1}{3} < \frac{8}{21} - \varepsilon < 0.4 - \varepsilon$, we can assume that $\alpha_2 + \alpha_3 > 0.4 - \varepsilon$. If $\alpha_2 + \alpha_3 \geqslant \frac{3}{7} + \varepsilon$, then we can assume that $\alpha_2 + \alpha_3 > \frac{4}{7} - \varepsilon$ and $\alpha_1 + \alpha_4 + \alpha_5 < \frac{3}{7} + \varepsilon$. But we also have $\alpha_1 + \alpha_4 + \alpha_5 > \frac{2}{7} - \frac{1}{2} \varepsilon + 2 \cdot \frac{1}{8} > \frac{3}{7} + \varepsilon$, making a contradiction. Now we assume that $\alpha_2 + \alpha_3 \in \left(0.4 - \varepsilon,\ \frac{3}{7} + \varepsilon\right)$. Since $\frac{1}{8} < \alpha_5 < \frac{1}{7} + \varepsilon$, suppose that $\frac{1}{8} < \alpha_5 < \frac{1}{7} - 2 \varepsilon$, we have
\begin{equation}
\alpha_2 + \alpha_3 + \alpha_5 \in \left(0.525 - \varepsilon,\ \frac{4}{7} - \varepsilon \right).
\end{equation}
Since $\theta + \varepsilon < 0.525 - \varepsilon$, we have (129) holds for $f(n)$. The remaining part is $\frac{1}{7} - 2 \varepsilon < \alpha_5 < \frac{1}{7} + \varepsilon$, leading to a ``loss'' of size $O(\varepsilon)$.

Hence (129) holds for $f(n)$ with an error $O(\varepsilon)$ in \textbf{Case 5}.

\textbf{Case 6: $\Omega(n) = 4$, a product of two variables lies in $\left(1 - \theta - \varepsilon,\ \theta + \varepsilon \right)$.}
Suppose that $n = p_1 p_2 p_3 p_4$ and $\alpha_1 + \alpha_2 \in \left(1 - \theta - \varepsilon,\ \theta + \varepsilon \right)$ (of course, $\alpha_3 + \alpha_4 \in \left(1 - \theta - \varepsilon,\ \theta + \varepsilon \right)$ too). Without loss of generality, we further assume that $\alpha_1 > \alpha_2$ and $\alpha_3 > \alpha_4$, while the remaining part $\alpha_i = \alpha_j$ gives a ``loss'' of size $O(\varepsilon)$. Now we have $\alpha_2, \alpha_4 < \frac{\theta}{2} \leqslant \frac{11}{42} - \varepsilon$, and we can assume that $\alpha_2, \alpha_4 < \frac{1}{5} + \varepsilon$. If $\alpha_4 < \frac{1}{5} - 2 \varepsilon$, then we have $\alpha_2 + \alpha_4 < 0.4 - \varepsilon$, and we can assume that $\alpha_2 + \alpha_4 < \frac{1}{3}$. If $\alpha_1 \geqslant \frac{1}{3}$, we can assume that $\alpha_1 > 0.4 - \varepsilon$. Thus we have $\alpha_1 + \alpha_2 > 0.4 - \varepsilon + \frac{1}{8} = 0.525 - \varepsilon > \theta + \varepsilon$, making a contradiction with the assumption $\alpha_1 + \alpha_2 \in \left(1 - \theta - \varepsilon,\ \theta + \varepsilon \right)$. Similar arguments also hold if $\alpha_3 \geqslant \frac{1}{3}$. Hence, we can assume that $\alpha_1, \alpha_3 < \frac{1}{3}$, and thus
\begin{equation}
1 = \alpha_1 + \alpha_2 + \alpha_3 + \alpha_4 = \alpha_1 + \alpha_3 + (\alpha_2 + \alpha_4) < \frac{1}{3} + \frac{1}{3} + \frac{1}{3} = 1,
\end{equation}
making a contradiction. The remaining part is $\frac{1}{5} - 2 \varepsilon < \alpha_4 < \frac{1}{5} + \varepsilon$, leading to a ``loss'' of size $O(\varepsilon)$. Hence (129) holds for $f(n)$ with an error $O(\varepsilon)$ in \textbf{Case 6}.

\textbf{Case 7: $\Omega(n) = 4$, one variable lies in $\left(1 - \theta - \varepsilon,\ \theta + \varepsilon \right)$.}
Suppose that $n = p_0 p_1 p_2 p_3$ and $\alpha_0 \geqslant \alpha_1 \geqslant \alpha_2 \geqslant \alpha_3$. Now we have $\alpha_0 \in \left(1 - \theta - \varepsilon,\ \theta + \varepsilon \right)$ and $\alpha_1 + \alpha_2 + \alpha_3 \in \left(1 - \theta - \varepsilon,\ \theta + \varepsilon \right)$. We also have $\alpha_1 + \alpha_3 \leqslant \frac{1}{2}$ and $\alpha_2 \leqslant \frac{1}{3}$. Note that $\Sigma_{42409}$ only counts numbers with all prime factors smaller than $\max(m_1, m_2, m_3) \leqslant x^{\frac{3}{7} + \varepsilon}$, one can easily find that this type of $n$ will not be counted in $\Sigma_{42409}$ by a simple observation. Hence, this type of $n$ will only be counted in one part of $\Sigma_{42405}$:
\begin{equation}
\sum_{\substack{n = p_1 p_2 p_3 \beta \\ 1 - \theta - \varepsilon < \alpha_1 + \alpha_2 + \alpha_3 < \theta + \varepsilon }} \psi\left(\beta, p_3 \right),
\end{equation}
where $\beta$ in (189) is a large prime. Using Buchstab's identity, we have
\begin{equation}
\sum_{\substack{n = p_1 p_2 p_3 \beta \\ 1 - \theta - \varepsilon < \alpha_1 + \alpha_2 + \alpha_3 < \theta + \varepsilon }} \psi\left(\beta, p_3 \right) = \sum_{\substack{n = p_1 p_2 p_3 \beta \\ 1 - \theta - \varepsilon < \alpha_1 + \alpha_2 + \alpha_3 < \theta + \varepsilon }} \psi\left(\beta, x^{\kappa} \right) - \sum_{\substack{n = p_1 p_2 p_3 p_4 \beta \\ 1 - \theta - \varepsilon < \alpha_1 + \alpha_2 + \alpha_3 < \theta + \varepsilon \\ \kappa \leqslant \alpha_4 < \min\left(\alpha_3, \frac{1}{2}(1 - \alpha_1 - \alpha_2 - \alpha_3) \right) }} \psi\left(\beta, p_4 \right).
\end{equation}
By Lemma~\ref{l217}, (129) holds for the first sum in (190) since $(\alpha_1 + \alpha_2, \alpha_3) \in A$, and we only need to deal with the second sum in (190) that counts numbers with $5$ or more prime factors, with a product of $3$ variables lies in $\left(1 - \theta - \varepsilon,\ \theta + \varepsilon \right)$. Now we reduce this case to \textbf{Case 3} and parts of \textbf{Cases 1 and 2} that a product of $3$ variables lies in $\left(1 - \theta - \varepsilon,\ \theta + \varepsilon \right)$. Since the second sum in (168) is a positive sum that cannot be dropped directly, we need (129) holds for $f(n)$ in \textbf{Case 3} and parts of \textbf{Cases 1 and 2} without an error $O(\varepsilon)$. By the discussions above, we know that (129) holds for $f(n)$ in \textbf{Case 3} without an error $O(\varepsilon)$, and we only need to consider the following 2 subcases.

\textbf{Case 7.1: $\Omega(n) = 7$, a product of three variables lies in $\left(1 - \theta - \varepsilon,\ \theta + \varepsilon \right)$.}
Suppose that $n = p_1 \cdots p_7$ and $\alpha_1 + \alpha_2 + \alpha_3 \in \left(1 - \theta - \varepsilon,\ \theta + \varepsilon \right)$ (of course, $\alpha_4 + \alpha_5 + \alpha_6 + \alpha_7 \in \left(1 - \theta - \varepsilon,\ \theta + \varepsilon \right)$ too). Without loss of generality, we further assume that $\alpha_1 \geqslant \alpha_2 \geqslant \alpha_3$ and $\alpha_4 \geqslant \alpha_5 \geqslant \alpha_6 \geqslant \alpha_7$. Since $\alpha_4 + \alpha_5 + \alpha_6 \geqslant \frac{3}{4}(1 - \theta - \varepsilon) \geqslant \frac{3}{4}\left(1 - \frac{11}{21} - \varepsilon \right) > 0.35 > \frac{1}{3}$, we can assume that $\alpha_4 + \alpha_5 + \alpha_6 > 0.4 - \varepsilon$. Now we have $1 = \alpha_1 + \cdots + \alpha_7 = (\alpha_1 + \alpha_2 + \alpha_3) + (\alpha_4 + \alpha_5 + \alpha_6) + \alpha_7 > (1 - \theta - \varepsilon) + (0.4 - \varepsilon) + \left(\frac{5 - 8 \theta}{6} - \varepsilon \right) > 1$, since $\theta < \frac{11}{21} < \frac{37}{70} - 100 \varepsilon$. This makes a contradiction, and we know that (129) holds for $f(n)$ in \textbf{Case 7.1}.

\textbf{Case 7.2: $\Omega(n) = 6$, a product of three variables lies in $\left(1 - \theta - \varepsilon,\ \theta + \varepsilon \right)$.}
Suppose that $n = p_1 \cdots p_6$ and $\alpha_1 + \alpha_2 + \alpha_3 \in \left(1 - \theta - \varepsilon,\ \theta + \varepsilon \right)$ (of course, $\alpha_4 + \alpha_5 + \alpha_6 \in \left(1 - \theta - \varepsilon,\ \theta + \varepsilon \right)$ too). Without loss of generality, we further assume that $\alpha_1 \geqslant \alpha_2 \geqslant \alpha_3$ and $\alpha_4 \geqslant \alpha_5 \geqslant \alpha_6$. Now we have $\alpha_1 + \alpha_4 \geqslant \frac{1}{3}$ and $\alpha_1 + \alpha_4 \geqslant \alpha_2 + \alpha_5 \geqslant \alpha_3 + \alpha_6$, and we can assume that $\alpha_1 + \alpha_4 > 0.4 - \varepsilon$. If $\alpha_2 + \alpha_5 \geqslant \frac{1}{3}$, we can assume that $\alpha_2 + \alpha_5 > 0.4 - \varepsilon$, and thus $1 = \alpha_1 + \cdots + \alpha_6 > (0.4 - \varepsilon) + (0.4 - \varepsilon) + 2 \cdot \frac{1}{8} = 1.05 - 2 \varepsilon > 1$, making a contradiction. Now we can assume that $\alpha_2 + \alpha_5 < \frac{1}{3}$, and thus $\alpha_2 + \alpha_6, \alpha_3 + \alpha_5 < \alpha_2 + \alpha_5 < \frac{1}{3}$. If $\alpha_4 + \alpha_5 < \frac{1}{3}$, we know that $\alpha_1 + \alpha_3 = 1 - (\alpha_2 + \alpha_6) - (\alpha_4 + \alpha_5) > \frac{1}{3}$, and we can assume that $\alpha_1 + \alpha_3 > 0.4 - \varepsilon$. But now we have $\alpha_1 + \alpha_2 + \alpha_3 > 0.4 - \varepsilon + \frac{1}{8} = 0.525 - \varepsilon > \theta + \varepsilon$, making a contradiction with the assumption $\alpha_1 + \alpha_2 + \alpha_3 \in \left(1 - \theta - \varepsilon,\ \theta + \varepsilon \right)$. Thus, we can assume that $\alpha_4 + \alpha_5 \geqslant \frac{1}{3}$. Similarly, we can assume that $\alpha_1 + \alpha_2 \geqslant \frac{1}{3}$. Now we can assume that $\alpha_1 + \alpha_2, \alpha_4 + \alpha_5 > 0.4 - \varepsilon$. But now we have $1 = \alpha_1 + \cdots + \alpha_6 > (0.4 - \varepsilon) + (0.4 - \varepsilon) + 2 \cdot \frac{1}{8} = 1.05 - 2 \varepsilon > 1$, making a contradiction. Hence (129) holds for $f(n)$ in \textbf{Case 7.2}.

\textbf{Case 8: $\Omega(n) = 4$, no product of variables lies in $\left[1 - \theta,\ \theta \right]$.}
Suppose that $n = p_1 p_2 p_3 p_4$ and $\alpha_1 > \alpha_2 > \alpha_3 > \alpha_4$, while the remaining part $\alpha_i = \alpha_j$ gives a ``loss'' of size $O(\varepsilon)$. In this case we can ``view'' $\left[\frac{3}{7} + \varepsilon,\ \frac{4}{7} - \varepsilon \right]$ as a ``fake'' Type-II range. Since we have $\theta < \frac{11}{21} < \frac{9}{17}$, the results in [\cite{MaynardLargeModuliI}, Chapter 9] can be applied here if any of the following conditions holds:

(1). $\alpha_4 \geqslant \frac{1}{7} + \varepsilon$;

(2). $\alpha_1 \leqslant \frac{3}{7} + \varepsilon$, $\alpha_1 + \alpha_4 > \frac{4}{7} - \varepsilon$.

By Condition (1), we can assume that $\frac{1}{8} < \alpha_4 < \frac{1}{7} + \varepsilon$. Suppose that $\alpha_1 > \frac{4}{7} - \varepsilon$. Note that $\Sigma_{42409}$ only counts numbers with all prime factors smaller than $\max(m_1, m_2, m_3) \leqslant x^{\frac{3}{7} + \varepsilon}$, one can easily find that this type of $n$ will not be counted in $\Sigma_{42409}$ by a simple observation. Hence, this type of $n$ will only be counted in one part of $\Sigma_{42405}$:
\begin{equation}
\sum_{\substack{n = p_2 p_3 p_4 \beta \\ \alpha_2 + \alpha_3 + \alpha_4 < \frac{3}{7} + \varepsilon}} \psi\left(\beta, p_3 \right),
\end{equation}
where $\beta$ in (191) is a large prime. Using Buchstab's identity twice, we have
\begin{align}
\nonumber \sum_{\substack{n = p_2 p_3 p_4 \beta \\ \alpha_2 + \alpha_3 + \alpha_4 < \frac{3}{7} + \varepsilon }} \psi\left(\beta, p_4 \right) =&\ \sum_{\substack{n = p_2 p_3 p_4 \beta \\ \alpha_2 + \alpha_3 + \alpha_4 < \frac{3}{7} + \varepsilon }} \psi\left(\beta, x^{\kappa} \right) - \sum_{\substack{n = p_2 p_3 p_4 p_5 \beta \\ \alpha_2 + \alpha_3 + \alpha_4 < \frac{3}{7} + \varepsilon \\ \kappa \leqslant \alpha_5 < \min\left(\alpha_4, \frac{1}{2}(1 - \alpha_2 - \alpha_3 - \alpha_4) \right) }} \psi\left(\beta, p_5 \right) \\
\nonumber =&\ \sum_{\substack{n = p_2 p_3 p_4 \beta \\ \alpha_2 + \alpha_3 + \alpha_4 < \frac{3}{7} + \varepsilon }} \psi\left(\beta, x^{\kappa} \right) - \sum_{\substack{n = p_2 p_3 p_4 p_5 \beta \\ \alpha_2 + \alpha_3 + \alpha_4 < \frac{3}{7} + \varepsilon \\ \kappa \leqslant \alpha_5 < \min\left(\alpha_4, \frac{1}{2}(1 - \alpha_2 - \alpha_3 - \alpha_4) \right) }} \psi\left(\beta, x^{\kappa} \right) \\
& + \sum_{\substack{n = p_2 p_3 p_4 p_5 p_6 \beta \\ \alpha_2 + \alpha_3 + \alpha_4 < \frac{3}{7} + \varepsilon \\ \kappa \leqslant \alpha_5 < \min\left(\alpha_3, \frac{1}{2}(1 - \alpha_2 - \alpha_3 - \alpha_4) \right) \\ \kappa \leqslant \alpha_6 < \min\left(\alpha_5, \frac{1}{2}(1 - \alpha_2 - \alpha_3 - \alpha_4 - \alpha_5) \right) }} \psi\left(\beta, p_6 \right).
\end{align}
Note that we have $\alpha_2 + \alpha_3 + \alpha_4 < \frac{3}{7} + \varepsilon$ and $\alpha_5 < \alpha_4 < \frac{1}{7} + \varepsilon$. Since we have $\left(\frac{3}{7} + \varepsilon, \frac{1}{7} + \varepsilon, \varepsilon \right) \in \boldsymbol{S}_{3}$ when $\theta < \frac{11}{20}$, by the Type-II range $\left[\varepsilon,\ \frac{17}{126} - \varepsilon \right]$, (129) holds for the first and the second sums in (192). Now we only need to deal with the last sum in (192) that counts numbers with $6$ or more prime factors. Since this sum is negative, we know that (129) holds for $f(n)$ with an error $O(\varepsilon)$ in \textbf{Case 8} if $\alpha_4 < \frac{1}{7} + \varepsilon$ and $\alpha_1 > \frac{4}{7} - \varepsilon$ by the discussions in \textbf{Cases 1 and 2}.

Now by Condition (2) and the ``fake'' Type-II range $\left[\frac{3}{7} + \varepsilon,\ \frac{4}{7} - \varepsilon \right]$, we only need to consider the case $\alpha_1 + \alpha_4 < \frac{3}{7} + \varepsilon$. Now we have $\alpha_2 + \alpha_3 > \frac{4}{7} - \varepsilon$, and thus $\alpha_1 + \alpha_2 > \alpha_1 + \alpha_3 > \alpha_2 + \alpha_3 > \frac{4}{7} - \varepsilon$. However, we cannot show that (129) holds for $f(n)$ in \textbf{Case 8}. An obvious counterexample is
\begin{equation}
\boldsymbol{\alpha}_4 = \left(\frac{2}{7} + 3 \varepsilon, \frac{2}{7} + 2 \varepsilon, \frac{2}{7} + \varepsilon, \frac{1}{7} - 6 \varepsilon \right).
\end{equation}
By a simple observation, one can easily find that such $f(n)$ that (129) does not hold will only be negative parts of $\Sigma_{42405}$ and $\Sigma_{42409}$ since $\Omega(n) = 4$. We can discard them in the upper bound case. We also have the following conditions on $n = p_1 p_2 p_3 p_4$ counted in such $f(n)$:
\begin{equation}
\frac{5 - 8 \theta}{6} - \varepsilon \leqslant \alpha_4 < \frac{1}{7} + \varepsilon, \quad \frac{3}{7} + \varepsilon > \alpha_1 + \alpha_4 > \alpha_2 + \alpha_4 > \alpha_3 + \alpha_4, \quad \frac{4}{7} - \varepsilon < \alpha_2 + \alpha_3 < \alpha_1 + \alpha_3 < \alpha_1 + \alpha_2.
\end{equation}
By the conditions (194), we can deduce an upper bound for $\alpha_1$, $\alpha_2$ and $\alpha_3$:
\begin{equation}
\alpha_3 < \alpha_2 < \alpha_1 \leqslant \frac{3}{7} - \alpha_4 \leqslant \frac{3}{7} - \frac{5 - 8 \theta}{6} + \varepsilon.
\end{equation}
Since $\alpha_1 + \alpha_2 + \alpha_3 + \alpha_4 = 1$, by (194) and (195), we can deduce a lower bound for $\alpha_1$, $\alpha_2$ and $\alpha_3$:
\begin{equation}
\alpha_1 > \alpha_2 > \alpha_3 = 1 - (\alpha_1 + \alpha_2 + \alpha_4) > 1 - 2 \left(\frac{3}{7} - \frac{5 - 8 \theta}{6} + \varepsilon \right) - \frac{1}{7} = \frac{5 - 8 \theta}{3} - 2 \varepsilon.
\end{equation}
Now, we know that such part of $\Sigma_{42405}$ that (129) does not hold is no more than
\begin{equation}
\sum_{\substack{n = p_1 p_2 p_3 \beta \\ \frac{5 - 8 \theta}{3} - 2 \varepsilon < \alpha_1 \leqslant \frac{3}{7} - \frac{5 - 8 \theta}{6} + \varepsilon \\ \frac{5 - 8 \theta}{3} - 2 \varepsilon < \alpha_2 \leqslant \frac{3}{7} - \frac{5 - 8 \theta}{6} + \varepsilon \\ \frac{5 - 8 \theta}{6} - \varepsilon \leqslant \alpha_3 < \frac{1}{7} + \varepsilon \\ \boldsymbol{\alpha}_3 \notin \boldsymbol{G}_{3} }} \psi\left(\beta, \left(\frac{2 x}{p_1 p_2 p_3} \right)^{\frac{1}{2}} \right) = \sum_{\substack{n = p_1 p_2 p_3 \beta \\ \frac{5 - 8 \theta}{3} - 2 \varepsilon < \alpha_1 \leqslant \frac{3}{7} - \frac{5 - 8 \theta}{6} + \varepsilon \\ \frac{5 - 8 \theta}{3} - 2 \varepsilon < \alpha_2 \leqslant \frac{3}{7} - \frac{5 - 8 \theta}{6} + \varepsilon \\ \frac{5 - 8 \theta}{6} - \varepsilon \leqslant \alpha_3 < \frac{1}{7} + \varepsilon \\ \boldsymbol{\alpha}_3 \notin \boldsymbol{G}_{3} }} 1.
\end{equation}
For such part of $\Sigma_{42409}$ that (129) does not hold, we have an extra condition $\alpha_1 \geqslant \tau = \frac{3(1-\theta)}{5} - \varepsilon$ since $\theta < \frac{11}{21}$. Note that we have $\frac{3(1-\theta)}{5} > \frac{5 - 8 \theta}{3}$ when $\theta > \frac{16}{31}$, and we can assume that $\theta \geqslant \frac{29}{56} > \frac{16}{31}$. Let $i, j \in \{1,2,3\}$ and $i \neq j$. By (194) and the condition $\boldsymbol{\alpha}_2 \in C$, we know that $m_1 \neq p_i p_j$ and $m_2 \neq p_i p_j$. Since we have $\Omega(m_1 m_2) \geqslant 3$, we know that $m_3 \neq p_i p_j$ and $m_3 \neq p_i p_4$. Now the only possible cases are $m_1 = p_i p_4$ or $m_2 = p_i p_4$, and such part of $\Sigma_{42409}$ that (129) does not hold is no more than
\begin{equation}
\sum_{\substack{n = p_1 p_2 p_3 p_4 \\ \tau \leqslant \alpha_1 \leqslant \frac{3}{7} - \frac{5 - 8 \theta}{6} + \varepsilon \\ \frac{5 - 8 \theta}{3} - 2 \varepsilon < \alpha_2 \leqslant \frac{3}{7} - \frac{5 - 8 \theta}{6} + \varepsilon \\ \frac{5 - 8 \theta}{6} - \varepsilon \leqslant \alpha_3 < \frac{1}{7} + \varepsilon \\ (\alpha_1 + \alpha_3, \alpha_2) \in C \text{ or } (\alpha_1, \alpha_2 + \alpha_3) \in C \\ \boldsymbol{\alpha}_3 \notin \boldsymbol{G}_{3} }} 1.
\end{equation}

Finally, we get the following theorem by combining all the 8 cases above:
\begin{theorem}\label{t410}
Define
\begin{align}
\nonumber \boldsymbol{D}_5 =&\ \boldsymbol{D}_5(\theta) = \left\{\boldsymbol{\alpha}_{3}: \frac{5 - 8 \theta}{3} \leqslant \alpha_1 \leqslant \frac{3}{7} - \frac{5 - 8 \theta}{6},\ \frac{5 - 8 \theta}{3} \leqslant \alpha_2 \leqslant \frac{3}{7} - \frac{5 - 8 \theta}{6},\ \frac{5 - 8 \theta}{6} \leqslant \alpha_3 \leqslant \frac{1}{7},\ \boldsymbol{\alpha}_3 \notin \boldsymbol{G}_{3} \right\} \\
\nonumber \boldsymbol{D}_6 =&\ \boldsymbol{D}_6(\theta) = \left\{\boldsymbol{\alpha}_{3}: \tau \leqslant \alpha_1 \leqslant \frac{3}{7} - \frac{5 - 8 \theta}{6},\ \frac{5 - 8 \theta}{3} \leqslant \alpha_2 \leqslant \frac{3}{7} - \frac{5 - 8 \theta}{6},\ \frac{5 - 8 \theta}{6} \leqslant \alpha_3 \leqslant \frac{1}{7},\ \boldsymbol{\alpha}_3 \notin \boldsymbol{G}_{3}, \right.\\
\nonumber & \left. \qquad \qquad \quad (\alpha_1 + \alpha_3, \alpha_2) \in C \text{ or } (\alpha_1, \alpha_2 + \alpha_3) \in C \right\}.
\end{align}
Let $(\theta_1, \theta_2) \in \boldsymbol{E}_{0210}$. Suppose that we have
$$
\theta_1 \leqslant \frac{1}{3} - \varepsilon, \quad \theta_2 \leqslant \frac{1}{5} \quad \text{and} \quad \frac{29}{56} \leqslant \theta_1 + \theta_2 \leqslant \frac{11}{21} - \varepsilon.
$$
Then we have
$$
C_1^{\prime}(\theta_1, \theta_2) \leqslant 1 + \int_{1 - \theta}^{\frac{1}{2}} \frac{\omega\left(\frac{1-t}{t} \right)}{t^2} d t + I_5 + I_6 + O(\varepsilon) = 1 + \log\left(\frac{\theta}{1 - \theta} \right) + I_5 + I_6 + O(\varepsilon),
$$
where
$$
I_5 = \int_{(t_1, t_2, t_3) \in \boldsymbol{D}_{5}} \frac{1}{t_1 t_2 t_3 (1 - t_1 - t_2 - t_3)} d t_3 d t_2 d t_1
$$
and
$$
I_6 = \int_{(t_1, t_2, t_3) \in \boldsymbol{D}_{6}} \frac{1}{t_1 t_2 t_3 (1 - t_1 - t_2 - t_3)} d t_3 d t_2 d t_1.
$$
Note that the ``loss integrals'' $I_5$ and $I_6$ correspond to the sums in (197) and (198) respectively.
\end{theorem}
\begin{remark*}
The above arguments also show that 
$$
C_1^{\prime}(\theta_1, \theta_2) \leqslant 1 + \int_{1 - \theta}^{\frac{1}{2}} \frac{\omega\left(\frac{1-t}{t} \right)}{t^2} d t + O(\varepsilon) = 1 + \log\left(\frac{\theta}{1 - \theta} \right) + O(\varepsilon)
$$
if $(\theta_1, \theta_2) \in \boldsymbol{E}_{0210}$, $\theta_1 \leqslant \frac{1}{3} - \varepsilon$, $\theta_2 \leqslant \frac{1}{5}$ and $\theta_1 + \theta_2 \leqslant \frac{29}{56} - \varepsilon$, since we have $\kappa > \frac{1}{7}$ in this case and (129) holds for $f(n)$ with an error $O(\varepsilon)$ in \textbf{Case 8}.

By Theorem~\ref{t410}, we get
$$
C_1^{\prime}(0.32, 0.20) \leqslant 1 + \log\left(\frac{13}{12} \right) + 10^{-5} + O(\varepsilon) < 1.0801
$$
and
$$
C_1^{\prime}(0.33, 0.19) \leqslant 1 + \log\left(\frac{13}{12} \right) + 10^{-5} + O(\varepsilon) < 1.0801.
$$
\end{remark*}

The Type-II range for $\boldsymbol{E}_{0211}$ is
\begin{equation}
\left[\varepsilon,\ \frac{1}{4}(2 - 2 \theta_1 - 5 \theta_2) - \varepsilon \right] \cup \left[2 \theta - 1 + \varepsilon,\ \frac{1}{6}(5 - 8 \theta) - \varepsilon \right] \cup \left[\theta_2 + \varepsilon,\ \frac{1}{6}(5 - 8 \theta_1 - 4 \theta_2) - \varepsilon \right] \cup \left[\theta_1 + \varepsilon,\ \frac{1}{2}(1 - \theta_2) - \varepsilon \right].
\end{equation}
The decompositions are similar to which in the case $\boldsymbol{E}_{0103}$.

The Type-II range for $\boldsymbol{E}_{0212}$ is
\begin{equation}
\left[\varepsilon,\ \frac{1}{4}(2 - 2 \theta_1 - 5 \theta_2) - \varepsilon \right] \cup \left[2 \theta - 1 + \varepsilon,\ \frac{1}{6}(5 - 8 \theta) - \varepsilon \right] \cup \left[\theta_2 + \varepsilon,\ \frac{1}{6}(5 - 8 \theta_1 - 4 \theta_2) - \varepsilon \right] \cup \left[\theta_1 + \varepsilon,\ \frac{1}{2}(2 - \theta_1 - 4 \theta_2) - \varepsilon \right].
\end{equation}
The decompositions are similar to which in the case $\boldsymbol{E}_{0103}$.

The Type-II range for $\boldsymbol{E}_{0213}$ is
\begin{equation}
\left[\varepsilon,\ \frac{1}{4}(2 - 2 \theta_1 - 5 \theta_2) - \varepsilon \right] \cup \left[2 \theta - 1 + \varepsilon,\ \frac{1}{6}(5 - 8 \theta) - \varepsilon \right] \cup \left[\theta_2 + \varepsilon,\ \frac{1}{4}(2 - 3 \theta_1) - \varepsilon \right] \cup \left[\theta_1 + \varepsilon,\ \frac{1}{2}(1 - \theta_2) - \varepsilon \right].
\end{equation}
The decompositions are similar to which in the case $\boldsymbol{E}_{0103}$.

The Type-II range for $\boldsymbol{E}_{0214}$ is
\begin{equation}
\left[\varepsilon,\ \frac{1}{4}(2 - 2 \theta_1 - 5 \theta_2) - \varepsilon \right] \cup \left[2 \theta - 1 + \varepsilon,\ \frac{1}{6}(5 - 8 \theta) - \varepsilon \right] \cup \left[\theta_2 + \varepsilon,\ \frac{1}{4}(2 - 3 \theta_1) - \varepsilon \right] \cup \left[\theta_1 + \varepsilon,\ \frac{1}{2}(2 - \theta_1 - 4 \theta_2) - \varepsilon \right].
\end{equation}
The decompositions are similar to which in the case $\boldsymbol{E}_{0103}$.

The Type-II range for $\boldsymbol{E}_{0215}$ is
\begin{equation}
\left[\varepsilon,\ \frac{1}{4}(2 - 2 \theta_1 - 5 \theta_2) - \varepsilon \right] \cup \left[\theta_2 + \varepsilon,\ \frac{1}{4}(2 - 3 \theta_1) - \varepsilon \right] \cup \left[\theta_1 + \varepsilon,\ \frac{1}{2}(2 - \theta_1 - 4 \theta_2) - \varepsilon \right].
\end{equation}
The decompositions are similar to which in the case $\boldsymbol{E}_{0104}$.

\subsubsection{$\boldsymbol{E}_{03}$}
We divide $\boldsymbol{E}_{03}$ into 7 subregions:
$$
\boldsymbol{E}_{03} = \boldsymbol{E}_{0301} \cup \boldsymbol{E}_{0302} \cup \boldsymbol{E}_{0303} \cup \boldsymbol{E}_{0304} \cup \boldsymbol{E}_{0305} \cup \boldsymbol{E}_{0306} \cup \boldsymbol{E}_{0307},
$$
where
\begin{align}
\nonumber \boldsymbol{E}_{0301} =&\ \left\{ (\theta_1, \theta_2) : \frac{7}{20} < \theta_1 \leqslant \frac{5}{14},\ \frac{1}{7}(6 - 14 \theta_1) \leqslant \theta_2 < \frac{1}{14}(5 - 8 \theta_1) \right\}, \\
\nonumber \boldsymbol{E}_{0302} =&\ \left\{ (\theta_1, \theta_2) : \frac{9}{28} < \theta_1 \leqslant \frac{7}{20},\ \frac{1}{7}(6 - 14 \theta_1) \leqslant \theta_2 < \frac{1}{18}(9 - 16 \theta_1) \right. \\
\nonumber & \left. \qquad \qquad \quad \text{ or } \frac{7}{20} < \theta_1 \leqslant \frac{5}{14},\ \frac{1}{14}(5 - 8 \theta_1) \leqslant \theta_2 < \frac{1}{18}(9 - 16 \theta_1) \right\}, \\
\nonumber \boldsymbol{E}_{0303} =&\ \left\{ (\theta_1, \theta_2) : \frac{20}{63} < \theta_1 \leqslant \frac{9}{28},\ \frac{1}{7}(6 - 14 \theta_1) \leqslant \theta_2 \leqslant \frac{1}{8}(4 - 7 \theta_1) \right. \\
\nonumber & \left. \qquad \qquad \quad \text{ or } \frac{9}{28} < \theta_1 \leqslant \frac{5}{14},\ \frac{1}{18}(9 - 16 \theta_1) \leqslant \theta_2 \leqslant \frac{1}{8}(4 - 7 \theta_1) \right\}, \\
\nonumber \boldsymbol{E}_{0304} =&\ \left\{ (\theta_1, \theta_2) : \frac{11}{35} < \theta_1 \leqslant \frac{20}{63},\ \frac{1}{7}(6 - 14 \theta_1) \leqslant \theta_2 < \frac{1}{3}(1 - \theta_1) \right. \\
\nonumber & \qquad \qquad \quad \text{ or } \frac{20}{63} < \theta_1 \leqslant \frac{13}{40},\ \frac{1}{8}(4 - 7 \theta_1) < \theta_2 < \frac{1}{3}(1 - \theta_1) \\
\nonumber & \left. \qquad \qquad \quad \text{ or } \frac{13}{40} < \theta_1 \leqslant \frac{5}{14},\ \frac{1}{8}(4 - 7 \theta_1) < \theta_2 < \frac{1}{20}(11 - 20 \theta_1) \right\}, \\
\nonumber \boldsymbol{E}_{0305} =&\ \left\{ (\theta_1, \theta_2) : \frac{13}{40} < \theta_1 \leqslant \frac{5}{14},\ \frac{1}{20}(11 - 20 \theta_1) \leqslant \theta_2 < \frac{1}{3}(1 - \theta_1) \right\}, \\
\nonumber \boldsymbol{E}_{0306} =&\ \left\{ (\theta_1, \theta_2) : \frac{43}{140} < \theta_1 \leqslant \frac{11}{35},\ \frac{1}{7}(6 - 14 \theta_1) \leqslant \theta_2 < \frac{1}{20}(11 - 20 \theta_1) \right. \\
\nonumber & \left. \qquad \qquad \quad \text{ or } \frac{11}{35} < \theta_1 < \frac{13}{40},\ \frac{1}{3}(1 - \theta_1) \leqslant \theta_2 < \frac{1}{20}(11 - 20 \theta_1) \right\}, \\
\nonumber \boldsymbol{E}_{0307} =&\ \left\{ (\theta_1, \theta_2) : \frac{2}{7} < \theta_1 \leqslant \frac{43}{140},\ \frac{1}{7}(6 - 14 \theta_1) \leqslant \theta_2 < \frac{1}{4}(2 - 3 \theta_1) \right. \\
\nonumber & \qquad \qquad \quad \text{ or } \frac{43}{140} < \theta_1 \leqslant \frac{13}{40},\ \frac{1}{20}(11 - 20 \theta_1) \leqslant \theta_2 < \frac{1}{4}(2 - 3 \theta_1) \\
\nonumber & \left. \qquad \qquad \quad \text{ or } \frac{13}{40} < \theta_1 \leqslant \frac{5}{14},\ \frac{1}{3}(1 - \theta_1) \leqslant \theta_2 < \frac{1}{4}(2 - 3 \theta_1) \right\}.
\end{align}
Note that we have $\theta < \frac{11}{20}$ for $(\theta_1, \theta_2) \in \boldsymbol{E}_{0301} \cup \boldsymbol{E}_{0302} \cup \boldsymbol{E}_{0303} \cup \boldsymbol{E}_{0304} \cup \boldsymbol{E}_{0306}$.

The Type-II range for $\boldsymbol{E}_{0301}$ is
\begin{equation}
\left[\varepsilon,\ \frac{1}{6}(5 - 8 \theta_1 - 4 \theta_2) - \varepsilon \right] \cup \left[\theta_1 + \varepsilon,\ \frac{1}{2}(1 - \theta_2) - \varepsilon \right].
\end{equation}
The decompositions in this case will be discussed later together with the case $\boldsymbol{E}_{0501}$.

The Type-II range for $\boldsymbol{E}_{0302}$ is
\begin{equation}
\left[\varepsilon,\ \frac{1}{6}(5 - 8 \theta_1 - 8 \theta_2) - \varepsilon \right] \cup \left[\theta_2 + \varepsilon,\ \frac{1}{6}(5 - 8 \theta_1 - 4 \theta_2) - \varepsilon \right] \cup \left[\theta_1 + \varepsilon,\ \frac{1}{2}(1 - \theta_2) - \varepsilon \right].
\end{equation}
The decompositions are similar to which in the case $\boldsymbol{E}_{0210}$. Note that this case covers the remaining parts of the region (176) where the case $\boldsymbol{E}_{0210}$ does not cover.

The Type-II range for $\boldsymbol{E}_{0303}$ is
\begin{equation}
\left[\varepsilon,\ \frac{1}{5}(4 - 6 \theta_1 - 8 \theta_2) - \varepsilon \right] \cup \left[2 \theta - 1 + \varepsilon,\ \frac{1}{6}(5 - 8 \theta) - \varepsilon \right] \cup \left[\theta_2 + \varepsilon,\ \frac{1}{6}(5 - 8 \theta_1 - 4 \theta_2) - \varepsilon \right] \cup \left[\theta_1 + \varepsilon,\ \frac{1}{2}(1 - \theta_2) - \varepsilon \right].
\end{equation}
The decompositions are similar to which in the case $\boldsymbol{E}_{0211}$.

The Type-II range for $\boldsymbol{E}_{0304}$ is
\begin{equation}
\left[\varepsilon,\ \frac{1}{5}(4 - 6 \theta_1 - 8 \theta_2) - \varepsilon \right] \cup \left[2 \theta - 1 + \varepsilon,\ \frac{1}{6}(5 - 8 \theta) - \varepsilon \right] \cup \left[\theta_2 + \varepsilon,\ \frac{1}{4}(2 - 3 \theta_1) - \varepsilon \right] \cup \left[\theta_1 + \varepsilon,\ \frac{1}{2}(1 - \theta_2) - \varepsilon \right].
\end{equation}
The decompositions are similar to which in the case $\boldsymbol{E}_{0213}$.

The Type-II range for $\boldsymbol{E}_{0305}$ is
\begin{equation}
\left[\varepsilon,\ \frac{1}{5}(4 - 6 \theta_1 - 8 \theta_2) - \varepsilon \right] \cup \left[\theta_2 + \varepsilon,\ \frac{1}{4}(2 - 3 \theta_1) - \varepsilon \right] \cup \left[\theta_1 + \varepsilon,\ \frac{1}{2}(1 - \theta_2) - \varepsilon \right].
\end{equation}
The decompositions are similar to which in the case $\boldsymbol{E}_{0215}$.

The Type-II range for $\boldsymbol{E}_{0306}$ is
\begin{equation}
\left[\varepsilon,\ \frac{1}{5}(4 - 6 \theta_1 - 8 \theta_2) - \varepsilon \right] \cup \left[2 \theta - 1 + \varepsilon,\ \frac{1}{6}(5 - 8 \theta) - \varepsilon \right] \cup \left[\theta_2 + \varepsilon,\ \frac{1}{4}(2 - 3 \theta_1) - \varepsilon \right] \cup \left[\theta_1 + \varepsilon,\ \frac{1}{2}(2 - \theta_1 - 4 \theta_2) - \varepsilon \right].
\end{equation}
The decompositions are similar to which in the case $\boldsymbol{E}_{0214}$.

The Type-II range for $\boldsymbol{E}_{0307}$ is
\begin{equation}
\left[\varepsilon,\ \frac{1}{5}(4 - 6 \theta_1 - 8 \theta_2) - \varepsilon \right] \cup \left[\theta_2 + \varepsilon,\ \frac{1}{4}(2 - 3 \theta_1) - \varepsilon \right] \cup \left[\theta_1 + \varepsilon,\ \frac{1}{2}(2 - \theta_1 - 4 \theta_2) - \varepsilon \right].
\end{equation}
The decompositions are similar to which in the case $\boldsymbol{E}_{0215}$.

\subsubsection{$\boldsymbol{E}_{04}$}
We divide $\boldsymbol{E}_{04}$ into 2 subregions:
$$
\boldsymbol{E}_{04} = \boldsymbol{E}_{0401} \cup \boldsymbol{E}_{0402},
$$
where
\begin{align}
\nonumber \boldsymbol{E}_{0401} =&\ \left\{ (\theta_1, \theta_2) : \frac{5}{14} < \theta_1 \leqslant \frac{29}{80},\ \frac{1}{2}(1 - 2 \theta_1) < \theta_2 < \frac{1}{3}(-1 + 4 \theta_1) \right. \\
\nonumber & \left. \qquad \qquad \quad \text{ or } \frac{29}{80} < \theta_1 \leqslant \frac{3}{8},\ \frac{1}{2}(1 - 2 \theta_1) < \theta_2 < \frac{1}{14}(5 - 8 \theta_1) \right\}, \\
\nonumber \boldsymbol{E}_{0402} =&\ \left\{ (\theta_1, \theta_2) : \frac{29}{80} < \theta_1 \leqslant \frac{3}{8},\ \frac{1}{14}(5 - 8 \theta_1) \leqslant \theta_2 < \frac{1}{3}(-1 + 4 \theta_1) \right\}.
\end{align}
Note that we have $\theta < \frac{11}{20}$ for $(\theta_1, \theta_2) \in \boldsymbol{E}_{04}$.

The Type-II range for $\boldsymbol{E}_{0401}$ is
\begin{equation}
\left[\varepsilon,\ \frac{1}{6}(5 - 8 \theta_1 - 4 \theta_2) - \varepsilon \right] \cup \left[\theta_1 + \varepsilon,\ \frac{1}{2}(1 - \theta_2) - \varepsilon \right].
\end{equation}
The decompositions in this case will be discussed later together with the case $\boldsymbol{E}_{0501}$.

The Type-II range for $\boldsymbol{E}_{0402}$ is
\begin{equation}
\left[\varepsilon,\ \frac{1}{6}(5 - 8 \theta_1 - 8 \theta_2) - \varepsilon \right] \cup \left[\theta_2 + \varepsilon,\ \frac{1}{6}(5 - 8 \theta_1 - 4 \theta_2) - \varepsilon \right] \cup \left[\theta_1 + \varepsilon,\ \frac{1}{2}(1 - \theta_2) - \varepsilon \right].
\end{equation}
The decompositions are similar to which in the case $\boldsymbol{E}_{0210}$.

\subsubsection{$\boldsymbol{E}_{05}$}
We divide $\boldsymbol{E}_{05}$ into 9 subregions:
$$
\boldsymbol{E}_{05} = \boldsymbol{E}_{0501} \cup \boldsymbol{E}_{0502} \cup \boldsymbol{E}_{0503} \cup \boldsymbol{E}_{0504} \cup \boldsymbol{E}_{0505} \cup \boldsymbol{E}_{0506} \cup \boldsymbol{E}_{0507} \cup \boldsymbol{E}_{0508} \cup \boldsymbol{E}_{0509},
$$
where
\begin{align}
\nonumber \boldsymbol{E}_{0501} =&\ \left\{ (\theta_1, \theta_2) : \frac{5}{14} < \theta_1 < \frac{29}{80},\ \frac{1}{3}(-1 + 4 \theta_1) \leqslant \theta_2 < \frac{1}{14}(5 - 8 \theta_1) \right\}, \\
\nonumber \boldsymbol{E}_{0502} =&\ \left\{ (\theta_1, \theta_2) : \frac{5}{14} < \theta_1 \leqslant \frac{29}{80},\ \frac{1}{14}(5 - 8 \theta_1) \leqslant \theta_2 < \frac{1}{18}(9 - 16 \theta_1) \right. \\
\nonumber & \left. \qquad \qquad \quad \text{ or } \frac{29}{80} < \theta_1 < \frac{3}{8},\ \frac{1}{3}(-1 + 4 \theta_1) \leqslant \theta_2 < \frac{1}{18}(9 - 16 \theta_1) \right\}, \\
\nonumber \boldsymbol{E}_{0503} =&\ \left\{ (\theta_1, \theta_2) : \frac{5}{14} < \theta_1 \leqslant \frac{4}{11},\ \frac{1}{18}(9 - 16 \theta_1) \leqslant \theta_2 \leqslant \frac{1}{8}(4 - 7 \theta_1) \right. \\
\nonumber & \qquad \qquad \quad \text{ or } \frac{4}{11} < \theta_1 < \frac{37}{100},\ \frac{1}{18}(9 - 16 \theta_1) \leqslant \theta_2 \leqslant \frac{1}{7}(2 - 2 \theta_1) \\
\nonumber & \left. \qquad \qquad \quad \text{ or } \frac{37}{100} \leqslant \theta_1 \leqslant \frac{3}{8},\ \frac{1}{18}(9 - 16 \theta_1) \leqslant \theta_2 < \frac{1}{20}(11 - 20 \theta_1) \right\}, \\
\nonumber \boldsymbol{E}_{0504} =&\ \left\{ (\theta_1, \theta_2) : \frac{37}{100} \leqslant \theta_1 \leqslant \frac{3}{8},\ \frac{1}{20}(11 - 20 \theta_1) \leqslant \theta_2 \leqslant \frac{1}{7}(2 - 2 \theta_1) \right\}, \\
\nonumber \boldsymbol{E}_{0505} =&\ \left\{ (\theta_1, \theta_2) : \frac{4}{11} < \theta_1 < \frac{11}{30},\ \frac{1}{7}(2 - 2 \theta_1) < \theta_2 \leqslant \frac{1}{2} \theta_1 \right. \\
\nonumber & \left. \qquad \qquad \quad \text{ or } \frac{11}{30} \leqslant \theta_1 < \frac{37}{100},\ \frac{1}{7}(2 - 2 \theta_1) < \theta_2 < \frac{1}{20}(11 - 20 \theta_1) \right\}, \\
\nonumber \boldsymbol{E}_{0506} =&\ \left\{ (\theta_1, \theta_2) : \frac{11}{30} \leqslant \theta_1 < \frac{37}{100},\ \frac{1}{20}(11 - 20 \theta_1) \leqslant \theta_2 \leqslant \frac{1}{2} \theta_1 \right. \\
\nonumber & \left. \qquad \qquad \quad \text{ or } \frac{37}{100} \leqslant \theta_1 \leqslant \frac{3}{8},\ \frac{1}{7}(2 - 2 \theta_1) < \theta_2 \leqslant \frac{1}{2} \theta_1 \right\}, \\
\nonumber \boldsymbol{E}_{0507} =&\ \left\{ (\theta_1, \theta_2) : \frac{5}{14} < \theta_1 \leqslant \frac{4}{11},\ \frac{1}{8}(4 - 7 \theta_1) < \theta_2 < \frac{1}{20}(11 - 20 \theta_1) \right. \\
\nonumber & \left. \qquad \qquad \quad \text{ or } \frac{4}{11} < \theta_1 < \frac{11}{30},\ \frac{1}{2} \theta_1 < \theta_2 < \frac{1}{20}(11 - 20 \theta_1) \right\}, \\
\nonumber \boldsymbol{E}_{0508} =&\ \left\{ (\theta_1, \theta_2) : \frac{5}{14} < \theta_1 < \frac{11}{30},\ \frac{1}{20}(11 - 20 \theta_1) \leqslant \theta_2 < \frac{1}{3}(1 - \theta_1) \right. \\
\nonumber & \left. \qquad \qquad \quad \text{ or } \frac{11}{30} \leqslant \theta_1 \leqslant \frac{3}{8},\ \frac{1}{2} \theta_1 < \theta_2 < \frac{1}{3}(1 - \theta_1) \right\}, \\
\nonumber \boldsymbol{E}_{0509} =&\ \left\{ (\theta_1, \theta_2) : \frac{5}{14} < \theta_1 \leqslant \frac{3}{8},\ \frac{1}{3}(1 - \theta_1) \leqslant \theta_2 < \frac{1}{4}(2 - 3 \theta_1) \right\}.
\end{align}
Note that we have $\theta < \frac{11}{20}$ for $(\theta_1, \theta_2) \in \boldsymbol{E}_{0501} \cup \boldsymbol{E}_{0502} \cup \boldsymbol{E}_{0503} \cup \boldsymbol{E}_{0505} \cup \boldsymbol{E}_{0507}$.

The Type-II range for $\boldsymbol{E}_{0501}$ is
\begin{equation}
\left[\varepsilon,\ \frac{1}{6}(5 - 8 \theta_1 - 4 \theta_2) - \varepsilon \right] \cup \left[\theta_1 + \varepsilon,\ \frac{1}{2}(1 - \theta_2) - \varepsilon \right].
\end{equation}
In this case we discuss $\boldsymbol{E}_{0201}$, $\boldsymbol{E}_{0301}$, $\boldsymbol{E}_{0401}$ and $\boldsymbol{E}_{0501}$, since the Type-II ranges for them are same:
$$
\boldsymbol{\mathcal{Z}_3} = \boldsymbol{E}_{0201} \cup \boldsymbol{E}_{0301} \cup \boldsymbol{E}_{0401} \cup \boldsymbol{E}_{0501} = \left\{ (\theta_1, \theta_2) : \frac{1}{3} < \theta_1 \leqslant \frac{3}{8},\ \frac{1}{2}(1 - 2 \theta_1) \leqslant \theta_2 < \frac{1}{14}(5 - 8 \theta_1) \right\}.
$$

First, we shall prove the following theorem.
\begin{theorem}\label{case3}
Let $(\theta_1, \theta_2) \in \boldsymbol{\mathcal{Z}_3}$. Suppose that we have
$$
16 \theta_1 + 8 \theta_2 \leqslant 7 - 14 \varepsilon.
$$
Then (129) holds for
$$
f(n) = \mathbbm{1}_{p}(n) + \sum_{\substack{n = p_1 \beta \\ 1 - \theta - \varepsilon < \alpha_1 < \frac{1}{2} }} \psi\left(\beta, p_1 \right),
$$
and we have
$$
C_1^{\prime}(\theta_1, \theta_2) \leqslant 1 + \int_{1 - \theta}^{\frac{1}{2}} \frac{\omega\left(\frac{1-t}{t} \right)}{t^2} d t + O(\varepsilon) = 1 + \log\left(\frac{\theta}{1 - \theta} \right) + O(\varepsilon).
$$
\end{theorem}
\begin{proof}
Since we have $(\theta_1, \theta_2) \in \boldsymbol{\mathcal{Z}_3}$ and $16 \theta_1 + 8 \theta_2 \leqslant 7 - 14 \varepsilon$, all of the following conditions hold true:
$$
\kappa = \frac{5 - 8 \theta_1 - 4 \theta_2}{6} - \varepsilon > \frac{1}{4}, \quad \frac{1}{2} < \theta < \frac{11}{21} - \varepsilon.
$$

Now, we can decompose our $\mathbbm{1}_{n \sim x, n = p}(n) = \psi\left(n, (2 x)^{\frac{1}{2}}\right)$ in a way similar to the decompositions in the proof of Theorem~\ref{case1}. By the discussions in Theorem~\ref{case1}, we only need to show that (129) holds for $f(n) = $ sums that count numbers with $4$ or more prime factors. Since $\kappa > \frac{1}{4}$, any $n \sim x$ with $4$ or more prime factors must has at least one factor smaller than $x^{\kappa}$, and we can use our Type-II range $\left[\varepsilon,\ \frac{5 - 8 \theta_1 - 4 \theta_2}{6} - \varepsilon \right]$ to show that (129) holds for $f(n)$ that count numbers with $4$ or more prime factors. The proof of Theorem~\ref{case3} is now completed.
\end{proof}

For the remaining parts of $\boldsymbol{\mathcal{Z}_3}$, we use $\kappa = \frac{5 - 8 \theta_1 - 4 \theta_2}{6} - \varepsilon$ as the ``starting point'' and use the second Type-II range $\left[\theta_1 + \varepsilon,\ \frac{1}{2}(1 - \theta_2) - \varepsilon \right]$ to discard sums explicitly. The decompositions are similar to which in the case $E_{0203}$.

The Type-II range for $\boldsymbol{E}_{0502}$ is
\begin{equation}
\left[\varepsilon,\ \frac{1}{6}(5 - 8 \theta_1 - 8 \theta_2) - \varepsilon \right] \cup \left[\theta_2 + \varepsilon,\ \frac{1}{6}(5 - 8 \theta_1 - 4 \theta_2) - \varepsilon \right] \cup \left[\theta_1 + \varepsilon,\ \frac{1}{2}(1 - \theta_2) - \varepsilon \right].
\end{equation}
The decompositions are similar to which in the case $\boldsymbol{E}_{0210}$.

The Type-II range for $\boldsymbol{E}_{0503}$ is
\begin{equation}
\left[\varepsilon,\ \frac{1}{5}(4 - 6 \theta_1 - 8 \theta_2) - \varepsilon \right] \cup \left[2 \theta - 1 + \varepsilon,\ \frac{1}{6}(5 - 8 \theta) - \varepsilon \right] \cup \left[\theta_2 + \varepsilon,\ \frac{1}{6}(5 - 8 \theta_1 - 4 \theta_2) - \varepsilon \right] \cup \left[\theta_1 + \varepsilon,\ \frac{1}{2}(1 - \theta_2) - \varepsilon \right].
\end{equation}
The decompositions are similar to which in the case $\boldsymbol{E}_{0303}$.

The Type-II range for $\boldsymbol{E}_{0504}$ is
\begin{equation}
\left[\varepsilon,\ \frac{1}{5}(4 - 6 \theta_1 - 8 \theta_2) - \varepsilon \right] \cup \left[\theta_2 + \varepsilon,\ \frac{1}{6}(5 - 8 \theta_1 - 4 \theta_2) - \varepsilon \right] \cup \left[\theta_1 + \varepsilon,\ \frac{1}{2}(1 - \theta_2) - \varepsilon \right].
\end{equation}
The decompositions are similar to which in the case $\boldsymbol{E}_{0305}$.

The Type-II range for $\boldsymbol{E}_{0505}$ is
\begin{equation}
\left[\varepsilon,\ \frac{1}{5}(4 - 6 \theta_1 - 8 \theta_2) - \varepsilon \right] \cup \left[2 \theta - 1 + \varepsilon,\ \frac{1}{6}(5 - 8 \theta) - \varepsilon \right] \cup \left[\theta_2 + \varepsilon,\ \frac{1}{2}(1 - 2 \theta_1 + \theta_2) - \varepsilon \right] \cup \left[\theta_1 + \varepsilon,\ \frac{1}{2}(1 - \theta_2) - \varepsilon \right].
\end{equation}
The decompositions are similar to which in the case $\boldsymbol{E}_{0303}$.

The Type-II range for $\boldsymbol{E}_{0506}$ is
\begin{equation}
\left[\varepsilon,\ \frac{1}{5}(4 - 6 \theta_1 - 8 \theta_2) - \varepsilon \right] \cup \left[\theta_2 + \varepsilon,\ \frac{1}{2}(1 - 2 \theta_1 + \theta_2) - \varepsilon \right] \cup \left[\theta_1 + \varepsilon,\ \frac{1}{2}(1 - \theta_2) - \varepsilon \right].
\end{equation}
The decompositions are similar to which in the case $\boldsymbol{E}_{0305}$.

The Type-II range for $\boldsymbol{E}_{0507}$ is
\begin{equation}
\left[\varepsilon,\ \frac{1}{5}(4 - 6 \theta_1 - 8 \theta_2) - \varepsilon \right] \cup \left[2 \theta - 1 + \varepsilon,\ \frac{1}{6}(5 - 8 \theta) - \varepsilon \right] \cup \left[\theta_2 + \varepsilon,\ \frac{1}{4}(2 - 3 \theta_1) - \varepsilon \right] \cup \left[\theta_1 + \varepsilon,\ \frac{1}{2}(1 - \theta_2) - \varepsilon \right].
\end{equation}
The decompositions are similar to which in the case $\boldsymbol{E}_{0304}$.

The Type-II range for $\boldsymbol{E}_{0508}$ is
\begin{equation}
\left[\varepsilon,\ \frac{1}{5}(4 - 6 \theta_1 - 8 \theta_2) - \varepsilon \right] \cup \left[\theta_2 + \varepsilon,\ \frac{1}{4}(2 - 3 \theta_1) - \varepsilon \right] \cup \left[\theta_1 + \varepsilon,\ \frac{1}{2}(1 - \theta_2) - \varepsilon \right].
\end{equation}
The decompositions are similar to which in the case $\boldsymbol{E}_{0305}$.

The Type-II range for $\boldsymbol{E}_{0509}$ is
\begin{equation}
\left[\varepsilon,\ \frac{1}{5}(4 - 6 \theta_1 - 8 \theta_2) - \varepsilon \right] \cup \left[\theta_2 + \varepsilon,\ \frac{1}{4}(2 - 3 \theta_1) - \varepsilon \right] \cup \left[\theta_1 + \varepsilon,\ \frac{1}{2}(2 - \theta_1 - 4 \theta_2) - \varepsilon \right].
\end{equation}
The decompositions are similar to which in the case $\boldsymbol{E}_{0307}$.

\subsubsection{$\boldsymbol{E}_{06}$}
We divide $\boldsymbol{E}_{06}$ into 2 subregions:
$$
\boldsymbol{E}_{06} = \boldsymbol{E}_{0601} \cup \boldsymbol{E}_{0602},
$$
where
\begin{align}
\nonumber \boldsymbol{E}_{0601} =&\ \left\{ (\theta_1, \theta_2) : \frac{3}{8} < \theta_1 \leqslant \frac{2}{5},\ \frac{1}{2}(1 - 2 \theta_1) < \theta_2 < \frac{1}{14}(5 - 8 \theta_1) \right\}, \\
\nonumber \boldsymbol{E}_{0602} =&\ \left\{ (\theta_1, \theta_2) : \frac{3}{8} < \theta_1 \leqslant \frac{2}{5},\ \frac{1}{14}(5 - 8 \theta_1) \leqslant \theta_2 < \frac{1}{3}(2 - 4 \theta_1) \right\}.
\end{align}
Note that we have $\theta < \frac{11}{20}$ for $(\theta_1, \theta_2) \in \boldsymbol{E}_{06}$.

The Type-II range for $\boldsymbol{E}_{0601}$ is
\begin{equation}
\left[\varepsilon,\ \frac{1}{6}(5 - 8 \theta_1 - 4 \theta_2) - \varepsilon \right] \cup \left[\theta_1 + \varepsilon,\ \frac{1}{2}(1 - \theta_2) - \varepsilon \right].
\end{equation}
The decompositions are similar to which in the non-asymptotic parts in the case $\boldsymbol{\mathcal{Z}_3}$.

The Type-II range for $\boldsymbol{E}_{0602}$ is
\begin{equation}
\left[\varepsilon,\ \frac{1}{6}(5 - 8 \theta_1 - 8 \theta_2) - \varepsilon \right] \cup \left[\theta_2 + \varepsilon,\ \frac{1}{6}(5 - 8 \theta_1 - 4 \theta_2) - \varepsilon \right] \cup \left[\theta_1 + \varepsilon,\ \frac{1}{2}(1 - \theta_2) - \varepsilon \right].
\end{equation}
The decompositions are similar to which in the case $\boldsymbol{E}_{0210}$.

\subsubsection{$\boldsymbol{E}_{07}$}
We divide $\boldsymbol{E}_{07}$ into 6 subregions:
$$
\boldsymbol{E}_{07} = \boldsymbol{E}_{0701} \cup \boldsymbol{E}_{0702} \cup \boldsymbol{E}_{0703} \cup \boldsymbol{E}_{0704} \cup \boldsymbol{E}_{0705} \cup \boldsymbol{E}_{0706},
$$
where
\begin{align}
\nonumber \boldsymbol{E}_{0701} =&\ \left\{ (\theta_1, \theta_2) : \frac{3}{8} < \theta_1 \leqslant \frac{2}{5},\ \frac{1}{3}(2 - 4 \theta_1) \leqslant \theta_2 < \frac{1}{18}(9 - 16 \theta_1) \right\}, \\
\nonumber \boldsymbol{E}_{0702} =&\ \left\{ (\theta_1, \theta_2) : \frac{3}{8} < \theta_1 \leqslant \frac{2}{5},\ \frac{1}{18}(9 - 16 \theta_1) \leqslant \theta_2 < \frac{1}{20}(11 - 20 \theta_1) \right\}, \\
\nonumber \boldsymbol{E}_{0703} =&\ \left\{ (\theta_1, \theta_2) : \frac{3}{8} < \theta_1 \leqslant \frac{2}{5},\ \frac{1}{20}(11 - 20 \theta_1) \leqslant \theta_2 \leqslant \frac{1}{7}(2 - 2 \theta_1) \right\}, \\
\nonumber \boldsymbol{E}_{0704} =&\ \left\{ (\theta_1, \theta_2) : \frac{3}{8} < \theta_1 < \frac{2}{5},\ \frac{1}{7}(2 - 2 \theta_1) < \theta_2 < \frac{1}{2} \theta_1 \right\}, \\
\nonumber \boldsymbol{E}_{0705} =&\ \left\{ (\theta_1, \theta_2) : \frac{3}{8} < \theta_1 < \frac{2}{5},\ \frac{1}{2} \theta_1 \leqslant \theta_2 < \frac{1}{3}(1 - \theta_1) \right\}, \\
\nonumber \boldsymbol{E}_{0706} =&\ \left\{ (\theta_1, \theta_2) : \frac{3}{8} < \theta_1 < \frac{2}{5},\ \frac{1}{3}(1 - \theta_1) \leqslant \theta_2 < \frac{1}{4}(2 - 3 \theta_1) \right\}.
\end{align}
Note that we have $\theta < \frac{11}{20}$ for $(\theta_1, \theta_2) \in \boldsymbol{E}_{0701} \cup \boldsymbol{E}_{0702}$.

The Type-II range for $\boldsymbol{E}_{0701}$ is
\begin{equation}
\left[\varepsilon,\ \frac{1}{6}(5 - 8 \theta_1 - 8 \theta_2) - \varepsilon \right] \cup \left[\theta_2 + \varepsilon,\ \frac{1}{6}(5 - 8 \theta_1 - 4 \theta_2) - \varepsilon \right] \cup \left[\theta_1 + \varepsilon,\ \frac{1}{2}(1 - \theta_2) - \varepsilon \right].
\end{equation}
The decompositions are similar to which in the case $\boldsymbol{E}_{0210}$.

The Type-II range for $\boldsymbol{E}_{0702}$ is
\begin{equation}
\left[\varepsilon,\ \frac{1}{5}(4 - 6 \theta_1 - 8 \theta_2) - \varepsilon \right] \cup \left[2 \theta - 1 + \varepsilon,\ \frac{1}{6}(5 - 8 \theta) - \varepsilon \right] \cup \left[\theta_2 + \varepsilon,\ \frac{1}{6}(5 - 8 \theta_1 - 4 \theta_2) - \varepsilon \right] \cup \left[\theta_1 + \varepsilon,\ \frac{1}{2}(1 - \theta_2) - \varepsilon \right].
\end{equation}
The decompositions are similar to which in the case $\boldsymbol{E}_{0303}$.

The Type-II range for $\boldsymbol{E}_{0703}$ is
\begin{equation}
\left[\varepsilon,\ \frac{1}{5}(4 - 6 \theta_1 - 8 \theta_2) - \varepsilon \right] \cup \left[\theta_2 + \varepsilon,\ \frac{1}{6}(5 - 8 \theta_1 - 4 \theta_2) - \varepsilon \right] \cup \left[\theta_1 + \varepsilon,\ \frac{1}{2}(1 - \theta_2) - \varepsilon \right].
\end{equation}
The decompositions are similar to which in the case $\boldsymbol{E}_{0504}$.

The Type-II range for $\boldsymbol{E}_{0704}$ is
\begin{equation}
\left[\varepsilon,\ \frac{1}{5}(4 - 6 \theta_1 - 8 \theta_2) - \varepsilon \right] \cup \left[\theta_2 + \varepsilon,\ \frac{1}{2}(1 - 2 \theta_1 + \theta_2) - \varepsilon \right] \cup \left[\theta_1 + \varepsilon,\ \frac{1}{2}(1 - \theta_2) - \varepsilon \right].
\end{equation}
The decompositions are similar to which in the case $\boldsymbol{E}_{0506}$.

The Type-II range for $\boldsymbol{E}_{0705}$ is
\begin{equation}
\left[\varepsilon,\ \frac{1}{5}(4 - 6 \theta_1 - 8 \theta_2) - \varepsilon \right] \cup \left[\theta_2 + \varepsilon,\ \frac{1}{4}(2 - 3 \theta_1) - \varepsilon \right] \cup \left[\theta_1 + \varepsilon,\ \frac{1}{2}(1 - \theta_2) - \varepsilon \right].
\end{equation}
The decompositions are similar to which in the case $\boldsymbol{E}_{0305}$.

The Type-II range for $\boldsymbol{E}_{0706}$ is
\begin{equation}
\left[\varepsilon,\ \frac{1}{5}(4 - 6 \theta_1 - 8 \theta_2) - \varepsilon \right] \cup \left[\theta_2 + \varepsilon,\ \frac{1}{4}(2 - 3 \theta_1) - \varepsilon \right] \cup \left[\theta_1 + \varepsilon,\ \frac{1}{2}(2 - \theta_1 - 4 \theta_2) - \varepsilon \right].
\end{equation}
The decompositions are similar to which in the case $\boldsymbol{E}_{0307}$.

\subsubsection{$\boldsymbol{E}_{08}$}
By Statement (3) of Theorem~\ref{t47}, we know that
$$
C_1^{\prime}(\theta_1, \theta_2) = C_0^{\prime}(\theta_1, \theta_2) = 1
$$
for $(\theta_1, \theta_2) \in \boldsymbol{E}_{08}$ if $7 \theta_1 + 12 \theta_2 \leqslant 4 - \varepsilon$. We divide the remaining parts of $\boldsymbol{E}_{08}$ into 2 subregions:
\begin{align}
\nonumber \boldsymbol{E}_{0801} =&\ \left\{ (\theta_1, \theta_2) : \frac{2}{5} < \theta_1 \leqslant \frac{13}{32},\ \frac{1}{12}(4 - 7 \theta_1) \leqslant \theta_2 < \frac{1}{14}(5 - 8 \theta_1) \right. \\
\nonumber & \left. \qquad \qquad \quad \text{ or } \frac{13}{32} < \theta_1 < \frac{4}{9},\ \frac{1}{12}(4 - 7 \theta_1) \leqslant \theta_2 < \frac{1}{3}(2 - 4 \theta_1) \right\}, \\
\nonumber \boldsymbol{E}_{0802} =&\ \left\{ (\theta_1, \theta_2) : \frac{2}{5} < \theta_1 < \frac{13}{32},\ \frac{1}{14}(5 - 8 \theta_1) \leqslant \theta_2 < \frac{1}{3}(2 - 4 \theta_1) \right\}.
\end{align}
Note that we have $\theta < \frac{11}{20}$ for $(\theta_1, \theta_2) \in \boldsymbol{E}_{0801} \cup \boldsymbol{E}_{0802}$.

The Type-II range for $\boldsymbol{E}_{0801}$ is
\begin{equation}
\left[\varepsilon,\ \frac{1}{6}(5 - 8 \theta_1 - 4 \theta_2) - \varepsilon \right] \cup \left[\theta_1 + \varepsilon,\ \frac{1}{2}(1 - \theta_2) - \varepsilon \right].
\end{equation}
The decompositions are similar to which in the non-asymptotic parts in the case $\boldsymbol{\mathcal{Z}_3}$.

The Type-II range for $\boldsymbol{E}_{0802}$ is
\begin{equation}
\left[\varepsilon,\ \frac{1}{6}(5 - 8 \theta_1 - 8 \theta_2) - \varepsilon \right] \cup \left[\theta_2 + \varepsilon,\ \frac{1}{6}(5 - 8 \theta_1 - 4 \theta_2) - \varepsilon \right] \cup \left[\theta_1 + \varepsilon,\ \frac{1}{2}(1 - \theta_2) - \varepsilon \right].
\end{equation}
The decompositions are similar to which in the case $\boldsymbol{E}_{0210}$.

\subsubsection{$\boldsymbol{E}_{09}$}
By Statement (3) of Theorem~\ref{t47}, we know that
$$
C_1^{\prime}(\theta_1, \theta_2) = C_0^{\prime}(\theta_1, \theta_2) = 1
$$
for $(\theta_1, \theta_2) \in \boldsymbol{E}_{09}$ if $\theta_2 \leqslant \min\left(\frac{4 - 7 \theta_1}{12}, \frac{10 - 19 \theta_1}{20}\right) - \varepsilon$. We divide the remaining parts of $\boldsymbol{E}_{09}$ into 5 subregions:
\begin{align}
\nonumber \boldsymbol{E}_{0901} =&\ \left\{ (\theta_1, \theta_2) : \frac{13}{32} < \theta_1 \leqslant \frac{4}{9},\ \frac{1}{3}(2 - 4 \theta_1) \leqslant \theta_2 < \frac{1}{14}(5 - 8 \theta_1) \right. \\
\nonumber & \qquad \qquad \quad \text{ or } \frac{4}{9} < \theta_1 \leqslant \frac{9}{20},\ \frac{1}{12}(4 - 7 \theta_1) \leqslant \theta_2 < \frac{1}{14}(5 - 8 \theta_1) \\
\nonumber & \qquad \qquad \quad \text{ or } \frac{9}{20} < \theta_1 \leqslant \frac{5}{11},\ \frac{1}{12}(4 - 7 \theta_1) \leqslant \theta_2 < 1 - 2 \theta_1 \\
\nonumber & \left. \qquad \qquad \quad \text{ or } \frac{5}{11} < \theta_1 < \frac{10}{21},\ \frac{1}{20}(10 - 19 \theta_1) \leqslant \theta_2 < 1 - 2 \theta_1 \right\}, \\
\nonumber \boldsymbol{E}_{0902} =&\ \left\{ (\theta_1, \theta_2) : \frac{2}{5} < \theta_1 \leqslant \frac{13}{32},\ \frac{1}{3}(2 - 4 \theta_1) \leqslant \theta_2 < \frac{1}{18}(9 - 16 \theta_1) \right. \\
\nonumber & \left. \qquad \qquad \quad \text{ or } \frac{13}{32} < \theta_1 < \frac{9}{20},\ \frac{1}{14}(5 - 8 \theta_1) \leqslant \theta_2 < \frac{1}{18}(9 - 16 \theta_1) \right\}, \\
\nonumber \boldsymbol{E}_{0903} =&\ \left\{ (\theta_1, \theta_2) : \frac{2}{5} < \theta_1 < \frac{9}{20},\ \frac{1}{18}(9 - 16 \theta_1) \leqslant \theta_2 < \frac{1}{20}(11 - 20 \theta_1) \right\}, \\
\nonumber \boldsymbol{E}_{0904} =&\ \left\{ (\theta_1, \theta_2) : \frac{2}{5} < \theta_1 < \frac{5}{12},\ \frac{1}{20}(11 - 20 \theta_1) \leqslant \theta_2 \leqslant \frac{1}{7}(2 - 2 \theta_1) \right. \\
\nonumber & \left. \qquad \qquad \quad \text{ or } \frac{5}{12} \leqslant \theta_1 < \frac{9}{20},\ \frac{1}{20}(11 - 20 \theta_1) \leqslant \theta_2 < 1 - 2 \theta_1 \right\}, \\
\nonumber \boldsymbol{E}_{0905} =&\ \left\{ (\theta_1, \theta_2) : \frac{2}{5} < \theta_1 < \frac{5}{12},\ \frac{1}{7}(2 - 2 \theta_1) < \theta_2 < 1 - 2 \theta_1 \right\}.
\end{align}
Note that we have $\theta < \frac{11}{20}$ for $(\theta_1, \theta_2) \in \boldsymbol{E}_{0901} \cup \boldsymbol{E}_{0902} \cup \boldsymbol{E}_{0903}$.

The Type-II range for $\boldsymbol{E}_{0901}$ is
\begin{equation}
\left[\varepsilon,\ \frac{1}{6}(5 - 8 \theta_1 - 4 \theta_2) - \varepsilon \right] \cup \left[\theta_1 + \varepsilon,\ \frac{1}{2}(1 - \theta_2) - \varepsilon \right].
\end{equation}

First, we shall prove the following theorem, which, combined with [\cite{MaynardLargeModuliI}, Theorem 1.1], implies Theorem~\ref{t11}.
\begin{theorem}\label{case4}
Define
\begin{align}
\nonumber \boldsymbol{\mathcal{Z}_4} =&\ \left\{ (\theta_1, \theta_2) : \frac{5}{11} < \theta_1 \leqslant \frac{8}{17},\ \frac{1}{20}(10 - 19 \theta_1) \leqslant \theta_2 \leqslant \frac{1}{12}(4 - 7 \theta_1) - \varepsilon \right. \\
\nonumber & \left. \qquad \qquad \quad \text{ or } \frac{8}{17} < \theta_1 < \frac{10}{21},\ \frac{1}{20}(10 - 19 \theta_1) \leqslant \theta_2 \leqslant 1 - 2 \theta_1 - 2 \varepsilon \right\}.
\end{align}
Let $(\theta_1, \theta_2) \in \boldsymbol{\mathcal{Z}_4}$. Then (129) holds for
$$
f(n) = \mathbbm{1}_{p}(n),
$$
and we have
$$
C_1^{\prime}(\theta_1, \theta_2) = C_0^{\prime}(\theta_1, \theta_2) = 1.
$$
\end{theorem}
\begin{proof}
Since we have $(\theta_1, \theta_2) \in \boldsymbol{\mathcal{Z}_4}$, all of the following conditions hold true:
$$
\kappa = \frac{5 - 8 \theta_1 - 4 \theta_2}{6} - \varepsilon > \frac{1}{6}, \quad \frac{1}{2} < \theta \leqslant \frac{9}{17} - \varepsilon, \quad 11 \theta_1 + 12 \theta_2 < 6, \quad 3 \theta_1 + 2 \theta_2 < \frac{11}{7}.
$$

Now, we can decompose our $\mathbbm{1}_{n \sim x, n = p}(n) = \psi\left(n, (2 x)^{\frac{1}{2}}\right)$ in a way similar to the decompositions in the proof of Theorem~\ref{case1}. By the discussions of the three-dimensional sieves in this section, we need to apply the new three-dimensional Harman's sieve on both $B$ and $C$. By Buchstab's identity, we have
\begin{align}
\nonumber \psi\left(n, (2 x)^{\frac{1}{2}}\right) =&\ \psi\left(n, x^{\kappa}\right) - \sum_{\substack{n = p_1 \beta \\ \kappa \leqslant \alpha_1 < \frac{3}{7} + \varepsilon }} \psi\left(\beta, x^{\kappa} \right) + \sum_{\substack{n = p_1 p_2 \beta \\ \kappa \leqslant \alpha_1 < \frac{3}{7} + \varepsilon \\ \kappa \leqslant \alpha_2 < \min\left(\alpha_1, \frac{1}{2}(1 - \alpha_1) \right) \\ \boldsymbol{\alpha}_2 \in \boldsymbol{G}_2 }} \psi\left(\beta, p_2 \right) \\
\nonumber & + \sum_{\substack{n = p_1 p_2 \beta \\ \kappa \leqslant \alpha_1 < \frac{3}{7} + \varepsilon \\ \kappa \leqslant \alpha_2 < \min\left(\alpha_1, \frac{1}{2}(1 - \alpha_1) \right) \\ \boldsymbol{\alpha}_2 \in A }} \psi\left(\beta, x^{\kappa} \right) - \sum_{\substack{n = p_1 p_2 p_3 \beta \\ \kappa \leqslant \alpha_1 < \frac{3}{7} + \varepsilon \\ \kappa \leqslant \alpha_2 < \min\left(\alpha_1, \frac{1}{2}(1 - \alpha_1) \right) \\ \boldsymbol{\alpha}_2 \in A \cup B \cup C \\ \kappa \leqslant \alpha_3 < \min\left(\alpha_2, \frac{1}{2}(1 - \alpha_1 - \alpha_2) \right) }} \psi\left(\beta, p_3 \right) \\
\nonumber & + \sum_{\substack{n = p_1 p_2 \beta \\ \kappa \leqslant \alpha_1 < \frac{3}{7} + \varepsilon \\ \kappa \leqslant \alpha_2 < \min\left(\alpha_1, \frac{1}{2}(1 - \alpha_1) \right) \\ \boldsymbol{\alpha}_2 \in B \cup C }} \psi\left(\beta, x^{\kappa} \right) - \sum_{\substack{n = p_1 \beta \\ \frac{3}{7} + \varepsilon \leqslant \alpha_1 \leqslant 1 - \theta - \varepsilon }} \psi\left(\beta, p_1 \right) - \sum_{\substack{n = p_1 \beta \\ 1 - \theta - \varepsilon < \alpha_1 < \frac{1}{2} }} \psi\left(\beta, p_1 \right) \\
=&\ \Sigma_{42901} - \Sigma_{42902} + \Sigma_{42903} + \Sigma_{42904} - \Sigma_{42905} + \Sigma_{42906} - \Sigma_{42907} - \Sigma_{42908}.
\end{align}
By Lemma~\ref{l212} and a Type-II range $\left[\varepsilon,\ \frac{5 - 8 \theta_1 - 4 \theta_2}{6} - \varepsilon \right]$, (129) holds for $f(n) = \Sigma_{42901}$ and $f(n) = \Sigma_{42902}$. By Lemma~\ref{l211} and Lemma~\ref{l216}, (129) holds for $f(n) = \Sigma_{42903}$ and $f(n) = \Sigma_{42907}$. By an application of Lemma~\ref{l35}, (129) holds for $f(n) = \Sigma_{42908}$. By the Type-II range $\left[\varepsilon,\ \frac{5 - 8 \theta_1 - 4 \theta_2}{6} - \varepsilon \right]$ and the discussions in the end of the three-dimensional sieves (133)--(138) in this section, (129) holds for $f(n) = \Sigma_{42904}$. For the remaining sums, $\Sigma_{42905}$ only counts numbers with $4$ or more prime factors.

For $\Sigma_{42906}$, since we have $3 \theta_1 + 2 \theta_2 < \frac{11}{7}$, $11 \theta_1 + 12 \theta_2 < 6$, $\theta < \frac{9}{17} < \frac{8}{15}$ and a Type-II range $\left[\varepsilon,\ \frac{5 - 8 \theta_1 - 4 \theta_2}{6} - \varepsilon \right]$, we can use Lemma~\ref{l38} and a three-dimensional Harman's sieve to get a ``loss term''
\begin{equation}
\Sigma_{42909} = \sum_{\substack{n = m_1 m_2 m_3 \\ \kappa \leqslant \alpha_1 < \frac{3}{7} + \varepsilon \\ \kappa \leqslant \alpha_2 < \min\left(\alpha_1, \frac{1}{2}(1 - \alpha_1) \right) \\ \boldsymbol{\alpha}_2 \in B \cup C \\ \Omega(m_1 m_2) \geqslant 3 }} \psi\left(m_1 m_2 m_3, x^{\kappa} \right).
\end{equation}
Since $\Omega(m_1 m_2 m_3) \geqslant \Omega(m_1 m_2) + 1 \geqslant 4$, $\Sigma_{42909}$ only counts numbers with $4$ or more prime factors.

Now, the proof of Theorem~\ref{case4} reduces to showing that (129) holds for $f(n) = $ sums that count numbers with $4$ or more prime factors. The proof of Theorem~\ref{case4} is thus completed by applying Lemma~\ref{l39} on those sums, since we have $\kappa > \frac{1}{6} > \frac{1}{7}$.
\end{proof}

For the remaining parts of $\boldsymbol{E}_{0901}$, The decompositions are similar to which in the non-asymptotic parts in the case $\boldsymbol{\mathcal{Z}_3}$.

The Type-II range for $\boldsymbol{E}_{0902}$ is
\begin{equation}
\left[\varepsilon,\ \frac{1}{6}(5 - 8 \theta_1 - 8 \theta_2) - \varepsilon \right] \cup \left[\theta_2 + \varepsilon,\ \frac{1}{6}(5 - 8 \theta_1 - 4 \theta_2) - \varepsilon \right] \cup \left[\theta_1 + \varepsilon,\ \frac{1}{2}(1 - \theta_2) - \varepsilon \right].
\end{equation}
The decompositions are similar to which in the case $\boldsymbol{E}_{0210}$.

The Type-II range for $\boldsymbol{E}_{0903}$ is
\begin{equation}
\left[\varepsilon,\ \frac{1}{5}(4 - 6 \theta_1 - 8 \theta_2) - \varepsilon \right] \cup \left[2 \theta - 1 + \varepsilon,\ \frac{1}{6}(5 - 8 \theta) - \varepsilon \right] \cup \left[\theta_2 + \varepsilon,\ \frac{1}{6}(5 - 8 \theta_1 - 4 \theta_2) - \varepsilon \right] \cup \left[\theta_1 + \varepsilon,\ \frac{1}{2}(1 - \theta_2) - \varepsilon \right].
\end{equation}
The decompositions are similar to which in the case $\boldsymbol{E}_{0303}$.

The Type-II range for $\boldsymbol{E}_{0904}$ is
\begin{equation}
\left[\varepsilon,\ \frac{1}{5}(4 - 6 \theta_1 - 8 \theta_2) - \varepsilon \right] \cup \left[\theta_2 + \varepsilon,\ \frac{1}{6}(5 - 8 \theta_1 - 4 \theta_2) - \varepsilon \right] \cup \left[\theta_1 + \varepsilon,\ \frac{1}{2}(1 - \theta_2) - \varepsilon \right].
\end{equation}
The decompositions are similar to which in the case $\boldsymbol{E}_{0504}$.

The Type-II range for $\boldsymbol{E}_{0905}$ is
\begin{equation}
\left[\varepsilon,\ \frac{1}{5}(4 - 6 \theta_1 - 8 \theta_2) - \varepsilon \right] \cup \left[\theta_2 + \varepsilon,\ \frac{1}{2}(1 - 2 \theta_1 + \theta_2) - \varepsilon \right] \cup \left[\theta_1 + \varepsilon,\ \frac{1}{2}(1 - \theta_2) - \varepsilon \right].
\end{equation}
The decompositions are similar to which in the case $\boldsymbol{E}_{0506}$.

\subsubsection{$\boldsymbol{E}_{10}$}
From here, we have $2 \theta_1 + \theta_2 \geqslant 1$ and Lemmas~\ref{l44}--\ref{l45} become trivial. We divide $\boldsymbol{E}_{10}$ into 3 subregions:
$$
\boldsymbol{E}_{10} = \boldsymbol{E}_{1001} \cup \boldsymbol{E}_{1002} \cup \boldsymbol{E}_{1003},
$$
where
\begin{align}
\nonumber \boldsymbol{E}_{1001} =&\ \left\{ (\theta_1, \theta_2) : \frac{5}{11} < \theta_1 < \frac{1}{2},\ 1 - 2 \theta_1 \leqslant \theta_2 < \frac{1}{15}(5 - 8 \theta_1) \right\}, \\
\nonumber \boldsymbol{E}_{1002} =&\ \left\{ (\theta_1, \theta_2) : \frac{5}{12} < \theta_1 \leqslant \frac{7}{16},\ 1 - 2 \theta_1 \leqslant \theta_2 < \frac{1}{10}(5 - 8 \theta_1) \right. \\
\nonumber & \qquad \qquad \quad \text{ or } \frac{7}{16} < \theta_1 < \frac{5}{11},\ 1 - 2 \theta_1 \leqslant \theta_2 < \frac{1}{15}(11 - 20 \theta_1) \\
\nonumber & \left. \qquad \qquad \quad \text{ or } \frac{5}{11} \leqslant \theta_1 < \frac{1}{2},\ \frac{1}{15}(5 - 8 \theta_1) \leqslant \theta_2 < \frac{1}{15}(11 - 20 \theta_1) \right\}, \\
\nonumber \boldsymbol{E}_{1003} =&\ \left\{ (\theta_1, \theta_2) : \frac{2}{5} < \theta_1 \leqslant \frac{5}{12},\ 1 - 2 \theta_1 \leqslant \theta_2 < \frac{1}{15}(11 - 20 \theta_1) \right. \\
\nonumber & \left. \qquad \qquad \quad \text{ or } \frac{5}{12} < \theta_1 < \frac{7}{16},\ \frac{1}{10}(5 - 8 \theta_1) \leqslant \theta_2 < \frac{1}{15}(11 - 20 \theta_1) \right\}.
\end{align}
Note that we have $\theta > \frac{17}{32}$ for $(\theta_1, \theta_2) \in \boldsymbol{E}_{1002} \cup \boldsymbol{E}_{1003}$.

The Type-II range for $\boldsymbol{E}_{1001}$ is
\begin{equation}
\left[2 \theta_1 + \theta_2 - 1 + \varepsilon,\ \frac{1}{6}(5 - 8 \theta_1 - 4 \theta_2) - \varepsilon \right].
\end{equation}
We first assume that $\theta < \frac{17}{32}$. Similar to the case $(\theta_1, \theta_2) \in \boldsymbol{A}_{1101}$, we want to replace $\boldsymbol{U}_{j}$ with $\boldsymbol{U}_{j}^{\prime \prime}$ in the decompositions. Let $\boldsymbol{\alpha}_{k} \in \boldsymbol{A}_{k}$. By a decomposing process similar to (107), we only need to prove that if
$$
\alpha_{j+1} + \cdots + \alpha_{k-1} < 2 \theta_1 + \theta_2 - 1 + \varepsilon \leqslant \alpha_{j+1} + \cdots + \alpha_{k}
$$
and
$$
\alpha_{j+1} < \frac{5 - 8 \theta_1 - 4 \theta_2}{6} - \varepsilon
$$
hold for some $j$ ($0 \leqslant j \leqslant k - 1$), then $\boldsymbol{\alpha}_{k} \in \boldsymbol{G}_{k}$. 

When $\alpha_{k} < \frac{11 - 20 \theta_1 - 10 \theta_2}{6} - 2 \varepsilon$, then
$$
2 \theta_1 + \theta_2 - 1 + \varepsilon \leqslant \alpha_{j+1} + \cdots + \alpha_{k} < (2 \theta_1 + \theta_2 - 1 + \varepsilon) + \frac{11 - 20 \theta_1 - 10 \theta_2}{6} - 2 \varepsilon = \frac{5 - 8 \theta_1 - 4 \theta_2}{6} - \varepsilon
$$
and $\boldsymbol{\alpha}_{k} \in \boldsymbol{G}_{k}$.

Suppose that $\alpha_{k} \geqslant \frac{11 - 20 \theta_1 - 10 \theta_2}{6} - 2 \varepsilon$. Since $\boldsymbol{\alpha}_{k} \in \boldsymbol{A}_{k}$, we have $\alpha_{k} < \alpha_{j+1} < \frac{5 - 8 \theta_1 - 4 \theta_2}{6} - \varepsilon$. Now we only need to prove that
$$
\frac{11 - 20 \theta_1 - 10 \theta_2}{6} - 2 \varepsilon \geqslant 2 \theta_1 + \theta_2 - 1 + \varepsilon,
$$
or
$$
32 \theta_1 + 16 \theta_2 \leqslant 17 - 18 \varepsilon
$$
when $\theta < \frac{17}{32}$ and $(\theta_1, \theta_2) \in \boldsymbol{E}_{1001}$. A simple verification then completes the proof. The remaining decompositions are similar to which in the case $\boldsymbol{A}_{1101}$ in Section 3.

Now we assume that $\theta \geqslant \frac{17}{32}$. In this case we use the Type-II range $\left[2 \theta_1 + \theta_2 - 1 + \varepsilon,\ \frac{1}{6}(5 - 8 \theta_1 - 4 \theta_2) - \varepsilon \right]$ to remove sums explicitly.

The Type-II range for $\boldsymbol{E}_{1002}$ is
\begin{equation}
\left[2 \theta_1 + \theta_2 - 1 + \varepsilon,\ \frac{1}{6}(5 - 8 \theta_1 - 9 \theta_2) - \varepsilon \right] \cup \left[\theta_2 + \varepsilon,\ \frac{1}{6}(5 - 8 \theta_1 - 4 \theta_2) - \varepsilon \right].
\end{equation}
The decompositions are similar to which in the case $\boldsymbol{E}_{1001}$, where the Type-II range $\left[\theta_2 + \varepsilon,\ \frac{1}{6}(5 - 8 \theta_1 - 4 \theta_2) - \varepsilon \right]$ is used to remove sums explicitly.

The Type-II range for $\boldsymbol{E}_{1003}$ is
\begin{equation}
\left[2 \theta_1 + \theta_2 - 1 + \varepsilon,\ \frac{1}{6}(5 - 8 \theta_1 - 9 \theta_2) - \varepsilon \right].
\end{equation}
Since we have $\theta > 0.57$ in this region, we do not discuss any further decompositions here.

\subsubsection{$\boldsymbol{E}_{11}$}
Since we have $\frac{1}{15}(11 - 20 \theta_1) < \frac{1}{10}(11 - 20 \theta_1)$ when $\theta_1 < \frac{11}{20}$, by Lemma~\ref{l43} and the definition of the region $\boldsymbol{F}_{1}$, we can extend the region $\boldsymbol{E}_{11}$ to a larger one:
$$
\boldsymbol{E}_{11}^{\prime} = \left\{ (\theta_1, \theta_2) : \frac{1}{2} \leqslant \theta_1 < \frac{11}{20},\ 0 < \theta_2 < \frac{1}{10}(11 - 20 \theta_1) \right\}.
$$
The Type-II range for $\boldsymbol{E}_{11}^{\prime}$ is
\begin{equation}
\left[2 \theta_1 + \theta_2 - 1 + \varepsilon,\ \frac{1}{6}(5 - 8 \theta_1 - 4 \theta_2) - \varepsilon \right].
\end{equation}
The decompositions are similar to which in the case $\boldsymbol{E}_{1001}$. When $\theta < \frac{17}{32}$, we can replace $\boldsymbol{U}_{j}$ with $\boldsymbol{U}_{j}^{\prime \prime}$ in the decompositions. When $\theta \geqslant \frac{17}{32}$, we use the Type-II range $\left[2 \theta_1 + \theta_2 - 1 + \varepsilon,\ \frac{1}{6}(5 - 8 \theta_1 - 4 \theta_2) - \varepsilon \right]$ to remove sums explicitly.

The decompositions in other parts of $\boldsymbol{U} \backslash \boldsymbol{J}$ stay the same as in Section 3. Working on each case above, we can get the following upper bounds for $C_1^{\prime}(\theta_1, \theta_2)$ ($0.5 < \theta \leqslant 0.56$):
\begin{center}
\begin{tabular}{|>{\centering\arraybackslash}p{0.6cm}|>{\centering\arraybackslash}p{0.6cm}|>{\centering\arraybackslash}p{0.6cm}|>{\centering\arraybackslash}p{0.6cm}|>{\centering\arraybackslash}p{0.6cm}|>{\centering\arraybackslash}p{0.6cm}|>{\centering\arraybackslash}p{0.6cm}|>{\centering\arraybackslash}p{0.6cm}|>{\centering\arraybackslash}p{0.6cm}|>{\centering\arraybackslash}p{0.6cm}|>{\centering\arraybackslash}p{0.6cm}|>{\centering\arraybackslash}p{0.6cm}|>{\centering\arraybackslash}p{0.6cm}|>{\centering\arraybackslash}p{0.6cm}|>{\centering\arraybackslash}p{0.6cm}|>{\centering\arraybackslash}p{0.6cm}|}
\hline \boldmath{$0.28$} & \tiny $1.5641$ & \tiny $1.9602$ & \tiny $2.3301$ & --- & --- & --- & --- & --- & --- & --- & --- & --- & --- & --- & --- \\
\hline \boldmath{$0.27$} & \tiny $1.4833$ & \tiny $1.6632$ & \tiny $1.9602$ & \tiny $2.1636$ & --- & --- & --- & --- & --- & --- & --- & --- & --- & --- & --- \\
\hline \boldmath{$0.26$} & \tiny $1.0830$ & \tiny $1.4833$ & \tiny $1.5641$ & \tiny $1.7620$ & \tiny $2.1185$ & --- & --- & --- & --- & --- & --- & --- & --- & --- & --- \\
\hline \boldmath{$0.25$} & \tiny $1.0401$ & \tiny $1.0805$ & \tiny $1.4115$ & \tiny $1.4795$ & \tiny $1.8120$ & \tiny $2.0732$ & --- & --- & --- & --- & --- & --- & --- & --- & --- \\
\hline \boldmath{$0.24$} & $1$ & $1$ & \tiny $1.0801$ & \tiny $1.3433$ & \tiny $1.5084$ & \tiny $1.7924$ & \tiny $2.0309$ & --- & --- & --- & --- & --- & --- & --- & --- \\
\hline \boldmath{$0.23$} & $1$ & $1$ & $1$ & \tiny $1.0801$ & \tiny $1.3938$ & \tiny $1.4975$ & \tiny $1.7977$ & \tiny $2.0268$ & --- & --- & --- & --- & --- & --- & --- \\
\hline \boldmath{$0.22$} & $1$ & $1$ & $1$ & $1$ & \tiny $1.0810$ & \tiny $1.4315$ & \tiny $1.5398$ & \tiny $1.8781$ & \tiny $2.0672$ & --- & --- & --- & --- & --- & --- \\
\hline \boldmath{$0.21$} & $1$ & $1$ & $1$ & $1$ & $1$ & \tiny $1.0817$ & \tiny $1.4644$ & \tiny $1.5605$ & \tiny $1.9836$ & \tiny $2.1666$ & --- & --- & --- & --- & --- \\
\hline \boldmath{$0.20$} & $1$ & $1$ & $1$ & $1$ & $1$ & $1$ & \tiny $1.0801$ & \tiny $1.4925$ & \tiny $1.5939$ & \tiny $2.0468$ & \tiny $2.2182$ & --- & --- & --- & --- \\
\hline \boldmath{$0.19$} & $1$ & $1$ & $1$ & $1$ & $1$ & $1$ & $1$ & \tiny $1.0801$ & \tiny $1.5168$ & \tiny $1.6144$ & \tiny $2.0702$ & \tiny $2.2490$ & --- & --- & --- \\
\hline \boldmath{$0.18$} & $1$ & $1$ & $1$ & $1$ & $1$ & $1$ & $1$ & $1$ & \tiny $1.0805$ & \tiny $1.5407$ & \tiny $1.6360$ & \tiny $2.0714$ & \tiny $2.2632$ & --- & --- \\
\hline \boldmath{$0.17$} & $1$ & $1$ & $1$ & $1$ & $1$ & $1$ & $1$ & $1$ & \tiny $1.0401$ & \tiny $1.0825$ & \tiny $1.5635$ & \tiny $1.6513$ & \tiny $2.1058$ & \tiny $2.2974$ & --- \\
\hline \boldmath{$0.16$} & $1$ & $1$ & $1$ & $1$ & $1$ & $1$ & $1$ & $1$ & $1$ & \tiny $1.0405$ & \tiny $1.0837$ & \tiny $1.1274$ & \tiny $1.6749$ & \tiny $2.1300$ & \tiny $2.3159$ \\
\hline \boldmath{$0.15$} & $1$ & $1$ & $1$ & $1$ & $1$ & $1$ & $1$ & $1$ & $1$ & $1$ & \tiny $1.0401$ & \tiny $1.0827$ & \tiny $1.1549$ & \tiny $1.6998$ & \tiny $2.1540$ \\
\hline \boldmath{$0.14$} & $1$ & $1$ & $1$ & $1$ & $1$ & $1$ & $1$ & $1$ & $1$ & $1$ & $1$ & \tiny $1.0401$ & \tiny $1.0810$ & \tiny $1.1368$ & \tiny $1.7389$ \\
\hline \boldmath{$0.13$} & $1$ & $1$ & $1$ & $1$ & $1$ & $1$ & $1$ & $1$ & $1$ & $1$ & $1$ & $1$ & \tiny $1.0407$ & \tiny $1.0811$ & \tiny $1.1274$ \\
\hline \boldmath{$0.12$} & $1$ & $1$ & $1$ & $1$ & $1$ & $1$ & $1$ & $1$ & $1$ & $1$ & $1$ & $1$ & $1 + \varepsilon$ & \tiny $1.0418$ & \tiny $1.0808$ \\
\hline \boldmath{$0.11$} & $1$ & $1$ & $1$ & $1$ & $1$ & $1$ & $1$ & $1$ & $1$ & $1$ & $1$ & $1$ & $1$ & $1 + \varepsilon$ & \tiny $1.0401$ \\
\hline \small \boldmath{$\theta_2 \backslash \theta_1$} & \boldmath{$0.26$} & \boldmath{$0.27$} & \boldmath{$0.28$} & \boldmath{$0.29$} & \boldmath{$0.30$} & \boldmath{$0.31$} & \boldmath{$0.32$} & \boldmath{$0.33$} & \boldmath{$0.34$} & \boldmath{$0.35$} & \boldmath{$0.36$} & \boldmath{$0.37$} & \boldmath{$0.38$} & \boldmath{$0.39$} & \boldmath{$0.40$} \\
\hline
\end{tabular} \\
\textbf{Table 4.1: Upper Bounds for }\boldmath{$C_1^{\prime}(\theta_1, \theta_2)$} \textbf{(}\boldmath{$0.5 < \theta \leqslant 0.56$}\textbf{) 1/2}
\end{center}
\begin{center}
\begin{tabular}{|>{\centering\arraybackslash}p{0.6cm}|>{\centering\arraybackslash}p{0.6cm}|>{\centering\arraybackslash}p{0.6cm}|>{\centering\arraybackslash}p{0.6cm}|>{\centering\arraybackslash}p{0.6cm}|>{\centering\arraybackslash}p{0.6cm}|>{\centering\arraybackslash}p{0.6cm}|>{\centering\arraybackslash}p{0.6cm}|>{\centering\arraybackslash}p{0.6cm}|>{\centering\arraybackslash}p{0.6cm}|>{\centering\arraybackslash}p{0.6cm}|>{\centering\arraybackslash}p{0.6cm}|>{\centering\arraybackslash}p{0.6cm}|>{\centering\arraybackslash}p{0.6cm}|>{\centering\arraybackslash}p{0.6cm}|>{\centering\arraybackslash}p{0.6cm}|}
\hline \boldmath{$0.15$} & \tiny $2.3347$ & --- & --- & --- & --- & --- & --- & --- & --- & --- & --- & --- & --- & --- & --- \\
\hline \boldmath{$0.14$} & \tiny $2.1645$ & \tiny $2.3628$ & --- & --- & --- & --- & --- & --- & --- & --- & --- & --- & --- & --- & --- \\
\hline \boldmath{$0.13$} & \tiny $1.7782$ & \tiny $2.2026$ & \tiny $2.3900$ & --- & --- & --- & --- & --- & --- & --- & --- & --- & --- & --- & --- \\
\hline \boldmath{$0.12$} & \tiny $1.1249$ & \tiny $1.8132$ & \tiny $2.2388$ & \tiny $2.3930$ & --- & --- & --- & --- & --- & --- & --- & --- & --- & --- & --- \\
\hline \boldmath{$0.11$} & \tiny $1.0807$ & \tiny $1.1266$ & \tiny $1.8455$ & \tiny $2.2375$ & \tiny $2.3891$ & --- & --- & --- & --- & --- & --- & --- & --- & --- & --- \\
\hline \boldmath{$0.10$} & \tiny $1.0408$ & \tiny $1.0812$ & \tiny $1.1320$ & \tiny $1.8464$ & \tiny $2.2577$ & \tiny $2.3885$ & --- & --- & --- & --- & --- & --- & --- & --- & --- \\
\hline \boldmath{$0.09$} & $1$ & \tiny $1.0401$ & \tiny $1.0828$ & \tiny $1.1386$ & \tiny $1.8463$ & \tiny $2.2082$ & \tiny $2.3888$ & --- & --- & --- & --- & --- & --- & --- & --- \\
\hline \boldmath{$0.08$} & $1$ & $1$ & $1$ & \tiny $1.0847$ & \tiny $1.1467$ & \tiny $1.8476$ & \tiny $2.2383$ & \tiny $2.3907$ & --- & --- & --- & --- & --- & --- & --- \\
\hline \boldmath{$0.07$} & $1$ & $1$ & $1$ & $1$ & $1$ & \tiny $1.1491$ & \tiny $1.8498$ & \tiny $2.2535$ & \tiny $2.4108$ & --- & --- & --- & --- & --- & --- \\
\hline \boldmath{$0.06$} & $1$ & $1$ & $1$ & $1$ & $1$ & $1$ & \tiny $1.7993$ & \tiny $1.8561$ & \tiny $2.2675$ & \tiny $2.4313$ & --- & --- & --- & --- & --- \\
\hline \boldmath{$0.05$} & $1$ & $1$ & $1$ & $1$ & $1$ & $1$ & $1$ & \tiny $1.7183$ & \tiny $1.8642$ & \tiny $2.2987$ & \tiny $2.4608$ & --- & --- & --- & --- \\
\hline \boldmath{$0.04$} & $1$ & $1$ & $1$ & $1$ & $1$ & $1$ & $1$ & $1 + \varepsilon$ & \tiny $1.7618$ & \tiny $1.8806$ & \tiny $2.3027$ & \tiny $2.4918$ & --- & --- & --- \\
\hline \boldmath{$0.03$} & $1$ & $1$ & $1$ & $1$ & $1$ & $1$ & $1$ & $1$ & \tiny $1.0452$ & \tiny $1.8057$ & \tiny $1.8963$ & \tiny $2.3271$ & \tiny $2.5324$ & --- & --- \\
\hline \boldmath{$0.02$} & $1$ & $1$ & $1$ & $1$ & $1$ & $1$ & $1$ & $1$ & $1 + \varepsilon$ & \tiny $1.1022$ & \tiny $1.8111$ & \tiny $1.9316$ & \tiny $2.3483$ & \tiny $2.5345$ & --- \\
\hline \boldmath{$0.01$} & $1$ & $1$ & $1$ & $1$ & $1$ & $1$ & $1$ & $1$ & $1$ & \tiny $1.0458$ & \tiny $1.7453$ & \tiny $1.8181$ & \tiny $1.9658$ & \tiny $2.3633$ & \tiny $2.5345$ \\
\hline \small \boldmath{$\theta_2 \backslash \theta_1$} & \boldmath{$0.41$} & \boldmath{$0.42$} & \boldmath{$0.43$} & \boldmath{$0.44$} & \boldmath{$0.45$} & \boldmath{$0.46$} & \boldmath{$0.47$} & \boldmath{$0.48$} & \boldmath{$0.49$} & \boldmath{$0.50$} & \boldmath{$0.51$} & \boldmath{$0.52$} & \boldmath{$0.53$} & \boldmath{$0.54$} & \boldmath{$0.55$} \\
\hline
\end{tabular} \\
\textbf{Table 4.2: Upper Bounds for }\boldmath{$C_1^{\prime}(\theta_1, \theta_2)$} \textbf{(}\boldmath{$0.5 < \theta \leqslant 0.56$}\textbf{) 2/2}
\end{center}

\subsection{Lower Bounds}
We shall construct the minorant $\rho_0(n)$ in this subsection. Before constructing, we first mention some existing results of $C_0^{\prime}(\theta_1, \theta_2)$.
\begin{theorem}\label{t413}
The function $C_0(\theta_1, \theta_2)$ satisfies the following conditions:

(1). $C_0^{\prime}(\theta_1, \theta_2) = C_0^{\prime}(\theta_2, \theta_1)$;

(2). $C_0^{\prime}(\theta_1, \theta_2) = 1$ for all $\theta_1, \theta_2$ satisfy $\theta_1 + \theta_2 \leqslant 0.5 - \varepsilon$; 

(3). $C_0^{\prime}(\theta_1, \theta_2) = 1$ for all $\theta_1, \theta_2$ satisfy $2 \theta_1 + \theta_2 \leqslant 1 - \varepsilon$ and $7 \theta_1 + 12 \theta_2 \leqslant 4 - \varepsilon$;

(4). $C_0^{\prime}(\theta_1, \theta_2) = 1$ for all $\theta_1, \theta_2$ satisfy $\theta_1 < \frac{1}{3}$, $\theta_2 < \frac{1}{5}$ and $\theta_1 + \theta_2 < \frac{29}{56}$;

(5). $C_0^{\prime}(\theta_1, \theta_2) = 1$ for all $\theta_1, \theta_2$ satisfy $\theta_1 + 3 \theta_2 < 1$, $\theta_1 + \theta_2 < \frac{29}{56}$ and $\theta_2 < \max\left(\frac{1 - 2 \theta_1}{2}, \frac{2 - 2 \theta_1}{5}\right)$;

(6). $C_0^{\prime}(\theta_1, \theta_2) = 1$ for all $\theta_1, \theta_2$ satisfy $\theta_1 + 3 \theta_2 < 1$, $\theta_1 + \theta_2 < \frac{29}{56}$, $4 \theta_1 + \theta_2 < \frac{403}{266}$ and $\frac{7}{4} \theta_1 + \theta_2 < \frac{403}{532}$;

(7). $C_0^{\prime}(\theta_1, \theta_2) \geqslant C_0(\theta_1, \theta_2) \geqslant C_0(\theta_1 + \theta_2)$ for $0.5 \leqslant \theta_1 + \theta_2 \leqslant 1$.
\end{theorem}
\begin{proof}
Statement (1) is obvious. Statements (2)--(6) follow easily from the Bombieri--Vinogradov Theorem, Theorem~\ref{t11}, [\cite{677}, Theorem 3] and [\cite{FouvryA2}, Page 621 and Corollaire 5]. Statement (7) holds trivially by the work done in Section 2 and Section 3. When there are no new arithmetic information inputs outside of those in previous sections, we use $C_0(\theta_1, \theta_2)$ as a lower bound for $C_0^{\prime}(\theta_1, \theta_2)$. Statement (8) holds from Statement (7) together with Statement (5) of Theorem~\ref{t311}.
\end{proof}

By the discussions in Subsection 2.5, we need to show that (129) holds for
\begin{equation}
f(n) = \sum_{\substack{n = p_1 \beta \\ 1 - \theta - \varepsilon < \alpha_1 < \frac{1}{2} }} \psi\left(\beta, p_1 \right).
\end{equation}
The only result that gives arithmetic information for this sum is Lemma~\ref{l35}, hence we need (108) and (109) to show that (129) holds for (243). From the Statement (3) of Theorem~\ref{t413}, we have $C_0^{\prime}(\theta_1, \theta_2) = 1$ in this case. By the discussions in Subsection 3.4, we cannot get any nontrivial lower bound for $C_0^{\prime}(\theta_1, \theta_2)$ if either (108) or (109) is not fulfilled.

\section{$3$-factored Moduli}
In this section we focus on the $3$-factored case. Since we almost do not have any new arithmetic information inputs, the $3$-factored case with absolute values is almost same as the first $2$-factored case, with only one additional Type-II information input [\cite{MaynardLargeModuliIII}, Proposition 5.2] that is only applicable for $Q_1 Q_2 Q_3 < x^{0.501}$. Hence we only discuss the trilinear case, or $3$-factored case with divisor-bounded coefficient weights. The initial setups on the sieves are similar to the bilinear case. We want to get the following result with some $0 < C_0^{\prime}(\theta_1, \theta_2, \theta_3) \leqslant 1$ and $C_1^{\prime}(\theta_1, \theta_2, \theta_3) \geqslant 1$:
\begin{theorem}\label{t51}
There exist functions $\rho_0$ and $\rho_1$ which satisfies the following properties:

(Majorant / Minorant). $\rho_0(n)$ is a minorant for the prime indicator function $\mathbbm{1}_{p}(n)$, and $\rho_1(n)$ is a majorant for the prime indicator function $\mathbbm{1}_{p}(n)$. That is, we have
$$
\rho_0(n) \leqslant \mathbbm{1}_{p}(n) \leqslant \rho_1(n).
$$

(Upper and Lower bounds). We have
$$
\sum_{n \leqslant x} \rho_0(n) \geqslant (1+o(1))\frac{C_0^{\prime}(\theta_1, \theta_2, \theta_3) x}{\log x} \quad \text{and} \quad \sum_{n \leqslant x} \rho_1(n) \leqslant (1+o(1))\frac{C_1^{\prime}(\theta_1, \theta_2, \theta_3) x}{\log x}
$$
for two functions $C_0^{\prime}(\theta_1, \theta_2, \theta_3)$ and $C_1^{\prime}(\theta_1, \theta_2, \theta_3)$ satisfy $0 < C_0^{\prime}(\theta_1, \theta_2, \theta_3) \leqslant 1$ and $C_1^{\prime}(\theta_1, \theta_2, \theta_3) \geqslant 1$.

(Distributions in Arithmetic Progressions). Let $\lambda_{j, q_j}$ ($j = 1, 2, 3$) be divisor-bounded complex sequences. For any $a \in \mathbb{Z} \backslash \{0\}$ and any $A>0$, we have
$$
\sum_{\substack{q_1 \sim Q_1 \\ q_2 \sim Q_2 \\ q_3 \sim Q_3 \\ (q_1 q_2 q_3, a) = 1}} \lambda_{1, q_1} \lambda_{2, q_2} \lambda_{3, q_3} \left( \sum_{\substack{n \leqslant x \\ n \equiv a (\bmod q_1 q_2 q_3)}} \rho_j(n) - \frac{1}{\varphi(q_1 q_2 q_3)} \sum_{\substack{n \leqslant x \\ (n, q_1 q_2 q_3) = 1}} \rho_j(n) \right) \ll \frac{x}{(\log x)^A}
$$
for $j = 0, 1$.
\end{theorem}

In order to prove Theorem~\ref{t51} with suitable $C_0^{\prime}(\theta_1, \theta_2, \theta_3)$ and $C_1^{\prime}(\theta_1, \theta_2, \theta_3)$, we need results of the form
\begin{equation}
\sum_{\substack{q_1 \sim Q_1 \\ q_2 \sim Q_2 \\ q_3 \sim Q_3 \\ (q_1 q_2 q_3, a) = 1}} \lambda_{1, q_1} \lambda_{2, q_2} \lambda_{3, q_3} \left( \sum_{\substack{n \sim x \\ n \equiv a (\bmod q_1 q_2 q_3)}} f(n) - \frac{1}{\varphi(q_1 q_2 q_3)} \sum_{\substack{n \sim x \\ (n, q_1 q_2 q_3) = 1}} f(n) \right) \ll \frac{x}{(\log x)^A}.
\end{equation}
Again, we may want the coefficients to satisfy \textbf{Conditions A and B} mentioned in Section 2.

\subsection{Preliminary Lemmas}
Before constructing the majorant and minorant, we need estimate results of the form (244). The results from Sections 2--4 (except for the results in Subsection 2.3) are still applicable in the final decompositions here.

\subsubsection{Type-II estimate}
The next lemma comes from \cite{Lichtman2}, and the readers can compare it with Lemma~\ref{l35} to see a difference. It is the only new arithmetic information input in the trilinear case.
\begin{lemma}\label{l52} ([\cite{Lichtman2}, Proposition 8.3]).
Let $M_1 M_2 \asymp x$. Let $a_{1, m_1}$, $a_{2, m_2}$, $\lambda_{1, q_1}$, $\lambda_{2, q_2}$ and $\lambda_{3, q_3}$ be divisor-bounded complex sequences. Suppose that $a_{2, m_2}$ satisfies \textbf{Conditions A and B}. If we have
$$
Q_1^7 Q_2^{12} Q_3^{10} < x^{4 - 20 \varepsilon},\ Q_2 < Q_1 Q_3^3,\ Q_1 x^{\varepsilon} < M_2 < Q_1^{-1} Q_3^{-2} x^{1 - 6 \varepsilon},
$$
then
$$
\sum_{\substack{q_1 \sim Q_1 \\ q_2 \sim Q_2 \\ q_3 \sim Q_3 \\ (q_1 q_2 q_3, a) = 1}} \lambda_{1, q_1} \lambda_{2, q_2} \lambda_{3, q_3} \left( \sum_{\substack{m_1 \sim M_1 \\ m_2 \sim M_2 \\ m_1 m_2 \equiv a (\bmod q_1 q_2 q_3)}} a_{1, m_1} a_{2, m_2} - \frac{1}{\varphi(q_1 q_2 q_3)} \sum_{\substack{m_1 \sim M_1 \\ m_2 \sim M_2 \\ (m_1 m_2, q_1 q_2 q_3) = 1}} a_{1, m_1} a_{2, m_2} \right) \ll \frac{x}{(\log x)^A}.
$$
\end{lemma}

\subsection{Sieve Asymptotic Formulas}
Using Lemma~\ref{l52}, we can get the following lemma.
\begin{lemma}\label{l53} ([\cite{Lichtman2}, Lemma 9.2]).
Define
\begin{align}
\nonumber \boldsymbol{k} =&\ \boldsymbol{k}(\theta_1, \theta_2, \theta_3) = \left\{(s, t): 1 - \theta_1 - \theta_2 - \theta_3 - \varepsilon < s < \theta_1 + \theta_2 + \theta_3 + \varepsilon \right\}, \\
\nonumber \boldsymbol{\mathcal{K}} =&\ \boldsymbol{\mathcal{K}}(\theta_1, \theta_2, \theta_3) = \left\{\boldsymbol{\alpha}_{j}: \boldsymbol{\alpha}_{j} \text{ partitions into } \boldsymbol{k} \right\}.
\end{align}
Let $P_1 P_2 \cdots P_j \asymp x$ and $P_1 \geqslant P_2 \geqslant \cdots \geqslant P_j \geqslant x^{\frac{1}{7} + 10 \varepsilon}$. Suppose that
$$
\theta_2 < \theta_1 + 3 \theta_3, \quad \theta_1 + 2 \theta_3 < \frac{1}{2} - 20 \varepsilon, \quad 2 \theta_1 + \theta_2 + \theta_3 < 1 - 20 \varepsilon \quad \text{and} \quad 7 \theta_1 + 12 \theta_2 + 10 \theta_3 < 4 - 20 \varepsilon.
$$
Then (244) holds for
$$
f(n) = \sum_{\substack{n = p_1 \cdots p_j \\ p_i \sim P_i,\ 1 \leqslant i \leqslant j \\ \boldsymbol{\alpha}_{j} \in \boldsymbol{\mathcal{K}} }} 1.
$$
\end{lemma}
Again, many asymptotic formulas used in the decompositions in this section will be adopted from previous sections.

\subsection{Upper and Lower Bounds}
We shall construct the majorant $\rho_1(n)$ and the minorant $\rho_0(n)$ in this subsection. Before constructing, we first mention some existing results of $C_1^{\prime}(\theta_1, \theta_2, \theta_3)$ and $C_0^{\prime}(\theta_1, \theta_2, \theta_3)$.
\begin{theorem}\label{t54}
The functions $C_1^{\prime}(\theta_1, \theta_2, \theta_3)$ and $C_0^{\prime}(\theta_1, \theta_2, \theta_3)$ satisfy the following conditions:

(1.1). $C_1^{\prime}(\theta_1, \theta_2, \theta_3) = C_1^{\prime}(\theta_1, \theta_3, \theta_2) = C_1^{\prime}(\theta_2, \theta_1, \theta_3) = C_1^{\prime}(\theta_2, \theta_3, \theta_1) = C_1^{\prime}(\theta_3, \theta_1, \theta_2) = C_1^{\prime}(\theta_3, \theta_2, \theta_1)$;

(1.2). $C_0^{\prime}(\theta_1, \theta_2, \theta_3) = C_0^{\prime}(\theta_1, \theta_3, \theta_2) = C_0^{\prime}(\theta_2, \theta_1, \theta_3) = C_0^{\prime}(\theta_2, \theta_3, \theta_1) = C_0^{\prime}(\theta_3, \theta_1, \theta_2) = C_0^{\prime}(\theta_3, \theta_2, \theta_1)$;

(2). $C_1^{\prime}(\theta_1, \theta_2, \theta_3) = C_0^{\prime}(\theta_1, \theta_2, \theta_3) = 1$ for all $\theta_1, \theta_2, \theta_3$ satisfy $\theta_1 + \theta_2 + \theta_3 \leqslant 0.5 - \varepsilon$; 

(3). $C_1^{\prime}(\theta_1, \theta_2, \theta_3) = C_0^{\prime}(\theta_1, \theta_2, \theta_3) = 1$ for all $\theta_1, \theta_2, \theta_3$ satisfy $\theta_1 + \theta_2 + \theta_3 = 0.5 + r$, $0 < r < 0.001$, $40r < \theta_2 < \frac{1}{20} - 7r$ and $\frac{1}{10} - \theta_2 + 12r < \theta_3 < \frac{1}{10} - \frac{3}{5} \theta_2 - 4r$; 

(4.1). $C_1^{\prime}(\theta_1, \theta_2, \theta_3) = C_0^{\prime}(\theta_1, \theta_2, \theta_3) = 1$ for all $\theta_1, \theta_2, \theta_3$ satisfy $\theta_1 < \frac{1}{3}$, $\theta_2 + \theta_3 < \frac{1}{5}$ and $\theta_1 + \theta_2 + \theta_3 < \frac{29}{56}$;

(4.2). $C_1^{\prime}(\theta_1, \theta_2, \theta_3) = C_0^{\prime}(\theta_1, \theta_2, \theta_3) = 1$ for all $\theta_1, \theta_2, \theta_3$ satisfy $\theta_1 + \theta_3 < \frac{1}{3}$, $\theta_2 < \frac{1}{5}$ and $\theta_1 + \theta_2 + \theta_3 < \frac{29}{56}$;

(5.1). $C_1^{\prime}(\theta_1, \theta_2, \theta_3) = C_0^{\prime}(\theta_1, \theta_2, \theta_3) = 1$ for all $\theta_1, \theta_2, \theta_3$ satisfy $\theta_1 + 3 \theta_2 + 3 \theta_3 < 1$, $\theta_1 + \theta_2 + \theta_3 < \frac{29}{56}$ and $\theta_2 + \theta_3 < \max\left(\frac{1 - 2 \theta_1}{2}, \frac{2 - 2 \theta_1}{5}\right)$;

(5.2). $C_1^{\prime}(\theta_1, \theta_2, \theta_3) = C_0^{\prime}(\theta_1, \theta_2, \theta_3) = 1$ for all $\theta_1, \theta_2, \theta_3$ satisfy $\theta_1 + 3 \theta_2 + \theta_3 < 1$, $\theta_1 + \theta_2 + \theta_3 < \frac{29}{56}$ and $\theta_2 < \max\left(\frac{1 - 2 \theta_1 - 2 \theta_3}{2}, \frac{2 - 2 \theta_1 - 2 \theta_3}{5}\right)$;

(6.1). $C_1^{\prime}(\theta_1, \theta_2, \theta_3) = C_0^{\prime}(\theta_1, \theta_2, \theta_3) = 1$ for all $\theta_1, \theta_2, \theta_3$ satisfy $\theta_1 + 3 \theta_2 + 3 \theta_3 < 1$, $\theta_1 + \theta_2 + \theta_3 < \frac{29}{56}$, $4 \theta_1 + \theta_2 + \theta_3 < \frac{403}{266}$ and $\frac{7}{4} \theta_1 + \theta_2 + \theta_3 < \frac{403}{532}$;

(6.2). $C_1^{\prime}(\theta_1, \theta_2, \theta_3) = C_0^{\prime}(\theta_1, \theta_2, \theta_3) = 1$ for all $\theta_1, \theta_2, \theta_3$ satisfy $\theta_1 + 3 \theta_2 + \theta_3 < 1$, $\theta_1 + \theta_2 + \theta_3 < \frac{29}{56}$, $4 \theta_1 + \theta_2 + 4 \theta_3 < \frac{403}{266}$ and $\frac{7}{4} \theta_1 + \theta_2 + \frac{7}{4} \theta_3 < \frac{403}{532}$;

(7.1). $C_1^{\prime}(\theta_1, \theta_2, \theta_3) = C_0^{\prime}(\theta_1, \theta_2, \theta_3) = 1$ for all $\theta_1, \theta_2, \theta_3$ satisfy $2 \theta_1 + \theta_2 + \theta_3 \leqslant 1 - \varepsilon$ and $7 \theta_1 + 12 \theta_2 + 12 \theta_3 \leqslant 4 - \varepsilon$;

(7.2). $C_1^{\prime}(\theta_1, \theta_2, \theta_3) = C_0^{\prime}(\theta_1, \theta_2, \theta_3) = 1$ for all $\theta_1, \theta_2, \theta_3$ satisfy $2 \theta_1 + \theta_2 + 2 \theta_3 \leqslant 1 - \varepsilon$ and $7 \theta_1 + 12 \theta_2 + 7 \theta_3 \leqslant 4 - \varepsilon$;

(8.1). $C_1^{\prime}(\theta_1, \theta_2, \theta_3) \leqslant \min\left( C_1^{\prime}(\theta_1, \theta_2 + \theta_3), C_1^{\prime}(\theta_2, \theta_1 + \theta_3), C_1^{\prime}(\theta_3, \theta_1 + \theta_2), C_1^{\prime}(\theta_1 + \theta_2, \theta_3), C_1^{\prime}(\theta_1 + \theta_3, \theta_2), C_1^{\prime}(\theta_2 + \theta_3, \theta_1) \right)$ for $0.5 \leqslant \theta_1 + \theta_2 \leqslant 1$;

(8.2). $C_0^{\prime}(\theta_1, \theta_2, \theta_3) \geqslant \max\left( C_0^{\prime}(\theta_1, \theta_2 + \theta_3), C_0^{\prime}(\theta_2, \theta_1 + \theta_3), C_0^{\prime}(\theta_3, \theta_1 + \theta_2), C_0^{\prime}(\theta_1 + \theta_2, \theta_3), C_0^{\prime}(\theta_1 + \theta_3, \theta_2), C_0^{\prime}(\theta_2 + \theta_3, \theta_1) \right)$ for $0.5 \leqslant \theta_1 + \theta_2 \leqslant 1$;

(9) $C_1^{\prime}(\theta_1, \theta_2, \theta_3) \leqslant 1 + \varepsilon$ for all $\theta_1, \theta_2, \theta_3$ satisfy $\theta_1 + \theta_2 + \theta_3 = 0.5$;

(10.1). $C_1^{\prime}(\theta_1, \theta_2, \theta_3) \leqslant 1 + \varepsilon$ for all $\theta_1, \theta_2, \theta_3$ satisfy $\theta_1 \leqslant 0.5$ and $\theta_2 + \theta_3 = \min\left(1 - 2 \theta_1, \frac{4 - 7 \theta_1}{12} \right)$;

(10.2). $C_1^{\prime}(\theta_1, \theta_2, \theta_3) \leqslant 1 + \varepsilon$ for all $\theta_1, \theta_2, \theta_3$ satisfy $\theta_1 + \theta_3 \leqslant 0.5$ and $\theta_2 = \min\left(1 - 2 \theta_1 - 2 \theta_3, \frac{4 - 7 \theta_1 - 7 \theta_3}{12} \right)$.
\end{theorem}
\begin{proof}
This theorem follows from Theorems~\ref{t11}, ~\ref{t222}, ~\ref{t223}, ~\ref{t310}, ~\ref{t311}, ~\ref{t47}, ~\ref{t413} and [\cite{MaynardLargeModuliIII}, Theorem 1.2].
\end{proof}

Next, We shall prove the following theorem.
\begin{theorem}\label{t55}
Let $Q_1 = x^{\theta_1}$, $Q_2 = x^{\theta_2}$ and $Q_3 = x^{\theta_3}$. Define
\begin{align}
\nonumber \boldsymbol{Y} =&\ \boldsymbol{E}_{0201} \cup \boldsymbol{E}_{0301} \cup \boldsymbol{E}_{0401} \cup \boldsymbol{E}_{0501} \cup \boldsymbol{E}_{0601} \cup \boldsymbol{E}_{0801} \cup \boldsymbol{E}_{0901}, \\
\nonumber \boldsymbol{J}^{\prime} =&\ \left\{ (\theta_1, \theta_2) : (\theta_1, \theta_2) \in \boldsymbol{U};\ \theta_1 + \theta_2 \leqslant \frac{1}{2} - \varepsilon \right. \\
\nonumber & \qquad \qquad \quad \text{ or } 2 \theta_1 + \theta_2 \leqslant 1 - \varepsilon,\ 7 \theta_1 + 12 \theta_2 \leqslant 4 - \varepsilon \\
\nonumber & \qquad \qquad \quad \text{ or } \theta_1 < \frac{1}{3},\ \theta_2 < \frac{1}{5},\ \theta_1 + \theta_2 < \frac{29}{56} \\
\nonumber & \qquad \qquad \quad \text{ or } \theta_1 + 3 \theta_2 < 1,\ \theta_1 + \theta_2 < \frac{29}{56},\ \theta_2 < \max\left(\frac{1 - 2 \theta_1}{2}, \frac{2 - 2 \theta_1}{5}\right) \\
\nonumber & \left. \qquad \qquad \quad \text{ or } \theta_1 + 3 \theta_2 < 1,\ \theta_1 + \theta_2 < \frac{29}{56},\ 4 \theta_1 + \theta_2 < \frac{403}{266},\ \frac{7}{4} \theta_1 + \theta_2 < \frac{403}{532} \right\}, \\
\nonumber \boldsymbol{K} =&\ \left\{ (\theta_1, \theta_2, \theta_3) : \theta_2 < \theta_1 + 3 \theta_3,\ \theta_1 + 2 \theta_3 < \frac{1}{2} - \varepsilon,\ 2 \theta_1 + \theta_2 + \theta_3 < 1 - \varepsilon,\ 7 \theta_1 + 12 \theta_2 + 10 \theta_3 < 4 - \varepsilon \right\}, \\
\nonumber \boldsymbol{\mathcal{W}} =&\ \left\{ (\theta_1, \theta_2, \theta_3) : (\theta_1, \theta_2, \theta_3) \text{ partitions into } \boldsymbol{Y},\ (\theta_1, \theta_2, \theta_3) \text{ partitions into } \boldsymbol{T}, \right. \\
\nonumber & \left. \qquad \qquad \qquad (\theta_1, \theta_2, \theta_3) \text{ does not partition into } \boldsymbol{J}^{\prime},\ (\theta_1, \theta_2, \theta_3) \in \boldsymbol{K} \right\}.
\end{align}
Suppose that $(\theta_1, \theta_2, \theta_3) \in \boldsymbol{\mathcal{W}}$. Let $\lambda_{1, q_1}$, $\lambda_{2, q_2}$ and $\lambda_{3, q_3}$ be divisor-bounded complex sequences. Then, for any fixed $a \in \mathbb{Z} \backslash \{0\}$ and any $A>0$, we have
$$
\sum_{\substack{q_1 \leqslant Q_1 \\ q_2 \leqslant Q_2 \\ q_3 \leqslant Q_3 \\ (q_1 q_2 q_3, a) = 1}} \lambda_{1, q_1} \lambda_{2, q_2} \lambda_{3, q_3} \left( \pi(x; q_1 q_2 q_3, a) - \frac{\pi(x)}{\varphi(q_1 q_2 q_3)} \right) \ll \frac{x}{(\log x)^A}.
$$
The explicit definition of the region $\boldsymbol{\mathcal{W}}$ (with $\varepsilon$ ignored) is available at \textit{https://runbolicarey.com/assets/downloads/W.pdf}.
\end{theorem}
\begin{proof}
When $(\theta_1, \theta_2, \theta_3)$ partitions into $\boldsymbol{Y}$, then at least one of $(\theta_1 + \theta_2, \theta_3)$, $(\theta_1 + \theta_3, \theta_2)$, $(\theta_2 + \theta_3, \theta_1)$, $(\theta_1, \theta_2 + \theta_3)$, $(\theta_2, \theta_1 + \theta_3)$ and $(\theta_3, \theta_1 + \theta_2)$ is in $\boldsymbol{Y}$. Thus, we can assume that $(\theta_1, \theta_2, \theta_3)$ partitions into $(s, t) \in \boldsymbol{Y}$. By the definition of the region $\boldsymbol{Y}$, we know that $2 s + t \leqslant 1 - 2 \varepsilon$. Also, for $(s, t) \in \boldsymbol{Y}$ we have a Type-II range $\left[\varepsilon,\ \frac{1}{6}(5 - 8 s - 4 t) - \varepsilon \right]$. Since the available Type-II range starts from $\varepsilon$ and $(\theta_1, \theta_2, \theta_3)$ partitions into $\boldsymbol{T}$, we can use the new three-dimensional Harman's sieve. Moreover, we have
\begin{equation}
\frac{5 - 8 s - 4 t}{6} - \varepsilon \geqslant \frac{5 - 4(1 - 2 \varepsilon)}{6} - \varepsilon = \frac{1}{6} + \frac{1}{3} \varepsilon,
\end{equation}
and the Type-II range can be simplified to $\left[\varepsilon,\ \frac{1}{6} + \frac{1}{3} \varepsilon \right]$. Now we can write $\kappa = \frac{1}{6} + \frac{1}{3} \varepsilon$. Since $(\theta_1, \theta_2, \theta_3) \in \boldsymbol{K}$ implies $\theta = \theta_1 + \theta_2 + \theta_3 \leqslant \frac{17}{32} - \varepsilon$, we can decompose our $\mathbbm{1}_{n \sim x, n = p}(n) = \psi\left(n, (2 x)^{\frac{1}{2}}\right)$ in a way similar to the decompositions in the proof of Theorem~\ref{case4}. By the discussions of the three-dimensional sieves in Section 4, we need to apply the new three-dimensional Harman's sieve on both $B$ and $C$. By Buchstab's identity, we have
\begin{align}
\nonumber \psi\left(n, (2 x)^{\frac{1}{2}}\right) =&\ \psi\left(n, x^{\kappa}\right) - \sum_{\substack{n = p_1 \beta \\ \kappa \leqslant \alpha_1 < \frac{3}{7} + \varepsilon }} \psi\left(\beta, x^{\kappa} \right) + \sum_{\substack{n = p_1 p_2 \beta \\ \kappa \leqslant \alpha_1 < \frac{3}{7} + \varepsilon \\ \kappa \leqslant \alpha_2 < \min\left(\alpha_1, \frac{1}{2}(1 - \alpha_1) \right) \\ \boldsymbol{\alpha}_2 \in \boldsymbol{G}_2 }} \psi\left(\beta, p_2 \right) \\
\nonumber & + \sum_{\substack{n = p_1 p_2 \beta \\ \kappa \leqslant \alpha_1 < \frac{3}{7} + \varepsilon \\ \kappa \leqslant \alpha_2 < \min\left(\alpha_1, \frac{1}{2}(1 - \alpha_1) \right) \\ \boldsymbol{\alpha}_2 \in A }} \psi\left(\beta, x^{\kappa} \right) - \sum_{\substack{n = p_1 p_2 p_3 \beta \\ \kappa \leqslant \alpha_1 < \frac{3}{7} + \varepsilon \\ \kappa \leqslant \alpha_2 < \min\left(\alpha_1, \frac{1}{2}(1 - \alpha_1) \right) \\ \boldsymbol{\alpha}_2 \in A \cup B \cup C \\ \kappa \leqslant \alpha_3 < \min\left(\alpha_2, \frac{1}{2}(1 - \alpha_1 - \alpha_2) \right) }} \psi\left(\beta, p_3 \right) \\
\nonumber & + \sum_{\substack{n = p_1 p_2 \beta \\ \kappa \leqslant \alpha_1 < \frac{3}{7} + \varepsilon \\ \kappa \leqslant \alpha_2 < \min\left(\alpha_1, \frac{1}{2}(1 - \alpha_1) \right) \\ \boldsymbol{\alpha}_2 \in B \cup C }} \psi\left(\beta, x^{\kappa} \right) - \sum_{\substack{n = p_1 \beta \\ \frac{3}{7} + \varepsilon \leqslant \alpha_1 \leqslant 1 - \theta - \varepsilon }} \psi\left(\beta, p_1 \right) - \sum_{\substack{n = p_1 \beta \\ 1 - \theta - \varepsilon < \alpha_1 < \frac{1}{2} }} \psi\left(\beta, p_1 \right) \\
=&\ \Sigma_{5301} - \Sigma_{5302} + \Sigma_{5303} + \Sigma_{5304} - \Sigma_{5305} + \Sigma_{5306} - \Sigma_{5307} - \Sigma_{5308}.
\end{align}
By Lemma~\ref{l212} and a Type-II range $\left[\varepsilon,\ \kappa \right]$, (244) holds for $f(n) = \Sigma_{5301}$ and $f(n) = \Sigma_{5302}$. By Lemma~\ref{l211} and Lemma~\ref{l216}, (244) holds for $f(n) = \Sigma_{5303}$ and $f(n) = \Sigma_{5307}$. By Lemma~\ref{l53}, (244) holds for $f(n) = \Sigma_{5308}$. By the Type-II range $\left[\varepsilon,\ \kappa \right]$ and the discussions in the end of the three-dimensional sieves (133)--(138) in Section 4, (244) holds for $f(n) = \Sigma_{5304}$. For the remaining sums, $\Sigma_{5305}$ only counts numbers with $4$ or more prime factors.

For $\Sigma_{5306}$, since we have a Type-II range $\left[\varepsilon,\ \kappa \right]$ and $(\theta_1, \theta_2, \theta_3)$ partitions into $\boldsymbol{T}$, we can use Lemma~\ref{l37}, Lemma~\ref{l38} and a three-dimensional Harman's sieve to get a ``loss term''
\begin{equation}
\Sigma_{5309} = \sum_{\substack{n = m_1 m_2 m_3 \\ \kappa \leqslant \alpha_1 < \frac{3}{7} + \varepsilon \\ \kappa \leqslant \alpha_2 < \min\left(\alpha_1, \frac{1}{2}(1 - \alpha_1) \right) \\ \boldsymbol{\alpha}_2 \in B \cup C \\ \Omega(m_1 m_2) \geqslant 3 }} \psi\left(m_1 m_2 m_3, x^{\kappa} \right).
\end{equation}
Since $\Omega(m_1 m_2 m_3) \geqslant \Omega(m_1 m_2) + 1 \geqslant 4$, $\Sigma_{5309}$ only counts numbers with $4$ or more prime factors. 

Now, the proof of Theorem~\ref{t55} reduces to showing that (244) holds for $f(n) = $ sums that count numbers with $4$ or more prime factors. Since $\kappa > \frac{1}{6}$, the only two cases are $\Omega(n) = 4$ and $\Omega(n) = 5$. Now we can use the arguments in [\cite{MaynardLargeModuliI}, Chapter 9] to complete the proof of Theorem~\ref{t55}. We only need the following two modifications: In \textit{Case 1} ($J = 4$) we have
\begin{equation}
Q_1 Q_2 K L^4 < x^{\frac{17}{32} + \varepsilon} P_1 P_4 (P_2 P_3)^{\frac{3}{2}} \ll x^{\frac{17}{32} + \varepsilon + 1 + \frac{3}{14} + \varepsilon} = x^{\frac{391}{224} + 2 \varepsilon} < x^{\frac{57}{32} - 10 \varepsilon}.
\end{equation}
In \textit{Case 2} ($J = 5$) we have
\begin{equation}
Q_1 Q_2 K L^4 < x^{\frac{17}{32} + \varepsilon} \frac{(P_1 P_2 P_3 P_4 P_5)^{\frac{4}{3}}}{(P_4 P_5)^{\frac{1}{3}}} \ll x^{\frac{17}{32} + \varepsilon + \frac{4}{3} - \frac{2}{7} \cdot \frac{1}{3}} = x^{\frac{1189}{672} + \varepsilon} < x^{\frac{57}{32} - 10 \varepsilon}.
\end{equation}
\end{proof}
\begin{remark*}
Since
$$
\frac{17}{32} > \frac{9}{17}, \quad \left(\frac{15}{32} + \varepsilon, \frac{1}{16} - 3 \varepsilon \right) \in \boldsymbol{E}_{0901} \cap \boldsymbol{T} \quad \text{and} \quad \left(\frac{15}{32} + \varepsilon, \frac{3}{64} - 2 \varepsilon, \frac{1}{64} - \varepsilon \right) \in \boldsymbol{K},
$$
we have
$$
\left(\frac{15}{32} + \varepsilon, \frac{3}{64} - 2 \varepsilon, \frac{1}{64} - \varepsilon \right) \in \boldsymbol{\mathcal{W}}.
$$
Hence, Theorem~\ref{t55} extends the range of $\theta$ in Theorem~\ref{t11} to $\frac{17}{32}$ under a trilinear form of moduli. Theorem~\ref{t55} also generalizes the main result of Lichtman \cite{Lichtman2}.
\end{remark*}

\section{Smooth Moduli}
In this section we focus on the smooth case, where the moduli $q \leqslant Q = x^{\theta}$ and $P^{+}(q) < x^{\delta}$. We also suppose that $q$ is square-free, hence $q \mid P(x^{\delta})$. Similarly, we want to get the following result with some $0 < C_0^{\delta}(\theta) \leqslant 1$ and $C_1^{\delta}(\theta) \geqslant 1$:
\begin{theorem}\label{t61}
There exist functions $\rho_0$ and $\rho_1$ which satisfies the following properties:

(Majorant / Minorant). $\rho_0(n)$ is a minorant for the prime indicator function $\mathbbm{1}_{p}(n)$, and $\rho_1(n)$ is a majorant for the prime indicator function $\mathbbm{1}_{p}(n)$. That is, we have
$$
\rho_0(n) \leqslant \mathbbm{1}_{p}(n) \leqslant \rho_1(n).
$$

(Upper and Lower bounds). We have
$$
\sum_{n \leqslant x} \rho_0(n) \geqslant (1+o(1))\frac{C_0^{\delta}(\theta) x}{\log x} \quad \text{and} \quad \sum_{n \leqslant x} \rho_1(n) \leqslant (1+o(1))\frac{C_1^{\delta}(\theta) x}{\log x}
$$
for two functions $C_0^{\delta}(\theta)$ and $C_1^{\delta}(\theta)$ satisfy $0 < C_0^{\delta}(\theta) \leqslant 1$ and $C_1^{\delta}(\theta) \geqslant 1$.

(Distributions in Arithmetic Progressions). For any $a \in \mathbb{Z} \backslash \{0\}$ and any $A>0$, we have
$$
\sum_{\substack{q \leqslant Q \\ q \mid P(x^{\delta}) \\ (q, a) = 1}} \left| \sum_{\substack{n \leqslant x \\ n \equiv a (\bmod q)}} \rho_j(n) - \frac{1}{\varphi(q)} \sum_{\substack{n \leqslant x \\ (n, q) = 1}} \rho_j(n) \right| \ll \frac{x}{(\log x)^A}
$$
for $j = 0, 1$.
\end{theorem}

In order to prove Theorem~\ref{t61} with suitable $C_0^{\delta}(\theta)$ and $C_1^{\delta}(\theta)$, we need results of the form
\begin{equation}
\sum_{\substack{q \leqslant Q \\ q \mid P(x^{\delta}) \\ (q, a) = 1}} \left| \sum_{\substack{n \asymp x \\ n \equiv a (\bmod q)}} f(n) - \frac{1}{\varphi(q)} \sum_{\substack{n \asymp x \\ (n, q) = 1}} f(n) \right| \ll \frac{x}{(\log x)^A}.
\end{equation}

As in previous sections, we sometimes want the coefficients to satisfy the Siegel-Wafisz condition (or \textbf{Condition A}). In our Type-I estimate, we also want the coefficients to satisfy an extra condition. Again, we use $\lambda_l$ as an example.

(\textbf{Condition C($L$)}: Smooth at scale $L$) $\lambda_l$ has the form of $\eta(\frac{l}{L})$ for some smooth function $\eta: \mathbb{R} \to \mathbb{C}$ supported on $[c_1, c_2]$ for fixed $0 < c_1 < c_2$. The function $\eta$ also satisfies the bound
$$
\left|\eta^{(j)}(x)\right| \ll (\log x)^{A}
$$
for all fixed $j \geqslant 0$, where $\eta^{(j)}$ denote the $j$-th derivative of $\eta$.

\subsection{Preliminary Lemmas}

\subsubsection{Type-II estimate}
The first estimate comes from \cite{Stadlmann525}, and it is nontrivial when $\theta < \frac{19}{36} \approx 0.5278$. In the proof of [\cite{Stadlmann525}, Theorem 2], Stadlmann used this lemma as the Type-II information input to construct a minorant for $\theta = 0.5253$.
\begin{lemma}\label{l62} ([\cite{Stadlmann525}, Proposition 1]).
Let $M_1 M_2 \asymp x$. Let $a_{1, m_1}$ and $a_{2, m_2}$ be divisor-bounded complex sequences. Suppose that $a_{2, m_2}$ satisfies \textbf{Condition A}. If we have
$$
x^{\frac{1}{2} - \sigma} \leqslant M_2 \leqslant x^{\frac{1}{2} + \sigma},
$$
where $\sigma$ and $\theta$ satisfy
$$
\sigma > 0,\ \theta > \frac{1}{2},\ 36 \theta + 24 \delta < 19,\ 24 \theta + 4 \sigma + 16 \delta < 13,\ 32 \theta + 2 \sigma + 20 \delta < 17,
$$
then
$$
\sum_{\substack{q \leqslant Q \\ q \mid P(x^{\delta}) \\ (q, a) = 1}} \left| \sum_{\substack{m_1 \asymp M_1 \\ m_2 \asymp M_2 \\ m_1 m_2 \equiv a (\bmod q)}} a_{1, m_1} a_{2, m_2} - \frac{1}{\varphi(q)} \sum_{\substack{m_1 \asymp M_1 \\ m_2 \asymp M_2 \\ (m_1 m_2, q) = 1}} a_{1, m_1} a_{2, m_2} \right| \ll \frac{x}{(\log x)^A}.
$$
\end{lemma}

The first estimate comes from \cite{Polymath8a}, and it is nontrivial when $\theta < \frac{9}{17} \approx 0.5294$. In the proof of [\cite{BakerIrving}, Theorem 1.1], Baker and Irving used this lemma as the Type-II information input to construct a minorant for $\theta \approx 0.5242$. Of course, we have $C_0^{\delta}(0.5242) = C_1^{\delta}(0.5242) = 1$ now by Stadlmann's result \cite{Stadlmann525}.
\begin{lemma}\label{l63} ([\cite{Polymath8a}, Theorem 2.8(iii)]).
Let $M_1 M_2 \asymp x$. Let $a_{1, m_1}$ and $a_{2, m_2}$ be divisor-bounded complex sequences. Suppose that $a_{2, m_2}$ satisfies \textbf{Condition A}. If we have
$$
x^{\frac{1}{2} - \sigma} \leqslant M_2 \leqslant x^{\frac{1}{2} + \sigma},
$$
where $\sigma$ and $\theta$ satisfy
$$
\sigma > 0,\ \theta > \frac{1}{2},\ 17 \theta + 7 \delta < 9,\ \frac{80}{3} \theta + \frac{34}{9} \sigma + 16 \delta < \frac{43}{3},\ 32 \theta + 2 \sigma + 18 \delta < 17,
$$
then
$$
\sum_{\substack{q \leqslant Q \\ q \mid P(x^{\delta}) \\ (q, a) = 1}} \left| \sum_{\substack{m_1 \asymp M_1 \\ m_2 \asymp M_2 \\ m_1 m_2 \equiv a (\bmod q)}} a_{1, m_1} a_{2, m_2} - \frac{1}{\varphi(q)} \sum_{\substack{m_1 \asymp M_1 \\ m_2 \asymp M_2 \\ (m_1 m_2, q) = 1}} a_{1, m_1} a_{2, m_2} \right| \ll \frac{x}{(\log x)^A}.
$$
\end{lemma}
When $\theta < \frac{19}{36}$, Lemma~\ref{l62} gives more Type-II information than Lemma~\ref{l63}. However, Lemma~\ref{l62} is not applicable when $\theta \geqslant \frac{19}{36}$. Lemma~\ref{l63} is the only Type-II information input when $\frac{19}{36} \leqslant \theta < \frac{9}{17}$.

\subsubsection{Type-I estimate}
The first Type-I estimate was proved by Baker and Irving \cite{BakerIrving} by combining 3 cases depending on the size of $M_2$, and it plays an important role in both \cite{BakerIrving} and \cite{Stadlmann525}. In \cite{LRBMinorant1019} the author used this as Type-I information input to construct a minorant for $\theta = \frac{10}{19} - \varepsilon$.
\begin{lemma}\label{l64} ([\cite{BakerIrving}, Lemma 5]).
Let $M_1 M_2 \asymp x$. Let $a_{1, m_1}$ and $a_{2, m_2}$ be divisor-bounded complex sequences. Suppose that $a_{1, m_1}$ satisfies \textbf{Condition A} and $a_{2, m_2}$ satisfies \textbf{Condition C$(M_2)$}. If we have
$$
\frac{1}{2} < \theta < \frac{10}{19},\ M_1 \leqslant x^{\frac{1}{2}(5 - 7 \theta) - 3 \delta},
$$
then
$$
\sum_{\substack{q \leqslant Q \\ q \mid P(x^{\delta}) \\ (q, a) = 1}} \left| \sum_{\substack{m_1 \asymp M_1 \\ m_2 \asymp M_2 \\ m_1 m_2 \equiv a (\bmod q)}} a_{1, m_1} a_{2, m_2} - \frac{1}{\varphi(q)} \sum_{\substack{m_1 \asymp M_1 \\ m_2 \asymp M_2 \\ (m_1 m_2, q) = 1}} a_{1, m_1} a_{2, m_2} \right| \ll \frac{x}{(\log x)^A}.
$$
\end{lemma}

For $\theta \geqslant \frac{10}{19}$, we cannot apply Lemma~\ref{l64}. Fortunately, we still have two valid Type-I information ranges. The first one can be proved by the method used in the discussions of the ``Polymath Type-0 sums'' in \cite{Polymath8a}. Readers can see the end of [\cite{Polymath8a}, Section 3] for more details.
\begin{lemma}\label{l65}
Let $M_1 M_2 \asymp x$. Let $a_{1, m_1}$ and $a_{2, m_2}$ be divisor-bounded complex sequences. Suppose that $a_{1, m_1}$ satisfies \textbf{Condition A} and $a_{2, m_2}$ satisfies \textbf{Condition C$(M_2)$}. If we have
$$
\theta > \frac{1}{2},\ \sigma > \theta - \frac{1}{2},\ M_1 \leqslant x^{\frac{1}{2} - \sigma},
$$
then
$$
\sum_{\substack{q \leqslant Q \\ q \mid P(x^{\delta}) \\ (q, a) = 1}} \left| \sum_{\substack{m_1 \asymp M_1 \\ m_2 \asymp M_2 \\ m_1 m_2 \equiv a (\bmod q)}} a_{1, m_1} a_{2, m_2} - \frac{1}{\varphi(q)} \sum_{\substack{m_1 \asymp M_1 \\ m_2 \asymp M_2 \\ (m_1 m_2, q) = 1}} a_{1, m_1} a_{2, m_2} \right| \ll \frac{x}{(\log x)^A}.
$$
\end{lemma}

The second one is [\cite{BakerIrving}, Lemma 3].
\begin{lemma}\label{l66} ([\cite{BakerIrving}, Lemma 3]).
Let $M_1 M_2 \asymp x$. Let $a_{1, m_1}$ and $a_{2, m_2}$ be divisor-bounded complex sequences. Suppose that $a_{1, m_1}$ satisfies \textbf{Condition A} and $a_{2, m_2}$ satisfies \textbf{Condition C$(M_2)$}. If we have
$$
\sigma > 0,\ \theta > \frac{1}{2},\ 14 \theta + 4 \sigma + 8 \delta < 8,\ x^{\frac{1}{2}} \leqslant M_1 \leqslant x^{\frac{1}{2} + \sigma},
$$
then
$$
\sum_{\substack{q \leqslant Q \\ q \mid P(x^{\delta}) \\ (q, a) = 1}} \left| \sum_{\substack{m_1 \asymp M_1 \\ m_2 \asymp M_2 \\ m_1 m_2 \equiv a (\bmod q)}} a_{1, m_1} a_{2, m_2} - \frac{1}{\varphi(q)} \sum_{\substack{m_1 \asymp M_1 \\ m_2 \asymp M_2 \\ (m_1 m_2, q) = 1}} a_{1, m_1} a_{2, m_2} \right| \ll \frac{x}{(\log x)^A}.
$$
\end{lemma}

Lemmas~\ref{l65}, \ref{l66} and \ref{l63} help us to obtain an analog of Lemma~\ref{l64} when $\theta \geqslant \frac{10}{19}$.
\begin{lemma}\label{l67}
Let $M_1 M_2 \asymp x$. Let $a_{1, m_1}$ and $a_{2, m_2}$ be divisor-bounded complex sequences. Suppose that $a_{1, m_1}$ satisfies \textbf{Condition A} and $a_{2, m_2}$ satisfies \textbf{Condition C$(M_2)$}. If we have
$$
\frac{1}{2} < \theta < \frac{9}{17} - \delta,\ M_1 \leqslant x^{\frac{1}{2}(5 - 7 \theta) - 3 \delta},
$$
then
$$
\sum_{\substack{q \leqslant Q \\ q \mid P(x^{\delta}) \\ (q, a) = 1}} \left| \sum_{\substack{m_1 \asymp M_1 \\ m_2 \asymp M_2 \\ m_1 m_2 \equiv a (\bmod q)}} a_{1, m_1} a_{2, m_2} - \frac{1}{\varphi(q)} \sum_{\substack{m_1 \asymp M_1 \\ m_2 \asymp M_2 \\ (m_1 m_2, q) = 1}} a_{1, m_1} a_{2, m_2} \right| \ll \frac{x}{(\log x)^A}.
$$
\end{lemma}
\begin{proof}
When $\frac{1}{2} < \theta < \frac{10}{19}$, this is Lemma~\ref{l64}.

When $\frac{10}{19} \leqslant \theta < \frac{9}{17} - \delta$, we follow the proof of Lemma~\ref{l64} by splitting the range of $M_1$ into 3 subranges:

(1). $M_1 \leqslant x^{1 - \theta - \delta}$: We take $\sigma = \theta - \frac{1}{2} + \delta$ and apply Lemma~\ref{l65};

(2). $x^{1 - \theta - \delta} < M_1 \leqslant x^{\frac{1}{2}}$: We take $\sigma = \theta - \frac{1}{2} + \delta$ and apply Lemma~\ref{l63};

(3). $x^{\frac{1}{2}} < M_1 \leqslant x^{\frac{1}{2}(5 - 7 \theta)}$: We take $\sigma = 2 - \frac{7}{2} \theta - 3 \delta$ and apply Lemma~\ref{l66}.

Combining the above 3 cases, Lemma~\ref{l67} is proved.
\end{proof}

\subsubsection{Type-I$_3$ estimate}
The next lemma is used together with Heath-Brown's identity to handle sums that count products of three large primes, and it can only give asymptotic formulas of the form (250) for such sums when $\theta < \frac{1}{2} + \frac{2}{79} \approx 0.525314$. Because of the lack of the Type-II information with a small variable, we cannot construct a three-dimensional Harman's sieve as previous sections based on this lemma.
\begin{lemma}\label{l68} ([\cite{Polymath8a}, Theorem 2.8(v)]).
Let $M_0 M_1 M_2 M_3 \asymp x$. Let $a_{0, m_0}$, $a_{1, m_1}$, $a_{2, m_2}$ and $a_{3, m_3}$ be divisor-bounded complex sequences. Suppose that $a_{i, m_i}$ satisfies \textbf{Condition C$(M_i)$} for $1 \leqslant i \leqslant 3$. If we have
$$
\theta < \frac{2}{3},\ \min\left(M_1 M_2, M_1 M_3, M_2 M_3 \right) > x^{\frac{14 \theta - 2}{9}},\ x^{\frac{28 \theta - 13}{9}} < M_1, M_2, M_3 < x^{\frac{11 - 14 \theta}{9}},
$$
then
$$
\sum_{\substack{q \leqslant Q \\ q \mid P(x^{\delta}) \\ (q, a) = 1}} \left| \sum_{\substack{m_0 \asymp M_0 \\ m_1 \asymp M_1 \\ m_2 \asymp M_2 \\ m_3 \asymp M_3 \\ m_0 m_1 m_2 m_3 \equiv a (\bmod q)}} a_{0, m_0} a_{1, m_1} a_{2, m_2} a_{3, m_3} - \frac{1}{\varphi(q)} \sum_{\substack{m_0 \asymp M_0 \\ m_1 \asymp M_1 \\ m_2 \asymp M_2 \\ m_3 \asymp M_3 \\ (m_0 m_1 m_2 m_3, q) = 1}} a_{0, m_0} a_{1, m_1} a_{2, m_2} a_{3, m_3} \right| \ll \frac{x}{(\log x)^A}.
$$
\end{lemma}

\subsection{Sieve Asymptotic Formulas}
In this subsection we prove results of the form (250) for some functions $f(n)$. We write
$$
\nu = \nu(\theta) = 1 - 2 \max\left(6 \theta - \frac{11}{4}, 16 \theta - 8 \right)
$$
and
$$
\nu^{\prime} = \nu^{\prime}(\theta) = 1 - 2 \max\left(\frac{120}{17} \theta - \frac{56}{17}, 16 \theta - 8 \right).
$$

\begin{lemma}\label{l69}
Let $\frac{1}{2} < \theta < \frac{19}{36}$ and $\alpha_1, \alpha_2, \ldots, \alpha_k > \varepsilon$. Suppose that we can partition $\{1, \ldots, k\}$ into $I$ and $J$ such that
$$
\frac{1 - \nu}{2} < \sum_{i \in I} \alpha_i < \frac{1 + \nu}{2},
$$
then (250) holds for
$$
f(n) = \sum_{n = p_1 p_2 \cdots p_k} 1.
$$

Let $\frac{1}{2} < \theta < \frac{9}{17}$ and $\alpha_1, \alpha_2, \ldots, \alpha_k > \varepsilon$. Suppose that we can partition $\{1, \ldots, k\}$ into $I$ and $J$ such that
$$
\frac{1 - \nu^{\prime}}{2} < \sum_{i \in I} \alpha_i < 1 - \frac{1 - \nu^{\prime}}{2},
$$
then (250) holds for
$$
f(n) = \sum_{n = p_1 p_2 \cdots p_k} 1.
$$
\end{lemma}
\begin{proof}
This lemma follows easily from Lemma~\ref{l62} and Lemma~\ref{l63}, taking $\sigma = \frac{\nu}{2} - 11 \delta$ and $\sigma = \frac{\nu^{\prime}}{2} - 11 \delta$ respectively.
\end{proof}

\begin{lemma}\label{l610}
Let $\frac{1}{2} < \theta < \frac{19}{36}$ and $\alpha_1, \alpha_2, \ldots, \alpha_k > \varepsilon$. Suppose that we can partition $\{1, \ldots, k\}$ into $I$ and $J$ such that
$$
\sum_{i \in I} \alpha_i < \frac{1 - \nu}{2},\ \sum_{j \in J} \alpha_j < \frac{4 - 7 \theta + \nu}{2},
$$
then (250) holds for
$$
f(n) = \sum_{n = p_1 p_2 \cdots p_k \beta} \psi\left(\beta, x^{\nu}\right) \quad \text{and} \quad f(n) = \psi\left(n, x^{\nu}\right).
$$

Let $\frac{1}{2} < \theta < \frac{9}{17}$ and $\alpha_1, \alpha_2, \ldots, \alpha_k > \varepsilon$. Suppose that we can partition $\{1, \ldots, k\}$ into $I$ and $J$ such that
$$
\sum_{i \in I} \alpha_i < \frac{1 - \nu^{\prime}}{2},\ \sum_{j \in J} \alpha_j < \frac{4 - 7 \theta + \nu^{\prime}}{2},
$$
then (250) holds for
$$
f(n) = \sum_{n = p_1 p_2 \cdots p_k \beta} \psi\left(\beta, x^{\nu^{\prime}}\right) \quad \text{and} \quad f(n) = \psi\left(n, x^{\nu^{\prime}}\right).
$$
\end{lemma}
\begin{proof}
This lemma can be proved by applying [\cite{BakerWeingartner}, Lemma 14], with Lemma~\ref{l67} as Type-I information and Lemmas~\ref{l62}--\ref{l63} as Type-II information.
\end{proof}

\begin{lemma}\label{l611}
Let $\frac{1}{2} < \theta < \frac{1}{2} + \frac{1}{79}$. We have (250) holds for
$$
f(n) = \sum_{\substack{n = p_1 p_2 p_3 \\ \nu \leqslant \alpha_2 < \alpha_1 < \frac{1 - \nu}{2} \\ \alpha_1 + \alpha_2 \geqslant \frac{1 + \nu}{2} \\ \alpha_3 > \alpha_2 \geqslant \frac{4 - 7 \theta + \nu}{2} }} 1.
$$
\end{lemma}
\begin{proof}
See [\cite{Stadlmann525}, Section 4.5]. Lemma~\ref{l68} was used in the proof.
\end{proof}

\subsection{Lower Bounds}
We shall construct the minorant $\rho_0(n)$ in this subsection. Before constructing, we first mention existing results of $C_0^{\delta}(\theta)$ proved by Stadlmann \cite{Stadlmann525}.
\begin{theorem}\label{t612} ([\cite{Stadlmann525}, Theorem 1]).
The function $C_1^{\delta}(\theta)$ satisfies the following conditions:

(1). $C_0^{\delta}(\theta) = 1$ for all $\theta < 0.525$.

(2). $C_0^{\delta}(\theta) \geqslant 0.9999$ for all $0.525 \leqslant \theta < 0.5253$.
\end{theorem}

Recalling that our aim is to decompose $\psi\left(n, (2 x)^{\frac{1}{2}}\right)$ using Buchstab's identity and show that (250) holds for most of the sums after the decomposition. For the remaining sums that we cannot ensure (250) holds, we must make them positive so that we can drop them in order to get a lower bound. Now we split the range $\theta \in \left[0.525,\ \frac{9}{17}\right)$ to 2 subranges.

\subsubsection{Case 1. $0.525 \leqslant \theta < \frac{19}{36}$}
Using Buchstab's identity twice, we have
\begin{align}
\nonumber \psi\left(n, (2 x)^{\frac{1}{2}}\right) =&\ \psi\left(n, x^{\nu}\right) - \sum_{\substack{n = p_1 \beta \\ \nu \leqslant \alpha_1 < \frac{1}{2} }} \psi\left(\beta, p_1 \right) \\
\nonumber =&\ \psi\left(n, x^{\nu}\right) - \sum_{\substack{n = p_1 \beta \\ \nu \leqslant \alpha_1 < \frac{1 - \nu}{2} }} \psi\left(\beta, p_1 \right) - \sum_{\substack{n = p_1 \beta \\ \frac{1 - \nu}{2} \leqslant \alpha_1 < \frac{1}{2} }} \psi\left(\beta, p_1 \right) \\
\nonumber =&\ \psi\left(n, x^{\nu}\right) - \sum_{\substack{n = p_1 \beta \\ \nu \leqslant \alpha_1 < \frac{1 - \nu}{2} }} \psi\left(\beta, x^{\nu} \right) - \sum_{\substack{n = p_1 \beta \\ \frac{1 - \nu}{2} \leqslant \alpha_1 < \frac{1}{2} }} \psi\left(\beta, p_1 \right) + \sum_{\substack{n = p_1 p_2 \beta \\ \nu \leqslant \alpha_1 < \frac{1 - \nu}{2} \\ \nu \leqslant \alpha_2 < \min\left(\alpha_1, \frac{1}{2}(1 - \alpha_1) \right) }} \psi\left(\beta, p_2 \right) \\
\nonumber =&\ \psi\left(n, x^{\nu}\right) - \sum_{\substack{n = p_1 \beta \\ \nu \leqslant \alpha_1 < \frac{1 - \nu}{2} }} \psi\left(\beta, x^{\nu} \right) - \sum_{\substack{n = p_1 \beta \\ \frac{1 - \nu}{2} \leqslant \alpha_1 < \frac{1}{2} }} \psi\left(\beta, p_1 \right) \\
\nonumber & + \sum_{\substack{n = p_1 p_2 \beta \\ \nu \leqslant \alpha_1 < \frac{1 - \nu}{2} \\ \nu \leqslant \alpha_2 < \min\left(\alpha_1, \frac{1}{2}(1 - \alpha_1) \right) \\ \alpha_1 + \alpha_2 \in \left(\frac{1 - \nu}{2},\ \frac{1 + \nu}{2} \right) }} \psi\left(\beta, p_2 \right) + \sum_{\substack{n = p_1 p_2 \beta \\ \nu \leqslant \alpha_1 < \frac{1 - \nu}{2} \\ \nu \leqslant \alpha_2 < \min\left(\alpha_1, \frac{1}{2}(1 - \alpha_1) \right) \\ \alpha_1 + \alpha_2 \notin \left(\frac{1 - \nu}{2},\ \frac{1 + \nu}{2} \right) }} \psi\left(\beta, p_2 \right) \\
=&\ S_{61} - S_{62} - S_{63} - S_{64} - S_{65}.
\end{align}
By Lemma~\ref{l610} we know that (250) holds for $f(n) = S_{61}$ and $f(n) = S_{62}$, and by Lemma~\ref{l69} we know that (250) holds for $f(n) = S_{63}$ and $f(n) = S_{64}$. 

We divide the sum $S_{65}$ into 4 subsums that summing over 4 regions $A$, $B$, $C$ and $D$ (only defined in this subsection) respectively. These regions are defined as
\begin{align}
\nonumber \nonumber A =&\ A(\theta) = \left\{\boldsymbol{\alpha}_{2}: \nu \leqslant \alpha_1 < \frac{1 - \nu}{2},\ \nu \leqslant \alpha_2 < \min\left(\alpha_1, \frac{1}{2}(1 - \alpha_1) \right),\ \alpha_1 + \alpha_2 \leqslant \frac{1 - \nu}{2},\ \alpha_2 < \frac{4 - 7 \theta + \nu}{2} \right\}, \\
\nonumber \nonumber B =&\ B(\theta) = \left\{\boldsymbol{\alpha}_{2}: \nu \leqslant \alpha_1 < \frac{1 - \nu}{2},\ \nu \leqslant \alpha_2 < \min\left(\alpha_1, \frac{1}{2}(1 - \alpha_1) \right),\ \alpha_1 + \alpha_2 \geqslant \frac{1 + \nu}{2},\ \alpha_2 < \frac{4 - 7 \theta + \nu}{2} \right\}, \\
\nonumber \nonumber C =&\ C(\theta) = \left\{\boldsymbol{\alpha}_{2}: \nu \leqslant \alpha_1 < \frac{1 - \nu}{2},\ \nu \leqslant \alpha_2 < \min\left(\alpha_1, \frac{1}{2}(1 - \alpha_1) \right),\ \alpha_1 + \alpha_2 \geqslant \frac{1 + \nu}{2},\ \alpha_2 \geqslant \frac{4 - 7 \theta + \nu}{2} \right\}, \\
\nonumber \nonumber D =&\ D(\theta) = \left\{\boldsymbol{\alpha}_{2}: \nu \leqslant \alpha_1 < \frac{1 - \nu}{2},\ \nu \leqslant \alpha_2 < \min\left(\alpha_1, \frac{1}{2}(1 - \alpha_1) \right),\ \alpha_1 + \alpha_2 \leqslant \frac{1 - \nu}{2},\ \alpha_2 \geqslant \frac{4 - 7 \theta + \nu}{2} \right\}.
\end{align}
Now, we have
\begin{align}
\nonumber S_{65} =&\ \sum_{\substack{n = p_1 p_2 \beta \\ \nu \leqslant \alpha_1 < \frac{1 - \nu}{2} \\ \nu \leqslant \alpha_2 < \min\left(\alpha_1, \frac{1}{2}(1 - \alpha_1) \right) \\ \alpha_1 + \alpha_2 \notin \left(\frac{1 - \nu}{2},\ \frac{1 + \nu}{2} \right) }} \psi\left(\beta, p_2 \right) \\
\nonumber =&\ \sum_{\substack{n = p_1 p_2 \beta \\ \boldsymbol{\alpha}_2 \in A }} \psi\left(\beta, p_2 \right) + \sum_{\substack{n = p_1 p_2 \beta \\ \boldsymbol{\alpha}_2 \in B }} \psi\left(\beta, p_2 \right) + \sum_{\substack{n = p_1 p_2 \beta \\ \boldsymbol{\alpha}_2 \in C }} \psi\left(\beta, p_2 \right) + \sum_{\substack{n = p_1 p_2 \beta \\ \boldsymbol{\alpha}_2 \in D }} \psi\left(\beta, p_2 \right) \\
\nonumber =&\ S_{65A} + S_{65B} + S_{65C} + S_{65D}.
\end{align}

For $S_{65A}$, since $\boldsymbol{\alpha}_3$ satisfies Lemma~\ref{l610}, we can use Buchstab's identity twice again to reach a four-dimensional sum
\begin{equation}
\sum_{\substack{n = p_1 p_2 p_3 p_4 \beta \\ \boldsymbol{\alpha}_2 \in A_1 \\ \nu \leqslant \alpha_3 < \min\left(\alpha_2, \frac{1}{2}(1 - \alpha_1 - \alpha_2) \right) \\ \boldsymbol{\alpha}_3 \text{ does not satisfy Lemma~\ref{l69}} \\ \nu \leqslant \alpha_4 < \min\left(\alpha_3, \frac{1}{2}(1 - \alpha_1 - \alpha_2 - \alpha_3) \right) \\ \boldsymbol{\alpha}_4 \text{ does not satisfy Lemma~\ref{l69}} }} \psi\left(\beta, p_4 \right).
\end{equation}
In the above sum, we can still perform Buchstab's identity twice when $(\alpha_1, \alpha_2, \alpha_3, \alpha_4, \alpha_4)$ satisfies Lemma~\ref{l610}. This process leads to a six-dimensional loss
\begin{equation}
\sum_{\substack{n = p_1 p_2 p_3 p_4 p_5 p_6 \beta \\ \boldsymbol{\alpha}_2 \in A_1 \\ \nu \leqslant \alpha_3 < \min\left(\alpha_2, \frac{1}{2}(1 - \alpha_1 - \alpha_2) \right) \\ \boldsymbol{\alpha}_3 \text{ does not satisfy Lemma~\ref{l69}} \\ \nu \leqslant \alpha_4 < \min\left(\alpha_3, \frac{1}{2}(1 - \alpha_1 - \alpha_2 - \alpha_3) \right) \\ \boldsymbol{\alpha}_4 \text{ does not satisfy Lemma~\ref{l69}} \\ (\alpha_1, \alpha_2, \alpha_3, \alpha_4, \alpha_4) \text{ satisfies Lemma~\ref{l610}} \\ \nu \leqslant \alpha_5 < \min\left(\alpha_4, \frac{1}{2}(1 - \alpha_1 - \alpha_2 - \alpha_3 - \alpha_4) \right) \\ \boldsymbol{\alpha}_5 \text{ does not satisfy Lemma~\ref{l69}} \\ \nu \leqslant \alpha_6 < \min\left(\alpha_5, \frac{1}{2}(1 - \alpha_1 - \alpha_2 - \alpha_3 - \alpha_4 - \alpha_5) \right) \\ \boldsymbol{\alpha}_6 \text{ does not satisfy Lemma~\ref{l69}} }} \psi\left(\beta, p_6 \right).
\end{equation}
Note that we have $\nu > \frac{1}{9}$ when $\theta < \frac{19}{36}$, only $\Omega(n) \leqslant 8$ will be counted in the sums, and further straightforward decompositions are not applicable since $n = p_1 \cdots p_8 m$ implies $\Omega(n) \geqslant 9$. For (252) and (253), we can also use reversed Buchstab's identity to gain possible savings.

For $S_{65B}$ we cannot perform a straightforward decomposition. However, since $\alpha_1 + \alpha_2 \geqslant \frac{1 + \nu}{2}$ implies $1 - \alpha_1 - \alpha_2 \leqslant \frac{1 - \nu}{2}$, a role-reversal can be applied. We first use Buchstab's identity once to get
\begin{align}
\nonumber S_{65B} =&\ \sum_{\substack{n = p_1 p_2 \beta \\ \boldsymbol{\alpha}_2 \in B_1 }} \psi\left(\beta, x^{\nu} \right) - \sum_{\substack{n = p_1 p_2 p_3 \beta \\ \boldsymbol{\alpha}_2 \in B_1 \\ \nu \leqslant \alpha_3 < \min\left(\alpha_2, \frac{1}{2}(1 - \alpha_1 - \alpha_2) \right) \\ \boldsymbol{\alpha}_3 \text{ satisfies Lemma~\ref{l69}} }} \psi\left(\beta, p_3 \right) - \sum_{\substack{n = p_1 p_2 p_3 \beta \\ \boldsymbol{\alpha}_2 \in B_1 \\ \nu \leqslant \alpha_3 < \min\left(\alpha_2, \frac{1}{2}(1 - \alpha_1 - \alpha_2) \right) \\ \boldsymbol{\alpha}_3 \text{ does not satisfy Lemma~\ref{l69}} }} \psi\left(\beta, p_3 \right) \\
=&\ S_{65B1} - S_{65B2} - S_{65B3}.
\end{align}
We know that (250) holds for $S_{65B1}$ and $S_{65B2}$ by Lemma~\ref{l610} (since we have $\alpha_1 < \frac{1 - \nu}{2}$ and $\alpha_2 < \frac{4 - 7 \theta + \nu}{2}$ in $B$) and Lemma~\ref{l69} respectively. For $S_{65B3}$, we change the roles of $p_1$ and $m$ in $S_{65B}$, and use Buchstab's identity on $p_1$ to get
\begin{align}
\nonumber S_{65B3} =&\ \sum_{\substack{n = p_1 p_2 p_3 \beta \\ \boldsymbol{\alpha}_2 \in B_1 \\ \nu \leqslant \alpha_3 < \min\left(\alpha_2, \frac{1}{2}(1 - \alpha_1 - \alpha_2) \right) \\ \boldsymbol{\alpha}_3 \text{ does not satisfy Lemma~\ref{l69}} }} \psi\left(\beta, p_3 \right) \\
\nonumber =&\ \sum_{\substack{n = \beta_1 p_2 p_3 \beta \\ \boldsymbol{\alpha}_2 \in B_1 \\ \nu \leqslant \alpha_3 < \min\left(\alpha_2, \frac{1}{2}(1 - \alpha_1 - \alpha_2) \right) \\ \boldsymbol{\alpha}_3 \text{ does not satisfy Lemma~\ref{l69}} }} \psi\left(\beta, p_3 \right) \psi\left(\beta_1, \left(\frac{2x}{p_2 p_3 m} \right)^{\frac{1}{2}} \right) \\
\nonumber =&\ \sum_{\substack{n = \beta_1 p_2 p_3 \beta \\ \boldsymbol{\alpha}_2 \in B_1 \\ \nu \leqslant \alpha_3 < \min\left(\alpha_2, \frac{1}{2}(1 - \alpha_1 - \alpha_2) \right) \\ \boldsymbol{\alpha}_3 \text{ does not satisfy Lemma~\ref{l69}} }} \psi\left(\beta, p_3 \right) \psi\left(\beta_1, x^{\nu} \right) \\
\nonumber & - \sum_{\substack{n = \beta_1 p_2 p_3 p_4 \beta \\ \boldsymbol{\alpha}_2 \in B_1 \\ \nu \leqslant \alpha_3 < \min\left(\alpha_2, \frac{1}{2}(1 - \alpha_1 - \alpha_2) \right) \\ \boldsymbol{\alpha}_3 \text{ does not satisfy Lemma~\ref{l69}} \\ \nu \leqslant \alpha_4 < \frac{1}{2} \alpha_1 \\ (1 - \alpha_1  - \alpha_2 - \alpha_3, \alpha_2, \alpha_3, \alpha_4) \text{ satisfies Lemma~\ref{l69}} }} \psi\left(\beta, p_3 \right) \psi\left(\beta_1, p_4 \right) \\
\nonumber & - \sum_{\substack{n = \beta_1 p_2 p_3 p_4 \beta \\ \boldsymbol{\alpha}_2 \in B_1 \\ \nu \leqslant \alpha_3 < \min\left(\alpha_2, \frac{1}{2}(1 - \alpha_1 - \alpha_2) \right) \\ \boldsymbol{\alpha}_3 \text{ does not satisfy Lemma~\ref{l69}} \\ \nu \leqslant \alpha_4 < \frac{1}{2} \alpha_1 \\ (1 - \alpha_1  - \alpha_2 - \alpha_3, \alpha_2, \alpha_3, \alpha_4) \text{ does not satisfy Lemma~\ref{l69}} }} \psi\left(\beta, p_3 \right) \psi\left(\beta_1, p_4 \right) \\
=&\ S_{65B31} - S_{65B32} - S_{65B33}.
\end{align}
In the above 3 sums, we have $\beta \asymp x^{1 - \alpha_1 - \alpha_2 - \alpha_3}$. We have $\beta_1 \asymp x^{\alpha_1}$ in $S_{65B31}$ and $\beta_1 \asymp x^{\alpha_1 - \alpha_4}$ in $S_{65B32}$ and $S_{65B33}$. Here, (250) holds for $S_{65B31}$ and $S_{65B32}$ by Lemma~\ref{l610} (since we have $(1 - \alpha_1 - \alpha_2 - \alpha_3) + \alpha_2 = 1 - \alpha_1 - \alpha_2 < \frac{1 - \nu}{2}$ and $\alpha_2 < \frac{4 - 7 \theta + \nu}{2}$ in $B$) and Lemma~\ref{l69} respectively. We discard $S_{65B33}$ which gives a four-dimensional loss. Note that further decompositions and the reversed Buchstab's identity can still be performed to gain possible savings.

For $S_{65C}$ and For $S_{65D}$ we cannot perform either straightforward decompositions or decompositions with role-reversals, and we can only discard the whole of them. However, Lemma~\ref{l611} is applicable to show that (250) holds for $S_{65C}$ when $\theta \leqslant 0.5253$. Since $\alpha_2 \geqslant \frac{4 - 7 \cdot 0.5253 + \nu(0.5253)}{2} > 0.256 > \frac{1}{4}$ when $0.525 \leqslant \theta \leqslant 0.5253$, we know that $S_{65C}$ only counts products of 3 primes, and we do not need to discard it in this range of $\theta$ by an application of Lemma~\ref{l611}.

Combining the loss from all 4 subsums, we can get the total loss and the lower bounds for $C_0^{\delta}(\theta)$.
\begin{center}
\begin{tabular}{|c|c|c|c|}
\hline \boldmath{$\theta$} & \boldmath{$C_0^{\delta}(\theta)$} & \boldmath{$\theta$} & \boldmath{$C_0^{\delta}(\theta)$} \\
\hline $0.5250$ & $0.9999$ & $0.5264$ & $0.7407$ \\
\hline $0.5251$ & $0.9999$ & $0.5265$ & $0.7259$ \\
\hline $0.5252$ & $0.9999$ & $0.5266$ & $0.7141$ \\
\hline $0.5253$ & $0.9999$ & $0.5267$ & $0.6947$ \\
\hline $0.5254$ & $0.8416$ & $0.5268$ & $0.6807$ \\
\hline $0.5255$ & $0.8337$ & $0.5269$ & $0.6662$ \\
\hline $0.5256$ & $0.8251$ & $0.5270$ & $0.6509$ \\
\hline $0.5257$ & $0.8166$ & $0.5271$ & $0.6312$ \\
\hline $0.5258$ & $0.8075$ & $0.5272$ & $0.6107$ \\
\hline $0.5259$ & $0.7979$ & $0.5273$ & $0.5903$ \\
\hline $0.5260$ & $0.7868$ & $0.5274$ & $0.5654$ \\
\hline $0.5261$ & $0.7759$ & $0.5275$ & $0.5374$ \\
\hline $0.5262$ & $0.7641$ & $0.5276$ & $0.5135$ \\
\hline $0.5263$ & $0.7527$ & $0.5277$ & $0.4760$ \\
\hline
\end{tabular} \\
\textbf{Table 6.1: Lower Bounds for }\boldmath{$C_0^{\delta}(\theta)$} \textbf{(}\boldmath{$0.525 \leqslant \theta < \frac{19}{36}$}\textbf{)}
\end{center}

\subsubsection{Case 2. $\frac{19}{36} \leqslant \theta < \frac{9}{17}$} The decompositions in this case are very similar to the first case; one can just replace the parameter $\nu$ occurred above with $\nu^{\prime}$ and calculate the total loss. Note that Lemma~\ref{l611} is not applicable in this case. Working like the above case we get the following lower bounds for $C_0^{\delta}(\theta)$.
\begin{center}
\begin{tabular}{|c|c|}
\hline \boldmath{$\theta$} & \boldmath{$C_0^{\delta}(\theta)$} \\
\hline $0.5278$ & $0.4291$ \\
\hline $0.5279$ & $0.3856$ \\
\hline $0.5280$ & $0.3354$ \\
\hline $0.5281$ & $0.2833$ \\
\hline $0.5282$ & $0.2189$ \\
\hline $0.5283$ & $0.1489$ \\
\hline $0.5284$ & $0.0648$ \\
\hline $0.5285$ & $-0.023$ \\
\hline
\end{tabular} \\
\textbf{Table 6.2: Lower Bounds for }\boldmath{$C_0^{\delta}(\theta)$} \textbf{(}\boldmath{$\frac{19}{36} \leqslant \theta < \frac{9}{17}$}\textbf{)}
\end{center}
Note that the lower bound becomes trivial when $\theta \geqslant 0.5285$.

\subsection{Upper Bounds}
We shall construct the majorant $\rho_1(n)$ in this subsection. Before constructing, we first mention an existing result of $C_1^{\delta}(\theta)$ proved by Stadlmann \cite{Stadlmann525}.
\begin{theorem}\label{t613} ([\cite{Stadlmann525}, Theorem 1]).
The function $C_1^{\delta}(\theta)$ satisfies the following condition:

$C_1^{\delta}(\theta) = 1$ for all $\theta < 0.525$.
\end{theorem}

Recalling that our aim is to decompose $\psi\left(n, (2 x)^{\frac{1}{2}}\right)$ using Buchstab's identity and show that (250) holds for most of the sums after the decomposition. For the remaining sums that we cannot ensure (250) holds, we must make them negative so that we can drop them in order to get an upper bound. Now we split the range $\theta \in \left[0.525,\ \frac{9}{17}\right)$ to 2 subranges.

\subsubsection{Case 1. $0.525 \leqslant \theta < \frac{19}{36}$}
Using Buchstab's identity, we get
\begin{align}
\nonumber \psi\left(n, (2 x)^{\frac{1}{2}}\right) =&\ \psi\left(n, x^{\nu}\right) - \sum_{\substack{n = p_1 \beta \\ \nu \leqslant \alpha_1 < \frac{1}{2} }} \psi\left(\beta, p_1 \right) \\
\nonumber =&\ \psi\left(n, x^{\nu}\right) - \sum_{\substack{n = p_1 \beta \\ \nu \leqslant \alpha_1 < \frac{1 + \nu}{4} }} \psi\left(\beta, p_1 \right) - \sum_{\substack{n = p_1 \beta \\ \frac{1 + \nu}{4} \leqslant \alpha_1 < \frac{1 - \nu}{2} }} \psi\left(\beta, p_1 \right) - \sum_{\substack{n = p_1 \beta \\ \frac{1 - \nu}{2} \leqslant \alpha_1 < \frac{1}{2} }} \psi\left(\beta, p_1 \right) \\
\nonumber =&\ \psi\left(n, x^{\nu}\right) - \sum_{\substack{n = p_1 \beta \\ \frac{1 + \nu}{4} \leqslant \alpha_1 < \frac{1 - \nu}{2} }} \psi\left(\beta, p_1 \right) - \sum_{\substack{n = p_1 \beta \\ \frac{1 - \nu}{2} \leqslant \alpha_1 < \frac{1}{2} }} \psi\left(\beta, p_1 \right) \\
\nonumber & - \sum_{\substack{n = p_1 \beta \\ \nu \leqslant \alpha_1 < \frac{1 + \nu}{4} }} \psi\left(\beta, x^{\nu} \right) + \sum_{\substack{n = p_1 p_2 \beta \\ \nu \leqslant \alpha_1 < \frac{1 + \nu}{4} \\ \nu \leqslant \alpha_2 < \min\left(\alpha_1, \frac{1}{2}(1 - \alpha_1) \right) }} \psi\left(\beta, p_2 \right) \\
\nonumber =&\ \psi\left(n, x^{\nu}\right) - \sum_{\substack{n = p_1 \beta \\ \frac{1 + \nu}{4} \leqslant \alpha_1 < \frac{1 - \nu}{2} }} \psi\left(\beta, p_1 \right) - \sum_{\substack{n = p_1 \beta \\ \frac{1 - \nu}{2} \leqslant \alpha_1 < \frac{1}{2} }} \psi\left(\beta, p_1 \right) \\
\nonumber & - \sum_{\substack{n = p_1 \beta \\ \nu \leqslant \alpha_1 < \frac{1 + \nu}{4} }} \psi\left(\beta, x^{\nu} \right) + \sum_{\substack{n = p_1 p_2 \beta \\ \nu \leqslant \alpha_1 < \frac{1 + \nu}{4} \\ \nu \leqslant \alpha_2 < \min\left(\alpha_1, \frac{1}{2}(1 - \alpha_1) \right) \\ \alpha_1 + \alpha_2 < \frac{1 - \nu}{2} }} \psi\left(\beta, p_2 \right) + \sum_{\substack{n = p_1 p_2 \beta \\ \nu \leqslant \alpha_1 < \frac{1 + \nu}{4} \\ \nu \leqslant \alpha_2 < \min\left(\alpha_1, \frac{1}{2}(1 - \alpha_1) \right) \\ \frac{1 - \nu}{2} \leqslant \alpha_1 + \alpha_2 < \frac{1 + \nu}{2} }} \psi\left(\beta, p_2 \right) \\
\nonumber =&\ \psi\left(n, x^{\nu}\right) - \sum_{\substack{n = p_1 \beta \\ \frac{1 + \nu}{4} \leqslant \alpha_1 < \frac{1 - \nu}{2} }} \psi\left(\beta, p_1 \right) - \sum_{\substack{n = p_1 \beta \\ \frac{1 - \nu}{2} \leqslant \alpha_1 < \frac{1}{2} }} \psi\left(\beta, p_1 \right) \\
\nonumber & - \sum_{\substack{n = p_1 \beta \\ \nu \leqslant \alpha_1 < \frac{1 + \nu}{4} }} \psi\left(\beta, x^{\nu} \right) + \sum_{\substack{n = p_1 p_2 \beta \\ \nu \leqslant \alpha_1 < \frac{1 + \nu}{4} \\ \nu \leqslant \alpha_2 < \min\left(\alpha_1, \frac{1}{2}(1 - \alpha_1) \right) \\ \frac{1 - \nu}{2} \leqslant \alpha_1 + \alpha_2 < \frac{1 + \nu}{2} }} \psi\left(\beta, p_2 \right) \\
\nonumber & + \sum_{\substack{n = p_1 p_2 \beta \\ \nu \leqslant \alpha_1 < \frac{1 + \nu}{4} \\ \nu \leqslant \alpha_2 < \min\left(\alpha_1, \frac{1}{2}(1 - \alpha_1) \right) \\ \alpha_1 + \alpha_2 < \frac{1 - \nu}{2} }} \psi\left(\beta, x^{\nu} \right) - \sum_{\substack{n = p_1 p_2 p_3 \beta \\ \nu \leqslant \alpha_1 < \frac{1 + \nu}{4} \\ \nu \leqslant \alpha_2 < \min\left(\alpha_1, \frac{1}{2}(1 - \alpha_1) \right) \\ \alpha_1 + \alpha_2 < \frac{1 - \nu}{2} \\ \nu \leqslant \alpha_3 < \min\left(\alpha_2, \frac{1}{2}(1 - \alpha_1 - \alpha_2) \right) \\ \boldsymbol{\alpha}_3 \text{ satisfies Lemma~\ref{l69}} }} \psi\left(\beta, p_3 \right) \\
\nonumber & - \sum_{\substack{n = p_1 p_2 p_3 \beta \\ \nu \leqslant \alpha_1 < \frac{1 + \nu}{4} \\ \nu \leqslant \alpha_2 < \min\left(\alpha_1, \frac{1}{2}(1 - \alpha_1) \right) \\ \alpha_1 + \alpha_2 < \frac{1 - \nu}{2} \\ \nu \leqslant \alpha_3 < \min\left(\alpha_2, \frac{1}{2}(1 - \alpha_1 - \alpha_2) \right) \\ \boldsymbol{\alpha}_3 \text{ does not satisfy Lemma~\ref{l69}} }} \psi\left(\beta, p_3 \right) \\
=&\ S_{661} - S_{662} - S_{663} - S_{664} + S_{665} + S_{666} - S_{667} - S_{668}.
\end{align}
We know that (250) holds for $S_{661}$, $S_{664}$, $S_{666}$ (by Lemma~\ref{l610}) and $S_{663}$, $S_{665}$, $S_{667}$ (by Lemma~\ref{l69}). We can perform a further decomposition on $S_{668}$ if $(\alpha_1, \alpha_2, \alpha_3, \alpha_3)$ satisfies Lemma~\ref{l610}, leading to a five-dimensional sum similar to (253). We discard this sum and remaining parts of $S_{668}$ where $(\alpha_1, \alpha_2, \alpha_3, \alpha_3)$ does not satisfy Lemma~\ref{l610}. Again, reversed Buchstab's identity and role-reversals can be applied.

When $\theta \geqslant \frac{1}{2} + \frac{1}{79}$, we discard the whole of $S_{662}$. When $\theta < \frac{1}{2} + \frac{1}{79}$, we perform a further decomposition on $S_{662}$ to get
\begin{align}
\nonumber S_{662} =&\ \sum_{\substack{n = p_1 \beta \\ \frac{1 + \nu}{4} \leqslant \alpha_1 < \frac{1 - \nu}{2} }} \psi\left(\beta, p_1 \right) \\
\nonumber =&\ \sum_{\substack{n = p_1 \beta \\ \frac{1 + \nu}{4} \leqslant \alpha_1 < \frac{1 - \nu}{2} }} \psi\left(\beta, x^{\nu} \right) - \sum_{\substack{n = p_1 p_2 \beta \\ \frac{1 + \nu}{4} \leqslant \alpha_1 < \frac{1 - \nu}{2} \\ \nu \leqslant \alpha_2 < \min\left(\alpha_1, \frac{1}{2}(1 - \alpha_1) \right) }} \psi\left(\beta, p_2 \right) \\
\nonumber =&\ \sum_{\substack{n = p_1 \beta \\ \frac{1 + \nu}{4} \leqslant \alpha_1 < \frac{1 - \nu}{2} }} \psi\left(\beta, x^{\nu} \right) - \sum_{\substack{n = p_1 p_2 \beta \\ \frac{1 + \nu}{4} \leqslant \alpha_1 < \frac{1 - \nu}{2} \\ \nu \leqslant \alpha_2 < \min\left(\alpha_1, \frac{1}{2}(1 - \alpha_1) \right) \\ \alpha_1 + \alpha_2 \in \left(\frac{1 - \nu}{2},\ \frac{1 + \nu}{2}\right) }} \psi\left(\beta, p_2 \right) - \sum_{\substack{n = p_1 p_2 \beta \\ \frac{1 + \nu}{4} \leqslant \alpha_1 < \frac{1 - \nu}{2} \\ \nu \leqslant \alpha_2 < \min\left(\alpha_1, \frac{1}{2}(1 - \alpha_1) \right) \\ \alpha_1 + \alpha_2 \notin \left(\frac{1 - \nu}{2},\ \frac{1 + \nu}{2}\right) \\ \alpha_2 \geqslant \frac{4 - 7 \theta + \nu}{2} }} \psi\left(\beta, p_2 \right) \\
\nonumber & - \sum_{\substack{n = p_1 p_2 \beta \\ \frac{1 + \nu}{4} \leqslant \alpha_1 < \frac{1 - \nu}{2} \\ \nu \leqslant \alpha_2 < \min\left(\alpha_1, \frac{1}{2}(1 - \alpha_1) \right) \\ \alpha_1 + \alpha_2 \notin \left(\frac{1 - \nu}{2},\ \frac{1 + \nu}{2}\right) \\ \alpha_2 < \frac{4 - 7 \theta + \nu}{2} }} \psi\left(\beta, p_2 \right) \\
\nonumber =&\ \sum_{\substack{n = p_1 \beta \\ \frac{1 + \nu}{4} \leqslant \alpha_1 < \frac{1 - \nu}{2} }} \psi\left(\beta, x^{\nu} \right) - \sum_{\substack{n = p_1 p_2 \beta \\ \frac{1 + \nu}{4} \leqslant \alpha_1 < \frac{1 - \nu}{2} \\ \nu \leqslant \alpha_2 < \min\left(\alpha_1, \frac{1}{2}(1 - \alpha_1) \right) \\ \alpha_1 + \alpha_2 \in \left(\frac{1 - \nu}{2},\ \frac{1 + \nu}{2}\right) }} \psi\left(\beta, p_2 \right) - \sum_{\substack{n = p_1 p_2 \beta \\ \frac{1 + \nu}{4} \leqslant \alpha_1 < \frac{1 - \nu}{2} \\ \nu \leqslant \alpha_2 < \min\left(\alpha_1, \frac{1}{2}(1 - \alpha_1) \right) \\ \alpha_1 + \alpha_2 \notin \left(\frac{1 - \nu}{2},\ \frac{1 + \nu}{2}\right) \\ \alpha_2 \geqslant \frac{4 - 7 \theta + \nu}{2} }} \psi\left(\beta, p_2 \right) \\
\nonumber & - \sum_{\substack{n = p_1 p_2 \beta \\ \frac{1 + \nu}{4} \leqslant \alpha_1 < \frac{1 - \nu}{2} \\ \nu \leqslant \alpha_2 < \min\left(\alpha_1, \frac{1}{2}(1 - \alpha_1) \right) \\ \alpha_1 + \alpha_2 \notin \left(\frac{1 - \nu}{2},\ \frac{1 + \nu}{2}\right) \\ \alpha_2 < \frac{4 - 7 \theta + \nu}{2} }} \psi\left(\beta, x^{\nu} \right) + \sum_{\substack{n = p_1 p_2 p_3 \beta \\ \frac{1 + \nu}{4} \leqslant \alpha_1 < \frac{1 - \nu}{2} \\ \nu \leqslant \alpha_2 < \min\left(\alpha_1, \frac{1}{2}(1 - \alpha_1) \right) \\ \alpha_1 + \alpha_2 \notin \left(\frac{1 - \nu}{2},\ \frac{1 + \nu}{2}\right) \\ \alpha_2 < \frac{4 - 7 \theta + \nu}{2} \\ \nu \leqslant \alpha_3 < \min\left(\alpha_2, \frac{1}{2}(1 - \alpha_1 - \alpha_2) \right) \\ \boldsymbol{\alpha}_3 \text{ satisfies Lemma~\ref{l69}} }} \psi\left(\beta, p_3 \right) + \sum_{\substack{n = p_1 p_2 p_3 \beta \\ \frac{1 + \nu}{4} \leqslant \alpha_1 < \frac{1 - \nu}{2} \\ \nu \leqslant \alpha_2 < \min\left(\alpha_1, \frac{1}{2}(1 - \alpha_1) \right) \\ \alpha_1 + \alpha_2 \notin \left(\frac{1 - \nu}{2},\ \frac{1 + \nu}{2}\right) \\ \alpha_2 < \frac{4 - 7 \theta + \nu}{2} \\ \nu \leqslant \alpha_3 < \min\left(\alpha_2, \frac{1}{2}(1 - \alpha_1 - \alpha_2) \right) \\ \boldsymbol{\alpha}_3 \text{ does not satisfy Lemma~\ref{l69}} }} \psi\left(\beta, p_3 \right) \\
=&\ S_{6621} - S_{6622} - S_{6623} - S_{6624} + S_{6625} + S_{6626}.
\end{align}
We have (250) holds for $S_{6621}$, $S_{6624}$ (by Lemma~\ref{l610}) and $S_{6622}$, $S_{6625}$ (by Lemma~\ref{l69}). For $S_{6623}$ we note that $\alpha_1 > \alpha_2 \geqslant \frac{4 - 7 \theta + \nu}{2} > \frac{1 - \nu}{4}$ when $0.525 \leqslant \theta \leqslant 0.5253$, hence the condition $\alpha_1 + \alpha_2 \notin \left(\frac{1 - \nu}{2},\ \frac{1 + \nu}{2}\right)$ is equivalent to $\alpha_1 + \alpha_2 \geqslant \frac{1 + \nu}{2}$ in this sum. Then we have
\begin{equation}
S_{6623} = \sum_{\substack{n = p_1 p_2 \beta \\ \frac{1 + \nu}{4} \leqslant \alpha_1 < \frac{1 - \nu}{2} \\ \nu \leqslant \alpha_2 < \min\left(\alpha_1, \frac{1}{2}(1 - \alpha_1) \right) \\ \alpha_1 + \alpha_2 \notin \left(\frac{1 - \nu}{2},\ \frac{1 + \nu}{2}\right) \\ \alpha_2 \geqslant \frac{4 - 7 \theta + \nu}{2} }} \psi\left(\beta, p_2 \right) = \sum_{\substack{n = p_1 p_2 \beta \\ \frac{1 + \nu}{4} \leqslant \alpha_1 < \frac{1 - \nu}{2} \\ \nu \leqslant \alpha_2 < \min\left(\alpha_1, \frac{1}{2}(1 - \alpha_1) \right) \\ \alpha_1 + \alpha_2 \geqslant \frac{1 + \nu}{2} \\ \alpha_2 \geqslant \frac{4 - 7 \theta + \nu}{2} }} \psi\left(\beta, p_2 \right) = S_{65C}
\end{equation}
where $S_{65C}$ is defined in Subsection 6.3. By Lemma~\ref{l611}, we know that (250) holds for this sum. We discard the remaining sum $S_{6626}$. Note that this process replaces a ``larger'' one-dimensional loss from $S_{662}$ with a ``smaller'' three-dimensional loss $S_{6626}$.

Numerical calculations show the following upper bounds for $C_1^{\delta}(\theta)$.
\begin{center}
\begin{tabular}{|c|c|c|c|}
\hline \boldmath{$\theta$} & \boldmath{$C_1^{\delta}(\theta)$} & \boldmath{$\theta$} & \boldmath{$C_1^{\delta}(\theta)$} \\
\hline $0.5250$ & $1.0001$ & $0.5264$ & $1.6367$ \\
\hline $0.5251$ & $1.0004$ & $0.5265$ & $1.6515$ \\
\hline $0.5252$ & $1.0011$ & $0.5266$ & $1.6660$ \\
\hline $0.5253$ & $1.0028$ & $0.5267$ & $1.6811$ \\
\hline $0.5254$ & $1.5112$ & $0.5268$ & $1.6949$ \\
\hline $0.5255$ & $1.5231$ & $0.5269$ & $1.7110$ \\
\hline $0.5256$ & $1.5349$ & $0.5270$ & $1.7269$ \\
\hline $0.5257$ & $1.5471$ & $0.5271$ & $1.7452$ \\
\hline $0.5258$ & $1.5594$ & $0.5272$ & $1.7584$ \\
\hline $0.5259$ & $1.5719$ & $0.5273$ & $1.7788$ \\
\hline $0.5260$ & $1.5843$ & $0.5274$ & $1.7961$ \\
\hline $0.5261$ & $1.5972$ & $0.5275$ & $1.8139$ \\
\hline $0.5262$ & $1.6103$ & $0.5276$ & $1.8332$ \\
\hline $0.5263$ & $1.6238$ & $0.5277$ & $1.8542$ \\
\hline
\end{tabular} \\
\textbf{Table 6.3: Upper Bounds for }\boldmath{$C_1^{\delta}(\theta)$} \textbf{(}\boldmath{$0.525 \leqslant \theta < \frac{19}{36}$}\textbf{)}
\end{center}

\subsubsection{Case 2. $\frac{19}{36} \leqslant \theta < \frac{9}{17}$} Again, one can just replace the parameter $\nu$ occurred above with $\nu^{\prime}$ and calculate the total loss. Working like the first case we get the following upper bounds for $C_1^{\delta}(\theta)$.
\begin{center}
\begin{tabular}{|c|c|}
\hline \boldmath{$\theta$} & \boldmath{$C_1^{\delta}(\theta)$} \\
\hline $0.5278$ & $1.8722$ \\
\hline $0.5279$ & $1.8885$ \\
\hline $0.5280$ & $1.9126$ \\
\hline $0.5281$ & $1.9315$ \\
\hline $0.5282$ & $1.9563$ \\
\hline $0.5283$ & $1.9751$ \\
\hline $0.5284$ & $1.9971$ \\
\hline $0.5285$ & $2.0156$ \\
\hline $0.5286$ & $2.0401$ \\
\hline $0.5287$ & $2.0649$ \\
\hline $0.5288$ & $2.0859$ \\
\hline $0.5289$ & $2.1117$ \\
\hline $0.5290$ & $2.1353$ \\
\hline $0.5291$ & $2.1628$ \\
\hline $0.5292$ & $2.1954$ \\
\hline $0.5293$ & $2.2446$ \\
\hline $0.5294$ & $2.2963$ \\
\hline
\end{tabular} \\
\textbf{Table 6.4: Upper Bounds for }\boldmath{$C_1^{\delta}(\theta)$} \textbf{(}\boldmath{$\frac{19}{36} \leqslant \theta < \frac{9}{17}$}\textbf{)}
\end{center}

\section{Lower Bounds: A General Case}
In this section, we focus on a general form of Mikawa's modified sieve \cite{Mikawa}. The sets $\mathcal{A}$ and $\mathcal{B}$ in this section can be other ``comparison'' sets, and the ``proportion'' may not be $\frac{1}{\varphi(q)}$ as in previous sections. Assume that variants of (34) and (38), with different proportion, hold true. Assume further that there are lots of Type-II information inputs so that the corresponding loss integrals (48) and (49) are zero. We want to find the minimum value of $\kappa$ such that the sum of other loss integrals that cannot be reduced using Type-II information is less than $1$. Let $\kappa \leqslant \frac{1}{8}$. For simplicity, we ignore the new integrals corresponding to sums that count numbers with $7$ or more prime factors, and only consider the 3 loss integrals corresponding to (43)--(45). Define
\begin{align}
\nonumber L_{7}(\kappa) =&\ 2 \int_{(t_1, t_2, t_3) \in U_{71}} \frac{1}{t_1 t_2 t_3 (1 - t_1 - t_2 - t_3)} d t_3 d t_2 d t_1 \\
\nonumber & + 2 \int_{(t_1, t_2, t_3, t_4) \in U_{72}} \frac{1}{t_1 t_2 t_3 t_4 (1 - t_1 - t_2 - t_3 - t_4)} d t_4 d t_3 d t_2 d t_1 \\
& + 20 \int_{(t_1, t_2, t_3, t_4, t_5) \in U_{73}} \frac{1}{t_1 t_2 t_3 t_4 t_5 (1 - t_1 - t_2 - t_3 - t_4 - t_5)} d t_5 d t_4 d t_3 d t_2 d t_1,
\end{align}
where
\begin{align}
\nonumber U_{71}(\boldsymbol{\alpha}_{3}) :=&\ \left\{ \kappa < \alpha_3 < \alpha_2 < \alpha_1,\ \alpha_1 + \alpha_2 + \alpha_3 > \frac{1}{2},\ 2 \alpha_1 + 2 \alpha_2 + \alpha_3 < 1 \right\}, \\
\nonumber U_{72}(\boldsymbol{\alpha}_{4}) :=&\ \left\{ \kappa < \alpha_4 < \alpha_3 < \alpha_2 < \alpha_1,\ 2 \alpha_1 + 2 \alpha_2 + \alpha_3 < 1 \right\}, \\
\nonumber U_{73}(\boldsymbol{\alpha}_{5}) :=&\ \left\{ \kappa < \alpha_5 < \alpha_4 < \alpha_3 < \alpha_2 < \alpha_1,\ 2 \alpha_1 + \alpha_2 + \alpha_3 + \alpha_4 + \alpha_5 < 1 \right\}.
\end{align}
Numerical calculations show that $L_7(\frac{1}{11}) < 0.84$ and $L_7(\frac{1}{12}) > 1.2$. This means that we cannot get a nontrivial lower bound using Mikawa's sieve if we only have a Type-II range $\left[\varepsilon,\ \frac{1}{12} \right]$ or even ``narrower'' ranges.

Another interesting question about Mikawa's sieve is: Can we apply this sieve on other sieve problems? In another paper, Mikawa \cite{MikawaGoldbach} used his method to study the distribution of Goldbach numbers in arithmetic progressions. However, we cannot apply Mikawa's sieve on many other classical sieve problems, such as primes in all short intervals and primes in almost all short intervals. The most important reason of that is the arithmetic information inputs in those problems often have extra restrictions, such as \textbf{Condition B} on one coefficient and the ``prime-factored'' condition (see [\cite{HarmanBOOK}, Chapter 7]). We cannot use those ``restricted'' arithmetic information inputs to prove that variants of (34) and (38) hold true.

\bibliographystyle{plain}
\bibliography{bib}
\end{document}